\magnification=1050
\hoffset=.0cm
\voffset=.0cm
\baselineskip=.55cm plus .55mm minus .55mm
\input eplain
\input epsf
\let\footnote=\numberedfootnote
\everyfootnote={\sevenrm}
 
\input amssym.def\input amssym.tex

\font\grassettogreco=cmmib10
\font\scriptgrassettogreco=cmmib7
\font\scriptscriptgrassettogreco=cmmib10 at 5 truept
\textfont13=\grassettogreco
\scriptfont13=\scriptgrassettogreco
\scriptscriptfont13=\scriptscriptgrassettogreco


\font\sansserif=cmss10
\font\scriptsansserif=cmss10 at 7 truept
\font\scriptscriptsansserif=cmss10 at 5 truept
\textfont14=\sansserif
\scriptfont14=\scriptsansserif
\scriptscriptfont14=\scriptscriptsansserif


\font\capital=rsfs10
\font\scriptcapital=rsfs10 at 7 truept
\font\scriptscriptcapital=rsfs10 at 5 truept
\textfont15=\capital
\scriptfont15=\scriptcapital
\scriptscriptfont15=\scriptscriptcapital


\font\euler=eusm10
\font\scripteuler=eusm7
\font\scriptscripteuler=eusm5 
\textfont12=\euler
\scriptfont12=\scripteuler
\scriptscriptfont12=\scriptscripteuler

\font\cps=cmcsc10
%

%

 %
%
%
%
%
\def\cite{\ref}
\def\ref#1{\lbrack#1\rbrack}
\def\Cal{\cal}
\def\paper{\it }
%
%
%
%
%
\def\dim{{\rm dim}\hskip 1pt}

\def\deg{{\rm deg}\hskip 1pt}

\def\ker{{\rm ker}\hskip 1pt}

\def\GL{{\rm GL}\hskip 1pt}

\def\Sp{{\rm Sp}\hskip 1pt}

\def\PGL{{\rm PGL}\hskip 1pt}
\def\Aut {{\rm Aut}\hskip 1 pt}
\def\Out {{\rm Out} \hskip 1 pt}

\def\Hom{{\rm Hom}\hskip 1pt}
\def\Out{{\rm Out}\hskip 1pt}

\def\Res{{\rm Res}\hskip 1pt}
\def\Ind{{\bf Ind}\hskip 1pt}
\def\Cores{{\rm Cores}\hskip 1pt}
\def\Card{{\rm Card}\hskip 1pt}
\def\Gal{{\bf Gal}\hskip 1pt}
\def\ppcm{{\rm ppcm}\hskip 1pt}
\def\pgcd{{\rm pgcd}\hskip 1pt}
\def\qed{{\makeblankbox{.2pt}{.2pt}} }
\def\Ni{{\bf Ni}\hskip 1pt}
\def\Spec{{\rm Spec}\hskip 1pt}
\def\Div{{\rm Div}\hskip 1pt}
\def\Irrep{{\rm Irrep}\hskip 1pt}
\def\Ext{{\rm Ext}\hskip 1pt}
\def\Isom{{\rm Isom}\hskip 1pt}
\def\Pic{{\rm Pic}\hskip 1pt}
\def\st{{\rm st}\hskip 1pt}
\def\Ima{{\rm Im} \hskip 1pt}
\def\Proj{{\rm Proj}\hskip 1pt}
\def\hst1{\hskip 1pt}
\def\R{{\rm R}\hskip 1pt}
\def\B{{\bf B}\hskip 1pt}
\def\Irrep{{\rm Irrep}\hskip 1pt}

\def\Ž {\'e}
\def\{\`e}
\def\M{{\rm M}\hskip 1pt}
%
%
%
%
%

\def\sectspace{\par\removelastskip\penalty-200
\vskip 18pt plus 4pt minus 6pt}
\font\larbf=cmbx10 scaled \magstep1
\def\section#1#2{\sectspace\message{Section #1. #2}
\centerline{\larbf #1. #2.}\nobreak
\bigskip\noindent\ignorespaces}


\def\rema#1{\vskip.3cm${\hbox{\cps Remarque {#1}}\vphantom{f}}~~$}
\def\exa#1{\vskip.3cm${\hbox{\cps Exemple {#1}}\vphantom{f}}~~$}
\def\dem{\vskip.3cm${\hbox{\it Preuve:}}~~~$}
\def\titlebf#1{\vskip.5cm${\underline{\hbox{\bf #1}}}$\vskip.5cm}

\def\tag{\eqno}


\hrule\vskip.5cm
\hbox to 16.5 truecm{  {Janvier 2007}   {version 2}   }
 
\vskip.5cm\hrule
\vskip.9cm
 
 \centerline{ \bf  CHAMPS DE HURWITZ}   
 \vskip.4cm
\centerline{par}
\vskip.4cm
\centerline{\bf  Jos\'e Bertin  et  Matthieu Romagny}
\bigskip
\centerline{\it  Institut Fourier, Universit\'e de Grenoble 1, F}
\centerline{\it Institut de Math\'ematiques de Jussieu, Universit\'e de Paris 6, F}
 
\centerline{\it E--mail:  jose.bertin@ujf-grenoble.fr, romagny@math.jussieu.fr}
\vskip.9cm
\hrule
\vskip.6cm
\centerline{\bf Abstract} 
\vskip.4cm
\par\noindent
 
 In this work we give a thorough study of Hurwitz stacks and  associated Hurwitz moduli spaces, 
both in  Galois  and  non Galois case. We compare our construction to those proposed by Harris-Mumford, Abramovich-Corti-Vistoli and Mochizuki-Wewers. We apply our  results to revisit some classical examples, 
 particularly the  stacks of  stable curves equipped with an arbitrary level-structure, and the stack of  tamely ramified cyclic covers. 
In a second part  we  exhibit some tautological  bundles and cohomology classes  naturally living on  Hurwitz stacks, and give 
some universal relations, in particular  a higher analogue of 
the Riemann-Hurwitz formula, between  these classes. Applications are given to the stack of cyclic covers of the projective line, with special attention  to Cornalba-Harris type  relations and to   Hyperelliptic Hodge integrals.   
\vskip.6cm
\hrule
\vskip.6cm
\par\noindent
Subject class: 14D15,14D20,14D22,14H30.  Mots cl\'es:  algebraic stack,   covering of curves, formal deformation,  group action, moduli  of curves.  \vfill\eject

\titlebf{Contenu}  
 \noindent {1. Introduction, Contenu}

\noindent {2. Classification des rev\^etements}

\indent {2.1. Actions de groupes et rev\^etements}

\indent {2.2. Donn\'ee de Hurwitz}

\indent {2.3. Nombre de Nielsen }

 \noindent {3. Familles de $G$-courbes lisses }

\indent {3.1.  G\'eom\'etrie du diviseur de branchement}

\indent {3.2. Inversion de la formule de Chevalley - Weil }

\noindent {4. Familles de $G$-courbes stables}

\indent {4.1. Actions stables, rev\^etements stables }

\indent {4.2. Collision des points de ramification }

\indent {4.3. Courbes stables marqu\'ees et actions de groupes}

\noindent {5.  D\'eformations des rev\^etements mod\'er\'ement
ramifi\'es }

\indent {5.1. D\'eformations \'equivariantes des courbes}

\indent {5.2. D\'eformations des rev\^etements mod\'er\'ement
ramifi\'es}

\indent {5.3. Mod\`ele stable marqu\'e d'un rev\^etement}

\noindent {6. Champs de Hurwitz }

\indent {6.1. $G$-champs et champs quotients}

\indent {6.2. Champs de Hurwitz }

\indent {6.3. Compactification du  sch\'ema  de Hurwitz:  Gieseker-Mumford}

\indent {6.4. Compactification stable du champ de Hurwitz (II)}

\indent {6.5. Champs de rev\^etements stables: cas non galoisien}

\noindent {7. Graphes et rev\^etements}

\indent {7.1. Graphes modulaires de Hurwitz}

\indent {7.2. Graphes de groupes et rev\^etements de graphes}

\indent {7.3. Groupe de Picard et rev\^etements}

\indent {7.4. Stratification canonique du bord}

\indent {7.5. Type topologique d'un point du bord}

\noindent {8. Structures de niveau  sur les courbes stables }

\indent {8.1. Le champ des courbes stables avec structure de niveau $G$}

\indent {8.2. Le niveau ab\'elien ($n$)}

\noindent {9. Rev\^etements cycliques}

\indent {9.1. Rev\^etements cycliques versus racines  d'un faisceau inversible}

\indent {9.2.  Rev\^etements cycliques stables}

\indent {9.3. Le champ ouvert des rev\^etements
cycliques de $\Bbb P^1$}

\noindent {10. Groupe de Picard et classes tautologiques}

\indent {10.1. Fibr\'e de Hodge}

\indent {10.2. Les fibr\'es en droites $\psi_{i,\chi}$ et $\mu_{i,v}$}

 \indent {10.3. Relations  de Riemann-Hurwitz d'ordre sup\'erieur}

\indent {10.4. Rev\^etements cycliques}

\vfill\eject

\section1{Introduction}

L'objectif de ce travail est   l'\'etude
syst\'ematique des champs de rev\^etements mod\'er\'ement ramifi\'es entre courbes alg\'ebriques. Rappelons que classiquement
\ref {35},  par {\sl espace de Hurwitz}, on entend la vari\'et\'e
alg\'ebrique $H_{g,d}$ param\'etrant les rev\^etements
$\pi: C\longrightarrow {\Bbb P}^1$ de degr\'e $d$ de la droite
projective, par une courbe lisse de genre $g$, et d\'efinis sur un
corps alg\'ebriquement clos $k$. {\sl Mod\'er\'e} signifie que l'ordre du
groupe de monodromie, c'est-\`a-dire le groupe de Galois de la
cl\^oture galoisienne, a un ordre premier \`a la caract\'eristique $p$
de $k$, si $p>0$.
Ces espaces (ou mieux les champs associ\'es ${\Cal H}_{g,d}$),  et
leurs compactifications naturelles
$\overline H_{g,d}$ et $\overline {\Cal H}_{g,d},$ jouent un r\^ole
significatif dans de nombreuses
questions. Citons comme exemples  du c\^ot\'e arithm\'etique,
la formulation g\'eom\'etrique du probl\`eme de Galois inverse
par Fried \ref {32}, Fried-V\"{o}lklein \ref {34}, mais aussi la
g\'en\'eralisation des tours modulaires sugg\'er\'ee par Fried \cite
{33}. Du c\^ot\'e g\'eom\'etrique  les r\'ecents et
profonds d\'eveloppements autour de l'\'etude asymptotique des
nombres de Hurwitz, en d'autres termes la th\'eorie de Gromov-Witten
sur ${\Bbb P}^1$,   d\'evelopp\'ee par Okounkov-Pandharipande
\cite {59},  Graber-Vakil \cite {38}, et
Ekedahl-Lando-Shapiro-Vainshtein \cite {27} ont une expression en termes de champs de Hurwitz.  Signalons aussi
l'\'etude parall\`ele \`a certains \'egards de l'espace des
modules des courbes avec structure de spin (dans un sens
g\'en\'eralis\'e) par Jarvis \cite {43}, \cite {44}, et
Jarvis-Kimura \cite {45}.   Un point de vue diff\'erent,
celui des applications stables de Kontsevich  $C\longrightarrow {\Bbb
P}^1$ \cite {49}, \cite {53}, s'est   impos\'e comme une approche
 alternative pour aborder  ce genre de questions, du moins lorsque la base est rigide.  Il n'est pas utilis\'e dans notre travail.  
 
  Les  champs  de Hurwitz   seront  essentiellement vus
comme correspondances entre espaces de modules de courbes. C'est   
  cette position qui a \'et\'e   la  motivation principale pour
entreprendre le pr\'esent travail.  On notera que la construction
des espaces modulaires de Hurwitz  a \'et\'e abord\'ee  \`a
plusieurs reprises dans la litt\'erature, mais sous des
angles particuliers, le plus souvent pour les rev\^etements de base
$\Bbb P^1$, et avec un traitement  partiel ou imparfait de la
compactification stable marqu\'ee. On se reportera au rapport de D\`ebes \cite {16},  ou \`a Wewers \cite {68}, \cite {69}.

 Notre objectif  est d'abord de proposer une
construction uniforme, puis ensuite une \'etude  d\'etaill\'ee des
champs ${\Cal H}_{g,h,d}$,   en corollaire  des espaces de modules
$H_{g,h,d}$, classifiant les rev\^etements  $\pi: C \longrightarrow
D$ d\'efinis sur un corps alg\'ebriquement clos
$k$, de degr\'e $d$, d'une courbe de genre $h$ par une courbe de
genre $g$, la ramification \'etant suppos\'ee mod\'er\'ee. Si $h =
0$, c'est la situation du d\'ebut; si $d = 1$,
${\Cal H}_{g,g,d} = {\Cal M}_g$.  On observera qu'en g\'en\'eral
${\Cal H}_{g,h,d}$ n'est pas connexe, mais poss\`ede la propri\'et\'e
 que les composantes connexes sont irr\'eductibles; propri\'et\'e
partag\'ee avec tout champ localement noeth\'erien et normal. On 
s\'epare de mani\`ere grossi\`ere ces composantes en fixant le type
de la ramification, appel\'e \/  la donn\'ee de Hurwitz, ou la
donn\'ee de ramification, et  not\'ee \/ $\xi$  dans le texte
(section 2.2). L'espace de Hurwitz  $H_{g,h,\xi}$ qui en d\'ecoule  n'est  en g\'en\'eral pas connexe. 
Le  nombre de composantes connexes (ou irr\'eductibles) est le nombre de Nielsen (D\'efinition
2.4). La contribution  principale du pr\'esent travail est la
construction, puis l'\'etude,  de la compactification naturelle  par
adjonction des rev\^etements stables, $\overline {\Cal H}_{g,h,\xi}$
de ${\Cal H}_{g,h,\xi}$.  A ce stade les choses sont plus subtiles, du fait que plusieurs choix naturels sont possibles, marquage ou pas par le diviseur de ramification.

Il doit \^etre not\'e que des constructions partielles  de l'espace
de Hurwitz  sont pr\'esentes dans la litt\'erature,  d'abord par
Fulton \cite {35} dans le cas des rev\^etements simples de la droite
projective, puis dans une logique diff\'erente par Fried \cite {32}.
On notera cependant que  l'espace modulaire consid\'er\'e par
Fried, et qui concerne la base ${\Bbb P}^1$, n'est pas identique au
n\^otre, bien que directement li\'e;  le champ de Hurwitz est le
champ quotient de celui  de Fried par $\PGL (1)$. La comparaison
entre  ces deux constructions, la construction de Fried  et la pr\'esente  construction, est  analys\'ee par D\`ebes \cite {16}, voir aussi Emsalem \cite {29},\cite {30}. Une  approche plus syst\'ematique de la construction du sch\'ema de Hurwitz, ainsi que d'une
compactification naturelle par adjonction des rev\^etements
admissibles, \`a sa source dans le travail de Harris-Mumford \cite
{42}, et dans les raffinements propos\'es par Mochizuki \cite {55},
Wewers  \cite {68} \cite {69}.

Les champs de Hurwitz sont, comme on va le voir, d'une extr\^eme
souplesse, du fait qu'il est possible de jouer avec les param\`etres:
groupe, donn\'ee de ramification et marquage.  Le champ ${\Cal
H}_{g,h,\xi}$ (resp. $\overline {\Cal H}_{g,G,\xi}$) est
du fait de sa d\'efinition, la source d'un morphisme naturel  et
important (le discriminant):
$(C\rightarrow D) \mapsto (D , (\{Q_j\}),$ qui associe \`a un
rev\^etement sa base marqu\'ee par les
points de branchement, d\'efinissant, si   $r$  est le nombre de
points de branchement, un morphisme
${\Cal H}_{g,h,\xi} \longrightarrow {\Cal M}_{h,r}$.
Plus important est le fait que ce morphisme se prolonge en $\overline
{\Cal H}_{g,h,\xi} \longrightarrow \overline {\Cal M}_{h,r}$, c'est
\`a dire aux compactifications stables.  On observera qu'un rev\^etement
stable est essentiellement stable marqu\'e par les points de
ramification, avec une subtilit\'e li\'ee \`a la coalescence \'eventuelle des points de ramification d'indice deux. De cette mani\`ere les champs (espaces) de Hurwitz
apparaissent comme des rev\^etements (ramifi\'es) de  ${\Cal
M}_{h,r}$, (resp. $\overline {\Cal M}_{h,r}$), l'espace des modules
des courbes de genre $h$ avec $r$ points marqu\'es.  On notera que
ces espaces contiennent comme cas particuliers les espaces modulaires
de courbes avec structure de niveau (\S~8). C'est de mani\`ere
simplifi\'ee, un point de vue finalement  assez fructueux, exploit\'e
 dans \cite {3},  de voir un rev\^etement comme  la donn\'ee d'une
sorte de structure de niveau (non ab\'elien) sur la base, et r\'eciproquement.

La construction de la compactification  stable du champ de Hurwitz,
model\'ee sur celle $\overline
{\Cal M}_{h,r}$ pour  ${\Cal M}_{h,r}$, est le sujet de la section 6.
L'id\'ee, somme toute banale, revient \`a ajouter \`a ${\Cal
H}_{g,h,\xi}$ les rev\^etements dits "stables", qui doivent \^etre
distingu\'es des rev\^etements admissibles introduits dans
Harris-Mumford \cite {40}, \cite {42}, et des applications stables de
Kontsevich \cite {49}. Contrairement aux applications stables, un
rev\^etement stable
est un objet qui, dans le cadre des rev\^etements, provient par sp\'ecialisation d'un rev\^etement entre courbes lisses,  ou d'une autre
mani\`ere  est tel que le rev\^etement g\'en\'erique dans sa d\'eformation universelle est  un rev\^etement entre courbes lisses, i.e.
qui est lissifiable.  Un point du bord est repr\'esent\'e,
dans le cas galoisien, par une courbe stable marqu\'ee par les points
de ramification, ce qui interdit \`a ceux-ci \`a se localiser en des
points doubles, cela joint \`a une condition de stabilit\'e de l'action du groupe de Galois $G$ du
rev\^etement. On notera que les points doubles qui peuvent avoir un
stabilisateur non trivial, alors cyclique,  ne sont pas consid\'er\'es comme des points de ramification, par contre  ils interviennent dans la d\'eformation universelle de $C \rightarrow D$ (voir la section 5 pour des d\'etails).  
Il y a  une d\'efinition alternative au concept de stabilit\'e, qui revient \`a ignorer le marquage par le diviseur de ramification. En effet les points de ramification ne sont pas mobiles, ils sont assign\'es par l'action du groupe de Galois. Bien que sans importance pour les espaces de modules, les deux d\'efinitions peuvent conduire \`a des champs distincts (voir \S 6.4).   On notera qu'un mod\`ele de cette  construction \'etait d\'eja
visible dans la compactification des rev\^etements doubles utilis\'ee
par Beauville \cite {7}.

Dans le cas non galoisien, un rev\^etement stable sera d\'efini
(section 6) comme  un rev\^etement qui \'etale-localement  admet une
cl\^oture galoisienne (voir section 6.5).
Ces uniformisations galoisiennes locales   doivent satisfaire \`a
une condition de coh\'erence. Un rev\^etement stable est admissible,
mais la r\'eciproque n'est pas vraie : les champs correspondants
diff\`erent.  Le champ de Hurwitz est de mani\`ere pr\'ecise la
normalisation  (d\'esingularisation) du champ de Harris-Mumford. On
doit noter qu'on arrive essentiellement au m\^eme r\'esultat en
enrichissant un rev\^etement admissible d'une log-structure,
construction \'etudi\'ee  par  Mochizuki \cite {55}, et Wewers \cite
{68}.  Cette construction diff\`ere l\'eg\`erement de la notre. 
C'est  d'ailleurs une vertu
essentielle des champs de Hurwitz, ou de ceux qui en d\'erivent, que
de proposer des compactifications lisses aux divers probl\`emes de
modules pos\'es par les rev\^etements mod\'er\'ement ramifi\'es.
Cela sera v\'erifi\'e pour le champ ${\Cal M}_g (G)$ des
structures de niveau $G$ sur les courbes lisses de genre $g$ ($g\geq
2$), i.e. des $G$-rev\^etements principaux de base une courbe lisse
de genre $g$. Nous d\'ecrirons (section 8) une compactification lisse
$\overline {\Cal M}_g (G)$ de ${\Cal M}_g (G)$, \`a la diff\'erence
de Deligne-Mumford \cite {20}, mais  par contre non repr\'esentable,
m\^eme si le niveau domine le niveau ab\'elien $(n), n\geq 3$. Les
objets de ce champ seront les $G$-rev\^etements principaux d\'eg\'en\'er\'es stables. Un objectif similaire au n\^otre  a \'et\'e
atteint, cependant  par des m\'ethodes et  points de vue diff\'erents, par Abramovich-Vistoli et leurs coauteurs  \cite {1}, \cite {2},
\cite {3}. Les objets  qu'ils d\'ecrivent sous le vocable de
"balanced twisted stable maps", constituent les points d'un champ de
Deligne-Mumford lisse et propre, et qui n'est pas autre chose que le
champ (quotient) (\S 6.4)  :
$$ \overline {\Cal M}_{g',b} (BG) \,  = \,
\coprod_{[\xi] }\, \overline {\Cal H}_{g,G,\xi} // Z(G)$$la somme \'etant sur les donn\'ees de Hurwitz, et $Z(G)$ \'etant le centre de $G$. Le point de vue des rev\^etements adopt\'e i\c{c}i rend les constructions tr\`es naturelles. Cela s'applique \`a l'extension $G$-\'equivariante des th\'eories topologiques des champs par Turaev, Jarvis at al (voir \cite {45},\cite {46}).

D\'ecrivons bri\`evement le contenu  du pr\'esent
travail.  Dans le paragraphe 2, on
introduit les d\'efinitions de base, donn\'ee de Hurwitz, ou de
ramification, ou d'holonomie autour des points de branchement, et les
op\'erations naturelles support\'ees par ces donn\'ees, comme la
restriction et l'induction. La description topologique de l'espace de
Hurwitz est pr\'esent\'ee, pour clarification, et conduit \`a la
d\'efinition du nombre de Nielsen, introduit par Fried \cite {32}.
Dans le paragraphe  3, nous \'etudions les familles de courbes
lisses munies de l'action d'un groupe fini, \`a donn\'ee de Hurwitz
fix\'ee, de mani\`ere \'equivalente les familles de $G$-rev\^etements
galoisiens. 

Le r\'esultat principal de cette section est un th\'eor\`eme d'inversion de la formule bien connue de Chevalley-Weil (Th\'eor\`eme 3.1). On montre que la donn\'ee de Hurwitz est totalement d\'etermin\'ee par la connaissance des repr\'esentations de $G$
dans les espaces $H^0 (C ,\omega_C^{\otimes m})$, en fait  par un
nombre fini d'entre elles. Dans le paragraphe 4, on \'etudie de la
m\^eme mani\`ere les familles $G$-\'equivariantes de courbes stables
marqu\'ees. L'inversion des relations de Chevalley-Weil est l'outil combinatoire 
crucial pour analyser le probl\`eme de la collision des points de ramification,
ph\'enom\`ene qui peut intervenir
dans certains cas, et qui correspond \`a  la pr\'esence de points
doubles \`a isotropie di\'edrale
(Th\'eor\`eme 4.1). L'objectif principal du paragraphe 4 est de
clarifier le comportement des points de
ramification par sp\'ecialisation et d\'eformation. On montre en
particulier qu'une $G$-courbe stable est essentiellement  marqu\'ee
par les points de ramification, la seule complication venant \'eventuellement des points de ramification d'indice deux.

Dans le paragraphe 5, nous \'etudions la d\'eformation verselle d'un
rev\^etement galoisien
mod\'er\'ement ramifi\'e $\pi: C \longrightarrow D$ de diviseur de
branchement $B$. Dans le cas lisse il est bien connu que les d\'eformations d'un rev\^etement correspondent bijectivement \`a celles
de la base \'equip\'ee du diviseur de branchement. Dans le cas stable, le r\'esultat
essentiel  (th\'eor\`eme 5.2),
d\'ecrit l'\'ecart   entre d'une part la d\'eformation verselle \'equivariante de $C$, et d'autre part
la d\'eformation verselle de la base $D$ marqu\'ee par les points de
branchement. 

Dans le paragraphe 6, on construit les champs (espaces)
de Hurwitz et leurs compactifications stables.    On montre, et cette remarque peut avoir une certaine utilit\'e, que
la projectivit\'e des espaces $\overline H_{g,h,\xi}$  peut
s'obtenir par
preuve  similaire \`a celle  donn\'ee par  Gieseker et Mumford  pour
la projectivit\'e de
$\overline M_g$. Cette approche permet  en fait d'obtenir de mani\`ere uniforme la projectivit\'e des
$\overline M_{g,n}$,  par une preuve qui  r\'eduit essentiellement
tout au cas de $\overline M_g$. La raison est bien s\^ur  le fait
que les points marqu\'es sont   implicites dans une action, ce sont
les points de branchements, ou de ramification.   Le th\'eor\`eme 6.2 compare notre construction et celle de Harris-Mumford. On montre
d'abord, dans le cas galoisien, puis dans le cas des rev\^etements
non galoisiens, que  le champ de Hurwitz (ou l'espace modulaire  de
Hurwitz) \'etudi\'e dans le pr\'esent travail, est la d\'esingularisation de celui, non normal, introduit dans \cite {42}.

Dans le paragraphe 7 on \'etudie la combinatoire du bord, qui se r\'esume \`a la version
\'equivariante de la description habituelle du bord de $\overline
{\Cal M}_{g,n}$. Un r\^ole cl\'e est jou\'e dans cette section
par la th\'eorie de Bass \cite {6} des rev\^etements de graphes ;
voir en particulier le  th\'eor\`eme 7.2.
Dans le paragraphe 8, on applique  en guise de test, la
construction des espaces de
Hurwitz \`a la description modulaire du bord des espaces de modules
de courbes avec structure de niveau $G$, ab\'elien ou non. Rappelons
que la  construction usuelle de la  compactification
$\overline M_g(G),$  qui proc\`ede de mani\`ere indirecte, (voir
Deligne-Mumford \cite {21}), ne conduit pas \`a une description
modulaire pr\'ecise des points du bord. Ici, ils apparaissent
comme des rev\^etements principaux  stablement
d\'eg\'en\'er\'es, de groupe $G$. C'est aussi le point
de vue r\'ecent d\'evelopp\'e par  Abramovich-Vistoli \cite {2}.  On retrouve   et on
clarifie quelques r\'esultats de Boggi-Pikaart \cite {13}, Oort-van
Geemen \cite {36}. 

Dans les paragraphes 9 et 10, on  introduit  les
classes     "tautologiques" qui vivent dans le
groupe de  Chow  d'un champ de Hurwitz,  particuli\`erement dans   le cas du
champ classifiant les  rev\^etements cycliques (voir aussi \cite {5}).
On met en \'evidence quelques relations universelles  entre ces \'el\'ements, de mani\`ere analoque au travail de Jarvis sur les courbes
\`a spin \cite {43}. Ces relations contiennent les relations de Cornalba-Harris dans le cas hyperelliptique. On observe que la relation de Riemann-Hurwitz    s'\'etend \`a un ordre sup\'erieur,  sous la  forme d'un relation simple entre classes kappa.  Comme application on calcule certaines int\'egrales de Hodge sur le champ des courbes hyperelliptiques sugg\'er\'ees par le travail de A. Bene \cite{10}. 

Les conventions et notations seront
rappel\'ees si n\'ecessaire en d\'ebut de chaque section.
Les r\'esultats du pr\'esent travail ont \'et\'e annonc\'es
en partie, et sous une forme restreinte,  dans la note \cite {11}. Ils s'appuient  occasionnellement   (\S \  6 \`a \S \ 8) sur la th\`ese du second auteur (voir \cite {61}). 

Nous remercions les  organisateurs de la semaine th\'ematique "Th\'eorie de galois et espaces de modules" au CIRM en Mars 04 durant laquelle plusieurs expos\'es furent consacr\'es aux espaces de Hurwitz.

\vfill\eject

\section2{Classification des rev\^etements}

Dans ce paragraphe, les d\'efinitions
de base concernant les rev\^etements sont introduites. Nous d\'efinissons et \'etudions les invariants fondamentaux, donn\'ee de
Hurwitz, repr\'esentations de Hurwitz.
\medskip

\beginsection 2.1. Actions de groupes et rev\^etements.

\medskip

  Fixons tout d'abord quelques  notations  et conventions communes \`a
 la totalit\'e de l'article,
puis celles particuli\`eres \`a ce paragraphe.  Soit
un corps $k$   alg\'ebriquement clos, de caract\'eristique arbitraire.
Une courbe, du moins lorsqu'il sera fait allusion \`a un corps de
d\'efinition, aura toujours pour sens, une courbe projective r\'eduite, connexe, d\'efinie sur
$k$.  On  se limitera dans la totalit\'e de  ce paragraphe, ainsi
que dans le suivant,  \`a des courbes  non singuli\`eres. Notons
cependant que dans la d\'efinition  qui suit,  comme dans d'autres,
cette restriction est  inutile.  Dans les paragraphes ult\'erieurs
(\S \  4 \`a \S \  8), une courbe sera de mani\`ere plus g\'en\'erale
une courbe  singuli\`ere connexe, avec  pour seules singularit\'es
des points doubles ordinaires,  donc selon une terminologie bien \'etablie, une courbe nodale ou pr\'estable.  le genre d'une telle courbe est $g(C) = p_a(C)$. Les d\'efinitions correspondantes,  sous ces hypoth\`eses,
seront rappel\'ees et pr\'ecis\'ees dans le paragraphe 4 et les
paragraphes ult\'erieurs. 
 \bigskip
{\it 2.1.1.  $G$-rev\^etements} 
\bigskip 
 Fixons un groupe fini $G$. On supposera, et cela de mani\`ere
permanente que si la caract\'eristique d'un corps $k$ est $p > 0$,
alors $p$ ne divise pas l'ordre de $G$.  De mani\`ere plus g\'en\'erale, dans les manipulations de familles,  on supposera que l'ordre
de $G$ est inversible
dans les  faisceaux structuraux des sch\'emas consid\'er\'es.  Dit
d'une autre mani\`ere, les actions  du groupe $G$ seront  mod\'er\'ement ramifi\'ees. Nous retenons  par ailleurs les conventions
suivantes concernant le groupe $G$. Le groupe des caract\`eres
(multiplicatifs), Hom $(G , k^*)$ sera not\'e 
$\hat G$. Nous noterons  $\Irrep(G)$ l'ensemble des   repr\'esentations (ou caract\`eres) irr\'eductibles de $G$, une m\^eme lettre d\'esignant en g\'en\'eral une repr\'esentation $V$, ou bien sans pr\'ecision
suppl\'ementaire  son caract\`ere $\chi_V: G \to k$. On notera respectivement
$\Aut  (G)$ et $\Out (G)$ le groupe des automorphismes de $G$,
respectivement des
automorphismes ext\'erieurs. Pour $s\in G, \,[[s]]$ d\'esignera la
classe de conjugaison de $s$.

L'objet fondamental \'etudi\'e dans de ce travail est une courbe
$C$,   munie d'une
{\it action } du  groupe $G$,   action signifiant toujours une action
fid\`ele. Celle-ci est donc
d\'ecrite par un morphisme injectif  
$\phi: G\to \Aut (C)$. On d\'esignera alors par $(C ,\phi)$ la donn\'ee de $C$ munie de  l'action
d\'efinie par $\phi$, lorsqu'il sera utile de pr\'eciser celle-ci. La
d\'efinition suivante d'\'equivalence,  pour les courbes munies d'une
action de $G$, est celle qui pr\'evaut  dans tout ce travail. On
notera cependant, et cela justifie qu'elles soient mises en \'evidence, qu'il y a des variantes \`a cette d\'efinition, que nous
allons rappeler,  voir  aussi par exemple, Fried \cite {32}, D\`ebes
\cite {16}, et aussi la section 6.

\proclaim D\'efinition 2.1. Soient deux actions $(C,\phi)$ et
$(C',\phi')$ de $G$ sur des courbes $C$
et $C'$. Ces actions sont  dites $G$-\'equivalentes, en abr\'eg\'e,
\'equivalentes,  s'il existe un $G$-isomorphisme $f: C \to C'$,
c'est \`a dire   un isomorphisme tel que pour tout $x \in C$ et
tout $s \in G$, on ait $f(\phi(s)).x = \phi'(s).f(x)$.

  Dans la suite  on omettra le plus souvent la lettre $\phi$ et on
parlera pour simplifier d'une
$G$-courbe;  l'action sera alors not\'ee de mani\`ere abr\'eg\'ee $(s
, x) \mapsto sx$.
Notons que  les courbes munies d'une action de $G$ que nous auront \`a
 consid\'erer  v\'erifierons en g\'en\'eral  la condition que le
groupe des automorphismes  $\Aut (C)$  est fini; on pourrait
remplacer cette hypoth\`ese par la finitude de $\Aut_G(C)$, groupe
des automorphismes $G$-\'equivariants. C'est le cas  par exemple  si
$C$ est lisse de genre $ g \geq 2$, ou plus g\'en\'eralement si $C$
est stable de genre $g \geq 2$ (section 3, et \cite {20},\ \cite {53}).
  Dans ce paragraphe on se limite  de mani\`ere exclusive  \`a la
consid\'eration de courbes lisses; les r\'esultats seront \'etendus
au paragraphe suivant aux courbes nodales stables.   Sous les
hypoth\`eses  en vigueur, il est bien connu que la courbe quotient
$C/G$ est d\'efinie; nous parlerons du morphisme   quotient $\pi: C
\to C/G$ comme  \'etant le rev\^etement (galoisien) de groupe
$G$ d\'efini par $(C , G)$. 

De mani\`ere plus g\'en\'erale,  un rev\^etement
$\pi : C \to D$, $C$ et $D$ \'etant deux courbes lisses
connexes, est  un morphisme fini  et   g\'en\'eriquement \'etale.
L'hypoth\`ese faite sur le cardinal de $G$, \`a savoir  que  si la
caract\'eristique  de $k$ est $p > 0$,
alors $p$ ne divise pas $\vert G\vert$ a pour traduction que le rev\^etement
associ\'e est
mod\'er\'ement ramifi\'e.  En fait le rev\^etement $ \pi : C
\rightarrow  C/G$ associ\'e \`a $(C
, G)$ est, en utilisant une terminologie r\'epandue, un $G$-rev\^etement,
ce qui signifie qu'on
inclut \`a la donn\'ee de $\pi : C \to D,$  l'action de $G$ sur
$C$, ou ce qui revient au m\^eme,
on fixe une identification $G\cong \Aut(\pi)$. 

Dans cette situation,
donc pour un rev\^etement galoisien mod\'er\'ement ramifi\'e,
rappelons qu'un
  point de {\it ramification} est un point  (ferm\'e) $x\in C$ tel que  si
$G_x$ est le stabilisateur, alors
$G_x \ne 1$. L'orbite d'un tel point sera souvent appel\'ee une
orbite singuli\`ere; on parlera aussi par
opposition d'orbites r\'eguli\`eres. Notons que si $x$ est un point
de ramification le stabilisateur $G_x$ est   cyclique, d'ordre $e(x)$ l'indice de ramification. Les points
images dans $ C/G$ des points de ramification sont appel\'es  les
points de {\it branchement}. Ces d\'efinitions s'\'etendent  comme il est bien connu au cas non galoisien. Finalement la relation de Riemann-Hurwitz pour un rev\^etement non n\'ecessairement galoisien $\pi: C\to D$ de degr\'e $d$ est
$$2g(C) - 2 = d(2g(D) - 2) + \sum_{x\in C} \, (e(x) - 1)$$

\medskip

 {\it 2.1.2.  Classification des rev\^etements}

\bigskip

La classification des
couples $(C , G)$, ou ce qui revient au m\^eme du $G$-rev\^etement
correspondant $\pi: C \to C/G$, introduite ci-dessus,  ne doit pas
\^etre confondue avec ce que l'on entend usuellement par la
classification des rev\^etements \cite {17}, \cite {32}  \cite {55}, \cite
 {67}, \cite {68}. 
 
 A titre de comparaison avec la d\'efinition
initiale,  et pour utilisation ult\'erieure, rappelons que deux rev\^etements,
galoisiens ou non, $ \pi: C\to D$ et $\pi': C'\to D'$ sont
dits {\it \'equivalents}, s'il existe  des isomorphismes $f:
C\buildrel\sim\over\rightarrow C'$ et $h:
D\buildrel\sim\over\rightarrow D'$ tels
que $\pi' f\, =\, h \pi$, soit un diagramme commutatif
\vskip 5pt
$$\commdiag {C & \mapright^f& C' \cr
 \mapdown\lft{\pi}&&\mapdown\lft{\pi'}\cr
  D& \mapright^h& D' \cr} \tag (2.1)$$
  \vskip 5pt

Naturellement si $f: (C , G) \buildrel\sim\over\rightarrow (C' , G)$
est une $G$-\'equivalence,  donc si  les $G$-rev\^etements associ\'es
sont  $G$-\'equivalents, ils sont  alors \'equivalents. La $G$-\'equivalence exige en plus  que l'isomorphisme  $f$ soit un
$G$-isomorphisme. On notera aussi  que pour les rev\^etements de
base fix\'ee $D$,  on peut aussi  choisir comme relation d'\'equivalence,  la relation d'\'equivalence $G$-stricte (resp. stricte). Cela revient \`a imposer $h
= 1$ dans la d\'efinition de dessus.

Dans le cas  galoisien, donc pour des $G$-rev\^etements,  on a  ˆ notre disposition     deux mani\`eres
distinctes d'identifier deux
rev\^etements, selon que l'on exige  que  $f$  soit $G$-\'equivariant, c'est ˆ dire la $G$-\'equivalence, ou pas. Cela correspond classiquement \`a la distinction
  entre les courbes modulaires $X_0(N)$
et $X_1(N)$. Notons que si la base est fix\'ee, \'egale \`a $\Bbb P
^1$ (le sch\'ema de Hurwitz classique), la d\'efinition de l'\'equivalence retenue par Fried \cite {32}, Fulton \cite {35} , est l'\'equivalence stricte (ou rigide), pour laquelle deux rev\^etements  $\pi:
C\rightarrow \Bbb P ^1$ et $\pi': C'\rightarrow \Bbb P ^1$ sont \'equivalents, si et seulement si il existe un isomorphisme \'equivariant $f: C \buildrel\sim\over\longrightarrow C'$ avec $\pi' f =
\pi$; c'est \`a dire avec les notations de dessus, $h = 1$.

Le cas g\'en\'eral des rev\^etements non galoisien sera trait\'e dans un
paragraphe  ult\'erieur (Chap 6,  \S\  6.5). Posons cependant le principe de classification. La classification des rev\^etements non galoisiens   est  celle donn\'ee par un diagramme (2.1). Dans le Chap. 6, elle sera ramen\'ee au cas galoisien par passage \`a la cl\^oture galoisienne.

Dans le reste de ce
paragraphe, les rev\^etements seront galoisiens,  mieux des
$G$-rev\^etements, l'\'equivalence sera la $G$-\'equivalence.
\medskip
 
 \beginsection 2.2. Donn\'ee de Hurwitz.

\medskip

Dans  cette partie,  une courbe est  une courbe projective connexe,
lisse et  d\'efinie sur $k$, corps alg\'ebriquement clos.  On d\'esigne par $G$ un groupe fini d'ordre premier \`a la caract\'eristique de $k$.

\bigskip

 {\it 2.2.1. Donn\'ee de Hurwitz: D\'efinition} 
\bigskip
Soit donc $C$ une courbe  munie
d'une  action, effective rappelons le, du groupe fini $G$. Si $x\in
C$ est un point de ramification, c'est \`a dire \`a isotropie non
triviale $G_x \ne 1$,  il en d\'ecoule  une repr\'esentation
naturelle  $$\chi_x : G_x \to GL(  {{\Cal M}_x \over{\Cal M}_x^2})
= k^*   $$
du groupe cyclique $G_x$ dans l'espace cotangent en $x$ (et en
dualisant dans l'espace tangent).
  Le caract\`ere $\chi_x$ de $G_x$ est primitif, c'est \`a dire
d'ordre $e_x = \vert G_x\vert$.
  Ceci nous am\`ene \`a consid\'erer l'ensemble des couples $(H ,
\chi)$, o\`u $H$ est un sous-groupe cyclique et $\chi$ un caract\`ere primitif de $H$. Introduisons la d\'efinition suivante: nous
dirons que les couples  
$(H , \chi)$ et $(H' , \chi')$ sont conjugu\'es, nous \'ecrirons
alors  $(H , \chi) \sim (H' ,
\chi')$, s'il existe $s\in G$ tel que
$$H' = sHs^{-1} \,\hbox{et}\,\, \chi'(sts^{-1}) = \chi(t)
\,\,\hbox{pour} \,\,t\in H.$$
Ou encore $\chi' = \,{^s\chi}$. Il est manifeste que si $y = sx$ alors $(G_y , \chi_y) \sim (G_x ,
\chi_x)$. Ceci permet de
d\'efinir le $G$- {\it type} d'une orbite singuli\`ere, disons
l'orbite de $x$,  comme \'etant la classe de
conjugaison de $(G_x , \chi_x)$, not\'ee dans la suite  $[G_x ,
\chi_x]$. Ainsi nous parlerons d'une orbite de type $[H , \chi]$.
D'une mani\`ere plus pr\'ecise, nous dirons qu'un point fixe $x$ est
{\it d'holonomie} $(H,\chi)$ si $G_x = H$ et $\chi_x = \chi$.
D\'ecomposons l'ensemble des points de ramification $F$ en une r\'eunion d'orbites (les orbites
singuli\`eres  )
$F = F_1 \cup \cdots \cup F_b$, et notons $[H_i , \chi_i]$ le type de
$F_i$. La d\'efinition suivante
est sous une forme ou une autre classique \cite {32}, \cite {50}, \cite {65}:
\proclaim D\'efinition 2.2.   La donn\'ee de Hurwitz (ou de
ramification) de l'action de $G$ sur $C$ d\'esigne  l'ensemble $\{
[H_i ,\chi_i]\}$ des types distincts on non, des orbites singuli\`eres compt\'ees avec multiplicit\'e. 

Il s'av\`ere parfois commode de sp\'ecifier seulement les types
distints $[H_i , \chi_i] ,
\,1 \leq  i \leq  t$, et de pr\'eciser dans ce cas la multiplicit\'e
de $[H_i , \chi_i]$.
Pour une formulation plus pr\'ecise formons le groupe ab\'elien libre suivant:
$$ R_+(G) \,=\, \bigoplus_{[H,\chi]} \;{\Bbb Z}[H , \chi]  \tag (2.2)$$
la somme directe \'etant index\'ee par les classes de conjugaison de
couples $[H , \chi]$, o\`u
rappelons le, $H$  parcourt les sous groupes cycliques de $G$ et
$\chi : H \to k^*$ les caract\`eres
non n\'ecessairement primitifs   de $H$. Ainsi une donn\'ee de
ramification est un \'el\'ement $\xi$ de $R_+(G)$  qui a ses
coordonn\'ees positives. On peut donc \'ecrire une telle donn\'ee
$\xi = \sum b_i [H_i , \chi_i]$, avec pour tout $i, \, b_i \geq 1)$;
l'entier $b = \sum b_i$ est le degr\'e de $\xi$. Il
est clair que la donn\'ee de Hurwitz ne d\'epend que de la classe
d'\'equivalence de $(C , G)$.
\bigskip
{\it 2.2.2. Op\'erations sur les donn\'ees de
Hurwitz} 
\bigskip
Il y a plusieurs op\'erations naturelles  sur l'ensemble $R_+(G)$  donn\'ees de
Hurwitz. La premi\`ere est l'induction. Soit
$J < G$ un sous-groupe, alors toute donn\'ee relative \`a $J$ d\'efinit
de mani\`ere \'evidente une donn\'ee de $G$, soit une application
$\Ind_J^G: R_+(J) \rightarrow R_+(G)$:
$$ \xi = \sum \,b_i [H_i , \chi_i]_H \mapsto \xi = \sum \,b_i [H_i , \chi_i]_G \tag (2.3)$$
l'indice pr\'ecisant le groupe dans lequel la classe de conjugaison
est prise. Si $J_1, J_2$ sont
deux sous-groupes de $G$ et $\xi_1, \xi_2$ des donn\'ees relatives \`a
 ces deux sous-groupes, on peut ainsi effectuer l'op\'eration
d'induction $\Ind_{J_1}^G (\xi_1) + \Ind_{J_2}^G (\xi_2)$ conduisant \`a
un \'el\'ement de $R_+(G)$.

Il y a deux  autres op\'erations naturelles  concernant les donn\'ees
de ramification  qui seront
utilis\'ees dans la suite, la restriction et  le passage au quotient.
Soit $J < G$ un sous groupe; on
d\'efinit le morphisme de restriction,  $\Res_J^G : R_+(G) \to R_+(J)$ par
  $$\Res_J^G ([H , \chi])\, = \,\sum_{J\setminus G/H} [J\cap gHg^{-1}
, \chi^g\mid _{J\cap gHg^{-1}}]
\tag (2.4)$$
o\`u la sommation porte sur les classes doubles mod. $(J , H)$, et
$\chi^g$ d\'esigne le caract\`ere
$\chi^g (s) = \chi(g^{-1}sg)$.

  Si maintenant $J = \frac G /{K}$ on d\'efinit le morphisme de corestriction
$\Cores_J^G : R_+(G) \to R_+(J) $  par 

$$   \Cores_J^G ([H
, \chi]) =  [ H /{H\cap
K} ,  \chi^{\Card H\cap K\vert }] \tag (2.5)$$

L'interpr\'etation de ces deux op\'erations, au niveau des actions,
est expliqu\'ee dans  l'\'enonc\'e
 ci-dessous, dont la preuve \'el\'ementaire est omise. L'indice 1
signifie qu'on prend la projection sur
le sous groupe $\sum_{H\ne1} {\Bbb Z}[H , \chi]$.

\proclaim Proposition 2.3.
Soit $\xi$ la donn\'ee de Hurwitz de l'action de $G$ sur $C$;
alors  i) $\Res_J^G (\xi)_1$ est la donn\'ee associ\'ee \`a la
restriction de l'action \`a $J$.\hfill\break
ii) $\Cores_J^G (\xi)_1$ est la donn\'ee  associ\'ee \`a l'action
de $J = \frac G /{K}$ sur la
courbe $C/K$.    \hfill\break \qed

 Pour utilisation ult\'erieure, notons qu'il existe une action \'evidente de $\Out(G)$ sur
$R_+(G)$,  d\'efinie par  
 $$\theta [H , \chi] \,=\, [ \theta (H) , \chi \,{\theta}^{-1}]  \tag (2.6)$$
Si on consid\`ere l'action de $G$ tordue au moyen de
$\theta$, alors  $\theta (\xi)$ n'est
autre que la donn\'ee de ramification de cette action.

Notons pour finir que le groupe $ \Gamma = \Gal_{{\overline {\Bbb
Q}/{\Bbb Q}}}$ agit lui aussi de
mani\`ere naturelle sur le groupe  $R_+(G)$ selon la r\`egle
$\sigma [H , \chi]\, = \,[H , \sigma . \chi]. $
On peut d'une autre mani\`ere interpr\'eter la donn\'ee $\sigma(\xi)$
comme celle associ\'ee \`a l'action de $G$ sur la courbe $C^{\sigma} =
C \otimes_k k$ tordue au moyen de $\sigma$. Il peut \^etre utile de
signaler une autre mani\`ere de d\'ecrire une donn\'ee de
ramification, par exemple  voir (Fried \cite {32}, Serre \cite {65}).
Supposons avoir fix\'e une racine de l'unit\'e $\zeta \in k$
d'ordre $ n=\Card G$, 
$\zeta \in k$.  Si  $\chi$ est un caract\`ere primitif du sous
groupe cyclique  $H$, et si $\Card H = e$, il lui
correspond un unique g\'en\'erateur $s$ de $H$ tel que $\chi(s) =
\zeta^{n/e}$. Si $s'$ correspond
de la m\^eme mani\`ere \`a $(H' ,\chi')$, alors  $(H , \chi) \, \cong
\,(H' , \chi')$  \'equivaut \`a
$s'$ est conjugu\'e \`a $s$. De cette mani\`ere la donn\'ee de
Hurwitz devient une somme formelle de classes de conjugaison
$\sum_{[[s]]} b_{[[s]]} [[s]]$.  On notera que le stabilisateur de
$[H,\chi] \leftrightarrow [[s]]$, pour l'action adjointe de $G$,est
le commutant $C(s) = C_G (H)$ de $s$ (de $H$) dans $G$.

On va \'etendre le  formalisme du $G$-type aux diviseurs
$G$-invariants. Soit $D = \sum_{i=1}^d x_i$ un diviseur invariant par
$G$, tel que si  $i \ne j,\,\, x_i \ne x_j$  . Alors $D$ est une somme
d'orbites deux \`a deux distinctes , et si $[H_{\alpha} ,
\chi_{\alpha}]$ sont les $G$ - types  de ces
orbites, on appellera l'\'el\'ement $\sum_\alpha [H_\alpha ,
\chi_\alpha]$ le $G$-type  de $D$.
On prolonge la d\'efinition en disant qu'une orbite r\'eguli\`ere est
de type $[1]$ (holonomie locale triviale).

\proclaim D\'efinition 2.4.   Soit une action de $G$ sur la courbe
$C$ (lisse ou nodale). On dit
que l'action est marqu\'ee par $b+r$ orbites, dont $r$  r\'eguli\`eres, si on a fait le choix de $b+r$ orbites, contenant les $b$
orbites singuli\`eres,  augment\'ees de  $r$ orbites r\'eguli\`eres,
les orbites \'etant index\'ees ${\Cal O}_1,\cdots , {\Cal
O}_{b+r}$. 

  Le cas non index\'e peut \^etre consid\'er\'e de la m\^eme
mani\`ere. Il est commode pour la suite d'\'etendre  aux  courbes
marqu\'ees  dans le sens pr\'ec\'edent, la donn\'ee de Hurwitz qui
est attach\'ee \`a une action de $G$, en ajoutant \`a la donn\'ee
initialement d\'efinie le terme $r[1].$  On voit alors facilement que
l'\'enonc\'e (2.4)  reste valable, l'indice $1$ (c'est \`a
dire la restriction $H \ne 1$) \'etant maintenant supprim\'e.

\bigskip

\beginsection 2.3. Nombre de Nielsen 

\bigskip

Dans cette section $k = \Bbb C$. Une courbe est   identifi\'ee \`a
 une surface de Riemann compacte et connexe.

\bigskip

 {\it 2.3.1. Classification topologique} 
\bigskip

On s'int\'eresse \`a la classification \'equivariante des actions de
groupes finis sur les courbes ($G$-rev\^etements), mais maintenant
sous le seul aspect topologique.   On
fixe    une surface  $C^\infty$  compacte  orient\'ee
$C$, de genre $g \geq  2$,  munie d'une action du groupe fini $G$; on
suppose que l'action de $G$ respecte l'orientation. Nous pouvons
reformuler dans le cadre topologique ($C^\infty$) l'essentiel de
la section 2.3. La premi\`ere observation  est que la donn\'ee de
Hurwitz, d\'efinie alg\'ebriquement,  est en fait un invariant
topologique de
l'action de $G$ sur la surface $C$. En effet, soit $x \in C$ un point
fixe de stabilisateur $G_X = H$;
alors il en d\'ecoule une repr\'esentation lin\'eaire r\'eelle de
degr\'e deux dans l'espace tangent r\'eel $T_x$ de $C$ en $x$. La
repr\'esentation complexifi\'ee $T_x \otimes \Bbb C$ est alors somme
directe de deux
repr\'esentations conjugu\'ees de degr\'e un. On choisit celle
parmi les deux dont la direction $\xi$ est telle que la $\Bbb R$-base
$(Re(\xi) , Im(\xi)$ fournit l'orientation de $T_x$. Alors le
caract\`ere local $\chi_x$ de la section 2.2 est retrouv\'e comme
celui qui d\'ecrit l'action de $H$ sur la droite  de direction $\xi$,
soit  pour  $s \in H, \,\, s\xi = \chi_x(s) \xi$.

La d\'efinition de l'\'equivalence de deux $G$-rev\^etements (D\'efinition 2.1)  prend maintenant la forme suivante. Soit deux surfaces
$C$ et $C'$ munies d'une action (fid\`ele) de $G$. 

Les surfaces sont
  $G$-{\it  \'equivalentes\/}  (en bref \'equivalentes), s'il
existe un diff\'eomorphisme $G$-\'equivariant et pr\'eservant
l'orientation $f : C \buildrel\sim\over\rightarrow C'$. La encore
cette relation ne doit pas \^etre confondue avec
la relation d'\'equivalence des rev\^etements, dans la situation pr\'esente, le rev\^etement quotient $C \to C/G$. Comme dans le cadre
alg\'ebrique, la relation d'\'equivalence coincide avec l'\'equivalence   au sens des $G$-rev\^etements. Ceci \'etant, il est de
nouveau clair que deux actions \'equivalentes  de $G$ sur des
surfaces $C$ et $C'$ ont des donn\'ees de Hurwitz identiques. Il est
aussi imm\'ediat que la donn\'ee de Hurwitz \'etant fix\'ee, ainsi que
le genre de $C$, alors le genre de $D =   C/G$ est  de fait d\'etermin\'e; c'est une cons\'equence directe de la formule de Riemann - Hurwitz.

Notre objectif maintenant est d'identifier les classes d'\'equivalences , en d'autres termes les types
topologiques, d'actions de $G$, groupe fix\'e, sur les surfaces de
genre $g \,(g \geq  2$, fix\'e); on fixe  en plus de $g$, le genre
$g'$ de la surface quotient $C/G$. La proc\'edure  pour conduire
cette classification est la suivante.

On fixe tout d'abord une surface $B$, et $\beta\subset B$ un ensemble
fini de points (distincts). On suppose en outre avoir fix\'e une
partition
$\beta =\beta_1\sqcup\cdots\sqcup\beta_t$ de $\beta$ en t - parties.
On pose  
$$Diff^+(B,\beta) = \{ h\in Diff^+(B), \, h(\beta_i) =
\beta_i\,\,(i=1,\cdots,t)\} $$
 il est connu que le groupe $Diff^+(B)$ agit transitivement sur
les parties finies num\'erot\'ees de cardinal fix\'e, le choix de
$\beta$ est   donc sans  importance. On suppose enfin avoir fix\'e une
donn\'ee de Hurwitz  $\xi = \sum_{i=1}^t  b_i [H_i,\chi_i] $,
les classes  $[H_i , \chi_i]$ \'etant  deux \`a deux distinctes. On
fait alors
l'hypoth\`ese que $\Card {\beta_i} = b_i$. Fixons enfin un point de
base $\star \in U = B - \beta$, et posons $\pi = \pi_1(U,\star)$.

\proclaim D\'efinition  2.5.
Un marquage d'une action de $G$ sur $C$ consiste en le choix d'un
morphisme $C^\infty$,  %
$\phi : C \to B$ qui identifie $B$ avec la surface quotient $  C/G$,
et qui en plus est tel qu'une orbite singuli\`ere est du type $[H_i ,
\chi_i]$ si et seulement si c'est une fibre $\phi^{-1}(b)$ avec $b\in
\beta_i$. 

Une \'equivalence entre deux actions marqu\'ees $(C , \phi), (C'
,\phi')$ est  par une  d\'efinition identique \`a (2.1) donn\'ee par  un
diff\'eomorphisme $G$-\'equivariant $f: C\buildrel\sim\over\rightarrow
C'$ tel que $\phi'\circ f = \phi$ (\'equivalence stricte).
Il est clair que $Diff^+(B,\beta)$ agit transitivement sur l'ensemble
des classes de marquages de l'action donn\'ee de $G$ sur $C$, mais
pas simplement transitivement en g\'en\'eral, du fait de l'existence
possible d'automorphismes de la base $B$. On associe comme d'habitude
\`a une action marqu\'ee par $\phi : C \to B$ un morphisme de
monodromie de $\pi$ dans le groupe des permutations de la fibre   
$\phi^{-1}(\star)$, qui correspond \`a une action \`a droite de
$\pi$ sur $\phi^{-1}(\star)$. Supposons les points de cette fibre
num\'erot\'es
$$\phi^{-1}(\star) = (\star_1, \cdots,\star_n)$$
C'est un fait bien connu que le groupe de monodromie, c'est \`a
dire l'image de $\pi \to \Bbb S_{\phi^{-1}(\star)},$ s'identifie gr\^ace \`a la num\'erotation
des points de la fibre avec le commutant de $G$ dans le groupe des
permutations $\Bbb S_n = \Bbb S_{\phi^{-1}(\star)}$. De
cette construction d\'ecoule un morphisme surjectif
$\psi_\phi : \pi \to G  $
qui \`a conjugaison pr\`es ne d\'epend pas de la num\'erotation des
points de $\phi^{-1}(\star)$. Rappelons la d\'efinition explicite de
ce morphisme; on choisit  $\star_1 \in \phi^{-1}(\star)$. L'image par
$\psi_\phi$ d'un lacet $\alpha$ bas\'e en
$\star$ s'obtient en relevant dans $\phi^{-1}(U), \,\alpha$ en
$\tilde\alpha$ d'origine $\star_1$.
Alors  
$$\psi_\phi (\langle\alpha\rangle) = \sigma \in
G\,\,\Longleftrightarrow \tilde\alpha(1) =
\sigma(\star_1) \tag (2.7)$$
Si on fait le choix d' un autre point de base $\star'\in U$, et si
$l$ est un chemin de $U$ qui joint
$\star$ \`a $\star'$, alors l'isomorphisme $l_\star :
\pi_1(U,\star) \to \pi_1(U,\star')$ conduit
\`a  la relation \`a conjugaison pr\`es
$\psi_\phi \,=\,\psi'_\phi \circ  l_\star.$

Dans la suite la surjection caract\'eristique $\psi_\phi$ sera
 comprise \`a conjugaison pr\`es dans $G$, et  modulo
un changement \'eventuel du point de base. Le fait suivant
qui fait abstraction de la donn\'ee de Hurwitz est  essentiellement
le  contenu du r\'esultat \'el\'ementaire et bien connu sur la
classification topologique des rev\^etements \footnote {Deux
rev\^etements \'etales de degr\'e $n$,  $\pi: C\to D$ et $\pi': C' \to
D$, de base fix\'ee $D$,  sont  topologiquement \'equivalents si et
seulement si les actions de monodromie, $\phi, \phi':
\pi_1(D,\star)\to \Bbb S_n$ sont conjugu\'ees (voir par exemple
\cite {35}, Proposition 1.2). Le groupe $\Aut (\pi)$ est le commutant
dans $\Bbb S_n$ de l'action de monodromie.}:

\proclaim Lemme 2.6.  la correspondance $(C ,\phi) \mapsto
\psi_\phi$ \'etablit une bijection
entre d'une part les classes d'actions marqu\'ees de $G$ sur $C$, et d'autre part  les classes
de conjugaison de surjections
$\pi \to G$. \hfill\break
\qed 

\bigskip

 {\it 2.3.2. Nombre de Nielsen} 
\bigskip

Pour poursuivre la classification topologique,  fixons une
pr\'esentation   de $\pi =
\pi_1(U)$ (pr\'esentation canonique):
$$\pi = \langle A_1,\cdots,A_g ; B_1,\cdots,B_g ;
\gamma_1,\cdots,\gamma_b \quad/
\quad\prod_{j=1\cdots,g'}[A_j,B_j]\gamma_1\cdots\gamma_b = 1 \rangle$$
 Le morphisme $\psi_\phi$ est d\'etermin\'e
par les images $a_j,b_j\,\,
\hbox{et}\,\,\sigma_k$ des g\'en\'erateurs $A_j,B_j,\,\hbox{et}\,
\gamma_k$. Pour abr\'eger nous noterons $(\underline a , \underline
b , \underline \sigma )$ cette donn\'ee; c'est essentiellement la
{\it description du rev\^etement\/} au sens de Fried \cite {32}. On a
donc
$$G = \langle \underline a , \underline b , \underline \sigma
  \quad / \quad \prod_{j=1\cdots,g'}[a_j,b_j]\sigma_1\cdots\sigma_b =
1 \rangle  \tag (2.8)$$

 Le lien entre la donn\'ee du $(2g' + b)$-uple  $(\underline a ,
\underline b , \underline \sigma )$
et la donn\'ee de Hurwitz $\xi$ est le suivant. Soient
$C_1,\cdots,C_r$ les classes de conjugaison
distinctes d\'efinies par $\sigma_1,\cdots,\sigma_b$; on suppose que
la multiplicit\'e de $C_i$ dans cette liste est $b_i$.  Alors la
donn\'ee de Hurwitz dans son interpr\'etation au moyen de classes de
conjugaison (\S  2.2 ) est  d\'ecrite comme \'etant
$\xi = \sum_{i=1}^r b_i [C_i]. $
D'une autre mani\`ere, on peut voir la donn\'ee des $r$ classes
distinctes comme un ensemble de $r$ couleurs; le lacet $\sigma_j =
\psi_\phi (\gamma_j)$ \'etant de couleur $i$ ssi $\sigma_j \in C_i$.
En cons\'equence le point de branchement $Q_j$ appartient \`a $\beta_j$
ssi $\sigma_j = \psi_\phi
(\gamma_j)$ est de couleur $i$. Dans la suite nous noterons
$\Hom_\beta(\pi , G)$ l'ensemble des
surjections $\psi: \pi \to G$ qui satisfont aux conditions ci-dessus,
dict\'ees par la donn\'ee de
ramification, donc en r\'esum\'e:
$$Q_j \in \beta_i \Longleftrightarrow \psi(\gamma_j) \in C_i  \tag (2.9)$$
  La classification topologique des actions de $G$ sur les surfaces de
genre $g$, \`a donn\'ee
fix\'ee, s'obtient \`a ce stade par la construction classique
suivante. Soit une action de $G$ sur $C$; on choisit  un marquage
(D\'efinition 2.4) $\phi: C \to B$, ce qui produit une surjection
$\psi_\phi \in
\Hom_\beta(\pi , G)$, d\'efinie \`a conjugaison pr\`es (Lemme 2.6).
Soit maintenant le mapping class group $\Bbb M _{g',b}$, qui est,
rappelons-le,  le sous groupe de $\Out(\pi)$,  compos\'e des
automorphismes ext\'erieurs qui permutent les classes de conjugaison
des lacets $\langle \gamma_j \rangle\,\,(j=1,\cdots,b)$. Il y a un
morphisme naturel $$Diff^+(B,\beta) \to \Out(\pi)$$ qui factorise
par $\Bbb M _{g',b}$. Un r\'esultat fondamental de Nielsen d\'ecrit
l'image de ce morphisme, qui
s'identifie \`a l'ensemble des automorphismes de $\pi$ qui
permutent les classes de conjugaison des $\langle \gamma_j \rangle$ \`a
l'int\'erieur de chaque partie $\beta_i\, (i=1,\cdots,r)$. Si on
consid\`ere la surjection canonique $
\Bbb M_{g',b} \to\Bbb S_b$, l'image ci-dessus est le groupe not\'e
 dans la suite $\Bbb
M_{g',(b_1,\cdots,b_r)}$ image r\'eciproque du sous-groupe $\Bbb
S_{b_1}\times \cdots\times
\Bbb S _{b_r}$. Il en r\'esulte donc une action de $\Bbb
M_{g',(b_1,\cdots,b_r)}$ sur
$\Hom_\beta(\pi , G)/G$.  Le r\'esultat  suivant est bien connu: 

\proclaim Proposition 2.7.
   Fixons une donn\'ee de Hurwitz. Les classes topologiques d'actions de $G$
sur les surfaces de genre $g$ \`a donn\'ee fix\'ee, sont en
correspondance bijective avec les ''classes doubles''
$$ \Bbb M_{g',(b_1,\cdots,b_r)} \backslash \left(\Hom_\beta(\pi , G) / G \right)\tag (2.10)$$
 
 \dem La preuve est claire.  Il suffit d'\'eliminer le marquage,
ce qui compte tenu du lemme 2.6 revient \`a effectuer le
quotient de $\Hom_\beta(\pi , G)/G$ par l'action d\'ecrite ci-dessus de
  $\Bbb
M_{g',(b_1,\cdots,b_r)}$. \qed  

\exa {2.1}
 Examinons ce que donne la correspondance ci-dessus dans le
cas $g' = 0$  (rev\^etements de la sph\`ere). On a dans ce cas $B = S^2$, et
le mapping class group est ici le groupe $\Bbb M_{0,b}$ image du
groupe des tresses \`a $b$ brins $B(b)$. Le groupe $\pi$ est dans ce
cas le groupe libre de rang $b-1$  d\'ecrit  par la pr\'esentation
usuelle
$$\langle \gamma_1,\cdots,\gamma_b , \,\,/ \,\, \prod \gamma _i = 1
\rangle \tag (2.11)$$
Une surjection $\psi \in \Hom_\beta(\pi , G)$ est d\'etermin\'ee par
le $r$ - uple $\underline\sigma =
(\sigma_1, \cdots,\sigma_b)$ avec $\sigma_i =
\psi(\gamma_i)\,\,(1\leq i\leq b),$ ces g\'en\'erateurs \'etant
soumis \`a  la relation  $\prod\sigma_i = 1$. Noter que $G$ est
engendr\'e par les $\sigma_i \,\,(i=1,\cdots,b)$. L'action par
conjugaison sur les surjections devient dans cette description la conjugaison diagonale
$$\underline \sigma \mapsto g\underline\sigma g^{-1} =
(g\sigma_1g^{-1},\cdots,g\sigma_bg^{-1})
\tag (2.12)$$
Si on fixe une donn\'ee de Hurwitz, la restriction qui en d\'ecoule,
c'est \`a dire, la  couleur  assign\'ee aux  points de branchement
est: $Q_j \in \beta_i
\Longleftrightarrow \sigma_j\in C_i$. L'action de $\Bbb M_{0,b}$ sur
$ \Hom_\beta (\pi ,
G) /{G}$ se d\'eduit de l'action standard du groupe des tresses de la
sph\`ere,   elle
est en particulier engendr\'ee par les op\'erations
$$S_i(\underline \sigma) = (\sigma_1,\cdots,\sigma_{i-1}, \sigma_i
\sigma_{i+1} \sigma_i^{-1}, \sigma_i, \sigma_{i+2}, \cdots,\sigma_b) \,\, (1\leq i\leq b) \tag (2.13)$$
La condition (2.9) exige  de se  limiter  \`a l'action du sous-groupe des
tresses partiellement color\'ees
$B(b_1,\cdots,b_r)$, image r\'eciproque de $\Bbb S_{b_1} \times
\cdots \times \Bbb S _{b_r}$ par la
surjection canonique $B(b) \to \Bbb S_b$.  $\lozenge$ 
\smallskip

Dans le cas g\'en\'eral, on est conduit \`a consid\'erer l'ensemble
$\Ni(\xi)$ des classes de Nielsen
relatives \`a une donn\'ee de Hurwitz $\xi$ fix\'ee. La d\'efinition
est la suivante \cite {16},\cite {32}:
$$\Ni(\xi) = \{(\underline a , \underline b , \underline \sigma)\} \in
G^{2g'+b} / \hbox{ modulo
conjugaison}\}$$
les \'el\'ements $(\underline a , \underline b , \underline \sigma)$ \'etant soumis aux conditions explicit\'ees ci-dessus, \`a savoir:\hfill\break
 \indent i) $\prod_{i=1}^{g'} [a_i , b_i] \sigma_1\cdots\sigma_b = 1$\hfill\break
\indent ii) $G = \langle\underline a , \underline b , \underline
\sigma \rangle$\hfill\break
\indent iii) $Q_j \in \beta_i \Longleftrightarrow  \sigma_j \in
C_i $ 

En conclusion l'ensemble quotient  $\Bbb M_{g',(b_1,\cdots,b_r)}\,
\backslash \Ni(\xi) $
est le classifiant pour les types topologiques d'actions de $G$ sur
les surfaces de genre $g$, avec
une donn\'ee de Hurwitz fix\'ee $\xi$; c'est un ensemble fini, ce qui
conduit \`a poser \cite {25}:
\proclaim D\'efinition 2.8.  Fixons une donn\'ee de Hurwitz $\xi$.
Le nombre de Nielsen  $h(\xi)$ est d\'efini comme \'etant le nombre de types topologiques  distincts
d'actions de $G$ sur les surfaces de genres $g$, avec une donn\'ee de Hurwitz
fix\'ee \'egale \`a $\xi$, soit
$$ h(\xi) = \Card {\Bbb M_{g',(b_1,\cdots,b_r)} \backslash
\Hom_\beta(\pi , G) /
G} \tag (2.14)$$ 

Le fait le plus important pour la suite est que le nombre de Nielsen\footnote {Ce nombre ne doit pas \^etre confondu avec le nombre de Hurwitz  \cite {59}, qui est    $h_{g,G,\xi} = \#
\Hom_{\beta}(\pi,G)/G$. La formule de sommation de Burnside  donne $h_{g,G,\xi} = \sum_{[\pi]}  {1\over {\#\Aut (\pi)}}$, la sommation
portant sur les classes d'\'equivalence (strictes) de $G$-rev\^etements
de donn\'ee de ramification fix\'ee, d'une courbe de genre
$g'$ fix\'ee, ainsi que les points de branchement $\beta$.} 
a une signification g\'eom\'etrique pr\'ecise, il  repr\'esente le
nombre de composantes connexes, ou irr\'eductibles, de l'espace des
modules des  rev\^etements, espace qui sera construit dans les
sections 5 et 6.  Pour le prouver, il faut resituer la classification
des rev\^etements dans son
contexte  analytique,  la th\'eorie de Teichm\"{u}ller \cite {24},
\cite {25}. La discussion qui suit est succinte, pour plus de d\'etails nous renvoyons aux  r\'ef\'erences pr\'ec\'edentes.

On fixe une surface de r\'ef\'erence $\Sigma_g$ de genre $g\,\,(g
\geq  2)$, ainsi qu'une action fid\`ele du groupe fini $G$ sur
$\Sigma_g$. Le couple $(\Sigma_g , G)$ d\'efinit donc ce qu'on peut
appeler pour abr\'eger, un type topologique.
Soit $ Conf(\Sigma_g)$ l'ensemble des structures conformes sur
$\Sigma_g$. le groupe
$Diff^+(\Sigma_g)$ agit naturellement sur $Conf(\Sigma_g)$, et si
$Diff_\circ (\Sigma_g)$  d\'esigne
la composante  connexe de l'\'el\'ement neutre du groupe
$Diff^+(\Sigma_g)$, donc le sous-groupe des diff\'eomorphismes
homotopes (ou isotopes) \`a l'identit\'e, le quotient  $$T_g =
{Conf(\Sigma_g)}/{Diff_0 (\Sigma_g)}$$  est une des formes de
l'espace de Teichm\"{u}ller \cite {24}. En
suivant la construction de Earle (loc.cit), on est amen\'e \`a
consid\`erer l'ensemble des points fixes  de $G$ dans
$Conf(\Sigma_g)$. Soit $Z_0 (G)$ le centralisateur de $G$ dans
$Diff_0 (\Sigma_g)$. On peut montrer que, $N(G)$ \'etant le
normalisateur de $G$ dans $Diff^+(\Sigma_g)$, on a
$ Conf(\Sigma_g)^G = N (G) \cap Diff_0 (\Sigma_g).$ 
 L'espace de Teichm\"{u}ller \'equivariant, c'est \`a dire
celui qui classifie les actions \`a type topologique fix\'e est
$ T_{g,G} =  {Conf(\Sigma_g)^G} /{Z_0 (G)}$

On montre sous ces conditions le r\'esultat  de fondamental suivant loc.cit\footnote{Une partie de ce r\'esultat est red\'emontr\'ee dans la note de Natazon:  ''The topological structure of the space of holomorphic morphisms of Riemann surfaces'', Comm. of the Moscow Math. Soc. }:
\proclaim Th\'eor\`eme 2.9.  $T_{g,G}$ est une sous vari\'et\'e
ferm\'ee et contractile de  $T_g$. En fait $T_{g,G}$ est diff\'eomorphe \`a $\Bbb R ^{6g'-6+2b}$. 
Ce r\'esultat donne imm\'ediatement le r\'esultat \'enonc\'e au
dessus, \`a savoir que l'espace  de Teichm\"{u}ller \'equivariant,
classifiant les actions de $G$ sur les surfaces (de Teichm\"{u}ller)
de genre $g$ \`a type topologique fix\'e est connexe;  par suite
l'espace modulaire classifiant les actions \`a type
topologique fix\'e, \'etant un quotient du pr\'ec\'edent,  est donc
aussi connexe. En conclusion, si on fixe seulement  la donn\'ee de
Hurwitz $\xi$, l'espace modulaire a   exactement $h(\xi)$ composantes
connexes.\hfill\break \qed

 \rema {2.2}    
 La d\'efinition 2.8 du nombre de Nielsen $h(\xi)$
montre imm\'ediatement que dans le cas $g' = 0$,  et le groupe $G$ \'etant suppos\'e ab\'elien, alors $h(\xi)$ = 1.     La d\'etermination  du nombre de Nielsen est, dans le cas g\'en\'eral, une
question connue comme difficile,  conduisant \`a \'etudier l'action
d'un groupe de tresses sur des uplets d'\'el\'ements de $G$; \cite
{32},\cite {41} \cite  {67}. On sait  cependant que $h(\xi) = 1$ dans plusieurs
situations utiles \cite {25}, \cite {32}, \cite {35}, \cite {67}. Le
cas le plus simple, et  par ailleurs bien connu, est celui des
rev\^etements simples de degr\'e $n$ de ${\Bbb P}^1$, avec
$r$-points de branchement \cite {35}. Le r\'esultat est \'etendu dans
\cite {41} \`a une base de genre $g'$ arbitraire, sous l'hypoth\`ese le nombre de points de branchement soit tel que $b \geq  2n$. Il
est imm\'ediat que la cl\^oture galoisienne d'un tel rev\^etement a
pour groupe de Galois ${\Bbb S}_n$,  la donn\'ee de ramification \'etant
 $\xi = r (\hbox {\rm classe de conjugaison des transpositions})$.
On peut prouver que le nombre de Nielsen est encore \'egal \`a un,
dans le cas o\`u la ramification est simple sauf au plus en une ou
deux fibres,  base $\Bbb P^1$ (Wajnryb \cite {67}). 

Un autre cas pour
lequel la d\'etermination de nombre de Nielsen est ais\'ee, et par
ailleurs classique, et celui des structures de niveau ab\'elien
$n\geq 3$ \cite {20}. Cela signifie que $G =   \Bbb Z /{n\Bbb
Z})^{2g'}$, et $\xi = \emptyset$. Alors (voir paragraphe 8), $h $ est
dans ce cas  l'indice $ [GL_{2g'} ( \Bbb Z /{n\Bbb Z}) : Sp_{2g'}
( \Bbb Z /{n\Bbb Z})]$.
\hfill\break $\lozenge$
 
Supposons maintenant   $G =  \Bbb Z /{n\Bbb Z}$ et  
$\xi \ne 0$  . La surjection
$\psi : \pi \to G$  est dans ce
cas d\'etermin\'ee par les seules images des lacets d'homologie, du
fait que les classes de conjugaison sont r\'eduites \`a un seul \'el\'ement; donc la donn\'ee de Hurwitz d\'etermine  les \'el\'ements
$\sigma_1, \cdots,\sigma_b$, ceux-ci \'etant soumis \`a la seule
relation $\sigma_1. \cdots\sigma_b = 1$. Ne subsiste alors que la
condition $G = \langle \underline a , \underline b , \underline
\sigma \rangle$, soit
$\pgcd {(a_j , b_k , \sigma _l , n)} = 1. $

Le r\'esultat qui suit est bien connu, il remonte \`a  Nielsen. Un
r\'esultat plus g\'en\'eral sur le type
topologique dans le cas d'un groupe ab\'elien est aussi disponible
(Edmonds \cite {25}, theorem 3.1)
\proclaim Proposition 2.10.  Si $G$ est cyclique, pour toute donn\'ee de Hurwitz $\xi \ne 0$, on a
$h(\xi) = 1$; en d'autres termes il y a un seul type topologique, et
l'espace modulaire correspondant est connexe.  \hfill\break
\qed

Soit par exemple le cas $n=3$, et $r = g+2$. La donn\'ee de ramification s'identifie \`a une partition (non ordonn\'ee) de $[1,g+2]$ en deux parties $\Lambda_1,\, \Lambda_2$, selon que le caract\`ere local est $\sigma \mapsto j$, ou $j^2 \,(j^3=1)$. On a $\vert \Lambda_1\vert + 2\vert\Lambda_2\vert \equiv 0 \pmod 3$, soit 
$\vert \Lambda_1\vert  \equiv \vert\Lambda_2\vert   \pmod 3.$
Le nombre de types topologiques, c'est \`a dire de composantes connexes de l'espace de Hurwitz est donc
$${1\over 2}\sum_{2l\equiv g-1\pmod 3} \, {g+2 \choose l}$$

\section3 { Familles de $G$-courbes lisses}  
\medskip 
Les notations et conventions de la section 2 sont conserv\'ees.  Dans cette section on s'int\'eresse  de
mani\`ere essentielle  \`a la g\'eom\'etrie de l'action de $G$ dans une famille de courbes lisses, puis nodales.

\beginsection 3.1.  G\'eom\'etrie du diviseur de branchement

\bigskip

  {\it 3.1.1. Diviseurs de points fixes} 
\bigskip

La situation g\'en\'erique  est la
suivante: soit $\pi : C \to S$ une famille de courbes, qu'on suppose
tout d'abord lisse, et avec des fibres de genre $g \geq  2$
(restriction superficielle): le morphisme  $\pi$ est donc projectif,
lisse, \`a  fibres connexes.  On supposera   en g\'en\'eral la base $S$ est connexe. Il est
alors bien connu que le foncteur en groupes
$\Aut_{C/S} : T \mapsto \Aut_T (C\times_S T)$
    d\'efini sur la cat\'egorie des $S$-sch\'emas, est
repr\'esentable,  en fait  repr\'esent\'e  par un $S$-sch\'ema
fini et non ramifi\'e  $\Aut_S(C)$ \cite {20}. Fixons la convention
suivante: une {\it action} de $G$ sur $C/S$ est un morphisme injectif
de $G$ dans le groupe des sections de $\Aut_S(C)$. Notons alors que
si $\sigma \ne 1$ est un
$S$-automorphisme de $C$, l'automorphisme $\sigma_s$ induit sur une
quelconque fibre $C_s$, est distinct de l'identit\'e, comme il d\'ecoule de la propri\'et\'e de non ramification. Nous allons tout
d'abord  rassembler quelques propri\'et\'es \'el\'ementaires des
sous-sch\'emas de points fixes. Rappelons
que $\Card  {G}$ est inversible dans ${\Cal O}_S$. Comme les
stabilisateurs des points dans les
fibres g\'eom\'etriques sont cycliques, les sous-sch\'emas de points
fixes non vides sont ceux relatifs
aux sous-groupes cycliques. Notons le fait  facile suivant

\proclaim Proposition 3.1.  Pour tout sous-groupe cyclique $H
\subset G$, le sous-sch\'ema des
points fixes $C^H$ de $H$ est  un diviseur de Cartier relatif,   \'etale sur $S$.

\dem C'est la version relative, pour la propri\'et\'e \'etale,  du
fait bien connu que lorsqu'un groupe r\'eductif op\`ere sur une
vari\'et\'e \/ lisse sur un corps, le sous sch\'ema des points fixes
est lisse. Dit d'une autre mani\`ere, si $P$ est un point fixe de
$H$, d'image $s\in S$,   on peut alors lin\'eariser (formellement)
l'action  de $H$ en $P$, donc identifier l'anneau local compl\'et\'e
 $\hat {\Cal O}_P$ \`a $\hat {\Cal O}_s [[T]]$, l'action d'un g\'en\'erateur $\sigma$ de $H$ \'etant $\sigma (T) = \zeta T$, pour une
racine de l'unit\'e \/ convenable $\zeta$. Alors l'\'equation de
$C^H$ en $P$est $\sigma (T) - T = (\zeta - 1) T = 0$. Le r\'esultat
est donc visible.\hfill\break
\qed 

Supposons un instant que $S = \Spec(k)$, avec $k$ alg\'ebriquement clos. Si $P
\in C$ est un point de ramification, donc qui a un stabilisateur $H
\ne 1$, rappelons (2.2)  qu'on
attache alors \`a $H$ un caract\`ere primitif $\chi_P$ de $H$.

\proclaim Lemme 3.2.  Soit un sous-groupe $H \subset G$   
stabilisateur d'un point d'une fibre
g\'eom\'etrique $C_{s_0}$, il est  alors le stabilisateur d'un point pour
toute fibre g\'eom\'etrique $C_s$. \hfill\break
 Soit   $\{H_1,\cdots,H_q\}$
la liste des sous-groupes cycliques de $G$ qui sont les
stabilisateurs des points dans une fibre g\'eom\'etrique $C_s$, alors
cette liste est ind\'ependante de $s$, de m\^eme que le nombre de
points de $C_s$ ayant pour stabilisateur  $H_i$.

\dem  Il r\'esulte imm\'ediatement  de la proposition 3.1 que pour un
sous-groupe $H$ de $G$, le fait d'avoir un point fixe dans une fibre
g\'eom\'etrique $C_s$ entra\^{\i}ne que $H$ a un point fixe dans
chaque fibre g\'eom\'etrique. De plus le cardinal de $C_s^H$, qui est
le degr\'e de $C^H$, est constant sur les fibres g\'eom\'etriques.
Soit
$\{H_1,\cdots,H_q\}$  la liste des sous-groupes  $H$, tels que   $C^H
\ne \emptyset$, c'est \`a dire $C_s^H \ne \emptyset$ pour tout
$s$. Posons
$$\Delta_i(s) = \{ x \in C_s\,\hbox{ tel que } \,\, H_i \,\,\hbox{est
le stabilisateur de}\,\, x \}$$
et $r_i(s) = \Card (\Delta_i(s))$. D\'efinissons aussi $r_i^*(s) =
\Card C_s^{H_i} - r_i(s)$,
c'est \`a dire le nombre de points dont le stabilisateur contient
strictement $H_i$. Par la
proposition 3.1, la fonction $r_i + r_i^*$ est constante. Pour
prouver que $r_i$ est constante,
supposons d'abord   $H_i$ maximal dans la famille $\{H_j\}$; dans ce cas la conclusion
est claire. Si tel n'est plus le
cas pour $H_i$, et si la conclusion est admise pour tout sous-groupe
$H_j$ tel que $H_i \subsetneqq H_j$, l'\'egalit\'e \/ $r_i^* =
\sum_{j,H_i\subsetneqq H_j} r_j,$ montre que $r_i^*$ est constante,
et donc aussi $r_i$.  \hfill\break\qed 

Fixons de nouveau une fibre g\'eom\'etrique $C_{s_0}$, et soit
$\{P_1,\cdots,P_l\}$ la liste des
points de $C_{s_0}$ qui ont pour stabilisateur $H = H_i$; le nombre
$l$ de ces points est
par le Lemme 3.2 ind\'ependant de la fibre, ce qui donne un sens \`a
 la d\'efinition. On note
$\chi_j\,(1\leq j\leq l)$ le caract\`ere de $H$ attach\'e \`a $P_j$.

\proclaim Lemme 3.3.  L'ensemble $\{\chi_1, \cdots, \chi_l\}$
constitu\'e de $l$ caract\`eres
(distincts ou non) de $H$, c'est \`a dire compt\'es avec
multiplicit\'e, est ind\'ependant de la fibre.

\dem Soient $H_i \subsetneqq H_j$ deux sous-groupes de la liste du
lemme 3.2. Le sch\'ema des points fixes $C^{H_j}$ est un sous-sch\'ema ferm\'e de $C^{H_i}$, et il est \'etale sur $S$, par la proposition
3.1.  L'injection de $C^{H_j}$ dans $C^{H_i}$ est  \'etale finie sur
$S*$, et donc ouverte et ferm\'ee. Si on pose
$$Z_i = C^{H_i} - \bigcup_{H_i\subsetneqq H_j} C^{H_j} \tag (3.1)$$
il est clair que $Z_i$ est fini \'etale sur $S$, de degr\'e relatif
$l$. Consid\'erons sur $Z_i$ le
faisceau inversible ${ \Cal L}_i = \Omega_{C/S}^1 \otimes {\Cal
O}_{Z_i};\,\,{ \Cal L}_i$  est un
$({\Cal O}_{Z_i} , H_i)$ module et $\pi_* ({ \Cal L}_i )$ est un
$({\Cal O}_S , H_i)$-
module localement libre de rang $l$.  Il admet   une d\'ecomposition en facteurs isotypiques
$\pi_* ({ \Cal L}_i ) = \bigoplus_{\chi\in \hat {H_i}} E_\chi \otimes 1_\chi, $
$1_\chi$ d\'esignant la repr\'esentation de degr\'e un de caract\`ere $\chi$, le module $E_i$ \'etant
localement libre. Le nombre de fois que $\chi$ appara\^{\i}t dans la
liste  $\{ \chi_i \}$ de
l'\'enonc\'e est \'egal au rang de $E_\chi$; la connexit\'e de
$S$ entra\^\i ne le r\'esultat. \hfill\break\qed 

Il est tr\`es utile dans les manipulations de familles de courbes de
pouvoir disposer de la famille
quotient. Le r\'esultat bien connu, et par ailleurs \'el\'ementaire
du fait de l'hypoth\`ese de r\'eductivit\'e,  suivant r\'esume les
propri\'et\'es essentielles de cette op\'eration:

\proclaim Proposition 3.4.  Soit $\pi: C \to S$ une courbe propre et
lisse de base $S$. Soit une
action de  $G$ sur $C/S$. Alors la courbe quotient  $B =  C/G $
existe, et le morphisme $\pi$
factorise en  $ \pi: C \to B \to S,\,\, \hbox{ le morphisme } \,\,f:
C \to B\,\,\hbox {\'etant fini plat de rang}
\,\,n$;  de plus la formation de $B$ commute aux changements de
bases. En particulier, au
niveau des fibres  on a $ B_s =  {C_s}/G \,\, (s\in S)$. \hfill\break
\qed 

Une autre construction  importante pour la suite  est celle du diviseur de ramification (resp. le
diviseur de branchement). Fixons comme ci-dessus une action de $G$
sur la courbe $C/S$. Dans la
liste du Lemme 3.2, des sous-groupes de $G$ qui sont des
stabilisateurs de points dans les fibres
g\'eom\'etriques, choisissons dans chaque classe de conjugaison un \'el\'ement, et notons de nouveau $\{ H_1,\cdots,H_q\}$ cette liste.
Formons pour chaque indice $i$ le sous sch\'ema $Z_i$ (Lemme 3.3)
form\'e des points de stabilisateur exactement $H_i$. Rappelons que
$Z_i$ est fini et \'e tale sur $S$. Soit $N_i = N_G(H_i)$, le
normalisateur de $H_i$ dans $G$, et $\overline N_i =   {N_i}/{H_i}$. Le groupe $\overline N_i$ op\`ere   librement  sur $Z_i$,
de sorte que    $D_i =   {Z_i}/{\overline N_i}$ est
  fini \'etale sur $S$.   Alors $D_i \subset B$ est un diviseur
de cartier relatif de $B$ dont les points ferm\'es correspondent aux
orbites singuli\`eres de points de $C$ avec isotropie $H_i$. Cette
remarque \'el\'ementaire sera  souvent utilis\'ee.  En
conclusion,  on observe que la donn\'ee de ramification est
invariante par d\'eformation.

\bigskip

  {\it 3.1.2. Ramification et branchement} 
\bigskip

Nous appelerons le diviseur de cartier relatif $D = \sum _{i=1}^q
D_i$ de $B/S$ le
{\sl diviseur de branchement r\'eduit  relatif},   de $C/S$. Le
degr\'e de $D$, qui est \'etale sur $S$,  est  le nombre d'orbites
singuli\`eres contenues dans chaque fibre, donc le nombre de points
de branchement. Rappelons que le diviseur de branchement $B(\pi)$,
qui ne doit pas \^etre confondu avec $D$,  est l'image directe
$B(\pi) = \pi_* (R)$ du diviseur de ramification. Ce dernier \'etant
d\'efini par $R = \Div  (\Omega_{C/D})$, donc  tel que
${\Cal O}(R) = \Omega^1_{C/S} \otimes \pi^*
(\Omega^1_{B/S})^{-1}.$  Notant   plus g\'en\'eralement, et
pour utilisation ult\'erieure,  $F_H = C^H$ le diviseur de Cartier
relatif des points fixes du sous-groupe cyclique $H$, et $\Delta_H$
la composante de $F_H$ form\'ee des points de stabilisateur
exactement $H$, on a l'expression classique  du diviseur
de ramification $R_\pi$:
$$R_\pi = \sum_ {H\ne 1} \, \varphi (\vert
H\vert) F_H \,=\, \sum_{H\ne 1} \, (\vert H\vert - 1) \Delta_H \tag
(3.2)$$
$\varphi $ d\'esignant le fonction d'Euler.
  Il est facile de relier le diviseur de branchement $B(\pi)$  aux
diviseurs $D_i$ d\'efinis au dessus. En effet, partant de
l'expression de dessus  $R = \sum_{H\ne 1} \varphi (\vert H\vert )
C^H$, et prenant l'image par $\pi$ de cette relation dans $D$, on
obtient :
$$B(\pi) = \sum_{i=1}^q \left(\vert G\vert -  { {\vert G\vert}\over
{\vert H_i\vert}}\right) D_i \tag (3.3)$$

On peut raffiner la d\'ecomposition de $D$ en invoquant le type d'une
orbite singuli\`ere. La donn\'ee de Hurwitz \'etant $\xi$, et notant
$[H_1,\chi_1], \ldots, [H_k,\chi_k]$ les classes deux \`a deux
distinctes qui apparaissent dans $\xi$,  soit
$\xi = \sum _{i=1}^k b_i [H_i ,\chi_i] , \,\, (b_i > 0)$, 
 on a   la d\'ecomposition suivante:
 
\proclaim Proposition 3.5.  Le diviseur de branchement  $B(\pi)$, ou
de branchement r\'eduit $D(\pi)$, se d\'ecompose en une somme
$D(\pi) = \sum_{i=1}^k T_\alpha$ de diviseurs deux \`a deux
disjoints, et $\deg(T_\alpha) = b_\alpha$. Les points de  $D_\alpha$
sont les images des fibres singuli\`eres de type $[H_\alpha ,
\chi_\alpha]$.

\dem Fixons l'un des sous-groupes $H\in\{ H_i\}$, et consid\'erons
les diff\'erents caract\`eres
primitifs $\chi_1,\cdots,\chi_m$ tels que la classe de $(H , \chi_i)$
appara\^\i sse dans la donn\'ee
$\xi$. Notons que  $(H' , \chi') \equiv (H , \chi)$ \'equivaut (voir
2.2)  \`a l'existence de $t\in
N = N_G(H)$ tel que $\chi'(s) = \chi(t^{-1}st),\,\, \hbox{pour
tout}\,\, s\in H$. Par
r\'eciprocit\'e de Frobenius cela est aussi \'equivalent \`a
$ \Ind_H^N (\chi) \,\,\cong \,\,\Ind_H^N (\chi')$.
La repr\'esentation\quad $\Ind_H^n (\chi)$\quad d\'esigne la repr\'esentation obtenue par induction de $H$ \`a $N$ de la repr\'esentation
de degr\'e un de caract\`ere $\chi$. Soit $Z$ le sous sch\'ema  \'etale sur $S$ associ\'e \`a $H$ comme dans le Lemme 3.3. Le faisceau
conormal $\Cal N$ de $Z$ dans $C$ coincide avec  
$\Omega_{C/S}^1 \otimes {\Cal O}_Z$. Il est inversible et muni d'une
action de $\overline N =
  N/{H}$. L'image directe $\pi_\star ({\Cal N})$ est  alors un
faisceau muni d'une action de
$\overline N$. On peut ainsi le d\'ecomposer en facteurs isotypiques
$$  \pi_* ({\Cal N}) = {\bigoplus _{v\in  \Irrep(\overline N)}}\,\,
{\Cal L}_v \otimes v $$
o\`u dans cette somme directe, $v$ parcourt les repr\'esentations
irr\'eductibles de  $\overline N$, et
${\Cal L}_v$ est un module localement libre. Sur les composantes
connexes de $\pi(Z)$, qui sont
disjointes deux \`a deux, le rang de chacun de ces faisceaux est
constant et de cela d\'ecoule  que
la repr\'esentation $ \pi_* ({\Cal N})_y$  est constante le long des
fibres $y \in \pi(Z)$. Cela
prouve le r\'esultat. \hfill\break
\qed 

\rema  {3.1}   Soit $[H,\chi]$ une classe de
conjugaison, et notons $\Delta_{(H,\chi)}$ (resp.
$\Delta_{[H,\chi]}$) le diviseur de $C$, lieu des points
fixes"d'holonomie" exactement $(H,\chi)$ (resp. un conjugu\'e).
Soit $C(H,\chi)$ le stabilisateur de $(H,\chi)$  pour l'action
adjointe de $G$. On a $C(H,\chi) = C_G(H)$ le commutant de $H$ dans
$G$. On note  que $C_G(H)\over H$ agit librement sur
$\Delta_{(H,\chi)}$, de quotient la composante $B_{[H,\chi]}$ de
$B(\pi)$ (voir preuve de la Proposition 3.4). Par ailleurs notons la
d\'ecomposition en une somme disjointe $$\Delta_{[h,\chi]} =
\sum_{g\in G/C_G(H)}  \,e\Delta_{(H,\chi)}\tag (3.4)$$
$\lozenge$


\medskip

\beginsection {3.2. Inversion de la formule de Chevalley - Weil }

\bigskip

Soit comme dans la section pr\'ec\'edente, une courbe lisse $C$ d\'efinie
sur le corps $k$, de genre $g \geq 2$, munie d'une action de  $G$ de
donn\'ee de Hurwitz $\xi$. Le faisceau canonique $\omega_C= \Omega_C$
est un $G$-faisceau inversible, en
particulier les espaces vectoriels $H^0 (C , \omega_C^{\otimes m})$
sont naturellement des repr\'esentations de $G$,   appel\'ees dans la
suite  {\sl les repr\'esentations de Hurwitz} associ\'ees \`a la
$G$-courbe $C$, et not\'ees $Hurw_m (C)$, ou $Hurw_m$ s'il n'y a pas
de doute sur $C$.
\bigskip 
  {\it 3.2.1.  Formule de Chevalley-Weil
} \bigskip
Il y a une connexion \'etroite entre la structure de ces repr\'esentations et la donn\'ee de Hurwitz;
pour l'essentiel elle est donn\'ee par le th\'eor\`eme de Chevalley -
Weil, que nous allons rappeler.
Le r\'esultat principal de cette section (Th\'eor\`eme 3.6) est  une
r\'eciproque \`a ce th\'eor\`eme.
  Noter que Morrison et Pinkham ont \'etudi\'e en d\'etail la repr\'esentation $H^0 (C ,\omega_C)$; ils en donnent en particulier une
caract\'erisation lorsque $G$ est cyclique \cite {56}. Rappelons  tout
d'abord la formule de Chevalley-Weil sous la forme donn\'ee par
Ellingsrud et L\"{o}nsted \cite {28} (voir aussi \cite {56} pour le
cas cyclique). Soit $\Cal L$ un  $({\Cal O}, G)$ module inversible;
la trace de Lefschetz de $\Cal L$ est l'\'el\'ement de l'anneau des
repr\'esentations $R(G)$ donn\'e par:
$$L_G({\Cal L}) = \sum_{i=0}^1 (-1)^{i} [H^{i}(C , {\Cal L}] \tag (3.5)$$
o\`u si $v$ est un $kG$-module, $[v]$ d\'esigne sa classe dans
$R(G)$. Lorsque l'action de $G$ est
libre, on a la formule bien connue
$L_G({\Cal L}) = {\chi ({\Cal L}) \over \Card  G} [k[G]] . $

Indiquons bri\`evement comment s'obtient la formule dans le cas g\'en\'eral. Soit $\pi: C\to B$ le quotient par $G$. On consid\`ere la
d\'ecomposition du $({\Cal O}_B ,G)$-faisceau $\pi_* ({\Cal L})$ en
facteurs isotypiques
$$\pi_* ({\Cal L}) = \bigoplus_{v\in \hat G} E_v \otimes v$$
avec  $E_v = \pi_*^G ({\Cal L} \otimes v^\lor)$, faisceau localement
libre de rang $\dim (v)$. La
relation
$L_G({\Cal L}) = \sum _{v\in G^\lor} \chi(E_v) [ v ]  $
montre que pour obtenir la trace de Lefschetz, il suffit d'expliciter
$\chi(E_v)$,  du fait du th\'eor\`eme de Riemann-Roch qui donne  $
\chi(E_v) = \deg (E_v) + \dim (v)(1-g')$. 

Notons $\{y_1, \cdots,
y_r\}$ les points de branchements de $\pi$. Il est clair que le
conoyau ${\Cal Q}({\Cal L},v)$ du morphisme canonique  $\pi^* (E_v)
\to {\Cal L} \otimes v^\lor$ a pour support les points de
branchement. De plus on a visiblement en fixant un point $x_i \in
\pi^{-1}(y_i)$ de stabilisateur $H_i$:
$${\Cal Q}({\Cal L},v)_{y_i} = \bigoplus_{\pi(x)=y_i} {\Cal Q}({\Cal
L},v)_x = \Ind_{H_i}^G ({\Cal
Q}({\Cal L},v)_{x_i})  $$
En particulier, si $\Card  H_i = e_i$ on a $\dim \,{\Cal Q}({\Cal
L},v)_{y_i} = { n \over {e_i}}
\dim \,{\Cal Q}({\Cal L},v)_{x_i}$.
Pour d\'ecrire le $kH_i$-module ${\Cal Q}({\Cal L},v)_{x_i}$ on peut
noter que si $\nu_i \in R(H_i)$
est la repr\'esentation de degr\'e un de $H_i$ dans la fibre ${\Cal
L}(x_i)$, on a $({\Cal L} \otimes
v^\lor)_{x_i} =  {\Cal O}_{x_i} \otimes v^\lor (\nu_i)$. Le caract\`ere $\chi_i$ de la repr\'esentation
cotangente au point fixe $x_i$ de $H_i$ \'etant primitif,  il vient  dans
$R(H_i)$ une d\'ecomposition
$$[v^\lor (\nu_i)\vert_{H_i}] = \sum_{l=1}^{e_i} m_{il}
\chi_i^l$$
Pour expliciter la contribution de  ${\Cal Q}({\Cal
L},v)_{x_i}$, on est ramen\'e  \`a calculer la dimension du
conoyau de $\left (\chi_i^l \otimes {\Cal O}_{x_i}\right )^{H_i}
{\otimes_{\Cal O}}_{y_i}  {\Cal O}_{x_i} \rightarrow \chi_i^l \otimes
{\Cal O}_{x_i}$, qui est $e_i - l$. On en tire imm\'ediatement   $\dim
{\Cal Q}({\Cal L},v)_{x_i} = \sum_{l=1}^{e_i-1} (e_i - l)m_{il}$,
d'o\`u le degr\'e de $E_v$
$$\deg{E_v} = {\dim (v)\, \deg({\Cal L})\over n} - \sum_{i=1}^r
\sum_{l=1}^{e_i-1} {e_i - l\over
e_i} m_{il} \tag (3.6)$$Cette relation  fournit l'essentiel du contenu
num\'erique de la formule de Chevalley-Weil.
\exa{ 3.2} (Rev\^etements simples de $\Bbb P^1$)  
 Pour \'eclairer la relation (3.6), soit  le cas
d'un rev\^etement "simple"
$\pi: C \rightarrow \Bbb P^1$, simple signifiant que le groupe de
Galois est le groupe sym\'etrique $G =
{\bf S}_n$, et la donn\'ee de ramification est $r (12)$, o\`u par
$(12)$ on d\'esigne la classe de
conjugaison des transpositions. Un calcul facile conduit pour ${\Cal L}
= \Omega_C$ \`a
$$\deg (E_v) \;=\; \xi_v (1) ({ r\over 2} - 2) \;-\;{  r \over4} \xi_v ((12))$$
et vu que $\deg (E_v) = \chi_v (1)$, on trouve   
$\chi (E_v) \;=\; ({ r \over 4} - 1) \;\chi_v (1)  - { r\over
4}\;\chi_v ((12))$.
Notons que par les relations d'orthogonalit\'es des caract\`eres, on a
  $\sum_{v \in \Cal S_n} \;\chi_v ((12)) \chi_v (1)  = 0,$
ce qui assure la consistance des relations de dessus avec le genre
donn\'e par la relation de Riemann-Hurwitz $g_C = \dim H^0 (C , \Omega_C) = 1 + n! ({ r\over 4} - 1)$ $\lozenge$
\smallskip
Sp\'ecialisons maintenant la relation (3.6), pour utilisation ult\'erieure, dans le cas $G$ cyclique
(d'ordre $n$) et    ${\Cal L} = \omega_C^{\otimes m}$. Fixons un g\'en\'erateur $\sigma$ de $G$, et soit $\chi $ la repr\'esentation de
degr\'e un de $G$ donn\'ee par $\chi(\sigma) = \epsilon$, o\`u
$\epsilon$ d\'esigne une racine primitive n-i\`eme de l'unit\'e fix\'ee. Le stabilisateur $H_i$ est   alors engendr\'e par $\sigma^{
n \over {e_i}}$, et le caract\`ere local $\chi_i$ est donn\'e par
$\chi(\sigma^{ n\over  {e_i}}) = \epsilon^{nk_i \over e_i}$ pour un
certain entier $k_i$, avec $0 < k_i < e_i$ et $(k_i,e_i) = 1$. Soit
$\eta_i$ d\'efini par
$\eta_ik_i \equiv 1 \pmod  {e_i}$, alors la restriction de $\chi^l$
\`a $H_i$ co\"\i ncide avec
$\chi_i^{l\eta_i}$.

Avec le choix ${\Cal L} = \omega_C^m \;(m \geq  1)$, les repr\'esentations introduites ci-dessus sont
$\nu_i = \chi_i^m$, et si $v = \chi^l$, la multiplicit\'e de
$\chi_i^j$ dans la restriction \`a $H_i$ de
$v^\lor (\nu_i)$  est
$$   m_{ij} = \cases{ 0 & si $\, m -l\eta_i \not\equiv
j\,\pmod { e_i}$ \cr
   1&
 si  $\, m - l\eta_i \equiv j \,\pmod {e_i} $} \tag
(3.7)$$
 Adoptons pour la suite la notation suivante: si $\langle
x\rangle$ d\'esigne la partie fractionnaire de $x\in \Bbb R$,
$\langle\langle  x \rangle\rangle = 1 - \langle -x \rangle$
Soit $g'$ le genre de $B$; alors avec $v = \chi^l$ et en notant dans
ce cas $E_l$ au lieu de $E_v$,
la relation (3.6) se  sp\'ecialise en
$$\chi(E_l) = {m(2g - 2)\over n} + (1 - g') - \sum_{i=1}^r \left( 1 -
\langle\langle  {m -
l\eta_i\over e_i}\rangle\rangle \right)  \tag (3.8)$$
et tenant compte de la formule de Riemann - Roch
$$\chi(E_l) = (2m - 1)(g' - 1) + mr - \sum_{i=1}^r \left(
\langle{l\eta_i - m\over e_i} \rangle +
{m\over e_i} \right) \tag (3.9) $$
On obtient en particulier, si $m = 1$, que la multiplicit\'e
$\mu_l$ de $\chi^l$ dans  $H^0(C ,
\omega_C)$ est
$$ \mu_l = \cases{g'  \quad si \quad $l = 0 $\cr
    g' - 1 + r - \sum_{i=1}^r
\left(\langle {l\eta_i - 1\over e_i}\rangle  + {1\over e_i}\right) &
 si \quad
$l=1,\cdots,n-1$}   \tag (3.10)$$
Les formules (3.9) et (3.10) forment le contenu de ce qu'on appelle,
semble t-il usuellement, les (ou la) formules de  Chevalley-Weil.

\bigskip

{\it 3.2.2.  Inversion des relations de Chevalley-Weil}  

\bigskip

Inverser les relations de
Chevalley-Weil revient \`a montrer que la connaissance des
repr\'esentations de Hurwitz $H^0(C,\omega_C^{\otimes m})$, en fait  d'un nombre fini d'entre
elles, permet de retrouver la donn\'ee de Hurwitz. L'utilit\'e
d'une telle observation s'explique par le fait que les repr\'esentations de
Hurwitz se comportent "bien" dans des familles. De mani\`ere pr\'ecise, le r\'esultat, qui \`a notre
connaissance ne semble pas avoir \'et\'e remarqu\'e, est:

\proclaim Th\'eoreme 3.6. Soit une action du groupe fini $G$ sur la
courbe projective lisse $C/k$.
On suppose l'action mod\'er\'ee, donc  que $n = \Card  G$ est
inversible dans $k$. Les $kG$-modules $H^0(C , \omega_C^{\otimes m})
(m\geq  1)$
d\'eterminent de mani\`ere unique la donn\'ee de Hurwitz de l'action
(un nombre fini d'entre eux
en fait suffisent). 

\dem La preuve consiste \`a faire l'\'etude d'un point de vue
combinatoire de l'action du groupe $G$. Elle d\'ecoule plus pr\'ecis\'ement de l'utilisation it\'er\'ee de la formule d'inversion de
M\"{o}bius appliqu\'ee \`a l'ensemble ordonn\'e  des sous-groupes
cycliques de $G$. On proc\`ede par \'etapes:

\proclaim Lemme 3.7.  Supposons tout d'abord le groupe $G$ cyclique
d'ordre $n$. Soit $N_k$ le nombre d'orbites avec $ n /k$ \'el\'ements. Alors les nombres $N_k \,\,(k>1, k\vert n)$ sont totalement
d\'etermin\'es par les repr\'esentations de Hurwitz. 

\dem Soient $e_1\leq \cdots\leq e_r$ les ordres des stabilisateurs
des   orbites singuli\`eres rang\'es dans l'ordre croissant. Du fait de 3.10, la multiplicit\'e de la repr\'esentation triviale dans 
  $H^0(C , \omega_C^{\otimes m})$ est donn\'ee par    
  $$
\cases {g'&   si $\, m=1$ \cr
 3g' - 3 + r & si $\, m=2$}
 $$
Pour $m\geq 3$, le  r\'esultat est, en tenant compte de $r = N_2 + N_3 \;+\;\cdots$
$$5(g' - 1) + 3r - \sum_{i=1}^r \left(\langle -{3\over e_i}\rangle +
{3\over e_i}\right) \,=\, 5(g' -
1) + 3r - (2N_2+N_3+\cdots )$$
  
$$=\,  5(g' - 1) + 2r - N_2$$
expression  qui visiblement d\'etermine $N_2$. Supposons maintenant
conna\^{\i}tre  inductivement $N_2,\cdots,N_{k-1}$ avec
$k\geq 3$. La multiplicit\'e de la repr\'esentation triviale dans
le $kG$-module $H^0(C ,\omega_C^{k+1})$
est, par la relation (3.9),
$$(2k + 1)(g' - 1) + (k + 1)r - \sum_{i=1}^r\left( \langle -{k+1\over
e_i}\rangle  + {k+1\over
e_i}\right)$$
$$= (2k + 1)(g' - 1) + (k + 1)r - \sum_{2\leq d\leq k} N_d \left(
\langle - {k+1\over d}\rangle + {k+1\over
d}\right) + 2N_k + N_{k+1} + \cdots$$
ce qui permet de d\'eduire la valeur de $N_k$. \hfill\break\qed

  Revenons maintenant au cas g\'en\'eral. Si $H$ est un sous groupe
cyclique de $G$ de cardinal $h$, on notera pour tout diviseur $k$ de
$h$, $N_k(H)$ le nombre de $H$-orbites \`a $ h\over k$ \'el\'ements. Le lemme 3.3 appliqu\'e \`a la restriction de l'action \`a
 $H$, montre que le cardinal $\Card{(C^H)} = {N_h (H)}$ de
l'ensemble des points de $C$ fix\'es par $H$ est d\'etermin\'e par
la
seule connaissance des repr\'esentations de Hurwitz. Notons
$\Delta_H$ l'ensemble des points de $C$ dont le stabilisateur est
exactement $H$. On a alors (3.1) 
$$\Delta_H = C^H - \bigcup_{K,H\subsetneqq K} C^K  \tag (3.11)$$
Une orbite singuli\`ere qui contient un point de stabilisateur $H$
coupe $\Delta_H$ en $\beta =
[N_G(H) : H]$  points. La technique usuelle d'inversion de
M\"{o}bius, conduit \`a:

\proclaim Lemme 3.8.  Consid\'erons des op\'erations du groupe  $G$
sur les courbes lisses $C_1$ et $C_2$. On suppose que pour tout $m
\geq  1$ on a un isomorphisme de $kG$-modules
$H^0(C_1 , \omega_{C_1}^{\otimes m}) \cong H^0(C_2 ,
\omega_{C_2}^{\otimes m})$.
Alors un sous-groupe cyclique $H$ de $G$ est le stabilisateur d'un
point de $C_1$ si et seulement si $H$ est le stabilisateur d'un point
de $C_2$. 

\dem  On utilise implicitement un argument bien connu, souvent
appel\'e   le Lemme de Brauer (\cite {65},
p 67), qui rappelons le, dit que si $X,Y$ sont deux ensembles finis
munis d'une action du groupe fini $G$, alors les repr\'esentations
par permutations de $G$ d\'efinies par $X$ et $Y$ sont \'equivalentes
si et seulement si pour tout sous-groupe cyclique $H$ de $G$ on a
$\Card X^H = \Card Y^H$. Il d\'ecoule tout d'abord du lemme 3.7, que
si   
$H$ est un sous-groupe cyclique de $G$, alors $\Card {C_1^H} = \Card
{C_2^H}$. Posons si $i = 1,2$   
$$ f_i(H) = \Card C_i^H \quad\hbox{et}\quad g_i(H) = \Card {\Delta_i(H)}$$
de sorte que $f_1(H) =f_2(H)$ pour tout sous-groupe cyclique $H$ de
$G$. L'ensemble $\Cal C$ des sous-groupes cycliques de $G$ \'etant
ordonn\'e par l'inclusion, on peut d\'efinir la fonction de
M\"{o}bius correspondante
$\mu_G : {\Cal C} \times {\Cal C} \to \{0, \pm 1\}$. Si $\mu$ repr\'esente la fonction de M\"{o}bius usuelle, on a pour $H,K \in {\Cal
C},\,\,H\subset K, \,\, \mu_G( H , K) = \mu([K : H])$.  De (3.11) on
tire imm\'ediatement l'\'egalit\'e
$$f_i(H) = \sum_{K,H\subset K} \, g_i(K) \tag (3.12) $$
 La relation d'inversion  de M\"{o}bius
conduit  de suite \`a  la relation
$g_i(H) = \sum_{K,H\subset K} \mu_G(H , K) f_i(K).  $
Comme $f_1 = f_2$, on a bien $g_1 = g_2$, ce qui est en substance le
r\'esultat annonc\'e. \hfill\break
 \qed 
\smallskip
  {\sl Preuve du th\'eor\`eme 3.1:} 
 On suppose   avoir pour tout $m\geq 1$ un isomorphisme de $kG$ - modules:
$ H^0(C_1 , \omega_{C_1}^m) \cong H^0(C_2 , \omega_{C_2}^m).  $
Il s'agit de prouver que les donn\'ees de Hurwitz  $\{ \xi_i\}\,\,
(i=1,2)$ attach\'ees \`a $C_1$ et
$C_2$ sont identiques, c'est-\`a-dire que la multiplicit\'e  dans
$\xi_1$ et $\xi_2$ d'une quelconque
classe   $[H , \chi],$  est la m\^eme. La premi\`ere observation, qui
d\'ecoule du lemme 3.2, est qu'un
sous-groupe cyclique $H$ stabilise un point de $C_1$ si et seulement
si il stabilise un point de
$C_2$. Il faut donc maintenant prouver, en travaillant avec $C = C_1$ ou
$C_2$,  que la classe de
conjugaison de 
$H$ \'etant fix\'ee, on peut extraire des repr\'esentations de
Hurwitz  de $C$ la multiplicit\'e d'une
classe  $(H_i , \chi_i),$ avec $H_i$ conjugu\'e \`a $H$. Ces
classes sont en bijection avec les
$N_G(H)$-orbites dans $\Delta_H$, o\`u $N_G(H)$ est comme d\'eja
indiqu\'e le normalisateur de $H$, et rappelons le $\Delta_H =
\{x\in C , G_x = H\}.$ On est ainsi ramen\'e \/ par restriction de
l'action \`a  $H$, au cas $G = H $ cyclique. Ceci \'etant, notons
$$2\leq e_1\leq  \cdots \leq  e_r \leq n = \Card G \tag (3.13) $$
les ordres des stabilisateurs des orbites singuli\`eres, et
$H_1,\cdots,H_r$ ces stabilisateurs.
Extraire la donn\'ee de Hurwitz revient \`a expliciter pour chaque
$H_i$,  l'entier $\nu_i, 1\leq \nu_i<e_i, \,(\nu_i,e_i)= 1$. La
multiplicit\'e du caract\`ere $\chi^l$ dans la repr\'esentation
$Hurw_m$ \'etant connue (noter que les entiers $g'$ et $r$ sont d\'etermin\'es), on peut supposer que l'expression
$$\sum_{i=1}^r \left( \langle{l\nu_i - m\over e_i }\rangle + {m\over
e_i} \right)$$
est d\'etermin\'ee pour tout $m\geq 1$, et tout $l$. En particulier
si $ m=n$, l'expression $\sum_{i=1}^r \langle{l\nu_i\over
e_i}\rangle$ est d\'etermin\'ee pour tout $l$. Ainsi l'expression
$$\sum_{i=1}^r \left( \langle{l\nu_i - m\over e_i }\rangle \,-
\,\langle{l\nu_i\over
e_i}\rangle\,-\,\langle-{m\over e_i}\rangle \right)$$
peut donc \^etre consid\'er\'ee comme  acquise pour tout $m\geq 1$ et
tout $l$. Si $m = an - 1 \,\,(a\geq 1)$, cette derni\`ere expression
vaut
$$\sum_{i=1}^r \left( \langle{l\nu_i + 1\over e_i }\rangle \,-
\,\langle{l\nu_i\over
e_i}\rangle\,-\,\langle-{1\over e_i}\rangle \right) \tag (3.14)$$
mais pour tout $q\in \Bbb Z$, $\langle{q+1\over e}\rangle  - {1\over
e} - {q\over e} = -1$ si
$e \vert q+1,\; 0$ sinon. On obtient ainsi que l'expression (3.14)
compte le nombre d'indices $i\in
[1 , r], $ tels que
$ l\nu_i \equiv -1 \pmod { e_i}.$
Pour terminer la preuve, notons que si $n = \Card  G = p$ est
premier, la conclusion est claire,
car alors $H_1 = \cdots = H_r = G$ , et le nombre des $\nu_i$ tels
que $\nu_i \equiv q \pmod
{p}$ pour tout $q$ premier \`a $p$. Si maintenant $n$ est quelconque,
on peut  raisonner par
induction croissante sur l'ordre des sous-groupes, et supposer que la
donn\'ee relative \`a tout
sous-groupe strict  $H \subsetneqq G$  est d\'etermin\'ee. Si
$\Delta_G \ne \emptyset$, la preuve est
termin\'ee;  dans le cas contraire, il faut pouvoir compter pour tout
$q$ premier \`a $n$, le nombre des $\{\nu_i\}$ \'egaux \`a $q$
pour lesquels $e_i = n$. Supposons que la liste (3.14) soit telle que
$$e_1\leq \cdots \leq e_t<e_{t+1} = \cdots =e_r = n\quad(0\leq t<r)$$
Alors connaissant les entiers $\nu_1,\cdots,\nu_t,$ et pout tout $l$
premier \`a $n$ le nombre
d'indices   $i\in  [1 , r] $ tels que $l\nu_i \equiv -1 \pmod {e_i}$,
on peut finalement obtenir la
totalit\'e de la donn\'ee de Hurwitz, ce qui d\'emontre le th\'eor\`eme.\hfill\break
\qed 


Il est maintenant clair que si $C/S$ est une courbe (lisse) sur une
base connexe, munie d'une
action de $G$, alors la donn\'ee de Hurwitz de l'action de $G$ le
long des fibres g\'eom\'etriques est
constante. Cela d\'ecoule imm\'ediatement du th\'eor\`eme 3.1 joint
au r\'esultat  suivant, qui est une
extension au cas \'equivariant du r\'esultat classique qui \'etablit
que sous les hypoth\`eses du
pr\'esent \'enonc\'e, la caract\'eristique d'Euler-Poincar\'e
$\chi({\Cal L}_s)$ est constante. De mani\`ere pr\'ecise:
\proclaim Proposition 3.9.  Sous les hypoth\`eses ci-dessus, soit
${\Cal L}$ un $({\Cal O}_C , G)$
faisceau localement libre de rang fini. La trace de Lefschetz
$L_G({\Cal L}_s)$ est constante de
long des fibres g\'eom\'etriques. En particulier  la repr\'esentation  $ 
H^0(C_s,(\omega_{C/S})_s)   $ est ind\'ependante de
$s\in S$.\hfill\break
\qed

\section4{Familles de $G$-courbes stables}
\bigskip
Dans cette  section   on \'etend   les r\'esultats
de la section 3 aux familles de courbes pr\'estables.  On introduit
    la stabilit\'e   de l'action d'un groupe fini $G$, sur
une courbe pr\'estable, et  on propose   la d\'efinition
des rev\^etements galoisiens stables, respectivement  stables marqu\'es.

 \beginsection 4.1.  Actions stables, rev\^etements stables  

\bigskip

{\it 4.1.1. Courbes stables et  stables
marqu\'ees}  
\bigskip
 Les d\'efinitions g\'en\'erales   en ce qui concerne les courbes pr\'estables et stables marqu\'ees  sont dans \cite {40}, \cite {53},   ou dans les  articles de Deligne-Mumford \cite {20}, et
Knudsen \cite {48}.  Nous suivrons essentiellement la
terminologie de Manin (\cite {53}, Chap.5),  en particulier nous
utiliserons la terminologie pr\'estable, comme \'equivalente \`a
celle de courbe nodale ou semi-stable.  Les restrictions sur les
caract\'eristiques des corps r\'esiduels sont identiques \`a celles des
sections pr\'ec\'edentes. 

\proclaim D\'efinition 4.1. 1) Une courbe   pr\'estable,
de base un sch\'ema S, est un morphisme propre et plat $\pi: C \to S$
dont les fibres g\'eom\'etriques sont des courbes nodales connexes
(voir \S 2.1). Si $S$ est connexe, le genre (arithm\'etique) d'une
fibre est constant, c'est le genre de la courbe. \hfill\break
2)  Une courbe pr\'estable est dite stable si les fibres g\'eom\'etriques  $\{ C_s
\}$sont des courbes stables dans le sens de Deligne-Mumford, donc si
le groupe des automorphismes $\Aut(C_s)$ est fini.  

On sait que la stabilit\'e pour une courbe d\'efinie sur un corps
alg\'ebriquement clos, de genre $g \geq  2$, \'equivaut \`a
propri\'et\'e suivante : Si $E$ est une composante non singuli\`ere
rationnelle de $C,$ alors $E$ rencontre les autres composantes en au
moins trois points. Usuellement cette condition s'exprime de mani\`ere  particuli\`erement agr\'eable  sur le graphe dual classiquement
attach\'e \`a $C$; cela sera rappel\'e dans la suite (\S 7.1).
Du fait que pour une courbe pr\'estable le morphisme $\pi: C \to S$
est un morphisme localement d'intersection compl\`ete, le faisceau
dualisant relatif, d\'efini par (\cite {46}, \S 1)
$$\omega_{C/S} = \det (\Omega^1_{C/S}) \tag (4.1)$$
est localement libre de rang un.  On sait que sa formation est
compatible avec les changements de base. Dans loc.cit. il est prouv\'e\/  que la courbe pr\'estable  $\pi:C \to S$ est stable si et
seulement si $\omega_{C/S}$ est relativement ample.
  La d\'efinition (4.1) se g\'en\'eralise, comme il est
bien connu, \`a la consid\'eration de courbes stables
marqu\'ees (Harris-Morrison \cite {40}, Knudsen \cite {48},
Manin \cite {53}, Wewers \cite {67}):

\proclaim D\'efinition 4.2.  Soit $\pi: C \to S$ une courbe pr\'estable de base $S$. Soit $B\subset C$
un diviseur de Cartier relatif de degr\'e $d\geq 1$. 
i)  On dit que $B$ d\'efinit un marquage de $C/S$, ou que $C/S$
est marqu\'ee par $B$, si $B$ est \'etale sur $S$, et si de plus le
support de $B$ est contenu dans la partie lisse de $\pi$, donc si
pour tout point g\'eom\'etrique $s\in S$, les points de  $B_s\subset
C_s$ sont des points non singuliers. On parlera  dans la suite de
$(C/S , B)$ comme d'une courbe marqu\'ee. \hfill\break
ii)  La courbe marqu\'ee $(C/S , B)$ est dite stable, ou encore
$C/S$ est dite stable marqu\'ee par $B$, si pour tout point g\'eom\'etrique $s\in S$, l'une des propri\'et\'es \'equivalentes suivantes est
satisfaite:\hfill\break
$\bullet$ Le groupe $Aut(C_s , B_s)$ des automorphismes de $C_s$ pr\'eservant $B_s$ est
fini.\hfill\break
$\bullet$ Toute composante de la normalisation de $C_s$ qui est
rationnelle (resp. de genre 1),
contient au moins trois points exceptionnels (resp. un). On appelle
point exceptionnel un point qui
est soit dans le support de $B$, soit l'origine d'une branche issue
d'un point double.   

Quitte \`a  effectuer un changement de base \'etale, on peut supposer que  $B$ est   une somme
de sections disjointes $B = \sum_i P_i$, avec en tout point g\'eom\'etrique $s\in S$ la condition de stabilit\'e   satisfaite, $P_i(s) $ est un point non singulier, et si $i\ne j$, $ P_i(s) \ne
P_j(s)$. 
On montre alors, et cela sert de d\'efinition alternative, qu'une
courbe pr\'estable marqu\'ee $\pi:C \to S$ est stable si et seulement
si le faisceau $\omega_{C/S} (B)$ est  relativement ample (\cite {53}
Ch 5, \S 1). Rappelons qu'il  importe de distinguer  deux
situations, dans l'une  la courbe est marqu\'ee  par  $B$, et  dans l'autre $B$  est une somme de  
sections disjointes num\'erot\'ees,  la courbe est alors {\it
piqu\'ee}.  

 Dans  les deux cas  les isomorphismes   sont contraints  de   respecter  le marquage, respectivement les piq\^ures.  On rencontrera  aussi la situation interm\'ediaire dans
laquelle  les points marqu\'es sont affect\'es d'une couleur, mais
pas de mani\`ere univoque, les isomorphismes devant dans ce cas
simplement respecter la couleur. Cela  \'equivaut  \`a num\'eroter
les points par paquets. Dans ce cas   $B =
B_1+\cdots+B_r$ est d\'ecompos\'e en une somme de diviseurs
disjoints deux \`a deux, les points de
$B_\alpha$ \'etant vus comme les points d'une m\^eme couleur,  disons
la couleur $\alpha$. Posons  $d_\alpha = \deg {(B_\alpha)}$, et $d =
\deg (B) = \sum_\alpha d_\alpha$.  


On suppose maintenant,  comme dans le paragraphe pr\'ec\'edent  dont
on conserve les notations, que le groupe fini $G$ agit fid\`element
sur  $C/S$. Sous les hypoth\`eses fix\'ees au d\'ebut, le groupe  $G$
est r\'eductif, ce qui assure que l'op\'eration de passage au
quotient se comporte comme dans le cas lisse, en particulier on a
l'analogue de la proposition 3.2:

\proclaim Proposition 4.3.  Soit $\pi:C \to S$ une courbe  pr\'estable munie d'une action de $G$.
Alors la courbe quotient $D =  C/G$ est  pr\'estable, et sa formation
commute aux
changements de bases. En particulier, pour tout point g\'eom\'etrique
$s \in S$, $D_s =   {C_s}/G$.

\dem  Le point essentiel, qui est la commutation du quotient \`a
tout changement de base est ais\'e; voir la  proposition 3.2.\hfill\break
\qed

\bigskip

{\it 4.1.2. Actions stables} 

\bigskip

Il est maintenant tout \`a fait  clair que la conclusion de la
proposition 3.9 s'applique sans
changement \`a la situation de la proposition 4.3, avec dans ce cas
pour $G$-faisceau, le faisceau dualisant relatif $\omega_{C/S}$
\cite {20}. La
repr\'esentation de $G$ dans $H^0(C_s , (\omega_{C/S})_s)$ est ind\'ependante de $s\in S, S$ \'etant connexe. Revenons aux  conditions de
la proposition 4.1. Soit
$x\in C_s$ un point qui correspond \`a un point double d'une fibre
g\'eom\'etrique, de sorte que
$\hat{\Cal O}_{C_s,x} \cong   k[[X,Y]] /(XY)$.
Dans l'anneau local compl\'et\'e, $X = 0 $ et   $ Y = 0$
d\'efinissent les branches
formelles en $x$. Supposons que le stabilisateur $G_x$ soit non
trivial, alors $G_x$ op\`ere sur
l'ensemble \`a deux \'el\'ements compos\'e par les deux branches.
Remarquons aussi qu'on obtient  de mani\`ere analogue au cas lisse
une repr\'esentation fid\`ele de $G_x$ dans l'espace cotangent au
point $x$, soit un morphisme injectif
$G_x \hookrightarrow GL(2,k). $

La d\'efinition suivante trouvera sa motivation ult\'erieurement;
elle est r\'eserv\'ee  pour le moment aux courbes d\'epourvues de
marquage, et ne porte que sur l'action de $G$ au voisinage des points
doubles:
\proclaim D\'efinition 4.4.  L'action de $G$ sur la courbe pr\'estable $ C/S$ est dite stable
(certains disent Kummerienne \cite {63}), si pour tout point double 
$x$ d'une fibre g\'eom\'etrique, la repr\'esentation de $G_x$ dans
l'espace cotangent de la fibre en $x$ poss\`ede la propri\'et\'e
suivante:  pour tout $\sigma\in G_x$:  
$\det(\sigma) = 1$  si $\sigma $  fixe les branches, et  $ \det(\sigma) = -1 $ si $\sigma  $ \'echange les branches. 

Il est clair que la condition de stabilit\'e \'equivaut \`a dire
que le stabilisateur d'un point double
est soit cyclique, soit di\'edral, et dans ce dernier cas les \'el\'ements de d\'eterminant -1 dans la
repr\'esentation cotangente sont les \'el\'ements qui \'echangent les
branches. Plus pr\'ecis\'ement soit 
$ H = G_x^0 \ \triangleleft \ G_x$
le sous groupe d'indice au plus deux form\'e des \'el\'ements qui
fixent les branches. Ce groupe est cyclique, et  l'action de $H$ sur
l'espace cotangent en $x$
sur la branche $X = 0$ (resp. $Y = 0$), d\'efinit  un caract\`ere
primitif $\chi_X$ (resp.
$\chi_Y$) de $H$. La condition sur le d\'eterminant est alors  \'equivalente \`a
$$ \chi_X . \chi_Y = 1 \tag (4.2) $$
 Si $G_x$ est d'ordre deux, et s'il y a \'echange des branches, le groupe est di\'edral.  
  Soit $G$ agissant stablement sur la courbe pr\'estable
$C$, d\'efinie sur le corps alg\'ebriquement clos $k$. Les consid\'erations ci-dessous s'appliqueront donc \`a une fibre
g\'eom\'etrique.   On peut  classer les points \`a isotropie non triviale  en trois
familles disjointes deux \`a deux:\hfill\break
 \indent (I) points fixes non singuliers\hfill\break
 \indent   (II)  points doubles \`a isotropie cyclique \hfill\break
 \indent  (III) points doubles \`a isotropie di\'edrale

L'observation suivante, \'el\'ementaire, mais n\'eanmoins cruciale,
revient \`a dire que les points de type II ne doivent pas \^etre
consid\'er\'es comme des points de ramification. On verra cependant
qu'ils contribuent  de mani\`ere non triviale \`a la d\'eformation
universelle du rev\^etement, \`a  la diff\'erence des points de ramification.  Soit   $\omega_{C/S}$ la faisceau dualisant de la courbe
pr\'estable $\pi: C\to S$.  On a $\omega_{C/S} = \det (\Omega^1_{C/S})$
 avec  la description locale suivante \cite {48},\cite
{53}.   Supposons que localement   (pour la topologie \'etale)  sur
$U$, on r\'ealise $C$ comme une hypersurface d'\'equation ${f=0}$
dans un sch\'ema lisse de dimension relative 2 sur $S$. Alors  on a
une identification canonique, avec $I = (f)$:
$${\omega_{C/S}}_{\vert U} =  {\Hom_U (I/I^2 ,{ \wedge^2
\Omega^1}_{M/S} \otimes \Cal O_C}_{\vert U}) \tag (4.3)$$
 Il d\'ecoule de la suite exacte  (sur $U$)
$O \rightarrow I/I^2 \rightarrow
{\Omega^1}_{M/S}\otimes _{\Cal O_C}  \rightarrow {\Omega^1}_{C/S}
\rightarrow 0$ 
un morphisme canonique  ${{\Omega^1}_{C/S}}_{\vert U} \to {\omega_{C/S}}_{\vert U}  $
qui    est un isomorphisme en dehors du lieu singulier
\cite {48}. Rappelons que si $S = \Spec (k)$, et si $P$ un point
double, on a $\widehat {\omega}_{C,P} = \hat{\Cal O}_{C,P} . \omega$,
o\`u $\omega$ est la forme m\'eromorphe  qui sur la normalisation
de C, vaut  ${dX\over X}$ sur la branche $Y = 0$, et $-{dY\over Y}$
sur la branche  $X = 0$.  Dans le cas d'une base arbitraire, et pour
un quotient $\pi: C\to D=C/G$, d'une courbe pr\'estable $C$ par une
action stable de $G$ (D\'efinition 4.4), on va prouver que le
morphisme naturel $d\pi:  \pi^* (\Omega^1_{D/S}) \rightarrow
\Omega^1_{C/S}$, s'\'etend  canoniquement  en un morphisme (not\'e encore $d\pi$):
$$d\pi: \pi^*
(\omega_{D/S}) \rightarrow \omega_{C/S}\tag (4.4)   $$ 

 Ce  morphisme est un isomorphisme en dehors des points fixes de types I
et III. On  le construit d'abord localement, et par
naturalit\'e  on l'\'etend en un morphisme global.   On obtient ce morphisme en combinant le
morphisme  $d\pi: \pi^* \Omega^1_{D/S} \rightarrow \Omega^1_{C/S} $
avec la fl\`eche $\Omega^1_{C/S} \to \omega_{C/S}$, et la fl\^eche analogue. Au voisinage d'un
point double $P$ de type II, d'image $s\in S$, on peut utiliser un
syst\`eme de coordonn\'ees locales pour d\'ecrire le morphisme $\pi$
(par exemple Wewers \cite{67}, \S 2.1, Jarvis  \cite {44} \S 5.1). On
peut de la sorte choisir \'etale-localement des coordonn\'ees
$(x,y)$ en $P$ ,  telles que $C$ soit  sur ces coordonn\'ees  la
courbe relative
$$(xy - a =0) \subset {\Bbb A_S}^2$$
et   que l'action du stabilisateur $G_P$ (d'ordre $e\geq 2$)  est d\'ecrite
par $\sigma (x) = \zeta x,\quad  \sigma (y )= \zeta^{-1} y$,
$\sigma$ \'etant un g\'en\'erateur de $G_P$\footnote {Il est bien
connu que ces coordonn\'ees sont bien d\'efinies aux changements
suivants: image r\'eciproque par localisation \'etale, substitution
$(x,y,a) \mapsto (\alpha x , \beta y , \gamma a)$ o\`u $\alpha,
\beta \in {\Cal O}_P^*, \,\, \gamma = \alpha \beta \in {\Cal O}_s^*$
\cite {44}.}. On en d\'eduit imm\'ediatement une description
locale-\'etale  de $D$ par  l'interm\'ediaire des coordonn\'ees   $u
= x^{e}, \,v = y^{e}$  soumises \`a l'\'equation   $uv - a^{e} =
0$. Ces descriptions   fournissent des g\'en\'erateurs locaux
$\omega (x,y)$ et $\omega (u,v)$ pour respectivement $\omega_{C/S},
$ et $ \omega_{D/S}$. Le morphisme local requis est l'isomorphisme d\'efini par  $\pi^*(\omega (u,v)) \mapsto  \omega (x,y) $.
 Il est facile de voir que ce morphisme est canonique,  du
fait que tout autre syst\`eme de coordonn\'ees locales se d\'eduit
de $(x,y)$  par $x' = \alpha x, \, y' = \beta y$, et $a' = \alpha
\beta a$, apr\`es  localisation \'etale.   On peut d'une autre mani\`ere invoquer l'argument suivant.  Le faisceau $\omega_{C/S}$ est le faisceau dualisant, de sorte qu'on a un isomorphisme fonctoriel  $\Hom_C({\Cal F} , \omega_{C/S}) \cong \Hom_S(R^1p_*({\Cal F}) , {\Cal O}_S)$.
 Si ${\Cal F} = \pi^*(\omega_{D/S})$, alors le morphisme global cherch\'e provient par dualit\'e du morphisme trace 
$$Tr = {1\over \vert G\vert } \sum_{g\in G} \, g: \pi_*({\Cal O}_C) \to {\Cal O}_D$$  

On notera que  les faisceaux
$\omega_{C/S}$ et $ \Omega^1_{C/S}$ sont  de mani\`ere naturelle  des
$G$-faisceaux, et seront   en permanence pris comme tels. La
structure du morphisme 4.4 peut \^etre pr\'ecis\'ee:

\proclaim Lemme 4.5.  Soit $p: C \to S$ une courbe pr\'estable munie
d'une action de $G$. On suppose que $G$ agit fid\`element et
stablement sur les fibres g\'eom\'etriques. Soit $\pi: C \to B =
C/G$ le morphisme de $C$ sur la courbe quotient. Pour tout entier $l
\geq 1$, on a une suite exacte
$$0\longrightarrow \pi^*\omega^{\otimes l}_{ D/S} \longrightarrow
\omega^{\otimes
l}_{ C/S}\longrightarrow {\Cal V}_l \longrightarrow 0 \tag (4.5)$$
o\`u ${\Cal V}_l$ est un ${\Cal O}_S$-module  coh\'erent plat, dont
la formation est compatible avec les changements de base.  Le support
de  ${\Cal V}_l$ est concentr\'e aux points de type (I) et
(III). 

\dem   Partant  du morphisme   $\pi^* (\omega_{D/S}) \rightarrow
\omega_{C/S}$,  par passage \`a la puissance tensorielle $l$-i\`eme,
on obtient un morphisme   $\pi^*\omega^{\otimes l}_{ D/S}
\longrightarrow \omega^{\otimes
l}_{ C/S}$. Soit ${\Cal V}_l$ le conoyau. Pour d\'ecrire ce faisceau,  supposons
d'abord $S = {\Spec}(k)$, le corps $k$ \'etant alg\'ebriquement clos.
Il n'est pas difficile dans ce cas
de  d\'ecrire (localement) explicitement   le morphisme
$\pi^* (\omega_D^{\otimes l}) \longrightarrow \omega_C^{\otimes l} $,
et  de la sorte ${\Cal V}_l$.
  Traitons d'abord le cas $l = 1$. En un point non singulier $x\in C$,
d'image $y\in D$,  soit $t$ un param\`etre local, tel que le
stabilisateur $G_x$ agisse par
$\sigma. t = \chi(\sigma) t, \quad (\sigma \in G_x)$, o\`u  $\chi \in
{\widehat G}_x $ est un caract\`ere primitif.
Soit $e$ l'ordre de $G_x$. Alors $u = t^{e}$ est un
param\`etre local en $y$ et le
morphisme (4.5) a pour description locale en $x$
$${\Cal O}_x du \longrightarrow {\Cal O}_x dt ,\quad\hbox{\rm avec}\quad
du = et^{e-1}dt \tag (4.6)$$
On trouve  ainsi que $\dim({\Cal V}_x) = e-1$. Si maintenant $x$ est
un point double de type (II),
selon la nomenclature pr\'ec\'edente, on peut supposer que
formellement l'action de $G_x$ sur les
branches est donn\'ee par l'interm\'ediaire de deux caract\`eres
oppos\'es, $\chi$ et $\chi^{-1}$, d'ordre $e > 1$ (condition de
stabilit\'e (1)). On peut ainsi consid\'erer des param\`etres
locaux $X,Y$ le long des branches de sorte que
$${\widehat{\Cal O}}_x = k[[X,Y]]/(XY)
\,\,,\,\,{\widehat{\Cal O}}_y = 
k[[U,V]]/(UV),\quad \hbox {\bf avec}\quad U = X^{e} , V = Y^{e} $$
o\`u $U,V$ d\'efinissent des param\`etres locaux le long des
branches en $y$. Notons par une lettre
minuscule  $x,y,u,v,$ les classes respectives de 
$X,Y,U,V$; alors le morphisme (4.5) a pour description locale au point $x$
$$ \widehat{\Cal O}_x ({du\over u} , -{dv\over v}) \longrightarrow
\widehat{\Cal O}_x ({dx\over x} ,
-{dy\over y}) \tag (4.7) $$
avec ${du\over u} = e {dx\over x}$, et ${dv\over v} = e {dy\over y}$.
Le  morphisme (4.5) est  bien en un tel point un
isomorphisme, comme impos\'e dans la construction
du morphisme 4.4.

Supposons maintenant que $x$ est un point fixe de type (III), et que
$G_x \cong {\Bbb D}_e, (e\geq 1)$. Dans ce cas il est clair que $y$
est non singulier, et que les param\`etres  le long des branches \'etant  choisis comme au dessus, $t = x^{e} + y^{e}$ est un param\`etre
en $y$. La description locale de (4.5) est dans ce cas:
$$\widehat{\Cal O}_x dt \longrightarrow \widehat{\Cal O}_x ({dx\over x} ,-
{dy\over y}),\quad \hbox{\rm avec}\quad dt = (ex^{e} , - ey^{e})({dx\over
x} , - {dy\over y}) \tag (4.8)$$
En particulier cette description donne ${\Cal V}_x \ne 0$, en fait le
rang est  $\dim({\Cal V}_x) = 2e-1$. Revenons au morphisme (4.5); on
vient de voir qu'il est injectif sur les fibres, il est donc injectif
avec un  
conoyau ${\Cal V}_l$ plat sur $S$, comme il r\'esulte du crit\`ere
local de platitude, et donc de
formation compatible aux changements de base. Le support de ${\Cal
V}_l$ est concentr\'e
exclusivement en les  points fixes de type (I) ou (III).\hfill\break
\qed

Dans le cas $l = 1$, on posera dans la suite  ${\Cal V} = {\Cal V}_1$.
On conserve les notations du Lemme 4.5. La suite exacte 4.5 justifie la d\'efinition suivante, dans laquelle $\Div $ est le diviseur de cartier construit dans \cite {58}:
 \proclaim D\'efinition 4.6.   Le diviseur de Cartier relatif $R_\pi
= \Div {\Cal V} = \Div(\pi^*(\omega_{D/S} \to \omega_{C/S})$ est appel\'e le diviseur de ramification du
morphisme $\pi: C \to D$. 

On notera que le support de $R_\pi$ est l'ensemble des points de type
I ou de type III. Par ailleurs, de la d\'efinition, jointe \`a la
suite exacte 4.4, il vient l'\'egalit\'e (formule de ramification)
$$ \omega_{C/S} \cong \pi^* (\omega_{D/S}) \otimes \Cal O (R_\pi)
\tag (4.9)$$
On souhaite relier le diviseur $R_\pi$ aux sous-sch\'emas
de points fixes des sous-groupes cycliques de $G$, comme dans le cas
lisse. Soit  $H$ un sous-groupe cyclique de $G, \,\, H\ne 1$. Ce qui
pr\'ec\`ede justifie la d\'efinition d'apparence  peu naturelle
suivante:
\proclaim D\'efinition 4.7.   Soit une action stable de $G$
sur la courbe pr\'estable $C$
d\'efinie sur   $k$. Soit $p\in C$:\hfill\break
  On dit que $H$ fixe strictement $p$, si $p$ \'etant non
singulier, on a $Hp = p$,  et   $p$ \'etant un point double, si on a
d'une part $Hp = p$, et d'autre part  tout \'el\'ement $\sigma \in H, \sigma\ne 1,$ 
inverse les branches. 

On  observera que si $H$ fixe strictement le point double $p$, alors
$H$ est d'ordre deux, et
l'\'el\'ement distinct de l'identit\'e de $H$ est une reflexion,
i.e. inverse les branches. On \'etend la
d\'efinition  au cas relatif. La proposition suivante pr\'ecise la
structure du sous-sch\'ema des points
fixes d'un sous-groupe cyclique en pr\'esence de points doubles.

\proclaim Proposition 4.8.   Soit $\phi:C\to S$ une courbe pr\'estable sur laquelle le groupe $G$
agit stablement, et soit $\pi:C \to D =  C/G$  le rev\^etement
associ\'e. \hfill\break  1) Pour tout sous-groupe cyclique $H$ de
$G$, le sous-sch\'ema $Fix(H) = C^H$ des points fixes de $H$ est r\'eunion disjointe d'une partie horizontale et d'une partie verticale:
$$Fix(H) = Fix(H)_{hor} \sqcup Fix(H)_{ver}  \tag (4.10)$$
o\`u $Fix(H)_{hor}$ est un diviseur de Cartier relatif qui a pour
support les points fix\'es strictement
par $H$. Le support de $Fix(H)_{ver}$ est inclus dans les points
doubles de type (II) de $H$.  \hfill\break
2)   Comme dans le
cas lisse on a l'\'egalit\'e de diviseurs, la somme \'etant \'etendue
aux sous-groupes cycliques non triviaux 
 $$R_\pi = \sum_H \, \phi (\vert H\vert ) Fix (H)_{hor} \tag (4.11)$$     
 
\dem Fixons $\sigma$ un g\'en\'erateur de $H$, avec ${\Card H} = e$.
Soit $x\in Fix(H)$ un point fixe
de type (I),  donc non singulier dans sa fibre, et   soit $\phi(x) =
s\in S$. C'est un fait standard, observ\'e dans la proposition 4.4,
qu'au point  $x, \,Fix(H)$ est un diviseur de Cartier relatif \'etale
sur $S.$  Supposons maintenant que $x$ soit du type (II). Quitte \`a
 passer aux anneaux locaux compl\'et\'es, on peut supposer que
${\widehat{\Cal O}_x} \cong   \widehat{\Cal O}_s [[X,Y]]/
(XY-a),\quad (a\in\hat {\Cal
M}_s),$
l'action de $H$ \'etant d\'ecrite   par $\sigma(X) = \zeta X, \,\,
\sigma(Y) = {\zeta}^{-1} Y$,
pour une  racine  primitive d'ordre $e$ convenable de l'unit\'e. Il
est imm\'ediat que l'id\'eal du sous-sch\'ema  ferm\'e $Fix(H)$ en
$x$ est engendr\'e par les classes de $\sigma(X) - X$ et $\sigma(Y)
- Y$, il est donc \'egal \`a $(X,Y)$. On voit  ainsi qu'en ce point
$Fix(H)$ n'est pas un diviseur de Cartier. Supposons maintenant le
point $x$ de type (III), et que $x$ soit un point fixe strict de $H$;
alors  avec les notations pr\'ec\'edentes, $\sigma(X) = Y, \sigma(Y)
= X$, et
dans ce cas $H$ est d'ordre deux. L'id\'eal de $Fix(H)$ est engendr\'e
 par la classe de $X - Y$, en particulier c'est un diviseur de
Cartier en ce point. L'anneau local compl\'et\'e de $Fix(H)$ en $x$
est
${\widehat{\Cal O}}_{Fix(H),x} \cong  \widehat{\Cal O}_s
[[X]] /(X^2 - a), $
avec $a \in \hat {\Cal M}_s.$  Il n'est pas \'etale sur
${\widehat{\Cal O}}_s$.

L'ensemble (ouvert) des points de $Fix(H)$ en lesquels $Fix(H)$ est
un diviseur de Cartier, est \'egal au support du diviseur de
ramification  du morphisme quotient $C \to   D=C/H$ (d\'efinition
4.6). Finalement $Fix(H)$ se d\'ecompose en une somme disjointe
$$Fix(H) = Fix(H)_{hor} \sqcup Fix(H)_{ver}  $$
o\`u $Fix(H)_{hor}$ est un diviseur de Cartier relatif support\'e
par les points fixes stricts de $H,$ et
$Fix(H)_{ver}$ a un support  inclus dans l'ensemble des points fixes
de type (II) de $H$, ce qui prouve 1). \hfill\break  2) L'\'egalit\'e
 4.10  a \'et\'e observ\'ee pr\'ec\'edemment pour les points
fixes non singuliers dans leur fibre.   Reste a v\'erifier  l'\'egalit\'e en les points de type III. Soit $x$ un tel point, avec un
stabilisateur di\'edral d'ordre $2e$. Les notations \'etant comme
au-dessus, une \'equation locale de $R_\pi$ en $x$, est $X^{e} -
Y^{e}$.  Si le stabilisateur $G_x$ a pour g\'en\'erateurs $\sigma$ et
$\tau$, avec $\sigma^{e} = \tau^2 = (\tau\sigma)^2 = 1$, alors l'\'equation locale en $x$ de $Fix(\tau\sigma^i)_{hor} $ est comme obtenu
ci-dessus $X - \zeta_e^i Y$. Le r\'esultat se r\'eduit \`a l'\'egalit\'e  $X^{e} - Y^{e} = \prod_{i = 0}^{e-1} \, (X - \zeta_e^i
\,Y)$.\hfill\break
\qed  

\proclaim Corollaire 4.9.  Sous les hypoth\`eses de la proposition
4.2, la base $S$ \'etant suppos\'ee connexe, notons $d$ le degr\'e
de $Fix(H)_{hor}$ relativement \`a $S$. Pour tout point
g\'eom\'etrique $s \in S$, on a
$$d =  \#{\{ x\in Fix(H)_{hor} (s) , x\,\,\hbox{\bf de type (I)} \}}
   + \,\,2\#{\{x \in Fix(H)_{hor} (s), x\,\,  \hbox {\bf de type (III)} \}}
$$ 

\dem On vient de voir que $Fix(H)_{hor}$ est \'etale sur $S$ en les
points de type (I), et que si $x\in
Fix(H)_{hor}$ est un point de type (III) d'image $s$, alors
$\dim({\widehat{\Cal O}}_{Fix(H)_{hor},x} \otimes k(s)) = 2$.
  Le r\'esultat  en d\'ecoule.\hfill\break   \qed  
  
  \bigskip     {\it 4.1.3. G\'eom\'etrie du quotient par une
action stable}  
 \bigskip

Pour utilisation ult\'erieure, la courbe  pr\'estable $C$ \'etant
toujours d\'efinie sur le corps $k$, nous allons comparer les
$G$-faisceaux $\pi^*(\Omega^1_{D/k})$ et $\Omega^1_{C/k}$. Les
notations pr\'ec\'edentes \'etant toujours en vigueur, rappelons que
le foncteur    
$\pi_\star^G: (G-Coh\, C) \rightarrow (Coh \,D)$
transforme un $G$-faisceau (coh\'erent) ${\Cal F}$ sur $C$ en le
faisceau sur $D$ des sections
$G$-invariantes $\pi_\star ({\Cal F})^G$  de $\pi_\star ({\Cal F})$.
Le premier r\'esultat est:

\proclaim Proposition 4.10.  1)  Le morphisme
$\Omega^1_{D/k} \to \pi_*^G
(\Omega^1_{C/k})$ est injectif et le conoyau est concentr\'e aux
points images des points doubles de type (II).\hfill\break
  2)  Le morphisme $\pi_*^G (\theta_C) \to \theta_D$ est injectif
et le conoyau est concentr\'e aux points de $B$ images d'un point
de type (I), donc nonsingulier, ou un point double de type (III); en
ces points le conoyau est de dimension un.  \hfill\break
En particulier, s'il n'y a pas de point fixe de type (III), on a     
$\pi_*^G (\theta_C) =  \pi_*^G (\theta_C(-R)) = \theta_D(- \Delta), $ 
o\`u $\Delta = \sum_j Q_j$ (resp. $R = \sum_i
P_i$)  d\'esigne le diviseur de branchement  (resp. de ramification)
r\'eduit. 

\dem  Partons de la suite exacte
$\pi^*(\Omega^1_{D/k}) \longrightarrow\Omega^1_{C/k}
\longrightarrow\Omega^1_{C/D}\longrightarrow0.$
Appliquant le foncteur $\pi_*^G$, on obtient un morphisme
$\Omega^1_{D/k} \longrightarrow\pi_*^G (\Omega^1_{C/k}) $
qui est un isomorphisme en dehors des points de branchement (lisses),
et des points doubles. Soit
$y = \pi(x) \in B$ un point de branchement non singulier, et soit $e$
l'indice de ramification
correspondant. Si on choisit des param\`etres locaux $t$ et $u$ en
$x$ et $y$ respectivement, comme dans le lemme 4.5,  le morphisme
(4.11) admet  alors la description locale suivante:
$$ g(u)du \in (\Omega^1_{B/k})_y \longmapsto eug(u) {dt\over t} \tag (4.12)$$
Comme il est clair que $\pi_*^G (\Omega^1_{C/k})_y  =  u{dt\over
t}{\Cal O}_{D,y}$, on
trouve bien que le morphisme (4.11) est bijectif en un tel point.
Soit maintenant le cas de $x$
point double avec isotropie cyclique, d'ordre $e \geq 2$. En
utilisant les coordonn\'ees locales le
long des branches, introduites dans le lemme 4.5, on obtient la
description suivante pour les
deux termes du morphisme (4.11); pour le premier de ces termes:
$$(\Omega^1_{D/k})_y = {\widehat{\Cal O}_y\, du \oplus \widehat{\Cal
O}_y\, dv \over (udv + vdu)}$$
et pour l'autre:
$$ (\pi_*^G \widehat{(\Omega^1_{C/k})}_y =  \left({\widehat{\Cal
O}_x\, dx \oplus
\widehat{\Cal O}_x\, dx\over (xdy + ydx)} \right)^{G_x}$$
Soit $\sigma$ un g\'en\'erateur de $G_x$, agissant par $\sigma x =
\zeta x$, et $ \sigma y =
{\zeta}^{-1}y$. Un calcul \'el\'ementaire conduit \`a l'expression
suivante d'une forme $\omega$,  $G_x$ - invariante:
$$ \omega = ({A(u)\over x} , by)dx + (0 , {D(v)\over y}) dy, \quad(b
\in k,\quad A(0) = D(0) =
0)\tag (4.13)$$
On voit alors imm\'ediatement que le morphisme (4.11) est injectif en
$x$, et que son image est
l'ensemble des 1-formes qui dans l'expression (4.13) ont $b = 0$. En
particulier le conoyau est de
dimension un, comme indiqu\'e. Pour la  seconde assertion, en dualisant le morphisme 4.11, on obtient un morphisme $\pi_*^G(\theta_C) \to \theta_D$; la conclusion  rel\`eve des m\^emes calculs locaux  que nous  omettons.\hfill\break
  \qed  

\rema {4.1 }  Si l'action
stable du groupe $G$ sur $C/S$ n'a pas de points du type III, alors
le diviseur de ramification est de la forme habituelle. La partie
horizontale (Proposition 4.8) de $Fix (H)$, $H$ sous-groupe cyclique,
est une somme  de diviseurs de Cartier relatifs \`a supports
disjoints $Fix_{hor} (H) = \sum_{C\subset H}
\,\Delta_C$, en notant $\Delta_C$   le lieu des points de stabilisateur
exactement $C$. Il est   imm\'ediat de retrouver l'expression
 usuelle  du diviseur de ramification   $R = \sum_H \, (\vert H\vert - 1) \Delta_H. 
\quad\lozenge$

\beginsection 4.2. Collision des points de ramification 

\bigskip

  {\it 4.2.1. $G$-type d'un $G$-diviseur} 
\bigskip

 Soit une action stable de $G$ sur $C$.
On suppose   que
$D\subset C$ est un diviseur de Cartier relatif, \'etale de degr\'e
$N$ sur $S$, de support contenu dans la partie lisse de $\pi: C\to
S$, et  $G$-stable.  Cela permet de  parler du $G$-type de la fibre $D_s$ en un point g\'eom\'etrique
$s \in S$ (section 2.1). Tout d'abord une remarque \'el\'ementaire:
\proclaim Lemme 4.11.  Le $G$-type de $D_s$ en un point g\'eom\'etrique $s \in S$ est constant.

\dem Par changement de base \'etale, on peut supposer que $D$ est
somme de $N$ sections disjointes $\{P_i \}$, donc $D = \sum_{i=1}^N
P_i$. Un \'el\'ement $\sigma \in G$ agit par permutation sur les
$\{P_i \}$. Soit pour tout indice $i$, $H_i$ le stabilisateur de
$P_i$. Le faisceau  conormal le long de la section $P_i$, 
${\Cal N}_{P_i/S}^\lor \,\cong \, P_i^* ( \Omega^1_{C/S}) $
est un $G$-faisceau inversible sur $S$. On peut donc l'\'ecrire
${\Cal N}_{P_i/S}^\lor \cong {\Cal
L}_i \otimes V_i$, pour une certaine repr\'esentation  $V_i$ de
degr\'e un de $H_i$; soit $\chi_i$
le caract\`ere de $V_i$. Alors le type de $D_s$ est $\sum_{i=1}^N
[H_i , \chi_i]$,  donc constant le long des fibres.\hfill\break
\qed

Notre objectif est maintenant de pr\'eciser le comportement des
points de ramification lorsqu'un
rev\^etement (galoisien de groupe $G$) d\'eg\'en\`ere, c'est \`a
dire lorsque le rev\^etement se d\'eplace dans une famille de courbes
stables. Le r\'esultat principal de cette section est l'observation
que seuls les points de type (III) (4.3), donc ceux \`a isotropie
di\'edrale, sont responsables de la
collision \'eventuelle des points de ramification. Ce comportement
pathologique sera \'elimin\'e  par le marquage \`a priori de la
courbe au moyen d'un diviseur contenant les points de
ramification.
\bigskip

 {\it 4.2.2. Chevalley-Weil (bis
repetita)}
 \bigskip

 On suppose  maintenant que    $C/S$ est une courbe stable et que l'action de $G$ sur $C$ est stable.  Dans ce contexte si $D$  marque la courbe, on suppose que $D$ est $G$-invariant. Si la base est connexe    le  $G$-type de  $D$ (Lemme 4.11), ou de $(C,D)$ est d\'efini. Si  en plus  $\pi$ est g\'en\'eriquement lisse,   on peut attacher  une  donn\'ee de ramification \`a l'action de $G$, qui est celle attach\'ee \`a la fibre g\'en\'erique. Cette donn\'ee est  dor\'enavant  fix\'ee.    Nous allons   suivre les caract\'eristiques  combinatoires lorsque
la fibre $C_s$  d\'eg\'en\`ere. Par inversion des formules de Chevalley-Weil au point g\'en\'erique, les repr\'esentations de Hurwitz (th\'eor\`eme 3.1),  
$$ Hurw_m = H^0(C_s , \omega_{C_s} ^{\otimes m})\quad (m\geq 1,\,\,
s\in S) \tag (4.14)$$
sont acquises.    On fixe une
fibre singuli\`ere, donc une   $G$-courbe stable d\'efinie sur un corps alg\'ebriquement
clos $k$. Soit   $D =  C/G$ la courbe quotient. Les points
doubles de $D$ sont donc les images des points doubles de type (II)
de $C$.

 Installons quelques
notations et remarques pr\'eliminaires.
Rappelons  que la suite exacte (4.4)  d\'efinit un
faisceau ${\Cal V}_l$ concentr\'e en les points  fixes de types (I)
et (III). Pour $l\geq 2$, cette suite  entra\^{\i}ne, par l'hypoth\`ese de stabilit\'e de $C$, et par le fait que pour $l\geq 2$,
$H^1(C , \omega_C^{\otimes l}) = 0$ (\cite {19} Thm 1.2), donc en
notant provisoirement $B = C/G$ la courbe quotient, pour \'eviter des
confusions avec le diviseur de marquage $D$, on note que
$$\dim( H^0(C , \omega_C^{\otimes l})^G) = \chi(B , \omega_B^{\otimes
l}) + dim ({\Cal V}_l^G)
\tag (4.15)$$
Si $l = 1$,  des consid\'erations  similaires  conduisent \`a
$H^0(B , \omega_B) = H^0(C ,
\omega)^G$, en particulier
$p_a(B) = g'$. 
Pour $l\geq 2$, $\dim ({\Cal V}_l^G)$ est la
contribution significative \`a $\dim( H^0(C ,
\omega_C^{\otimes l})^G)$.   Nous sommes maintenant  en mesure de
reprendre les calculs \`a la Chevalley-Weil de la section 3.2 dans
le contexte de cette section, donc pour une courbe  nodale (stable).
Rappelons tout d'abord quelques faits \'el\'ementaires sur les repr\'esentations irr\'eductibles  des groupes di\'edraux.

On identifie le groupe di\'edral ${\Bbb D}_e \,\, (e\geq 1)$, d'ordre
$2e$ avec le sous-groupe de  $\GL(2,k)$ engendr\'e par les matrices
$$\tau = \pmatrix {0&1\cr1&0 }, \,\,\sigma = \pmatrix
{\epsilon&0\cr 0&\epsilon  }$$
o\`u $\epsilon$ est une racine primitive e-i\`eme de l'unit\'e.
Notons $\psi_1$ et $\psi'_1$ les
caract\`eres de degr\'e un de  ${\Bbb D}_e$ d\'efinis par
$ \psi_1 = 1, \quad\hbox{et}\quad \psi'_1(\sigma) = 1,
\;\;\psi'_1(\tau) = -1.$
Pour tout $j\geq 0$, soit le caract\`ere $\phi_j$ de degr\'e deux
donn\'e par
$$\phi_j(\tau) = 0,\quad \phi_j(\sigma) = \epsilon^j +
\epsilon^{-j}\quad (2\cos ({2\pi j\over e})
\quad\hbox{si}\quad k = {\Bbb C})$$
de sorte que $\phi_j$ est le caract\`ere de la repr\'esentation d\'ecrite par
$$\tau \mapsto \pmatrix{ 0&1\cr1&0 } , \quad \sigma \mapsto
\pmatrix {\epsilon^j&0\cr 0&\epsilon^{-j} }$$
  On observe  les relations de p\'eriodicit\'e $\phi_{ je } =
\phi_j$, et $\phi_j$ est irr\'eductible sauf si $j \equiv 0\pmod e$
ou bien $ e = 2m$ et $j \equiv m \pmod e$. Dans ces deux cas les
repr\'esentations se d\'ecomposent en
$$\phi_0 = \psi_1 + \psi'_1, \quad \hbox{ et si }\quad e = 2m,
\;\,\phi_m = \phi'_m + \phi''_m \tag
(4.16)$$
$\phi'_m$ et $\phi''_m$ \'etant deux caract\`eres de degr\'e un donn\'es par
$\phi'_m(\sigma) = \phi''_m(\sigma) = -1,\quad \hbox{et}\quad
\phi'_m (\tau) = 1, \,\phi''_m(\tau) =
-1 $.

On suppose maintenant que les orbites singuli\`eres de l'action de
$G$ sur $C$ sont, en dehors des
orbites de type (II), compos\'ees de $b$-orbites de type (I) et de
$s$-orbites de type (III).
Intuitivement, en un point double de type (III), si on regarde la
courbe comme sp\'ecialisation
\'equivariante d'une courbe lisse, alors il y a coalescence de deux
orbites de points fixes. Cela va
\^etre pr\'ecis\'e sous peu;  on va voir que dans cette hypoth\`ese
le stabilisateur est d'ordre deux, donc di\'edral, engendr\'e par
une reflexion qui \'echange les deux branches.

Rappelons, pour \'eviter toute confusion, la terminologie utilis\'ee
(voir en particulier la d\'efinition 4.7), lorsqu'on a affaire \`a
un point double \`a isotropie non triviale.

Un sous-groupe cyclique   $H\subset G$   {\sl fixe strictement} $x$,
si lorsque $x$ est non
singulier $H \subset G_x$, et si $x$ est un point double, $H$ fixe
$x$ dans le sens de la
d\'efinition 4.5, donc $H$ est d'ordre deux engendr\'e par une r\'eflexion \'echangeant les branches. On dira  aussi que $H$ est un {\sl
stabilisateur} de $x$, si $H = G_x$ lorsque $x$ est non singulier, et
si $x$ est un point double, comme ci-dessus, $H$ fixe strictement
$x$, donc est en particulier "di\'edral d'ordre deux". Dans ce
dernier cas, si $G_x = {\Bbb D}_l$, il y a donc $l$ sous-groupes de
$G_x$ qui sont des stabilisateurs de $x$ !. Si $l>1$, $G_x$ tout
entier n'est ainsi pas consid\'er\'e comme un sous-groupe
stabilisateur de $x$.

Rappelons que dans le probl\`eme de classification  qui nous occupe,
on a fix\'e au pr\'ealable une
donn\'ee de Hurwitz $\xi$, qui  sp\'ecifie $r$ classes de conjugaison
de caract\ res primitifs de sous-groupes cycliques,
compt\'ees avec leurs  multiplicit\'es respectives, donc  fournit en
particulier une liste de classes de
sous-groupes cycliques, non n\'ecessairement distinctes, not\'ees
$\Gamma_1,\cdots,\Gamma_r$. On suppose que les ordres respectifs sont
$\delta_1 \leq \delta_2\leq  \cdots \leq \delta_r$. Pour tout $e \geq
2,\,\, e\mid \Card G$, posons
$$N_e = \Card \{ i \in [1 , r] \,\,,\,\, \delta_i = e \} \tag (4.17)$$
de sorte que $N_e$ repr\'esente  le nombre d'orbites singuli\`eres de
cardinal $\Card G \over e$ dans toute courbe lisse $C$ supportant une
action de $G$ du type sp\'ecifi\'e. Les m\^emes conventions de
notations s'appliquent \`a un quelconque sous-groupe  $H$ de $G$,
et $N_e(H)$ aura donc la signification ci-dessus. Sous les m\^emes
conditions sur $e$, on notera $N'_e$ le nombre de $G$-orbites
singuli\`eres de cardinal $\Card G \over e$ qui sont contenues dans
la partie lisse de $C$, donc
$$r' = \sum_{2\leq e,e\mid\Card G} N'_e \tag (4.18)$$
Pour tout sous-groupe $H$, on peut  aussi d\'efinir $r'(H), s(H),
N'_e (H)$. Le r\'esultat qui suit est
crucial pour clarifier la coalescence  \'eventuelle des points de ramification:
\proclaim Th\'eor\`eme 4.12.  On a les \'egalit\'es:  
i)  $r = r' + 2s$, et \ \   
ii)  $ N_e = N'_e$ si \ $e \geq  3$, et $N_2 = N'_2 + 2s$.
 
\dem  L'assertion (2) est visiblement plus forte que (1) et donc
l'implique;  on commence cependant
par prouver (1). Choisissons dans chaque orbite singuli\`ere un point
$x_i\,\,(1\leq i\leq r')$, de stabilisateur $H_i = G_{x_i}$ pour
celles form\'ees de points non singuliers, et  
$y_j, \,\,(1\leq j\leq s)$ avec $K_j = G_{y_j}$ pour celles form\'ees
de points doubles. Posons
$\lambda_i = ({\Cal V}_l)_{x_i}, \,\,\mu_j = ({\Cal V}_l)_{y_j}\quad
(l=1,2,\cdots).$
On a  une d\'ecomposition  de $G$-modules
$$\Gamma (C , {\Cal V}_l) = \left( \bigoplus_i \Ind_{H_i}^G
(\lambda_i) \right) \bigoplus
\left(\bigoplus_j \Ind_{K_j}^G (\mu_j)\right)  \tag (4.19)$$
d'o\`u par r\'eciprocit\'e de Frobenius
$${\dim \Gamma (C , {\Cal V}_l)^G} = \sum
\dim(\lambda_i^{H_i})\,+\,\sum_j \dim( \mu_j^{K_j})
\tag (4.20)$$
L'analyse de la partie $H_i$-invariante de $\lambda_i$ rel\`eve de la
description faite
dans la section 3.2; on trouve de mani\`ere analogue:
$$\dim(\lambda_i^{H_i}) = l-1 + \left( \langle\langle {l\over
e}\rangle\rangle - {l\over e} \right)
\quad \hbox{si}\quad e = \Card {H_i}$$
Supposons maintenant $K_j = {\Bbb D}_e$. On peut lin\'eariser
l'action de $K_j$ relativement \`a des
coordonn\'ees $x,y$ le long des branches au point $y_j$, et ainsi
supposer que les g\'en\'erateurs
$\sigma$ et $\tau$ de $K_j$ agissent par
$\sigma(x) = \epsilon x,\,\,\sigma(y) = \epsilon^{-1} y;
\quad\tau(x) =  y ,\,\,\tau(y) = x$,
avec $\epsilon \in k^*$ une racine primitive d'ordre $e$. Le
morphisme quotient a pour forme locale
au point $y_j$:  $(x,y) \to u = x^{e} + y^{e}$, de sorte que  la
1-forme $\omega = ({dx\over x}, -{dy\over y})$ est une base locale de
$\omega_C$ en $y_j$. On a $\sigma^* (\omega) = \omega$ et $\tau^*
(\omega) = -\omega$. Le module $\mu_j$  admet la description
$$\mu_j =    { {\Cal O}_{C,y_j} \omega^{\otimes l}\over  {\Cal
O}_{C,y_j} (du)^{\otimes l}}
\,\cong\,\left({{\Cal O}_{C,y_j} \over  (x^{el} , (-1)^l
y^{el})}\right) \omega^{\otimes l} \tag
(4.21)$$
du fait que  $(du)^{\otimes l} =  (x^{el} , (-1)^l y^{el}) \omega^{\otimes l}$.

Il est ais\'e de d\'ecomposer le $K_j$-module $\mu_j$ en repr\'esentations irr\'eductibles  du groupe $K_j = {\Bbb D}_e$. Le r\'esultat est
  $$\mu_j =  \cases {\left( 2\psi_1 + \sum_{i=1}^{el-1} \phi_i \right)
\psi'_1 & si $ l$  impair  \cr
    \psi_1 + \psi'_1 + \sum_{i=1}^{el-1}
\phi_i & si $ l$   pair} \tag (4.22)$$
Si $l$ = 2, ceci conduit \`a ${\dim {\Cal V}_2^G} = r' + 2s$. Du
fait que $\chi (B
,\omega_B^{\otimes 2}) = 3g' - 3$,  l'\'egalit\'e (4.14) entraine  $r =
r' + 2s$, c'est \`a dire
l'assertion (1).

Maintenant si $l\geq 2$, on a par les m\^emes \'egalit\'es
$$\dim {\Cal V}_l^G = (l-1)r' + \sum_{i=1}^{r'} \left( \langle
\langle {l\over e_i}\rangle \rangle
- {l\over e_i} \right) + \cases  {s(l-1) & si $l$   impair\cr
  sl &  si $l$   pair}  $$
$$= (l-1)r' + \sum_{2\leq e\mid \Card G} N_e \left( \langle \langle
{l\over e}\rangle \rangle
- {l\over e} \right) + \cases { s(l-1)&  si $l$  impair\cr
  sl&    si $l$   pair   }$$
Ces relations   pour $l = 3$ donnent 
$2r - N_2 = 2r' - N'_2 + 2s$, 
ce qui tenant compte de l'assertion (1) conduit \`a $N_2 = N'_2 + 2s$.
Avec $l = 4$, on obtient de la
m\^eme mani\`ere
$3r - N_2 - N_3 = 3r' - N'_2 - N'_3 + 4s, $
d'o\`u on tire l'\'egalit\'e $N_3 = N'_3$.  Si on continue ce
raisonnement, on obtient  de proche en proche l'\'egalit\'e  $N_e =
N'_e$ pour tout $e \geq  3$, c'est \`a dire l'assertion (2). \hfill\break
\qed

Les arguments combinatoires utilis\'es ci-dessus  pour relier
la donn\'ee de Hurwitz aux
repr\'esentations $Hurw_m$  peuvent \^etre  repris sans changement
dans le cas stable. Dans cette situation on peut   encore lire
  la donn\'ee de Hurwitz.
  
\proclaim Proposition 4.13.   Sous les hypoth\`eses du th\'eor\`eme
4.12, les sous-groupes cycliques de $G$ qui sont des stabilisateurs
stricts de points (r\'eguliers ou points doubles) de $C$, sont les
conjugu\'es des sous-groupes  $\{ \Gamma_1,\cdots,\Gamma_r\}$ sp\'ecifi\'es par la donn\'ee de
Hurwitz. 

\dem On reprend la m\'ethode de d\'enombrement du paragraphe 3.2
(Lemme 3.4). Pour tout sous-groupe cyclique $H\subset G, H \ne 1$,
soit
$$C^H = \{ x\in C , H\,\,  \hbox {\rm  fixe  strictement}\,\, x\} \tag (4.23)$$
Posons:
$$f(H) = \cases {\Card {C^H}&   si $e = \Card H \geq 3$ \cr
N'_2(H) + 2 s(H)& si $e = 2$}  $$
Donc $f(H) = {\Card C^H }$,  chaque point double \'etant si $e = 2$
compt\'e deux fois. Il r\'esulte du th\'eor\`eme 4.1 que $f$ est
totalement d\'etermin\'e   par les seules repr\'esentations de Hurwitz. Soit encore  $\Delta_H = \{ x\in C
, H\,\, \hbox{\rm est  un stabilisateur de}\,\, x\}$. Rappelons
que si $x$ est un point double cela signifie  que $e = \Card H = 2$,
et que $H$ est engendr\'e par
une reflexion qui \'echange les branches en $x$. Il d\'ecoule des d\'efinitions que si $H \ne 1$, on a:
$$ C^H = \bigsqcup_{H\subseteq K} \Delta_K  $$
On peut poser $C^H = \Delta_H = \emptyset$ si $H$ n'est pas cyclique.
Si $x \in C^H$ est un point
double, alors $\Card H  = 2$ et $x \in \Delta_H$. Si on pose de
nouveau $g(H) = \Card
{\Delta_H}$, un point double \'etant  toujours compt\'e deux fois,
on retrouve  avec la convention $f(H) = g(H) = 0$ si $H$ n'est pas
cyclique, la relation famili\`ere  (voir (3.12))
$ f(H) = \sum_{H\subseteq K} g(K). $
Par un argument similaire \`a celui utilis\'e dans le lemme 3.4,
 on obtient  par inversion de M\"{o}bius 
$$g(H) = \sum_{H\subseteq K} \mu (H,K) f(K) \tag (4.24) $$
o\`u comme dans le \S \  3.2, $\mu (.,.)$ d\'esigne la fonction
de M\"{o}bius de l'ensemble ordonn\'e des sous-groupes de $G$. La
preuve se poursuit comme dans la section 3.2;  on obtient   que
$\Card \Delta_H$ ne d\'epend finalement que de la donn\'ee de Hurwitz
du probl\`eme modulaire. Le r\'esultat en d\'ecoule.\hfill\break\qed 

\bigskip

 {\it 4.3. Courbes stables marqu\'ees et actions de groupes}

\bigskip

 Parall\`element  \`a la stabilit\'e de la courbe $C$ qui supporte l'action de $G$, on aura  \`a  imposer   la stabilit\'e de la courbe marqu\'ee par un diviseur $G$-invariant\footnote{Dans la section 7, on aura besoin d'un marquage plus pr\'ecis que le seul marquage par un diviseur $G$-invariant. On peut souhaiter   num\'eroter les points du diviseur, i.e piquer la courbe.  Cela doit \^etre fait en coh\'erence avec l'action de $G$, donc on doit prendre en compte   l'holonomie en chaque point $P$ du diviseur, i.e le couple $(H,\chi)$, form\'e du stabilisateur de $P$, et du caract\`ere local de l'action de $H$ sur l'espace cotangent en $P$ (\S \ 2.2.1).  Dans une num\'erotation des points de $D$ les points de m\^eme rang  doivent avoir une m\^eme holonomie.} .  Le diviseur  sera  alors contraint de   contenir le diviseur de ramification.    
Lorsque  le diviseur se r\'eduit au diviseur de ramification, la stabilit\'e de la courbe  marqu\'ee  sera  en fait \'equivalente \`a la stabilit\'e de la courbe sans le marquage, sauf
 en  pr\'esence d'isotropie di\'edrale, ce dysfonctionnement
pouvant \^etre \'elimin\'e par ailleurs.   Stabiliser la courbe par les points de
ramification, reviendra de fait \`a \'eliminer les points fixes de
type III. 

 Dor\'enavant    $C/S$ est une $G$-courbe stable
marqu\'ee par    le diviseur   $G$-invariant  $D$ (D\'efinition 4.2). Pour
compl\'eter  l'\'etude pr\'ec\'edente   limit\'ee au cas stable non
marqu\'e, on doit maintenant
faire une \'etude  \`a la Chevalley-Weil du $G$-faisceau
$\omega_{C/S} (D)^{\otimes m} $,
substitut naturel \`a $\omega_{C/S}^{\otimes m}$. Il s'agit, $s$ \'etant un point g\'eom\'etrique,    d'extraire des informations des $kG$-modules
$H^0 (C_s ,\omega_{C_s} (D_s)^{\otimes m})$. On sait que   le faisceau dualisant modifi\'e $\omega_{C/S}(D)$ 
  poss\`ede la propri\'et\'e d'annulation    (Knudsen \cite {48})
$R^{1}\pi_\star (\omega_{C/S} (D)^{\otimes m}) = 0$ si
 $m\geq 2$.
 De plus si $m\geq 1$,  $\pi_\star (\omega_{C/S} (D)^{\otimes m})$ est localement
libre de rang $(2g - 2 + N)m - (g -
1)$. On montre  aussi que la formation de faisceau
$\omega_{C/S} (D)^{\otimes m}$ commute aux changements de base. Le lemme
suivant d\'ecoule d'arguments d\'eja utilis\'es, la preuve ne sera
donc pas r\'ep\'et\'ee.

\proclaim Lemme 4.14. Si $\pi: C \to S$ est une courbe stable marqu\'ee
par le diviseur $G$-invariant $D = \sum_{i=1}^N P_i$, avec $S$ connexe.
Les $kG$-modules $H^0 (C_s ,\omega_{C_s} (D_s)^{\otimes m})$ sont
ind\'ependants de la fibre g\'eom\'etrique $C_s$, $s\in S$ ($m\geq 1$).

 \qed

Supposons    d'abord   la fibre $C_s$   lisse. La suite
exacte canonique
$$0\longrightarrow \omega_{C_s}^{\otimes m} \longrightarrow
\omega_{C_s}(D_s)^{\otimes m}
\longrightarrow \bigoplus_i N_{P_i(s)}^{\otimes m-1} \longrightarrow 0$$
avec $N_{P_i(s)} = T_{C_s,P_i(s)}^*$, montre que le $G$-type (section
2.1) \'etant suppos\'e connu, alors les repr\'esentations de
Hurwitz $H^0(C_s , \omega_{C_s}^{\otimes m})$ sont  elles aussi
connues. Il en est donc de m\^eme pour la donn\'ee de ramification,
\'evalu\'ee le long des fibres non singuli\`eres. Pour tout
sous-groupe cyclique $H$ de $G$, consid\'erons le diviseur de Cartier
relatif $Fix(H)_{hor}$ (Proposition 4.8). Le corollaire de cette
proposition montre que le degr\'e $d$ de ce diviseur est \'egal \`a
$d = d_1 + 2d_2$ avec
$$ \cases{ d_1 & =  $\Card $   (nombre des points fixes de type $(I) $) \cr
d_2 &= $\Card $ (nombre des points fixes de type  $(III) $)
 \cr  }\tag (4.25)$$
Ce degr\'e coincide avec $f(H)$, fonction introduite dans la
proposition 4.3. De sorte que l'argument combinatoire invoqu\'e
dans cette proposition conduit clairement \`a un r\'esultat
identique. Pour l'\'enoncer, posons $g(H) = g_1(H) + 2g_2(H)$, avec
$$  \cases{ g_1(H)& = $\Card$   (points fixes de type $ (I)$ de
stabilisateur $H$)\cr
   g_2(H) &= $\Card $  (points fixes stricts de $ H $ de type
$(III)$)  \cr}\tag (4.26)$$
de sorte que $g(H)$  est  constant le long des fibres g\'eom\'etriques,  qui plus est,   totalement
d\'etermin\'e par la donn\'ee de Hurwitz, conjointement avec le
$G$-type de $D$. 

Introduisons quelques notations. Fixons une fibre g\'eom\'etrique  que nous noterons pour simplifier $C$,   les points marqu\'es seront not\'es $P_i$, de sorte que $D = \sum_i P_i$. Soit $b_1$  le nombre
d'orbites singuli\`eres (points de type I) contenues dans les points marqu\'es $\{ P_i
\}$ i.e. dans $D$; noter que ce nombre est fix\'epar le $G$-type de $D$, et ind\'ependant de la fibre.
Soit aussi $b_2$ le nombre d'orbites singuli\`eres (points de type $(I)$) disjointes de $D$. Ce nombre d\'epend de la fibre choisie;  on a $b = b_1 + b_2$ si la fibre est non
singuli\`ere. Soit enfin $s$ le nombre d'orbites constitu\'ees de
points de type $(III)$,   \`a isotropie di\'edrale.

 Le r\'esultat suivant   \'etend  le  th\'eor\`eme 4.12
aux courbes stables
marqu\'ees, fournissant ainsi une explication \`a  la coalescence
\'eventuelle des points de ramification, cela en fonction des
conditions initiales,  de marquage, et la donn\'ee de ramification.
 
 Soit $\pi: C \to S$   une courbe   de base connexe,  le lieu de
lissit\'e de $\pi$  \'etant suppos\'e  partout dense, de la sorte une quelconque fibre   est   sp\'ecialisation d'une
fibre non singuli\`ere.   La donn\'ee de ramification est  ainsi d\'efinie. Dire que le $G$-type de $D$ contient la donn\Ž e de ramification a  alors un sens. 
 
\proclaim Theoreme 4.15. Sous les hypoth\`eses pr\'ec\'edentes, donc
le groupe $G$ agissant stablement sur la courbe stable marqu\'ee
$\pi: C \to S$, on  a le long des fibres g\'eom\'etriques l'\'egalit\'e 
$$ b = b_1 + b_2 + 2s \tag (4.27)$$  

\dem Montrons en premier que pour toute fibre, on a l'\'egalit\'e
relative \`a  la fonction  $g$  (4.26):
$$ b = {1\over \Card G} \sum_{H \ne 1} (\Card H)\,\, g(H)  $$
La somme  est limit\'ee aux sous-groupes cycliques. Comme $g(H)$ ne
d\'epend pas de la fibre, il suffit de v\'erifier cette \'egalit\'e
pour une fibre non singuli\`ere, dans ce cas $b$ repr\'esente le
nombre d'orbites singuli\`eres de la fibre. On peut alors \'ecrire $b
= \sum_{H\ne 1} b_H$, avec $b_H$ \'egal au nombre de ces orbites
telles que le stabilisateur d'un de  ses points est dans la classe de
conjugaison de $H$. Soit $\omega$ une telle orbite; le nombre de
points de $\omega$ de stabilisateur $H$ est \'egal \`a ${\Card
{N_G(H)} \over \Card H}$, d'o\`u:
$$g(H) = \Card {\Delta_H}  =  b_H {\Card {N_G(H)} \over \Card H}$$
On peut ainsi \'ecrire (4.27), la somme portant maintenant sur les
classes de conjugaison, notation
$\langle -\rangle $
$$ b  = \sum_{\langle H \rangle \ne 1} g(H) {\Card {N_G(H)} \over
\Card H} = {1\over
\Card G} \sum_{H\ne 1} (\Card H)  g(H) \tag (4.28)$$
la derni\`ere somme portant  elle sur les sous-groupes, et non sur
les classes de conjugaison. Supposons maintenant la fibre singuli\`ere, avec $s$-orbites de type $(III)$. En substituant \`a $g(H)$
l'expression (4.26), on trouve
$$b = {1\over \Card G} \sum_{H \ne 1} (\Card H)  g_1(H) + {2\over
\Card G} \sum_{H\ne 1} (\Card H)  g_2(H) \tag (4.29)$$
Il est clair que le premier terme de la somme de droite est \'egal \`a
  $b_1 + b_2$. Reste \`a prouver que le second, qui ne fait
intervenir que les sous-groupes d'ordre deux, conduit \`a:
$$\Card G . s = 2 \left( \sum_{H, \Card H = 2} g_2(H) \right) $$
Dans ce but, on \'ecrit $s = \sum_{m\geq 1} s_m$, o\`u $s_m$ est le
nombre d'orbites \`a groupe d'isotropie ${\Bbb D}_m$. Soit $\omega$
une orbite  \`a groupe d'isotropie ${\Bbb D}_m$; notons
$g_{2,\omega} (H)$ le nombre de points fixes stricts de $H$ contenus
dans $\omega$. Dans la somme $\sum_{H,\Card H = 2} g_{2,\omega} (H)$
un point est en fait compt\'e $m$ fois, car  ${\Bbb D}_m$  contient
$m$ r\'eflexions. Cette somme vaut donc $ m \Card \omega = m {\Card G
\over 2m} = {\Card G \over 2}$; alors
$$\sum_{H, \Card H = 2} g_2(H) = \sum_{H , m\geq 1} s_m g_{2,\omega_m} (H) \tag (4.30)$$
$\omega_m$ d\'esignant une orbite avec isotropie ${\Bbb D}_m$. En
reportant  (4.30) dans la relation qui
pr\'ec\`ede, on obtient  ais\'emment  le r\'esultat. \hfill\break
\qed 
\proclaim Corollaire 4.16.  Les hypoth\`eses du th\'eor\`eme 4.15 \'etant conserv\'ees, supposons en outre que le $G$-type du marquage
$\{P_i\}$ contienne la donn\'ee de ramification.  Alors  pour toute
fibre g\'eom\'etrique    $s = 0$.   En d'autres termes, il n'y a pas
d'orbite \`a isotropie di\'edrale. 

\dem En effet on a $b = b_1$. \hfill\break
\qed

\section5 {  D\'eformations des rev\^etements mod\'er\'ement
ramifi\'es } 
\beginsection 5.1. D\'eformations \'equivariantes des courbes

\bigskip
{\it 5.2.1.  La d\'eformation universelle}
\bigskip

 On d\'ecrit la d\'eformation formelle  universelle
d'un rev\^etement galoisien
mod\'er\'ement ramifi\'e $\pi: C\to B=C/G$. On  notera
$B$ la courbe quotient pour \'eviter tout conflit avec le diviseur
$D$ qui marque la courbe). Dans ce contexte la th\'eorie des
d\'eformations  se simplifie de mani\`ere substantielle (comparer
avec \cite {12}).  Les seules contributions locales \`a
 la d\'eformation verselle viennent des points doubles et des
points marqu\'es,   points de ramification except\'es, car la
ramification \'etant mod\'er\'ee, ceux-ci  n'apportent aucune
contribution au foncteur des d\'eformations. On comparera  ensuite la
d\'eformation formelle verselle \'equivariante de la courbe marqu\'ee
$C$, \`a la d\'eformation formelle verselle de la base $B$ marqu\'ee
par les images des points exceptionnels de $C$.

On fixe une courbe stable marqu\'ee $(C , \{P_i\} _{i=1,\dots,N})$
munie d'une action du groupe
$G$. Il n'est pas n\'ecessaire \`a ce stade de pr\'eciser la mani\`ere
dont les points marqu\'es  $P_i$ sont trait\'es, c'est-\`a-dire,
ordonn\'es, ordonn\'es  par paquets, ou pas.  Le diviseur $D = \sum_i P_i$ est
suppos\'e  $G$-invariant, contenant     
  les points de ramification. Dans ces conditions on sait qu'une   
  d\'eformation formelle verselle existe,  et qu'elle est universelle du fait de
l'absence d'obstructions. Pour une discussion plus pr\'ecise   
voir Bertin-M\'ezard (\cite {12},\S 3).  

 On peut pr\'esenter le foncteur des d\'eformations de la mani\`ere
suivante. Consid\'erons un rel\`evement de la courbe  marqu\'ee
$(C,\{P_i\})$, avec oubli de l'action de $G$, \`a une
$k$-alg\`ebre locale artinienne $A$ de corps r\'esiduel $k$ (une
$W(k)$-alg\`ebre si la caract\'eristique $p$ de $k$ est positive).
Cela signifie qu'on a une courbe $\Cal C$ au dessus de
$A$,   plate sur $A$, marqu\'ee par des sections que nous noterons
encore $\{P_i\}_{i=1,\dots,N}$. La donn\'ee  contient aussi un
isomorphisme
$$j: C \to {\Cal C}\otimes k \tag (5.1)$$
tel que  $P_i \circ j$ soit le point marqu\'e \/ initial $P_i$ de
$C$, ceci pour tout indice $i, 1\leq i\leq N$. La courbe $C$ \'etant
suppos\'ee stable marqu\'ee, et donc sans automorphisme infinit\'esimal,
si l'action de $G$ se rel\`eve \`a $\Cal C$, ce rel\`evement est
alors unique \cite {20}. Deux rel\`evements sont dits \'equivalents,
s'ils sont isomorphes par un isomorphisme qui rel\`eve l'identit\'e.
On notera finalement $D(A)$ l'ensemble des d\'eformations, i.e.
des classes de rel\`evements de $(C ,\{P_i\})$ \`a $A$;  une classe
\'etant not\'ee
$[({\Cal C} , \{P_i\} , j)]$. On \'etend  comme d'habitude le
foncteur des d\'eformations  \`a tout anneau local noeth\'erien
complet qui est une $k$-alg\`ebre ( resp. une $W(k)$-alg\`ebre). Le
groupe $G$ op\`ere  sur le foncteur des d\'eformations de la mani\`ere
\'evidente suivante.  L'action de $\sigma \in G$ sur $[{\Cal C} ,
\{P_i\} , j]$ est
$$\sigma [{\Cal C} , \{P_i\} , j] =  [{\Cal C} , \{P_i\} , j\circ
\sigma^{-1}] \tag (5.2)$$
Avec ces d\'efinitions, la d\'eformation $[{\Cal C} , \{P_i\} , j]$ est un point fixe de
l'action de $G$ si et seulement si il
existe un  $A$-automorphisme $\Sigma$ de la courbe marqu\'ee  $({\Cal
C} , \{P_i\})$ qui rel\`eve
$\sigma$, c'est \`a dire tel que $\Sigma \circ j = j \circ \sigma$. Ce
rel\`evement est unique  et
d\'efinit ainsi un rel\`evement de l'action de $G$ \`a $({\Cal C} ,
\{P_i\})$. On observe donc que le
foncteur des d\'eformations \'equivariantes $D_G$ est exactement le
foncteur des points fixes  relativement \`a  l'action (5.2) de  
$G$ sur le foncteur des d\'eformations $D$. Il est bien connu  que le
foncteur des d\'eformations
$D$ de la courbe $C$ est effectivement pro-repr\'esentable   et
formellement lisse \cite {20}, \cite {40}.  La d\'eformation
universelle de la courbe point\'ee $({\Cal C} , \{P_i\})$ a donc pour
base le spectre d'une alg\`ebre de s\'eries formelles en $3g - 3 + N$
variables
$$ R_{ver} = W(k) [[ t_1,\dots t_{3g-3+N}]],\quad (g = \hbox {\rm  genre
de } C) \tag (5.3)$$
La th\'eorie cohomologique qui gouverne le foncteur $D$ est
$\Ext^\bullet_{{\Cal O}_C} (\Omega_C^1 , {\Cal O}_C ( - \sum_i P_i))
$. On sait que   $\Ext^2_{{\Cal O}_C} (\Omega_C^1 ,
{\Cal O}_C ( - \sum_i P_i)) = 0 $, et que l'espace tangent
s'identifie canoniquement  \`a
$$ D(k[\epsilon]) =  \Ext^1_{{\Cal O}_C} (\Omega_C^1 , {\Cal O}_C ( -
\sum_i P_i)) \tag (5.4)$$

Dans le cas \'equivariant, le foncteur $D_G$ des classes de d\'eformations \'equivariantes \'etant le
foncteur des points fixes, et comme par ailleurs  l'ordre de $G$ est
inversible dans $k$, la th\'eorie
cohomologique qui gouverne $D_G$ est particuli\`erement simple
(comparer avec   (\cite {12}, \S 3.1)), c'est le groupe des Ext \'equivariants
$$ \Ext^\bullet_{{\Cal O}_C , G} (\Omega_C^1 , {\Cal O}_C (- \sum_i
P_i)) \, \cong\,
\Ext^\bullet_{{\Cal O}_C} , (\Omega_C^1 , {\Cal O}_C (- \sum_i
P_i))^G \tag (5.5)$$
On a  ainsi $ \Ext^2_{{\Cal O}_C , G} (\Omega_C^1 , {\Cal O}_C (-
\sum_i P_i)) = 0$, il n'y a donc pas
d'obstructions au rel\`evement infinit\'esimal, et l'espace tangent est
$$D_G(k[\epsilon]) \, \cong\, \Ext^1_{{\Cal O}_C} (\Omega_C^1 , {\Cal
O}_C (- \sum_i P_i))^G
\tag (5.6)$$
En cons\'equence,   la d\'eformation universelle $G$-\'equivariante
se d\'eduit de la d\'eformation universelle   
  $ {\Cal C} \to \Spec (R) = {\Cal B}$ de $(C , \{P_i\})$, simplement
par restriction au sous sch\'ema
ferm\'e ${\Cal B}^G$ des points fixes, donc   
$ {\Cal C}^G = {\Cal C} \times_{\Cal B} {\Cal B}^G$.

L'absence d'obstruction pour le probl\`eme \'equivariant signifie
que ${\Cal B}^G$ est formellement
lisse, donc est aussi le spectre d'une alg\`ebre de s\'eries
formelles. Pour \^etre  plus pr\'ecis,
rappelons la forme du discriminant  dans la d\'eformation universelle
${\Cal C} \to {\Cal B}$ de $C\to B$ \cite  {20}. Notons pour cela
$x_1,\dots,x_k$ les points doubles de $C$. Par localisation en
$x_i$, on obtient un morphisme
$$  \Ext^1_{{\Cal O}_C}  (\Omega_C^1 , {\Cal O}_C (- \sum_i P_i)) \,
\to \,\prod_{i=1}^k
\Ext^1_{\widehat{\Cal O}_{x_i}} (\widehat{\Omega}_{C,x_i}^1 , \widehat{\Cal O}_{x_i}) \tag (5.7)$$
qui repr\'esente l'application lin\'eaire tangent au morphisme de
localisation  en le $k$-uplet des
points doubles $x_i$. Cette application est surjective, et pour tout $i$
$$ \dim \Ext^1_{\widehat{\Cal O}_{x_i}} (\widehat{\Omega}_{C,x_i}^1 ,
\widehat{\Cal O}_{x_i})  = 1 \tag (5.8)$$
Choisissons des coordonn\'ees $t_1,\dots,t_{3g-3+N}$ sur la base $\Cal B$,
num\'erot\'ees de telle sorte que si $1\leq i\leq k$, $t_i$ soit le
param\`etre induit par la d\'eformation
verselle du point double $x_i$. Le discriminant de la d\'eformation
universelle,  c'est \`a dire le lieu
param\'etrant les fibres singuli\`eres, a donc pour \'equation
$$t_1\dots t_k = 0 \tag (5.9)$$
Si on retourne \`a la situation \'equivariante, on voit que
l'action de $G$ est lissifiable, si et
seulement si la fibre g\'en\'erique de la d\'eformation universelle
${\Cal C}^G \to {\Cal B}^G$ est
lisse, donc si ${\Cal B}^G$ n'est pas inclus dans le lieu
discriminant (5.9). La traduction de
cette condition donne le r\'esultat  suivant (comparer avec Ekedahl
\cite {27}):
\proclaim Th\'eor\`eme 5.1.   L'action de $G$ sur la courbe stable
marqu\'ee $(C ,\{P_i\}_{1\leq i\leq N})$ admet une d\'eformation \'equivariante non singuli\`ere si et seulement si l'action de $G$ est stable
dans le sens de la d\'efinition 4.4. 

\dem Notons   que les points marqu\'es sont sans effet sur
la stabilit\'e de l'action, vu que
celle-ci ne porte que sur l'isotropie aux points doubles. On peut donc
les ignorer de mani\`ere provisoire.
Soit $\Cal I$ l'id\'eal de ${\Cal B}^G$ dans l'anneau $R =
k[[t_1,\dots,t_{3b-3+N}]] = {\Cal O}_{\Cal
B}$, et soit $\Cal M$ l'id\'eal maximal de $R$. Posons par ailleurs
$$W_i = \Ext^1_{\widehat{\Cal O}_{x_i}} (\widehat{\Omega_{C,x_i}^1}
, \widehat{\Cal O}_{x_i})^*\,\,,\,\, V =  \Ext^1_{{\Cal O}_C}
(\Omega_C^1 , {\Cal O}_C (- \sum_i
P_i))^* \tag (5.10)$$
L'application (5.7) transpos\'ee   donne un morphisme $G$-\'equivariant
$$ W = \prod_{i=1}^k W_i \to V \tag (5.11) $$
D\'ecomposons l'ensemble des points doubles en  $G$-orbites
$\{x_1, \dots , x_k\} = \coprod_{j=1}^s G y_j, \quad y_j = x_{i_j}$
de sorte qu'au niveau des $G$-modules,  on a
$\bigoplus_{x_i \in G y_j} W_i\, = \,\Ind_{H_j}^G (W_{i_j}) , \quad
(H_j = G_{y_j})$.
Soit $\chi_j \in \hat H_j$ le caract\`ere attach\'e \`a l'action
de $H_j$ sur $W_{i_j}$. Soit la
d\'ecomposition du $G$-module $V$ en $ V =  V^G \oplus V_G$, les
notations \'etant les notations usuelles, alors l'id\'eal $\Cal I$
est engendr\'e par $V_G$. D\`es lors la condition \`a v\'erifier
est la suivante: pour   
tout indice $i$, $t_i \not\in V_G$, condition \'equivalente \`a
$\Ind_{H_j}^G (W_{i_j})^G \ne 0$, ou
encore par r\'eciprocit\'e de Frobenius, $W_{i_j}^{H_j} \ne 0$.
Finalement la condition s'exprime par
$\chi_j = 1$ pour $1\leq j\leq s$.
Pour  en faire la traduction, il est commode de  passer \`a
l'anneau local compl\'et\'e ${\Cal O} = \widehat{\Cal O}_{y_j}$ du
point double $y_j$. On a la r\'esolution libre de $\Omega^1_{\Cal O}$:
$$0\longrightarrow I/I^2 \longrightarrow {\Cal O}dx \oplus {\Cal O}dy
\longrightarrow \Omega^1_{\Cal
O} \longrightarrow 0$$
o\`u $I = (xy))$, la fl\`eche de gauche \'etant d\'efinie par $xy
\mapsto xdy + ydx  \pmod I $. On
en tire l'identification standard
$$\Ext^1_{\Cal O} ( \Omega^1_{\Cal O} , {\Cal O}) =  \Hom (I/I^2 ,
{\Cal O}) /  \Hom ({\Cal
O}dx \oplus {\Cal O}dy , {\Cal O} ) \tag (5.12)$$
L'application $\phi : I/I^2 \to {\Cal O},\quad \phi(xy) = 1$ est un
\'el\'ement non nul de cet espace
vectoriel. Soit $\sigma \in H = H_j$; si $\sigma$ fixe les branches,
on peut supposer que les
coordonn\'ees sont choisies de telle sorte que
$$\sigma(x) = \alpha x + \dots , \quad \sigma(y) = \beta y + \dots $$
On a $\sigma (\phi) = \alpha \beta \phi$, de sorte que la condition
s'exprime alors par $\alpha \beta
=  \det (d\sigma) = 1$. S'il y a \'echange des branches, on peut dans
ce cas supposer que
$\sigma (x) = \alpha y + \dots, \quad \sigma (y) = \beta x + \dots$, de sorte que 
dans ce cas la condition est encore $\alpha \beta = 1$ , soit $ \det
(d\sigma) = -1$. On
retrouve bien de la sorte  la condition de stabilit\'e de la d\'efinition 4.4. \hfill\break\qed 

Dans la suite, nous noterons ${\Cal C} \to {\Cal B}$ la d\'eformation
universelle \'equivariante de la
courbe stable marqu\'ee $(C ,\{P_i\})$, oubliant \`a partir de
maintenant l'exposant $G$ dans les notations. On pourra supposer que les points $(P_1,\cdots,P_s), \,(s\leq N)$ sont les points de ramification. La base not\'ee ${\Cal
B}$ est  donc le spectre d'une alg\`ebre de s\'eries formelles sur
$k$, ou  $W(k)$ selon les cas. Si le nombre  d'orbites de points
doubles est $d$, alors on peut choisir des variables $t_1,\dots,t_M$,
de sorte  que $t_1,\dots,t_d$ correspondent aux param\`etres  de d\'eformation
des points doubles. L'\'equation du discriminant est alors
$t_1\dots t_d = 0$. Il est facile de
voir que l'espace tangent (5.6) s'ins\`ere dans une suite exacte
$$ \eqalign{ 0\longrightarrow  H^1_G (C , \Theta_C (-\sum_j P_j))&
\longrightarrow \Ext^1_{{\Cal O}_C
, G} (\Omega_C^1 , {\Cal O}_C (-\sum_j P_j)) \longrightarrow \cr
\qquad\longrightarrow\prod_{i=1}^d \Ext^1_{\widehat{\Cal O}_{x_i}}
(\widehat{\Omega}_{C,x_i}^1 ,
\widehat{\Cal O}_{x_i}) \rightarrow 0 \cr }\tag (5.13)$$

\noindent o\`u dans cette suite exacte, on a fait le choix d'un
point $x_i$ dans chaque  orbite
singuli\`ere, et $H^1_G(-)$ d\'esigne la cohomologie \'equivariante
(\cite {12} , \S 3). On note comme d'habitude
$\pi_\star^G$ le fonteur ${\Cal F} \mapsto \pi_\star ({\Cal F})^G$;
on prouve facilement que  
$H^1_G(C , \Theta_C(-\sum_j P_j)) \cong H^1(C/G , \pi_\star^G (\Theta_C(-\sum P_j))),$
notant toujours $\pi: C \to B=C/G$ le morphisme quotient.  On
observera que dans le terme
de droite de la suite exacte (5.13), le stabilisateur $G_{x_j}$ est
omis en exposant, du fait que par
la condition de stabilit\'e, il agit trivialement. On en tire la
dimension  \footnote { Une mani\`ere d'obtenir ce r\'esultat est de calculer $\dim  H^1(C/G , \pi_\star^G (\Theta_C(-\sum P_j)))$. Si $C$ est une courbe pr\Ž stable de composantes normalis\Ž es $C_\alpha \,(1\leq \alpha \leq d)$, et si $i_\alpha: C_\alpha \rightarrow C$ est le morphisme canonique, on  v\Ž rifie facilement  que 
$\theta_C = \bigoplus_\alpha \,i_{\alpha,*} (\theta_{C_\alpha} (-\sum p_j)$, les points $p_j$ \Ž tant les origines des branches sur la composante $C_\alpha$. Cela ram\ ne  le calcul \`a un simple calcul de dimension sur une courbe lisse. Appelons pour toute composante $D_\beta$ de $C/G$, $v_\beta$ le nombre des images des origines des branches de type III, et $l_\beta$ le nombre de points de branchement situ\'es sur $D_\beta$. La conclusion vient de la relation $b = \sum_\beta v_\beta + \sum_\beta l_\beta$, identique \`a  la relation $b = r + 2s$ (Th\'eor\`eme 4.15).}, soit $\dim (R) = M = 3g' - 3 + b + r$, 
o\`u g' est le genre de $C/G$, $b$ est le degr\'e de la donn\'ee
de ramification (d'une d\'eformation lisse), et $r$  est le nombre d'orbites de points
marqu\'es \`a isotropie triviale. Si $r=0$, donc si le diviseur de marquage est exactement le diviseur de ramification (r\'eduit), alors on voit que la dimension de l'espace tangent est identique \`a celle correspondante au probl\`eme sans marquage. Comme application du th\'eor\`eme 5.1,  notons  le r\'esultat suivant:

\proclaim Proposition 5.2.  Soit $A$ un anneau de valuation discr\`ete complet de corps r\'esiduel $k$  alg\'ebriquement clos, et de
corps des fractions $K$. Soit ${\Cal C} \to \Spec A$ une courbe
stable marqu\'ee par les points $\{P_i\}_{i=1}^N$, de fibre g\'en\'erique lisse. On suppose la fibre g\'en\'erique ${\Cal C}_K$ munie,
d'une action du groupe $G$. Alors l'action de $G$ s'\'etend de mani\`ere unique en une action stable sur $({\Cal C} , \{P_i\})$. 

\dem Du fait de l'unicit\'e du mod\`ele stable marqu\'e $\Cal C$
de ${\Cal C}_K$, l'action de $G$ se prolonge de mani\`ere unique,  et
donc $G$ s'\'etend en un groupe d'automorphismes de ${\Cal C} \to
\Spec (A)$. Le seul point non \'evident est que l'action de $G$ sur
la fibre sp\'eciale $C = {\Cal
C}\otimes k$ est stable. Pour le voir, soit $\Sigma \to \Spec (R)$ la
d\'eformation  universelle \'equivariante de la courbe marqu\'ee $C$.
L'absence d'automorphismes infinit\'esimaux justifie l'existence
d'un morphisme unique $u: R \to A$  tel que ${\Cal C} \cong \Sigma
\otimes_R A$. Soit $x = x_i$ un quelconque point double de $C$, alors
l'anneau local compl\'et\'e de $\Cal C$ en $x$ est de la forme
$$ \widehat{\Cal O}_{{\Cal C},x} \,\cong \,   A[[ X , Y]]/{(XY - a)} , \quad (a\in {\Cal
M}_A) \tag (5.14) $$
Les hypoth\`eses en vigueur i.e. ${\Cal C}_K$  lisse en particulier exigent  $a \ne 0$. Il en
d\'ecoule que $x$ doit se d\'eformer dans la d\'eformation universelle
de $C$, donc que le param\`etre $t =t_i$ qui  mesure la d\'eformation
de $x$ dans la d\'eformation universelle non \'equivariante, doit
avoir une restriction  non nulle \`a
$\Spec (R)$. La preuve du th\'eor\`eme 5.1 montre que cette
condition, exprim\'ee en chaque point
double, revient \`a dire que l'action de $G$ est stable.\hfill\break
\qed 

\bigskip

{\it 5.2.2. Stabilit\'e de la courbe quotient} 

\bigskip

On conserve les notations de la section pr\'ec\'edente. Soit   $(C
, \{P_i\}_{1\leq i\leq N})$  une
courbe stable marqu\'ee, \'equip\'ee d'une action stable du groupe
$G$.  En particulier  le
diviseur $D = \sum_i P_i$ est $G$-invariant. On suppose   que
les stabilisateurs des points doubles ne sont pas di\'edraux,  excluant les points
fixes du type (III) (voir 4.3). On verra comment  ramener la
situation g\'en\'erale \`a cette
condition.  Soit le  quotient $ \pi: C \rightarrow B = C/G$, de sorte 
  que $B$ est  
pr\'estable (proposition 4.3),  marqu\'ee par les images $\{Q_j\}_{1\leq j\leq r}$     
 des $\{P_i\}$.  On supposera 
  que les points de ramification, qui sont tous lisses par
hypoth\`ese, sont inclus dans les points marqu\'es.  Il n'y a
  aucune perte de g\'en\'eralit\'e car les points de
ramification sont des points
marqu\'es d'office.  On a donc $b \leq  r$, $b$
\'etant le nombre de points de branchement.  Notons la remarque \'el\'ementaire suivante:

\proclaim Proposition 5.3. \hfill\break
i) Sous les hypoth\`eses pr\'ec\'edentes, la courbe quotient
$(B ,\{Q_j\}_{1\leq j\leq r})$ est stable marqu\'ee.\hfill\break
ii) R\'eciproquement, si $P_1,\dots, P_m$ sont les
points de $C$ d'images $Q_1,\dots, Q_r$, alors la stabilit\'e de
$(B , \{Q_j\}_{1\leq j\leq r})$ implique la stabilit\'e de $(C ,
\{P_i\}_{1\leq i\leq m})$.

\dem i) Soit ${\Cal E} = \{ Q_j\}_{1\leq j\leq r}$, l'ensemble des
points exceptionnels. On doit v\'erifier la condition de stabilit\'e\/  (d\'efinition 4.4), donc si $\Delta$ est une composante de $B$,
$\Delta \cong {\Bbb P}^1$, alors
$$(\#  \,\,\hbox{\rm  points doubles de}\,\, B \in \Delta) \, + \,  \# (\Delta \cap {\Cal E}) \geq  3 \tag (5.15)$$
Soit $\Gamma$ une composante de $C$ d'image $\Delta$. Comme l'action
de $G$ est suppos\'ee sans inversion de branches, $\Gamma$ est non
singuli\`ere. Notons $H$ le stabilisateur de $\Gamma$. Si $Q \in
\delta$ est un point double de $\Sigma$, n\'ecessairement $Q$ est
l'image d'un point double $P \in C$, situ\'e sur $\Gamma$. Si $P'$
est un second point double  d'image $Q$, il y a un \'el\'ement $g\in
G$ tel que $gP = P'$. Mais $\Delta$ \'etant non singuli\`ere, on doit
avoir $g\in H$ et ainsi les points doubles de
$\Sigma$ situ\'es sur $\Gamma$ correspondent bijectivement  aux
$H$-orbites de points doubles de $C$ situ\'es sur $\Gamma$.
Distinguons plusieurs cas:\hfill\break
i) Le genre h de $\Gamma$ est non nul.  
Soit $t \geq 1$ le nombre de points de branchement du rev\^etement
$\Gamma \to   \Delta = \Gamma /H$; notons
$e_1,\dots,e_t$ les indices de ramification respectifs. Si $\Card H
= N$, la formule de
Riemann-Hurwitz donne
$2h - 2 = -2N + N\left( \sum_{i=1}^t (1 - {1\over e_i})\right),$
en particulier ${2h-2 \over N} = t-2-\sum_{i=1}^t {1\over e_i} \geq
0$, donc $t\geq 3$ et la condition de stabilit\'e est
satisfaite.\hfill\break
ii) Le genre $h$ est  nul. 
On peut supposer que, avec les notations pr\'ec\'edentes, $t<3$. Si
$t=2$, alors $H$ est cyclique et
op\`ere sur $\Gamma = {\Bbb P}^1$ avec deux points fixes. La courbe
$C$ \'etant stable marqu\'ee, $\Gamma$ contient au moins un point
parmis les points doubles de $C$, ou  parmis les points $\{P_i\}$
hors les points de ramification. La condition de stabilit\'e est
encore satisfaite. Reste pour conclure \`a examiner les deux  cas:
$\Delta$ est une courbe rationnelle avec un point double, et $\Delta$
est non singuli\`ere de genre un. Dans ces deux cas, la conclusion,
imm\'ediate, est laiss\'ee  sans v\'erification.
\hfill\break
ii) La preuve reprend les m\^emes arguments que i)
i), nous l'omettons.   \hfill\break
 \qed

 Il est n\'ecessaire \`a ce
stade de clarifier le statut des points de ramification parmi les
points
marqu\'es. Rappelons que par convention (voir section 4.2) la courbe
$C$ est marqu\'ee par un
diviseur $G$-invariant,   suppos\'e contenir les points de
ramification.  On peut donc
l\'egitimement demander de quelle mani\`ere ces points participent
effectivement \`a la stabilisation de la courbe. Dans la suite, on
fera syst\'ematiquement l'hypoth\`ese que dans un rev\^etement
$C\rightarrow B$, la courbe  $C$ \'etant marqu\'ee par les points
$P_1,\dots,P_N$,  alors  $P_1,\dots,P_m$ sont les points de
ramification (donc $m\leq N$).

\proclaim Lemme 5.4.   Soit comme ci-dessus une courbe stable marqu\'ee $(C ,\{P_i\}_{1\leq i\leq N})$ supportant une action stable de $G$.
La courbe $(C ,\{P_i\}_{m+1\leq i\leq N})$ (on omet les points de
ramification) est stable marqu\'ee (stable si $N = m$) si et
seulement si  $C$ n'a pas de composante non singuli\`ere rationnelle
$\Gamma$ telle que, notant $\eta$ le nombre des 
$P_i, \,  m+1\leq i\leq N$,  situ\'es sur $\Gamma$,   et $H$ le
stabilisateur de $\Gamma$, alors on n'est pas dans l'une ou l'autre
des  situations suivantes:\hfill\break
 \indent i)   La courbe $\Gamma$ contient un seul point double de $C$, $H
\ne 1$ et $\eta = 0$.
\indent ii)   La courbe $\Gamma$ contient deux points doubles de $C$, $H \ne
1$ est di\'edral, et $\eta =
0$.

\dem Il est clair que la non stabilit\'e est cons\'equence de
l'existence de composantes $\Gamma  =
{\Bbb P}^1,$ telles que si $\mu$ est le nombre de points doubles de
$C$ situ\'es sur $\Gamma$, alors $\mu + \eta \leq  2$. On peut
exclure le cas trivial $\mu = 0$, et supposer $\mu \geq  1$.
\hfill\break
Supposons  $\mu = 1$.  
La stabilit\'e initiale de la courbe permet d'exclure le cas $H =
1, \eta \leq  1$. Si maintenant $H \ne
1$, $H$ est n\'ecessairement cyclique et si $p\in \Gamma$ est le seul
point double de $C$ sur
$\Gamma$, on a $Hp = p$, et  $H$ doit avoir un second point fixe $p'
\in \Gamma$. La stabilit\'e est
mise en d\'efaut que si $\eta = 0$ car dans le cas contraire, $H$ \'etant non trivial, on a $\eta \geq  2$.

Supposons maintenant $\mu = 2$ et $\eta = 0$. Soient $p,q$ les deux
points doubles de $C$ situ\'es sur $\Gamma$. Si $p \not\in Hq$, alors
$Hp=p$, et $Hq=q$, et $H$ est alors cyclique. Ce cas est
exclu  par la stabilit\'e initiale. Si maintenant $q \in Hp$, le
sous-groupe $H$ a une orbite \`a
deux \'el\'ements sur $\Gamma$, il est donc en vertu de la nature des
groupes finis de
transformations homographiques, cyclique ou di\'edral. Si $H$ est
cyclique, alors $\Card H =
2$ par la stabilit\'e de l'action, et ce cas rentre dans le cas
di\'edral. Donc finalement, le
seul cas pour lequel les points de ramification contribuent \`a la
stabilit\'e intiale de la courbe
est $H = {\Bbb D}_l, \,\,l\geq 1, \,\,\mu = 2,$ et $\eta = 0$, ce qui
prouve le lemme. \hfill\break
\qed 


Dans le cas $N = m$, donc si le marquage initial $\{P_i\}$ est
exactement l'ensemble des points de
ramification, alors seul le point 2) du lemme 5.4 peut arriver; dans
ce cas, apr\`es
stabilisation de la courbe \cite {48},\cite {53}, voir aussi \S
6.4.1, on tombe sur un point double \`a isotropie di\'edrale.

\bigskip 
  {\it 5.2.3.  D\'eformations, versus d\'eformations
du diviseur de branchement}  \bigskip

 On souhaite comparer 
  les d\'eformations
universelles respectives de  la courbe $(C ,\{P_i\}_{1\leq i\leq N})$  munie de l'action de $G$,  avec celles de la base $(B
, \{Q_j\}_{1\leq j\leq r})$. Les points doubles de la
base $B$, qui correspondent aux
$G$-orbites de points doubles de $C$, contribuent individuellement
par un param\`etre ind\'ependant, cela dans les deux d\'eformations
universelles. Notons $\{y_1,\dots,y_d\}$ les points doubles de $B$;
soit alors  
$e_i \geq  1$ l'ordre du stabilisateur, cyclique par hypoth\`ese,  d'un point au-dessus de 
$y_i$.  Le passage au quotient par $G$ commute   avec un changement de base arbitraire
(Proposition 4.3), il  est alors clair que le morphisme $\pi: C \to B =
C/G$ induit un morphisme entre les foncteurs de d\'eformations
$$d\pi : D_{(C, G ,\{P_i\})} \longrightarrow D_{(B , \{Q_j\})}  \tag (5.16)$$
On a finalement le r\'esultat important suivant, qui peut \^etre
compar\'e  avec l'interpr\'etation en termes de
log-structures  donn\'ee dans \cite {44}, \cite {55}, \cite {67}. Il
dit en substance que pour avoir un espace de d\'eformations formelles
d'un rev\^etement non galoisien   formellement lisse, il
faut  d\'eformer simultan\'ement la cl\^oture galoisienne.

\proclaim Th\'eor\`eme 5.5. On peut choisir des syst\`emes de coordonn\'ees
  $(t_1,\dots,t_d,\dots \hfill\break \dots,t_{3g'-3+r})$  et
$(\tau_1,\dots,\tau_d, \dots,\tau_{3g'-3+r})$, respectivement sur
les bases des d\'eformations universelles de $(C, G
,\{P_i\})$ et $(B , \{Q_j\})$, de sorte que le morphisme  $d\pi$ est
donn\'e \/ relativement \`a ces coordonn\'ees par
$$d\pi^* : W(k) [[ \tau_1,\dots,\tau_d,\dots,\tau_{3g'-3+r}]]
\longrightarrow W(k)
[[t_1,\dots,t_d,\dots,t_{3g'-3+r})]] \tag (5.17)$$
avec $(d\pi)^*(\tau_i) = t^{e_i}$ si $1\leq i\leq d$, et
$(d\pi)^*(\tau_i) = t_i$ si $i>d$.

\dem On sait, voir la suite exacte (5.13), que la base de la d\'eformation \'equivariante universelle   
admet d'une part une  contribution venant des d\'eformations
localement triviales autour des points doubles, c'est-\`a-dire au
niveau de l'espace tangent le  sous-espace vectoriel  $H_G^1(C , \Theta_C(-\sum_i
P_i))$, et d'autre part une contribution de dimension un apport\'ee par  chaque orbite de points
doubles. L'action du groupe \'etant
stable et sans isotropie di\'edrale, la proposition 4.10 montre que
la premi\`ere contribution
appara\^{\i}t de mani\`ere identique dans le foncteur des d\'eformations de $(B ,\{Q_j\})$.  Seule la
contribution des points doubles diff\`ere dans les deux foncteurs. Au
niveau des espaces
tangents, cela correspond au diagramme suivant:
 $$ \halign {\hfill$#$\hfill & \hfill$#$\hfill& \hfill$#$\hfill& \hfill$#$\hfill&
\hfill$#$\hfill& \hfill$#$\hfill& \hfill$#$\hfill \cr
0 & \rightarrow & H_G^1(C ,\Theta_C(-\sum P_i)) & \rightarrow &
\Ext_{{\Cal O}_C , G}^1
(\Omega_C^1 , {\Cal O}_C (- \sum_i P_i)) & \rightarrow &
\prod_{i=1}^d \Ext_{\widehat {\Cal
O}_{x_i}} ^1( \widehat \Omega_{\widehat {\Cal O}_{x_i}} , \widehat
{{\Cal O}_{x_i}})\cr
& &\downarrow & & \downarrow & & \downarrow \cr
0 & \rightarrow & H_G^1(B ,\Theta_B(-\sum Q_j)) & \rightarrow & \Ext_
{{\Cal O}_B  }^1
(\Omega_B^1 , {\Cal O}_B (- \sum_j Q_j)) & \rightarrow &
\prod_{j=1}^d \Ext_{\widehat {\Cal
O}_{y_j}} ^1(\widehat \Omega_{\widehat {\Cal O}_{y_j}} , \widehat
{{\Cal O}_{y_i}})\cr}$$
\medskip
 Cela  ram\`ene \`a comparer  par le principe local-global (\cite
{12} thm 3.3.4)  les foncteurs de
d\'eformations locaux, c'est-\`a-dire pour un point double.
Supposons donc maintenant  que $x \in C$ un point double d'image $y
\in B$. On suppose que le stabilisateur de $x$  est $H \cong
{\Bbb Z} /{e\Bbb Z}$, et que $\sigma$ en est un g\'en\'erateur,
qui agit le long des branches au point $x$, relativement \`a des
param\`etres convenablement choisis, par
$$\sigma X = \zeta X\quad , \quad\sigma Y = \zeta^{-1} Y \quad $$
$\zeta $ \'etant une racine primitive d'ordre $e$. La d\'eformation
verselle  \'equivariante du point
$x$ est, du fait de la stabilit\'e de l'action   (voir la preuve du
th\'eor\`eme 5.1),  repr\'esent\'ee
par $ {k [[X,Y,\lambda]]}/ {(XY - \lambda)}$ (on remplace $k$
par l'anneau des vecteurs de Witt $W(k)$ dans le cas d'in\'egales
caract\'eristiques). De mani\`ere analogue, la
d\'eformation universelle du point double $y$ est repr\'esent\'ee
par $ {k [[U,V,\mu]]}/
{(UV - \mu)}$, o\`u $U = X^e , V = Y^{e}$. Le morphisme $d\pi$
correspond  alors \`a
$$  {k[[U,V,\mu]]} /{(UV - \lambda)} \longrightarrow  { k
[[X,Y,\lambda]]} /{(XY
- \lambda)} \tag (5.18)$$
avec $d\pi (U) = X^{e} , d\pi (V) = Y^{e} , d\pi (\mu ) =
\lambda^{e}$. Le r\'esultat en d\'ecoule.\hfill\break
 \qed 
 
 \bigskip

\beginsection {5.3. Mod\`ele stable marqu\'e d'un rev\^etement}

\bigskip

Pour un rev\^etement $C \rightarrow B = C/G$ on a essentiellement deux
d\'efinitions de la stabilit\'e. La
premi\`ere d\'efinition revient \`a imposer la stabilit\'e de
$C$, jointe \`a celle de l'action (D\'efinition 4.4). La seconde
plus flexible \`a certains \'egards,
revient \`a demander \`a c\^ot\'e de la stabilit\'e de
l'action, la stabilit\'e de $C$ marqu\'ee  maintenant par un
diviseur  $G$-invariant contenant les points de ramification.  Le cas non galoisien sera trait\'e dans la section 6.6. 

Avec l'une ou l'autre de ces deux d\'efinitions, on est en mesure de formuler un th\'eor\`eme de r\'eduction
stable. Fixons un anneau de valuation discr\`ete complet
$R$, de corps des fractions $K$, et de corps r\'esiduel $k$, suppos\'e
 alg\'ebriquement clos. Donnons nous un rev\^etement $\pi: C_K \to
B_K$, s\'eparable, les courbes \'etant suppos\'ees lisses connexes. Si
la caract\'eristique de $K$ est $p>0$, on suppose que $p$ ne divise pas l'ordre du groupe de monodromie.
Sous ces conditions on s'attend \`a ce que le rev\^etement donn\'e, d\'efini sur le corps
$K$, se prolonge de mani\`ere unique, apr\`es extension de degr\'e
fini \'eventuelle de $R$, en un
rev\^etement stable (resp. stable marqu\'e)\quad $\pi:\,\, {\Cal
C}\, \to \,{\Cal B}$.
On notera cependant,  qu'en g\'en\'eral, on ne peut pas prendre pour
$\Cal C$ et $\Cal B$
les mod\`eles stables respectifs. La clarification viendra de la
remarque que la seule
difficult\'e provient des isotropies di\'edrales (\'eventuelles)
dans la cl\^oture galoisienne.   La d\'efinition suivante  est
essentiellement la notion de stabilit\'e utilis\'ee par
Harris-Mumford pour compactifier le sch\'ema de Hurwitz "classique"
\cite {35}, \cite {40}, \cite {42}.
\proclaim D\'efinition 5.6.  Soit un rev\^etement entre courbes pr\'estables $\pi: C \to B =  C/G$,
galoisien de groupe $G$, au dessus de la base $S$. Il est dit stable,
mieux stable
marqu\'e  (ou HM-stable),   si d'une part l'action de $G$ est
stable (D\'efinition 4.3), et si d'autre part la
courbe  $C/S$ est stable marqu\'ee par les points de ramification. 

Noter que la courbe \'etant marqu\'ee par les points de ramification,
cela signifie en particulier
que les points de ramification sont toujours non singuliers,   il n'y
a donc pas de points \`a isotropie
di\'edrale  (points de type III) (cor 4.16). Sous ces hypoth\`eses on
sait aussi que la courbe quotient $B$
est stable marqu\'ee par les points de branchement (Proposition 5.3).
La variante du th\'eor\`eme de
r\'eduction stable s'\'enonce alors, la stabilit\'e \'etant prise
dans l'un ou l'autre des deux sens (Def 4.4 , Def 5.6):
\proclaim  Proposition 5.7.  Soit $S = \Spec (R)$ le spectre d'un
anneau de valuation discr\`ete complet, de corps des fractions $K$,
et de corps r\'esiduel $k$. Soit donn\'e  un rev\^etement galoisien
$\pi: C_K \to B_K$ entre courbes lisses connexes d\'efinies sur $K$,
et de groupe de Galois $G$. Si $\Card G$ est inversible dans $R$,
alors apr\`es extension finie \'eventuelle de $R$, le rev\^etement
admet un mod\`ele stable unique  $\pi: {\Cal C} \to {\Cal
B}$.

\dem  C'est simplement une reformulation de la  proposition 5.2. En
effet celle-ci  entra\^{\i}ne comme  cons\'equence directe du th\'eor\`eme de r\'eduction stable (resp. stable marqu\'e), que si
on consid\`ere le mod\`ele stable (resp. stable marqu\'e)  $\Cal C$
de $C_K$, alors l'action de $G$ s'\'etend stablement \`a $\Cal C$.\hfill\break
\qed  

\rema {5.1} On  notera qu'on
ne peut, sauf si la donn\'ee de Hurwitz l'interdit, exclure \`a
priori les points
doubles \`a isotropie di\'edrale dans la r\'eduction stable (non
marqu\'ee)  $\Cal C$ d'un rev\^etement. Pour \'eliminer ces points il
faut stabiliser la courbe $\Cal C $ au moyen des points de
ramification, c'est \`a dire d\'eployer ceux qui sont localis\'es
en les points doubles de type III,  et ainsi aboutir au mod\`ele
stable marqu\'e. Illustrons cette proc\'edure bien connue. Soit  $C
={\Cal C}_k$ la fibre sp\'eciale; soit  aussi $P \in C$ un point
double de $C$ d'\'epaisseur $k\geq 1$, donc
$\widehat {\Cal O}_{{\Cal C}_P} \, \cong \,  {R[[X,Y]]}/{(XY -
a)} ,$
avec $a \in {\Cal M}_R$ de valuation $k$. Si $k\geq 2, \,P$ est un
point singulier de la surface
normale ${\Cal C}$, de type $A_{k-1}$. Il est dans ce cas bien connu
qu'on r\'esoud cette singularit\'e  par $\left[ {k\over 2}\right]$
\'eclatements successifs de l'id\'eal maximal (par exemple \cite
{20}).  Le graphe de la fibre exceptionnelle correspondante \'etant
une chaine de $(k-1)$-droites projectives, chaque brin \'etant de
self-intersection $-2$.

\vskip 10pt
\epsfxsize=2truein
\centerline{\epsfbox{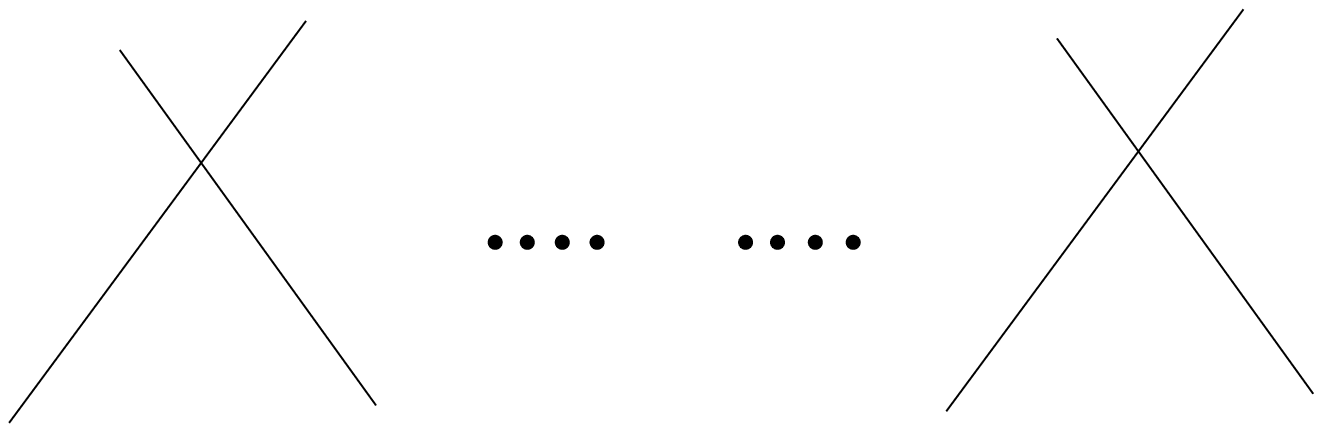}}
\vskip 10pt


Pour rendre  \`a chaque \'etape la proc\'edure d'\'eclatement \'equivariante on \'eclate  en fait chacun des points de l'orbite de
$P$. Soit $H \cong {\Bbb D}_m$ le stabilisateur de $P$, et supposons
ce groupe engendr\'e par $\sigma$ et $\tau$  agissant sur les
branches en $P$ par $\sigma x = \zeta x\,,\, \sigma y = \zeta^{-1} y$
et
$\tau x = y$, pour une certaine racine $m$-i\`eme de l'unit\'e.
Distinguons deux cas selon la parit\'e de $k$.\hfill\break
$\bullet\quad k $ pair\hfill\break
Si $k\geq 4$, l'\'eclatement des points de  l'orbite de $P$ produit
au dessus   de chacun de ces  points, une paire de courbes
rationnelles se coupant en un point de type $A_{k-3}$. Il est facile
de voir que ce point reste \`a isotropie ${\Bbb D}_m$, les
transform\'ees strictes des branches en $P$ \'etant par contre s\'epar\'ees et coupant respectivement les deux composantes du lieu
exceptionnel en deux points \`a isotropie cyclique. Noter que la
stabilit\'e de l'action est pr\'eserv\'ee.

\vskip 10pt
\epsfxsize=2truein
\centerline{\epsfbox{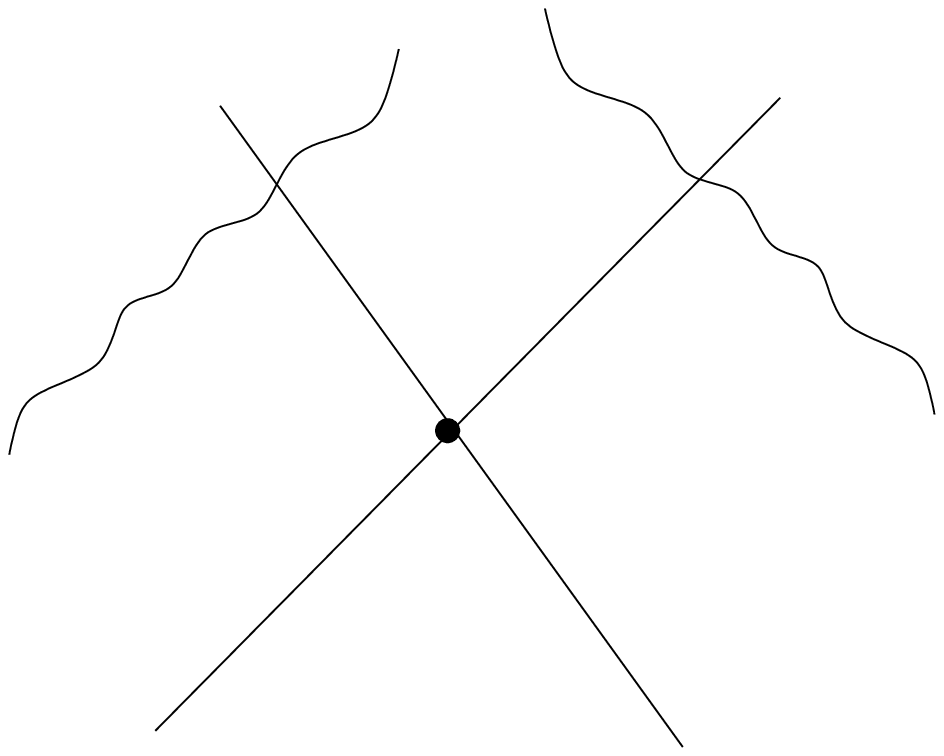}}
\vskip 10pt

Apr\`es ${k\over 2 }-1$ \'eclatements ponctuels on arrive finalement
\`a un point d'\'epaisseur deux, et d'isotropie di\'edrale ${\Bbb
D}_m$. Un dernier \'eclatement produit une courbe $E = \Bbb P^1$, et
r\'esoud la singularit\'e  $A_{k-1}$ initiale.

\vskip 10pt
\epsfxsize=2truein
\centerline{\epsfbox{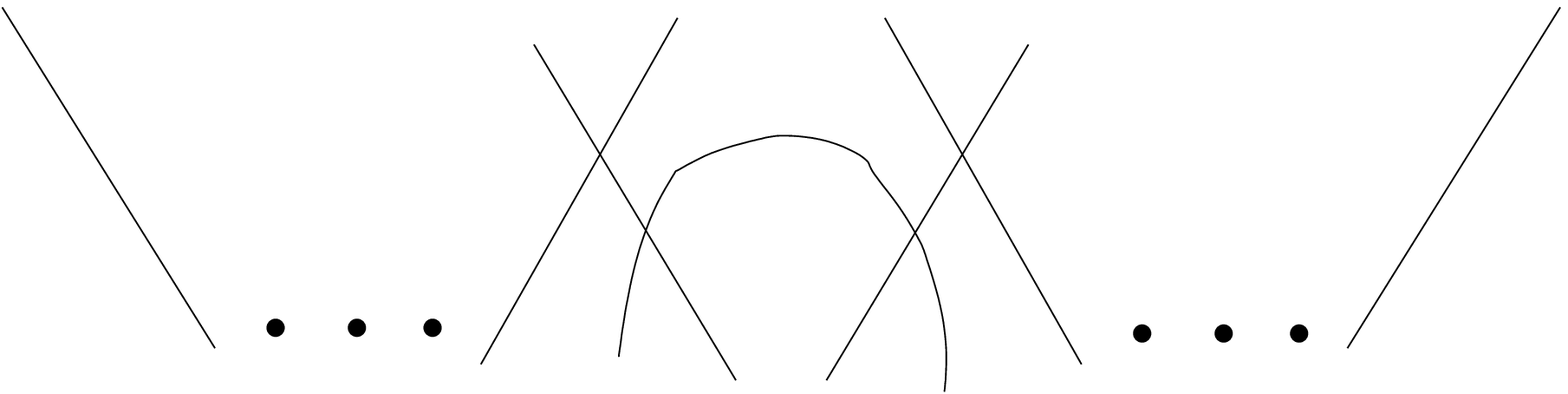}}
\vskip 10pt


Il est imm\'ediat  \`a ce stade que les stabilisateurs sont tous
cycliques, l'action de $H$ sur
la composante $E$, de param\`etre $t$, \'etant l'action standard
$\sigma t = \zeta t\,,\, \tau t =
t^{-1}$.  Notons que la courbe obtenue n'est pas stable marqu\'ee;
pour obtenir cette condition
on doit contracter les deux cha\^{\i}nes de longueur ${k\over 2}-1$
de part et d'autre de $E$,
conduisant \`a la situation

\vskip 10pt
\epsfxsize=2truein
\centerline{\epsfbox{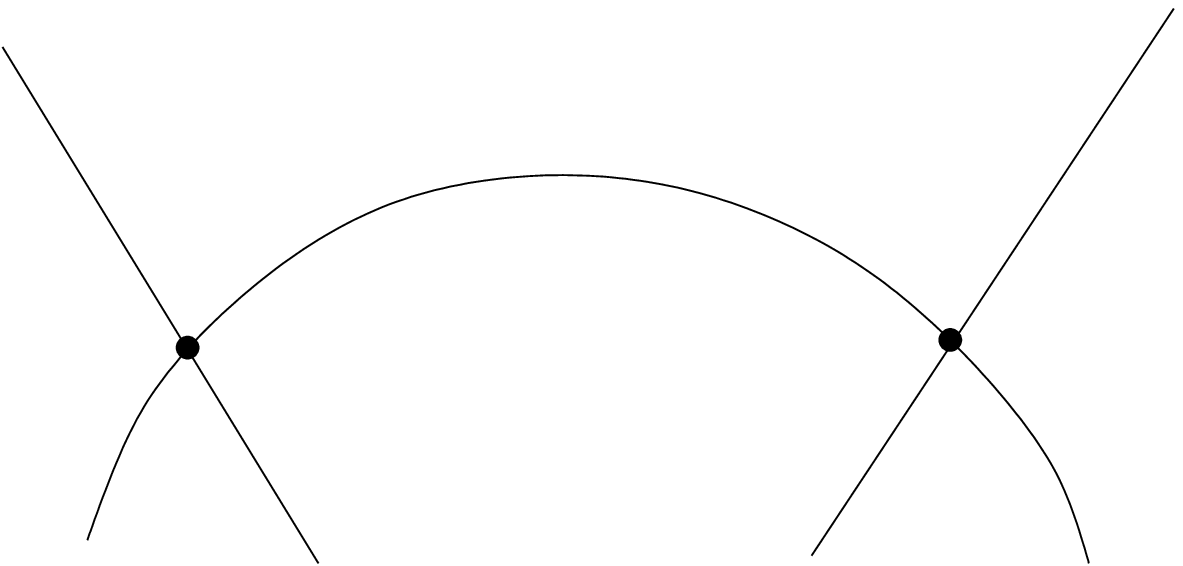}}
\vskip 10pt


o\`u les points $P_1$ et $P_2$ sont des points singuliers de type
$A_{k\over 2}$. Les points
fixes de $H$ sur la composante $E$  maintenant stabilisent la courbe.
\hfill\break
$\bullet \quad k$ impair \hfill\break

Dans ce cas, on r\'esoud le point singulier par une suite de
$[{k\over 2}]$ \'eclatements. On doit
noter que le point d'intersection des deux brins qui proviennent du
dernier \'eclatement est
encore \`a isotropie di\'edrale $H = {\Bbb D}_m$ au contraire du
cas $k$ pair.

\vskip 10pt
\epsfxsize=2truein
\centerline{\epsfbox{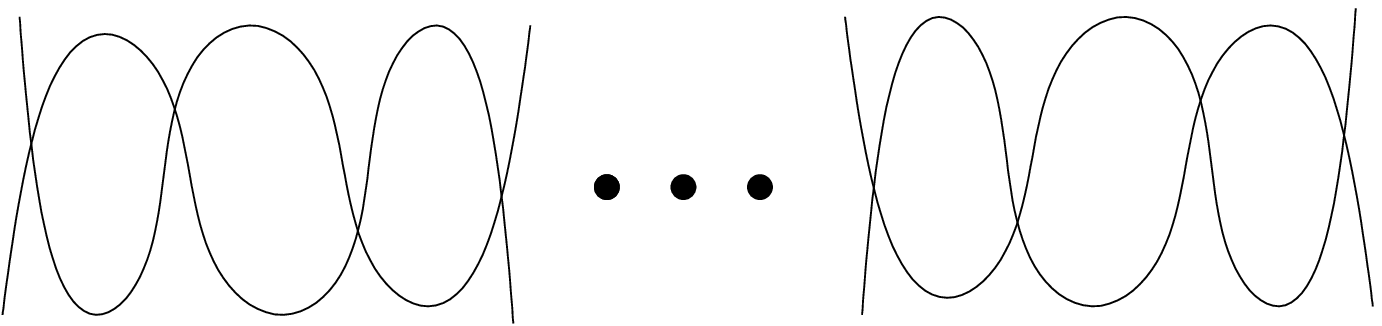}}
\vskip 10pt

  L'\'eclatement de ce dernier point (r\'egulier) conduit \`a une
composante exceptionnelle rationnelle $E$. Supposons le point donn\'e
  par l'\'equation $xy - z = 0$; on peut prendre pour param\`etre
sur $E$, $t = {x\over y}$, l'action de $H$ \'etant alors donn\'ee par
$\sigma t = \zeta^2 t\,,\, \tau t =
t^{-1}$.   La stabilit\'e de l'action impose $m$ impair, ce qui
nous ram\`ene \`a une situation  identique  \`a la pr\'ec\'edente. Par contraction des deux cha\^{\i}nes de longueur $[{k\over
2}]$ de part et d'autre de la courbe centrale, on arrive au mod\`ele
stable marqu\'e.   La construction montre  clairement l'unicit\'e
 du mod\`ele stable marqu\'e. $\lozenge$

\rema {5.2}     On conserve
les hypoth\`eses du th\'eor\`eme 5.7. Soit un rev\^etement galoisien
$ \pi:C_K \to B_K$, d\'efini sur $K$, de groupe de Galois $G$. \hfill\break 
$\bullet$ \quad Il est facile de voir que si $C$ admet un mod\`ele stable
$\Cal C$ sur
$R$, alors $B$ poss\`ede un mod\`ele stable $\Cal B$ sur $R$
et   $\pi$ se prolonge en un morphisme, non n\'ecessairement
fini,  ${\Cal C} \to {\Cal B}$. Soit en effet le mod\`ele stable
${\Cal C'} \to {\Cal B'}$ du rev\^etement  $ \pi:C_K \to B_K$ (Th\'eor\`eme 5.7); la preuve donn\'ee montre en fait qu'un tel mod\`ele
existe sur $R$. La courbe quotient ${\Cal B}'$ est alors stable
marqu\'ee par les points de branchement. Le proc\'ed\'e de
stabilisation de Knudsen \cite {48} conduit \`a un morphisme ${\Cal
B'} \to {\Cal B}$, o\`u ${\Cal B}$ est le mod\`ele stable  de $B$,
donc d\'efini sur $R$. Il est par ailleurs clair que le morphisme
compos\'e ${\Cal C'} \to {\Cal B'} \to {\Cal B}$ factorise par la
stabilisation  ${\Cal C}$ de ${\Cal C'}$. On notera que m\^eme si
${\Cal C'} = {\Cal C}$, en g\'en\'eral ${\Cal B'} \ne {\Cal B}$; l'\'egalit\'e n'est r\'ealis\'ee que si le rev\^etement g\'en\'erique
est \'etale.\hfill\break
$\bullet$ \quad  Sous les m\^emes hypoth\`eses,  
supposons maintenant que $B_K$ poss\`ede un mod\`ele lisse sur $R$
($B_K$ a bonne r\'eduction sur $R$), et soit $D\subset B$ un diviseur
de Cartier
relatif \'etale sur $R$. On suppose que $D_K$ est le diviseur de
branchement de $\pi: C_K \to B_K$.
Alors sous ces conditions, $C_K$ a bonne r\'eduction sur $R$, plus
pr\'ecis\'ement,  ${\Cal C}$ \'etant le mod\`ele lisse de $C_K$ sur
$R$, le rev\^etement g\'en\'erique se prolonge en un rev\^etement
$\pi: {\Cal C} \to {\Cal B}$ de lieu de branchement $D$. On prouve
cela de la mani\`ere suivante: soit une  extension finie galoisienne
$L/K$ de groupe de Galois $I$, relativement \`a laquelle $C_K$, \'equip\'ee de l'action de $G$,  acquiert un mod\`ele stable marqu\'e
${\Cal C}$. Alors la courbe quotient $B =
C/G$ est le mod\`ele stable sur la normalisation de $R$ de la courbe
$B_K$ marqu\'ee par les points de branchement. La lissit\'e de
$B_k$ implique celle de $C_k$. Reste \`a prouver que $C$ est en
fait d\'efinie sur $R$. Pour cela notons que le th\'eor\`eme 5.5
montre que le foncteur des d\'eformations  
$G$-\'equivariantes infinit\'esimales de $C_k$  est isomorphe \`a
celui de $(B_k , D_k)$. Donc si
$\sigma \in I$, le mod\`ele tordu ${\Cal C}^\sigma$ est une d\'eformation formelle de $C_k$ de
quotient ${\Cal B}$ et de discriminant $D$. On a donc un isomorphisme
unique ${\Cal C}^\sigma
\cong \Cal C$, donc en fait une donn\'ee de descente, ce qui permet
de descendre le mod\`ele $\Cal C$ \`a $R$.

On d\'eduit  facilement de cette remarque,  par le raisonnement
indiqu\'e dans \cite {29},  le th\'eor\`eme de  
S.Beckmann  (loc.cit)  qui montre que le corps des modules d'un
rev\^etement $C \to B$ d\'efini sur la cl\^oture alg\'ebrique de $K$, et
dont la base est d\'efinie sur $K$, a bonne r\'eduction sur $K$,  et
est non ramifi\'e en $R$  si les points de
branchement ne coalescent pas dans la fibre sp\'eciale.
 
\bigskip

\section6 { Champs de Hurwitz } 
\bigskip

On d\'efinit le champ de Hurwitz qui  classifie les rev\^etements, galoisiens ou non, entre
courbes  lisses de genres fix\'es, et  \`a monodromie fix\'ee. On  construit le champ compactifi\'e   
par addition des rev\^etements   stables.    Le cas   non galoisien est ramen\'e au cas galoisien par un argument de cl\^oture galoisienne.

Un $G$-rev\^etement galoisien stable $\pi: C\to D$ peut \^etre interpr\'et\'e comme
induisant  sur $D$,   courbe stable marqu\'ee par
les points de branchement, une structure additionnelle. Dans le
langage de Abramovich,Vistoli, et al \cite {1} cela \'equivaut \`a un morphisme $  D \longrightarrow BG$, la courbe $D$
\'etant  munie de sa structure d'orbifold, et  $BG$ \'etant le
champ classifiant de $G$ (\S \ 6.1). Par passage \`a la cl\^oture
galoisienne, on  r\'eduit la construction dans le cas g\'en\'eral
 \`a celle  classifiant les rev\^etements
galoisiens  \`a  groupe de monodromie $G$ fix\'e.   

Une  difficult\'e courante  dans  les  probl\`emes modulaires  est que les objets, 
c'est \`a dire dans le cas pr\'esent les rev\^etements,  ont en
g\'en\'eral des automorphismes permanents, par exemple ceux donn\'es
par le centre de $G$.  Il est  souhaitable de les \'eliminer pour
arriver \`a une d\'efinition raisonnable du champ de Hurwitz. De
cette mani\`ere, le champ de Hurwitz est reli\'e au champ
$\overline {\Cal M}_{g',n} (BG)$ introduit par Abramovich Corti et
Vistoli \cite {1},\cite {2},\cite {3}, par \footnote {On peut prendre
cette description comme une d\'efinition. Dans ce cas les propri\'et\'es de $\overline {\Cal M}_{g',n} (BG)$ d\'ecoulent  
des propri\'et\'es correspondantes du champ de Hurwitz, d\'evelopp\'ees dans cette section.}$$\coprod_{[\xi]} \overline {\Cal
H}_{g,G,\xi}/ / Z(G)  = \overline {\Cal M}_{g',n} (BG)$$ 
   Bien que  des r\'esultats g\'en\'eraux assurent l'existence des espaces
 grossiers  de modules, appel\'es   les
espaces   de Hurwitz,  nous donnerons   
  une construction projective directe de ces espaces grossiers  
  en suivant la construction de
$\overline {M_g}$ par   Gieseker \cite {37}.    La projectivit\'e des espaces de
Hurwitz conduit, si on le d\'esire,  \`a une  preuve uniforme de la
projectivit\'e des
$\overline {M}_{g,n}$, preuve qui se r\'eduit essentiellement \`a
  $\overline {M}_g.$
\bigskip

\beginsection 6.1. $G$-champs et champs quotients

\bigskip

   Pour la commodit\'e du lecteur  on rassemble une s\'erie de d\'efinitions
et constructions sur les champs alg\'ebriques   utilis\'ees
 dans la suite.  Les d\'etails, op\'erations dans une 2-cat\'egorie,  sont dans
\cite {61}.  
   Un champ ${\Cal M}$  est un
champ  de Deligne-Mumford au dessus d'un sch\'ema de base $S$.  On rappelle  que $BG$ d\'esigne  le champ classifiant du groupe  (fini) $G$
\cite {51}.  
\bigskip

  {\it 6.1.1. Quotient d'un champ par un groupe fini}
 
\bigskip

  La d\'efinition d'un $G-$champ   est comme suit  \cite{61}.       Notons
$m: G\times G \to G$  la multiplication de $G$, et  $  e: \Spec k
\to G$ l'\'el\'ement  neutre.
\proclaim  D\'efinition 6.1. i)  Une action de $G$ sur ${\cal M}$ est un morphisme $\mu: G\times {\cal M}$ rendant 2-commutatif le diagramme
$$\commdiag{G\times G\times  {\cal M}&\mapright^{m\times 1}&G\times  {\cal M}\cr
\mapdown\lft{1\times \mu} &&\mapdown\mu\cr
G\times {\cal M}&\mapright^\mu&{\cal M}\cr} \tag (6.1)$$
  En d'autres termes il existe un 2-morphisme (isomorphisme fonctoriel) $ \alpha: \mu . (m\times 1) \buildrel\sim\over \to \mu.(1\times \mu)$, satisfaisant \`a une relation d'associativit\'e d'ordre sup\'erieur, que nous n'\'ecrirons pas (voir \cite{61}, Definition 1.3 i)).  L'action est stricte si  $\alpha = 1$. \hfill\break
  ii) Un morphisme de $G$-champs $({\cal M},\mu,\alpha) \to ({\cal M'},\mu',\alpha')$ est une paire $(f,\sigma)$, avec $f: {\cal M}\to {\cal M'}$   un morphisme de champs, et $\sigma: \mu'.(1\times f) \to f.\mu$ un 2-morphisme v\'erifiant une relation de compatibilit\'e  avec $\alpha$ et $\alpha'$ (loc.cit iii).

On montre que tout $G$-champ est \'equivalent \`a un $G$-champ ''strict'' \cite{61}, ce qui permet de se limiter en principe  \`a des actions strictes. Dans la suite cependant certaines constructions conduiront   \`a des actions non strictes, qui seront rendues strictes par ''strictification'' (loc.cit proposition 1.5).
Dans la suite l'indice sera  souvent omis, on \'ecrira si $x\in {\Cal M} (U), gx$ pour l'image  $\mu(g,x)$, et aussi $gf: gx\buildrel\sim\over\rightarrow gx'$ pour
l'image de $f: x \buildrel\sim\over\rightarrow x'$ par $g$. Noter que
par  d\'efinition  $f^* (gx) = g f^* (x)  $.   
  
 
\proclaim D\'efinition 6.2. Soient  $\Cal M, \;\Cal N$  deux
$G$-champs.  Un $G$-isomorphisme
  $\phi: \Cal M \longrightarrow \Cal N$,
est  un $G$-morphisme $f: \Cal M \longrightarrow \Cal N$
qui est un isomorphisme de champs (\'equivalence de cat\'egories).  Il existe un quasi-inverse qui est un $G$-morphisme.

Une action triviale est une action \'equivalente \`a l'action  strictement triviale $\mu = {\rm pr_2}$.
On rencontrera la situation suivante que nous d\'ecrivons sans d\'etail. Soit ${\cal M}$ un $G$-champ, et soit $H \triangleleft G$ un sous-groupe distingu\'e. Si la restriction de l'action \`a $H$ est ''triviale'', la trivialisation \'etant compatible \`a l'action de $G$ sur $H$ par automorphismes int\'erieurs, il en d\'ecoule une action de $G/H$ sur ${\cal M}$.  On peut par exemple passer par une  transversale  de $G$ modulo $H$. On  notera   que l'action construite peut ne pas \^etre stricte.  

De mani\`ere  plus explicite  si $f: {\cal M}\to {\cal M'}$ est un morphisme de $G$-champs, le  2-morphisme $\sigma$ est contraint par la relation de compatibilit\'e  $g,g'\in G$:
$$ f(\alpha^x_{g,g'}) . \sigma^x_{g.g'} = \sigma^{g'x}_g . g\sigma^x_{g'}.{\alpha'}^x_{g,g'} \tag (6.2)$$  
 Si $X$ est un $G$-sch\'ema, un $G$-morphisme $X
\rightarrow \Cal M$ est  d\'etermin\'e par un objet $x \in \Cal M (X)$ $G$-lin\'earis\'e,  donc \'equip\'e pour tout
$g\in G$ d'un isomorphisme 
$\alpha_g (x): gx \buildrel\sim\over\longrightarrow g^* (x),  $
les $\alpha_g(x)$ \'etant contraints par une relation de cocycle
$$ \alpha_{hg}(x)  = g^* (\alpha_h (x)) . h(\alpha_g (x))\tag (6.3) $$
     Il est ais\'e de voir que si  $p: \Cal M
\rightarrow {\Cal Y}, \, q: {\Cal N}
\rightarrow {\Cal Y}$ sont $G$-\'equivariants, alors le produit fibr\'e
 ${\Cal M}\times_{p,{\Cal Y},q}
\Cal N$ supporte une action canonique de $G$. 
   Les $G$-champs forment une 2-cat\'egorie. 
   
    Le  champ quotient $[{\Cal M} / G]$ du $G$-champ $\Cal M$ par l'action de $G$ est caract\'eris\'e   par une propri\'et\'e universelle  dans la 2-cat\'egorie des $G$-champs (\cite{61}, D\'efinition  2.3).

\proclaim D\'efinition 6.3.  Le champ quotient $[{\Cal M} / G]$ est caract\'eris\'e par le fait qu'il existe un $G$-morphisme (pour l'action triviale de $G$ sur $[{\Cal M} / G]$)
 ${\Cal M} \rightarrow [{\Cal M} / G] $
qui au sens des 2-cat\'egories est universel pour les morphismes de but un champ avec action triviale de $G$.  On montre (sous des hypoth\`eses plus g\'en\'erales) que le quotient $[{\cal M}/G]$  existe en tant que champ alg\'ebrique.

Si on se r\'ef\`ere au cas   ${\Cal M}$   repr\'esentable, on peut interpr\'eter $[{\cal M}/G]$ comme un champ de $G$-torseurs. 
 Les  sections de   $[{\Cal M} / G]$   
au dessus de $U$ sont les diagrammes
 $$  {\Cal M}  \mapleft^f  E  \mapright^p  U   \tag (6.4) $$  
  o\`u  $p: E \rightarrow U$ est un $G$-torseur, et le morphisme $f: E
\rightarrow {\Cal M}$ est   $G$-\'equivariant (D\'efinition 6.2). Les morphismes sont les morphismes  de diagrammes,
avec comme stipul\'e, la 2-commutativit\'e.  Un morphisme est repr\'esent\'e par un carr\'e cart\'esien de $G$-torseurs
  $$ \commdiag{E  & \mapright^h &   E'   \cr
   \mapdown\lft{p}& &  \mapdown\lft{p'}\cr
  U  &\mapright^{u}& U'\cr }  $$
 Les $G$-morphismes companions $f: E \rightarrow \Cal
M$, et $f': E' \rightarrow \Cal
M$  \'etant augment\'es d'un 2-morphisme  $\Phi: f \buildrel \sim\over \longrightarrow f'h.  $
Les objets qui d\'efinisent $f$ et $f'$ \'etant $x\in \Cal M (E)$ et
$x' \in \Cal M (E')$, le
2-morphisme $\Phi$ est d\'etermin\'e par $\phi: x \buildrel
\sim\over \rightarrow h^* (x')$. La
compatibilit\'e aux actions de $G$ se traduit pour tout $g\in G$,
par la relation de cobord
$$ g^* (\phi) . \alpha_g (x) = h^* (\alpha_g (x')). g\phi \tag (6.5)$$
 Pour r\'esumer   (\cite{61} Th\'eor\`eme 4.1):
\proclaim  Th\'eor\`eme 6.4. i)  La cat\'egorie fibr\'ee
en groupo\"{\i}des $[{\Cal M} / G]$ est un champ alg\'ebrique, de
Deligne-Mumford si ${\Cal M}$ l'est.  \hfill\break
ii) Le morphisme  ${\Cal M}
\longrightarrow [{\Cal M} / G]$ est un $G$-torseur, qui g\'en\`ere 
  un carr\'e 2-cart\'esien  
  $$ \commdiag{
  S   &\mapright& BG  \cr
  \mapup&&\mapup \cr
    {\Cal M} & \mapright& [{\Cal M} / G]   } \tag (6.6)$$ 
  iii) Le morphisme ${\Cal M} \rightarrow [{\Cal M}/G]$ est un quotient cat\'egorique universel (au sens des 2-cat\'egories).\hfill\break
  iii) Si $\Cal M$ a un espace grossier des modules $M$, $G$ agit sur $M$, et
l'espace alg\'ebrique $M/G$ est un espace grossier des modules pour
$[{\Cal M}/G]$.  \hfill\break\qed

Dans   (6.6) le morphisme $[{\Cal M }/ G ]\rightarrow BG$
est celui qui \`a un objet
$(E\rightarrow U , f)$ du champ $[{\Cal M }/ G]$ associe le $G$-torseur
$E\rightarrow U$ objet de $BG$
au dessus de $U$. Par exemple  pour l'action triviale de $G$ sur
$\Cal M$, le champ quotient est
${\Cal M} \times BG$. Soit $H$ un sous-groupe distingu\'e de $G$;
il y a une action ''naturelle'' non stricte en g\'en\'eral
  de $G/H$ sur $[{\Cal M} / H]$, d'o\`u on tire  une \'equivalence \cite{61} 
$$\left[\left([{\Cal M }/ H]\right) / (G/H)\right] \buildrel\sim\over\rightarrow  [{\Cal M }/ G] \tag (6.7)$$
Du fait de la propri\'et\'e universelle dont jouit le quotient, si $\Phi: \Cal M \longrightarrow \Cal N$ est un $G$-morphisme, alors
il induit par passage au
quotient, un morphisme $\overline \Phi: [{\Cal M }/ G] \longrightarrow
{\Cal N }/ G$, le carr\'e suivant
\'etant 2-cart\'esien
$$ \commdiag{  \Cal M & \mapright^{\Phi}&  \Cal N  \cr
  \mapdown&&\mapdown\cr
[{\Cal M} / G]  &\mapright^{\overline \Phi} &   [{\Cal N} / G] }\tag (6.8) $$

Observons d'autre part que si $P \rightarrow \Cal M$ est un atlas,
alors le morphisme compos\'e
$P\rightarrow  \Cal M \rightarrow [{\Cal M} / G]$ d\'efinit un atlas de
${\Cal M} / G$. En sens inverse si $Q
\rightarrow [{\Cal M} / G]$ est un atlas de $[{\Cal M} / G]$, alors la projection
$q: Q \times_{{\Cal M} / G} \Cal M \rightarrow \Cal M$
  d\'efinit un $G$-atlas, i.e. un atlas  muni d'une action libre de
$G$, le morphisme
$q$ \'etant  $G$-\'equivariant.  

 La cat\'egorie des faisceaux coh\'erents sur le
champ quotient $[{\Cal M} / G]$  s'identifie canoniquement  \`a la cat\'egorie des
$G$-faisceaux coh\'erents sur $\Cal M$. Un $G$-faisceau coh\'erent
sur $\Cal M$  se d\'efinit comme la donn\'ee pour tout
$G$-atlas  $P \rightarrow \Cal M$, d'un $G$-faisceau $\Cal F$ sur
$P$, avec la compatibilit\'e usuelle si $P \rightarrow
P'\rightarrow \Cal M$ est un  $G$-atlas qui domine un autre
$G$-atlas.

\rema{6.1}  Supposons l'action de   $G$  sur $\Cal M$   libre, c'est \`a dire libre sur les objets ($gx
\ne x$ si $g\ne 1$),  et donc sur les fl\`eches. Sous cette hypoth\`ese  on d\'ecrit le champ
quotient comme suit.  Soit le pr\'echamp dont les objets sont les
classes modulo $G$ d'objets de $\Cal M$. Si $\overline x , \overline
y$ sont deux tels objets, un morphisme $\overline x \to \overline y$
est repr\'esent\'e par une  classe de $\sqcup_{s,t\in G} \Hom (sx ,
ty) / G$. La
composition se d\'efinit de mani\`ere  \'evidente. 
Cela   d\'efinit  un pr\'efaiseau, de sorte que le faisceau des morphismes de
$\overline x$ vers $\overline y$
est le faisceau  associ\'e, c'est \`a  dire le faisceau somme
$$\Hom (\overline x , \overline y) \,=\, \sqcup_{s,t\in G} \hom (sx ,
ty) / G = \sqcup \Hom_{s\in
G} \Hom (x , sy)\tag (6.9)$$
 Il y a  une  obstruction \`a l'effectivit\'e d'une  donn\'ee de descente d'objets
au dessus de $S$,  qui r\'eside dans
$H^1(S,G)$.   En effet, soit  $(\overline x_i , \sigma_{i,j})$ un
objet donn\'e localement; il en d\'ecoule un \'el\'ement
$\theta_{i,j} \in G$ tel que (en notation exponentielle)
$\sigma_{j,i}: x_i\buildrel\sim\over\rightarrow {x_j}^{
\theta_{i,j}}$. La condition de cocycle se r\'esumant \`a
$$\theta_{i,j}\sigma_{k,j} . \sigma_{j,i} = \sigma_{k,i}, \,\,\,\,
\theta_{i,j}\theta_{j,k} = \theta_{i,k}\tag (6.10)$$
La classe
d'obstruction est la classe $[\theta_{i,j}] \in H^1(S , G)$.  Le
champ quotient $[{\Cal M} /G]$ s'identifie au champ
associ\'e \`a ce
pr\'echamp \cite {51}. L'identification des deux constructions est
imm\'ediate.\bigskip

  {\it 6.1.3. 2-quotient  $\Cal M //{\rm  G}$ } 
\bigskip
  La   seconde op\'eration  que nous aurons \`a utiliser est    le 2-quotient  (ou $G$-rigidification)  ${\Cal M} // G$,  notation distinguant cette op\'eration  du 1-quotient d\'efini au-dessus \cite{2}, \cite{61}.      Dans cette construction on suppose   que
pour tout objet  $x$ de $\Cal M$ on a une inclusion ($G$ un groupe fini)
$$i_x: G \hookrightarrow \Aut (x) \tag (6.11)$$
o\`u par commodit\'e on notera $g$ \`a la place de $i_x(g)$. On
impose la relation de compatibilit\'e  aux diagrammes cart\'esiens
$f^* (g) = g $, cela pour tout morphisme $f: T \rightarrow S$, et
tout $x \in \Cal M (S)$. Pour rendre cette notion op\'erationnelle,
on doit exiger  que pour tout couple
d'objets $x,y \in \Cal M (U)$ l'ensemble $\Hom (x,y)$, sur lequel
$G$ op\`ere des deux c\^ot\'es, est un $G$-bi-ensemble normal, donc  si $u\in \Hom(x,y), \, g\in G, \,\, u^{-1}Gu  = G$.
 Il est alors clair que la composition des morphismes
 passe au quotient par $G$.  On a (\cite{61}, th\'eor\`eme 5.1):
\proclaim  Proposition 6.5.  Soit $\Cal M$ un $S$-champ alg\'ebrique,
et soit $G$ un groupe fini qui
admet un plongement $G \hookrightarrow \Aut (x)$ pour tout objet $x$,
avec les conditions de dessus. Alors il
existe un champ alg\'ebrique not\'e ${\Cal M} // G$ et un morphisme
$f: {\Cal M} \longrightarrow {\Cal M} // G$, avec la propri\'et\'e
universelle suivante:  tout morphisme de
source $\Cal M$  qui envoie les \'el\'ements de $G$ sur l'identit\'e factorise par ${\Cal M} // G$ .   \hfill\break
i)  Si $\Cal M$ est de Deligne-Mumford, alors ${\Cal M} // G$ l'est
aussi, et le    morphisme ${\Cal M} \rightarrow
{\Cal M} // G$ est une gerbe   \'etale de lien $G$,  de degr\'e $ 1\over {\vert G\vert }$ ($\Cal M$ \'etant   int\`egre).\hfill\break  
ii) Si $\Cal M$
est propre, ${\Cal M} // G$ l'est \'egalement, et si $M$ est un espace
des modules grossier de $\Cal M$, c'est aussi l'espace des modules  de ${\Cal M} // G$. \hfill\break
 \qed

\rema{6.2}  On peut penser \`a un  objet  $x\in ({\Cal M}//G)(S)$ comme la donn\'ee d'un recouvrement
$(S_i\to S)_i$, et pour tout $i$ d'un objet  $x_i\in {\Cal M}(S_i)$, avec en
sus une donn\'ee de descente 
$ \sigma_{j,i}:{ x_i}_{\vert S_{i,j}}
\buildrel\sim\over \longrightarrow {x_j}_{\vert S_{i,j}} $
 \`a un \'el\'ement de $G$-pr\`es, donc tel que  pour tout triplet $i,j,k$,
$${\sigma_{k,i}^{-1}\sigma_{k,j}\sigma_{j,i}}_{\vert S_{i,j,k}}  \in
G \tag (6.12)$$ 
 Si $G\ne 1$, le morphisme  $ {\Cal M} \rightarrow {\Cal M}
// G$    
n'est pas repr\'esentable. Par exemple ${\rm B}G//G = \Spec \, k$.  
Soit $H\triangleleft G$. Si $H$ agit trivialement sur $\cal M$, alors on peut former  d'une part le double quotient $[{\cal M}/G]//H$, et d'autre part $[{\cal M}/(G/H)]$. Il n'est pas difficile d'observer  l'identification
$$ [{\cal M}/G]//H \cong [{\cal M}/(G/H)]\tag (6.13)$$
En fait $[{\cal M}/G] = [[{\cal M}/H]/(G/H)]$ (\cite{61}, remark 2.3), mais $[{\cal M}/H] = {\cal M}\times {\rm B}H$, donc $[{\cal M}/G]//H = [{\cal M}/(G/H)]\times {\rm B}H//H = [{\cal M}/(G/H)]$.
\bigskip

\beginsection {6.2. Champs de Hurwitz }

\bigskip

  Un champ de Hurwitz classifie les
rev\^etements  $\pi: C \to D$ entre courbes
projectives lisses de genres respectifs, $g$ et $h$ fix\'es, et de
degr\'e fix\'e $d$. La relation de Riemann-Hurwitz relie $g,h$ et $d$
par
$$2g -2 = d(2h -2) + B\tag (6.14)$$
  o\`u $B$ est le degr\'e du diviseur de ramification.  Si on fixe
$D$,  et  les $b$ points de branchement, il n'y a qu'un nombre fini de rev\^etements
de degr\'e  fix\'e, donc  la dimension  du champ de Hurwitz est  $3h -3 + b$.  Si  $B > 0$,
on voit que cette dimension est maximum lorsque $B$ est minimum, donc
si la ramification est simple. Si $h = 0$, i.e. $D = {\Bbb P}_1$,
c'est la situation classique \'etudi\'ee par Fulton \cite {35}, et
plus rec\'emment par Harris-Mumford \cite {42},  Mochizuki
\cite {55} et Wewers \cite {68}, \cite {69}. \bigskip

  {\it 6.2.1.
Champs de rev\^etements: cas galoisien} 
\bigskip

  Supposons    les rev\^etements   galoisiens
de groupe $G$,  le cas g\'en\'eral non galoisien  sera abord\'e
 en fin de paragrahe. Comme cela a \'et\'e d\'efini dans les
sections 2.1 et 2.2,  on fixe une donn\'ee de
ramification $\xi  = \sum _{i=1}^r b_i [H_i , \chi_i]$
 de degr\'e $b = \sum
b_i$, $b$ est  le nombre de points de branchement.  Le degr\'e du diviseur de branchement est alors
$$B = {\Card G } \left( \sum_i b_i \left(1 - {1 \over \Card
{H_i}}\right) \right) \tag (6.15)$$
 Nous allons construire le champ de Hurwitz ${\Cal
H}_{h,G,\xi}$, puis d\'ecrire sa
compactification stable $\overline {\Cal H}_{h,G,\xi}$. Le r\'esultat
est un  champ non  connexe
en g\'en\'eral, mais lisse sur $\Bbb Z [  { 1\over {\Card G}} ]$. Il en
r\'esultera   que les composantes connexes  sont identiques aux  composantes irr\'eductibles. Le nombre de composantes
connexes de $\Cal H_{h,G,\xi},$ qui est aussi celui de
$\overline {\Cal H}_{h,G,\xi},$ est   le  nombre de
Nielsen $h(\xi)$ introduit dans  la d\'efinition 2.7. 

Rappelons que
la d\'efinition retenue pour classifier les rev\^etements en g\'en\'eral (galoisiens
ou non), par opposition avec les $G$-rev\^etements est  l'\'equivalence (2.1). 
  Fixons dor\'enavant
le groupe   $G$, ainsi qu'une donn\'ee de Hurwitz aff\'erente
$\xi\in R_+(G)$. La d\Ž finition suivante   sera ult\'erieurement renforc\Ž e par le marquage au moyen du diviseur de ramification.
\proclaim D\'efinition 6.7. i)  Le champ de Hurwitz
${\Cal H}_{g,G,\xi}$ a pour objets au dessus de $S$ les $G$-rev\^etements
$\pi: C\to D\to S, \quad S\in Sch[{1\over \Card G}]$, les
courbes $C$ et $D$ \'etant lisses, propres sur $S$, \`a fibres g\'eom\'etriques connexes de genre $g$, et pour morphismes les
isomorphismes $(f,h) : (C\to D) \longrightarrow
(C'\to D')$, i.e. les diagrammes cart\'esiens de la d\'efinition 6.1.\hfill\break
ii)   Le champ de Hurwitz  $\overline{\Cal H}_{g,G,\xi}$ a pour
objets au dessus de $S$
les $G$-rev\^etements $\pi: C\to D\to S,$ o\`u \/ $C\to S$ est une
$G$-courbe stable, donc une courbe stable munie d'une action stable
de $G$ (D\'efinition 4.4); les morphismes \'etant d\'efinis  de la
m\^eme mani\`ere que dans 1).

 Notons que  ${\Cal H}_{g,G,\xi}$ (resp.$\overline{\Cal
H}_{g,G,\xi})$
n'est a priori   qu'une cat\'egorie fibr\'ee en
groupo\"\i des  au dessus de
$Sch[{1\over \Card G}]$. Par des arguments standards    ce groupo\"\i de est  
un champ de Deligne-Mumford lisse \cite {19} , \cite {54}, et que  le
champ ${\Cal H}_{g,G,\xi}$  est un sous champ ouvert de
$\overline{\Cal H}_{g,G,\xi}$, ce dernier \'etant propre sur
$Sch[{1\over \Card G}]$.   Noter aussi que pour un $G$-rev\^etement
$\pi: C\to D$, le groupe des automorphismes $\Aut (\pi)$ est
\'egal au centre $Z(G)$. D'une mani\`ere un peu diff\'erente, on a
observ\'e que la consid\'eration d'un $G$-rev\^etement $\pi: C\to
D$, se r\'eduit \`a la seule donn\'ee de la courbe $C$ munie de
l'action de $G$. En d'autres termes, le foncteur  $C \mapsto (C \to
C/G)$ \'etablit une \'equivalence entre le champ dont les objets
au-dessus de $S$, sont les $S$- courbes  lisses (resp. stables)
munies d'une action de $G$ (resp. action stable) de donn\'ee $\xi$ le
long des fibre g\'eom\'etriques, et le champ ${\Cal H}_{g,G,\xi}$
(resp. $\overline {\Cal H}_{g,G,\xi}$). Noter que fixer la donn\'ee
de ramification \'equivaut, comme  cela a \'et\'e v\'erifi\'e
dans le th\'eor\`eme 3.1, \`a fixer les repr\'esentations de Hurwitz, d\'efinition plus commode que la d\'efinition initiale, du moins dans le
cadre des familles  de courbes stables. Les morphismes sont dans
cette interpr\'etation les isomorphismes $G$-\'equivariants, et
l'isomorphisme inverse est simplement l'oubli de la base $D \cong
C/G$. Dans la suite  ${\Cal H}_{g,G,\xi}$ d\'esignera l'une ou
l'autre des deux d\'efinitions \'equivalentes pr\'ec\'edentes.  
 
\bigskip
 {\it 6.2.2.  Fonctorialit\'e des champs de
Hurwitz} \bigskip

Un champ de Hurwitz est la source, ou le but, de foncteurs naturels
que nous allons pr\'eciser. Soit
$H\subset G$ un sous-groupe de $G$. Les op\'erations (2.2) et (2.3)
sur les donn\'ees de Hurwitz d\'efinies dans la section  2.2
apparaissent dans le pr\'esent contexte comme  des morphismes de
champs. Par exemple, la restriction de l'action de $G$ au sous-groupe
$H$, soit
$(C\to S, G) \mapsto (C\to S , H)$, d\'efinit un morphisme
$\overline{\Cal H}_{g,G,\xi} \longrightarrow \overline{\Cal
H}_{g,H,{\Res}_H^G (\xi)}$.
Un second morphisme   correspond au
passage au quotient par le
sous-groupe  $H$ (si $H \triangleleft\,  G$), op\'eration justifi\'ee
par la proposition 3.2, $\quad (C , G) \to (C/H  ,   G/H)$ d\'efinissant un morphisme (voir Proposition 2.1)
$${\Cal H}_{g,G,\xi} \longrightarrow {\Cal H}_ {g',G/H, Cores
_H^G (\xi)}
\tag (6.16)$$
Il n'y a pas a priori d'extension naturelle directe du morphisme (6.16)
\`a $\overline{\Cal H}_{g,G,\xi}$, car la courbe  pr\'estable
quotient  peut ne pas \^etre stable.  Cela  justifie  la n\'ecessit\'e d'un renforcement du concept de stabilit\'e. Notons   les deux cas particuliers, $H =1$, c'est \`a dire l'oubli de l'action de
$G$, et   $G = H$. Dans ce second cas, on associe \`a un rev\^etement
$\pi: C \to D$  la base $D$, marqu\'ee par les points de
branchement,  ou plus g\'en\'eralement un ensemble de points les
contenant.
\rema{6.3}  
  Le morphisme oubli de l'action de $G$, $ \overline {\Cal
H}_{g,G,\xi} \longrightarrow \overline {\Cal M}_g, \,\, (\pi: C\to D)
\mapsto C$
  est repr\'esentable,  m\^eme
fini.    Notons  l'action libre naturelle, par torsion de l'action,
du groupe $\Aut_\xi (G)$ (stabilisateur de $\xi$) sur le champ $\overline {\Cal
H}_{g,G,\xi}$.  $\lozenge$\hfill\break

Il est n\'ecessaire maintenant d'\'etendre la d\'efinition 6.2 aux
$G$-courbes stables marqu\'ees, ce qui  conduit \`a la d\'efinition
des champs de Hurwitz versus Harris-Mumford:
\proclaim D\'efinition 6.8.  Soit $\eta \in R_+(G)$ sp\'ecifiant le
$G$-type d'un diviseur ( \S \ 
2.2). On note    $\overline{\Cal H M}_{g,G,\xi,\eta}$ le champ
(de Harris-Mumford) dont les objets sont les courbes  munies d'une
action stable de $G$, et qui sont  stables marqu\'ees par un diviseur
$G$-invariant de $G$-type fix\'e $\eta$, la donn\'ee de
ramification \'etant fix\'ee \'egale \`a $\xi$. 

Dans la suite  on supposera que $\xi \leq  \eta$, en d'autres
termes que les points de ramification sont inclus dans les points
marqu\'es. Rappelons que cela implique  
(corollaire  4.6 du Th\'eor\`eme 4.5) que les stabilisateurs des
points doubles sont maintenant
cycliques, ces points ne devant pas \^etre  consid\'er\'es comme des
points de ramification. Si la base est  
${\Bbb P}^1$, et si la ramification est simple,  on retrouve
essentiellement  le champ introduit par
Harris et Mumford \cite {42}. Pour \^etre plus pr\'ecis,
$\overline{\Cal HM}_{g,G,\xi,\eta}$ en  est comme nous le verrons,
une d\'esingularisation.

\bigskip 
 {\it 6.2.3. Compactification stable
du champ de Hurwitz (I)}  
\bigskip

Le r\'esultat suivant  suivant est bien connu; comparer avec Jarvis \cite {43}, \cite {44}, ou  \cite
{1}, \cite {59}:
\proclaim Th\'eor\`eme 6.9.    $\overline{\Cal
H}_{g,G,\xi}$ (resp.  $\overline{\Cal H
M}_{g,G,\xi,\eta}$)  est un champ de Deligne-Mumford lisse et
propre sur  ${\Bbb Z}[{1\over
\Card G}]$. 

\dem Rappelons  
que le contenu de ce r\'esultat se
r\'eduit essentiellement \`a la v\'erification de deux choses:\hfill\break
\indent $\quad \bullet$ a) Si $C_1$ et $C_2$ sont deux objets de l'un des
champs consid\'er\'es, au dessus de $S \in Sch _{{\Bbb
Z}[{1\over\Card G}]}$, alors le foncteur $\Isom_{S,G} (C_1,C_2)$ qui
en $T$,
objet au dessus de $S$ est \'egal \`a $\Isom_{T,G} (C_1\times_S T ,
C_2\times_S T)$ est
repr\'esentable par un sch\'ema fini et non ramifi\'e sur $S$.\hfill\break
\indent $\quad \bullet$ b) Mettre en \'evidence un atlas $X \to
\overline{{\Cal H}_{g,G,\xi}}$, ce qui
signifie un morphisme \'etale surjectif du champ repr\'esent\'e
par $X$ sur le champ de Hurwitz
sus-mentionn\'e.

La preuve de a) est  claire car le foncteur $\Isom_{S,G}
(C_1,C_2)$ est visiblement un
sous-foncteur ferm\'e du foncteur $\Isom_S (C_1,C_2)$ (i.e. non \'equivariant); le r\'esultat d\'ecoule
alors du r\'esultat classique (\cite {20}, Theorem 1.11). Avant de
proc\'eder \`a une v\'erification plus
d\'etaill\'ee de l'assertion b), on peut observer que par le proc\'ed\'e d'adjonction de points marqu\'es (section 2.3), on peut se
limiter sans perte de g\'en\'eralit\'e au cas $\xi = \eta$.    Dans
ce
cas, pour simplifier les notations, les champs seront not\'es
$\overline{\Cal H}_{g,G,\xi}$ et
$\overline{\Cal H M}_{g,G,\xi}$ (on omet l'indice $\eta$). On obtient
l'atlas $X$ en associant \`a
un point de $\overline{\Cal H}_{g,G,\xi}$ un point de  Hilbert
convenable, selon la m\'ethode
de Gieseker-Mumford (Gieseker \cite {36}). 

Rappelons les points
essentiels de cette construction, en  traitant en premier le cas
lisse, donc  en se concentrant sur le champ ${\Cal H}_{g,G,\xi}$.
Fixons un entier $m\geq 3$, le choix de $m$ sera pr\'ecis\'e ult\'erieurement, et fixons le $G$-module (de Hurwitz) $V = V_m = H^0(C ,
\omega_c^{\otimes m})$; soit  $\dim (V) = n+1$. Consid\'erons le
sch\'ema de Hilbert not\'e pour abr\'eger Hilb, qui classifie les
courbes de genre $g$ et de degr\'e $m(2g - 2)$ dans ${\Bbb P}(V) =
{\Bbb P}^n$, c'est-\`a-dire de polyn\^ome de Hilbert $P(\nu) =
(2m\nu - 1)(g - 1)$; on peut fixer pour sch\'ema de base ${\Bbb
Z}[{1\over \Card G}, \zeta]$, la racine de l'unit\'e $\zeta$
convenablement choisie, assurant
l'existence de $V$ comme $G$-module sur ce sch\'ema de base. Il y a  une action
naturelle du groupe $PGL(n)$ sur Hilb ainsi que, par restriction, sur
la courbe universelle $Z
\subset {\Bbb P}^n \times Hilb$. Si $C$ est une telle courbe, d\'efinie par exemple sur le corps
alg\'ebriquement clos $k$, on a pour $\nu >> 0$, une surjection
$$H^0({\Bbb P}^n , {\Cal O}_{{\Bbb P}^n} (\nu)) \longrightarrow H^0(C
, {\Cal O}_C (\nu))  \tag (6.17)
$$
d\'efinissant le $\nu^{e}$ point de Hilbert,  qui est le point de l'espace
projectif  ${\Bbb P}^N = {\Bbb P} \left(\bigwedge^{P(\nu)} H^0({\Bbb
P}^n , {\Cal O}_{{\Bbb P}^n}
(\nu)) \right)$ fourni par la surjection
$$\bigwedge^{P(\nu)} H^0({\Bbb P}^n , {\Cal O}_{{\Bbb P}^n} (\nu))
\longrightarrow
\bigwedge^{P(\nu)}  H^0(C , {\Cal O}_C (\nu)) \tag (6.18)$$
 Il y a une version \'equivariante de cette construction, dans laquelle les courbes
plong\'ees $m$-canoniquement dans ${\Bbb P}^n$, non d\'eg\'en\'er\'ees  i.e  $H^0({\Bbb P}^n ,
{\Cal O}_{{\Bbb P}^n} (1)) \cong  H^0(C , {\Cal O}_C (1))$ et invariantes par
l'action de $G$. La courbe
$C$ h\'erite alors canoniquement d'une action de $G$. Ces courbes
correspondent dans l'interpr\'etation par les points de Hilbert,  aux
surjections (6.17)  \'equivariantes. Soit  $\pi: \Sigma
\to K$ la courbe universelle de Gieseker-Mumford \cite {37}; le
groupe $G$ agit sur $\Sigma$ et  $K$, soit alors $H = K^G$ le
sous-sch\'ema des points fixes. Comme $K$ est lisse, il en est   
de m\^eme de $H$; alors si $\Gamma = \pi^{-1} (H)$, le morphisme
$\Gamma \to H$ obtenu par
restriction, est la courbe universelle cherch\'ee. On notera que la
lissit\'e de $H$ se d\'eduit
aussi des r\'esultats sur la structure des d\'eformations \'equivariantes
du paragraphe 5. En
conclusion $H$ conduit \`a un morphisme repr\'esentable $H \to
{\Cal H}_{g,G,\xi}$, et qui par la
propri\'et\'e de lissit\'e assure que ${\Cal H}_{g,G,\xi}$ est un
champ alg\'ebrique de Deligne-Mumford.

Pour terminer la preuve, il suffit de montrer que la proc\'edure pr\'ec\'edente s'\'etend \`a
$\overline{\Cal H}_{g,G,\xi}$, resp. \`a $\overline{\Cal H
M}_{g,G,\xi}$. On donne quelques   d\'etails pour le premier cas, les
modifications de routine pour traiter le second sont omises; notons
pour rassurer le lecteur, qu'en corollaire, les deux champs seront
identifi\'es. Soit encore, en conservant les m\^emes notations
utilis\'ees au dessus,
$V = V_m$ la repr\'esentation de Hurwitz de degr\'e $m$, de rang
$n+1 = (2m - 1)(g - 1)$, et
${\Bbb P}^n = {\Bbb P}(V)$.  Maintenant  rappelons q'un r\'esultat
fondamental de Gieseker et Mumford, assure l'existence d'un $m>>0$,
puis d'un $\nu>>0$, de
sorte que le $\nu^{e}$ point de Hilbert de toute courbe connexe lisse
(d\'efinie sur un corps alg\'ebriquement clos $k$) soit G.I.T-stable
pour l'action du groupe lin\'eaire. Une fois $\nu$ fix\'e, on peut
consid\'erer la courbe universelle $Z \to H$ de base le sch\'ema de
Hilbert. Soit alors $W$  d\'efini comme \'etant l'ensemble des
points $h\in H$, tels que les deux conditions
de dessous sont satisfaites:\hfill\break
\indent $\quad \bullet$ i) $h\in H^{ss}$, l'indice sup\'erieur signifiant que
le $\nu^{e}$ point de  Hilbert
est semi-stable,\hfill\break
\indent $\quad \bullet$ ii) La courbe $Z_h$ est connexe et plong\'ee
$m$-canoniquement dans ${\Bbb P}^n$.

Il est prouv\'e dans \cite {37} que $W$ est un sous-sch\'ema ferm\'e
 de $H^{ss}$, et que si $h \in H$, la courbe $Z_h$ est r\'eduite
avec seulement des points doubles ordinaires. Rappelons que l'espace
modulaire grossier classifiant les courbes stables de genre $g\geq
2$, s'identifie au quotient
g\'eom\'etrique\footnote{Ne pas confondre avec le 2-quotient de la section 6.1.} $\overline{M_g} = W // SL(n)$.   De mani\`ere
identique a ce qui a \'et\'e fait au-dessus, on doit traiter  la
situation \'equivariante, ce qui am\`ene d'abord \`a former le
sous-sch\'ema des points fixes $T = W^G$; l'action de $G$ est celle
naturelle d\'ecrite dans la section 6.1. Le sch\'ema $T$ est lisse,
du fait de la lissit\'e de $W$, ou si on pr\'ef\`ere par la th\'eorie des d\'eformations \'equivariantes du paragraphe 5.  Les
composantes connexes  $\{T_\alpha\}$ de $T$ sont alors irr\'eductibles, et le long d'une telle composante, les repr\'esentations
de Hurwitz $H^0(C_h , \omega_{C_h}^{\otimes t})$ sont constantes (Lemme
4.3); noter que pour $t=m$, cela r\'esulte des hypoth\`eses. Si on
travaille avec
les objets du champ $\overline{\Cal H M}_{g,G,\xi}$, le m\^eme
argument que celui utilis\'e au-dessus nous assure que le $G$-type
du diviseur des points marqu\'es est constant le long des
$T_\alpha$. Parmi les composantes $\{T_\alpha\}$, on ne retient que
celles dont le point
g\'en\'erique d\'efinit une courbe lisse, et le long desquelles la
donn\'ee de Hurwitz, et
\'eventuellement le $G$-type du marquage  sont constants \'egaux \`a
$\xi$ et $\nu$; la r\'eunion de
ces composantes est encore not\'ee $T$. La restriction \`a $T$ de
la famille universelle fournit un morphisme repr\'esentable surjectif
$T\to \overline{\Cal H}_{g,G,\xi}$, prouvant que $\overline{\Cal
H}_{g,G,\xi}$ est un champ alg\'ebrique, de Deligne-Mumford. Notons
pour conclure que la propret\'e dans l'un ou l'autre des deux cas
r\'esulte du th\'eor\`eme de r\'eduction stable sous sa forme \'equivariante (Proposition 5.1).\hfill\break
\qed 

 \rema{6.4}  
 Si  $k = {\Bbb C}$, le type
topologique de l'action
s\'epare les composantes connexes. On invoque la version \'equivariante de la th\'eorie  de
Teichm\"{u}ller \cite {25}, de laquelle on d\'eduit que deux points de $H$ qui correspondent \`a
 des actions topologiquement
conjugu\'ees peuvent \^etre joints par un arc. On a donc (voir aussi
Th\'eor\`eme 2.1) $\Card {\pi_0(H)} = h(\xi) \,\, =   \hbox{\rm le  nombre  de  Nielsen}. \,\lozenge$ 

\bigskip

\beginsection 6.3. Compactification du  sch\'ema  de Hurwitz:  Gieseker-Mumford

\bigskip

   On se concentre  dor\'enavant sur les espaces de modules  
$\overline{H}_{g,G,\xi}$ (resp.$\overline{HM}_{g,G,\xi}$) associ\'es
aux champs $\overline{\Cal H}_{g,G,\xi}$ (resp.
$\overline{\Cal H M}_{g,G,\xi}$). La construction G.I.T esquiss\'ee
ci-dessus (\S \ 6.2)  conduit \`a  la
projectivit\'e des   sch\'emas modulaires de Hurwitz, et en corollaire \`a la   
   projectivit\'e des  $\overline{M}_{g,n}$.   On remarquera que finalement tout  est r\'eduit  \`a une et une seule construction,  la projectivit\'e
de $\overline{M_g}$.

Identifions $G$ \`a son image dans $GL(V)$, et notons $A$ (resp.
$B$) le centralisateur (resp.
le normalisateur) de $G$. Si $V = \bigoplus_i n_iV_i$ est la d\'ecomposition de  $V$ en
repr\'esentations irr\'eductibles deux \`a deux distinctes, alors
$A \cong \prod_i GL(n_i)$. Notons
que $A$ op\`ere sur $H$ en fixant individuellement  les composantes
connexes. La construction de $H_{g,G,\xi}$ (resp.
$\overline{H}_{g,G,\xi})$, r\'esulte de la proposition suivante:

\proclaim Proposition 6.10.  Si $m$ est assez grand, le quotient g\'eom\'etrique $ H/A$ existe, et
$$H_{g,G,\xi} =  H/A \quad (\hbox{\rm  resp}.\overline{H}_{g,G,\xi} =  {\overline H}/A.)\tag (6.19)$$ 

\dem C'est une simple adaptation \`a la situation \'equivariante
pr\'esente des arguments de Mumford et Gieseker \cite {37}. Traitons
en premier, pour plus de lisibilit\'e, l'espace modulaire
$H_{g,G,\xi}$. On observe   qu'une courbe  $C \subset {\Bbb
P}^n$,lisse, connexe, de genre
$g\geq 2$, est Chow-stable ou Hilbert-stable, si son degr\'e est
grand par rapport \`a $g$. La stabilit\'e  d'abord relative au groupe
lin\'eaire, entra\^\i ne a fortiori celle relative \`a $A$. Pour
que la preuve de la proposition  soit compl\`ete, on doit  tout
d'abord  noter la validit\'e \/ des deux points
suivants:\hfill\break
\indent $\quad \bullet$ i) Toute action de $G$ sur $C$, du type indiqu\'e,
appara\^\i t dans une fibre de la courbe
universelle \hfill\break
\indent $\quad \bullet$ ii) Si $C_1$ et $C_2$ sont deux fibres
$G$-isomorphes, alors l'isomorphisme est
induit par un \'el\'ement de $A$. En fait si $\pi: C \to S$ est un
objet de ${\Cal H}_{g,G,\xi}$, on prouve que localement pour la
topologie de Zariski sur $S$, $C/S$ appartient \`a ${{\Hom (S,H)}
\over {\Hom (S,A)}}$. Pour cela, on consid\`ere le faisceau ${\Cal E}
= \pi_* (\omega_{C/S}^{\otimes m})$; on a clairement $R^1\pi_*
(\omega_{C/S}^{\otimes m}) = 0$, et ${\Cal E}$ est localement libre,
muni d'une action de $G$. On peut alors  d\'ecomposer ${\Cal E}$ en
facteurs isotypiques
$${\Cal E} = \bigoplus_{V_i \in Irrep(G)}  {\Cal E}_i \otimes V_i \tag (6.20)$$
il suffit alors de noter que  les modules de covariants ${\Cal E}_i$
sont localement libres, et
que  $V = \oplus rg({\Cal E}_i) V_i$.

On peut r\'ealiser de la m\^eme mani\`ere l'espace modulaire grossier
$\overline {H}_{g,G,\xi}$,
prouvant ainsi que c'est un sch\'ema projectif, et que $H_{g,G,\xi}$
s'identifie \`a un ouvert
partout dense de $\overline {H}_{g,G,\xi}$; en particulier les deux
sch\'emas ont le m\^eme nombre
de composantes connexes, nombre \'egal au nombre de Nielsen $h(\xi)$
(Th\'eor\`eme 2.1).
Conservant les notations du d\'ebut de section, on prouve comme
ci-dessus que le quotient g\'eom\'etrique 
$T//A$ existe et que finalement $\overline {H}_{g,G,\xi} = T//A$.  On
doit  dans une seconde \'etape s'assurer que toute courbe stable $C$
munie d'une action stable de $G$ de donn\'ee de ramification
prescrite, appara\^\i t dans le bord. On sait, du fait de la
structure de la d\'eformation universelle de $(C,G)$ (section 5.1)
qu'on peut \'etendre de mani\`ere \'equivariante $C$ en une courbe
stable $\pi: X \to S = \Spec (R)$, l'anneau $R$ \'etant de valuation
discr\`ete et complet, de sorte que si on  note $a$ le point ferm\'e
, et $b$ le point g\'en\'erique, alors $\pi^{-1} (a) \cong C$ et la
courbe $\pi^{-1} (b)$ est lisse. Cette courbe d\'efinit un morphisme
classifiant $f : \Spec (R) \longrightarrow H^G$
avec $f(b) \in W$. L'argument de Gieseker (\cite {36}, thm 2.0.2)
montre qu'il existe une courbe
appartenant \`a la $A$-orbite de $X$ et qui correspond \`a un
morphisme classifiant $\Spec (R) \to
\sqcup T_\alpha$. Les fibres au dessus du point ferm\'e sont ainsi
dans une m\^eme $A$-orbite. Cela \'etant, la projectivit\'e de
$\overline {H}_{g,G,\xi}$  d\'ecoule bien \'evidemment de la
construction de cet espace modulaire comme quotient. On peut aussi,
la propret\'e \'etant cons\'equence du th\'eor\`eme de r\'eduction
stable sous sa forme \'equivariante, s'appuyer sur la projectivit\'e
 connue de
$\overline M_g$ et utiliser le morphisme oubli de l'action de $G$
$$ \overline {H}_{g,G,\xi} \longrightarrow \overline {M_g}\tag (6.21)$$
qui est alors propre et quasi-fini, donc fini.\hfill\break\qed

\proclaim Proposition 6.11.  Soit $\xi$ une donn\'ee de Hurwitz
attach\'ee au groupe $G$, et soit
$h(\xi)$ le nombre de Nielsen correspondant. Si $k$ est un corps
alg\'ebriquement clos de
caract\'eristique non facteur de  $\Card G$, alors
$$h(\xi)  = \Card {\pi_0(\overline{H}_{g,G,\xi})} = \Card {\pi_0
({H_{g,G,\xi}})} \tag (6.22)$$
 
\dem On reprend la construction de la proposition 6.1, maintenant
sur un  sch\'ema de  base  $S = \Spec (R)$, o\`u $R$ est un anneau
de valuation discr\`ete complet, de corps des fractions $K$ de
caract\'eristique z\'ero, et de corps r\'esiduel $k$. Si $G$ agit
stablement sur une courbe stable
d\'efinie sur $k$, on sait qu'on peut relever la courbe, ainsi que
l'action de $G$, \`a $R$, cela \`a
donn\'ee de ramification constante. Notons $a$ (resp. $b$) le point
ferm\'e (resp. le point
g\'en\'erique) de $S$. Consid\'erons le morphisme quotient
$\sqcup_\alpha T_\alpha \longrightarrow \overline{H}_{g,G,\xi} =
(\sqcup_\alpha T_\alpha)// A$. 
Rappelons que chaque composante connexe $T_\alpha$ est lisse sur $S$.
Si on prend la fibre  en  $a$ (resp. $b$), on a par un r\'esultat de Seshadri
(rappel\'e dans \cite {35}),
$$(T_\alpha // A)_a = (T_\alpha)_a // A \quad (\hbox{\rm resp.}) \quad
(T_\alpha // A)_b =
(T_\alpha)_b // A$$
Le th\'eor\`eme de connexion de Zariski montre  que  $(T_\alpha //
A)_a $ est connexe; il s'ensuit que  $(T_\alpha)_a$ est connexe, donc
irr\'eductible, car lisse. On en d\'eduit  que $(T_\alpha)_a // A =
(\overline{H}_{g,G,\xi})_\alpha$ est irr\'eductible, d'o\`u le r\'esultat annonc\'e. \hfill\break
  \qed

\proclaim Corollaire 6.12.  Sous les conditions pr\'ec\'edentes, et
si maintenant $G$ est un groupe
cyclique, alors pour tout corps alg\'ebriquement clos $k$ de caract\'eristique ne divisant pas $\Card
G$, l'espace de Hurwitz $\overline {H}_{g,G,\xi}$ est irr\'eductible.

\dem  On sait en effet que sous ces hypoth\`eses que $h(\xi) = 1$
(Proposition 2.3). \hfill\break\qed

\bigskip

\beginsection 6.4. Compactification stable du champ de Hurwitz (II)

\bigskip

  {\it 6.4.1. Le champ de Hurwitz versus Harris-Mumford} 

\bigskip

La construction du sch\'ema de Hurwitz au moyen d'un quotient
g\'eom\'etrique conduit   avec des modifications mineures,
\`a  la construction  des espaces modulaires
$HM_{g,G,\xi}$ (resp. $\overline{HM}_{g,G,\xi}$). Rappelons
que dans le cas Harris-Mumford les objets du champ sont les
$G$-courbes stables marqu\'ees par les points de
ramification,  ou plus g\'en\'eralement marqu\'ees par un diviseur
$G$-invariant les contenant.   On va prouver ci-apr\`es que
les espaces de modules $\overline{H}_{g,G,\xi}$ et $\overline{
HM}_{g,G,\xi}$ sont en fait  isomorphes, m\^eme si cela n'est pas  en g\'en\'eral le cas pour les champs.  Rappelons  au pr\'ealable  les constructions de Knudsen \cite {48}, \cite
{53}:

\indent $\bullet $ Soit $(C' , \{y_j\})$ est une courbe pr\'estable marqu\'ee, connexe. Parmi tous les morphismes
surjectifs  $ (C' , \{y_j\}) \longrightarrow (C ,\{x_i\})$
de but une courbe stable marqu\'ee, les $(y_j)$ s'envoyant
surjectivement sur les $(x_i)$, il y a si l'ensemble des $\{y_j\}$ est non vide
un morphisme {\it minimal}.  Il
s'obtient en contractant en un point les composantes instables (de
genre z\'ero) contenues dans
$C'$. On montre que ce morphisme est  d\'efini par une puissance
convenable du faisceau
inversible $\omega_{C'} (\sum y_j)$. Il est clair que si $G$ agit sur
$C'$ en permutant les $y_j$, alors l'action de $G$ descend \`a $C$,
de sorte que le  morphisme
commute \`a l'action de $G$. On applique  cette proc\'edure  lorsque  l'oubli de certains points  marqu\'es d'une
courbe stable $(C', \{y_j\})$ d\'etruit la stabilit\'e.    Les sections oubli\'ees induisent des sections de $C$,  eventuellement non disjointes et  passant par des points doubles.  Important pour la suite
est que cette op\'eration de contraction (ou de stabilisation) est compatible aux changements de base, donc d\'efinit apr\`es oubli du marquage par les points de
ramification,  un morphisme  
$$\st :  \overline{\Cal HM}_{g,G,\xi} \longrightarrow
\overline{\Cal H}_{g,G,\xi} \tag (6.23)$$
\indent $\bullet$ Il y a un foncteur   en sens inverse   (\cite {46},
cor 2.6). Ce morphisme  qui revient \`a \'eclater certains points doubles, transforme une courbe stable marqu\'ee par un ensemble de $n$ sections, \`a laquelle on ajoute une section
suppl\'ementaire, en une courbe stable marqu\'ee par $n+1$ sections, la $(n+1)$-i\`eme \'etant la section additionnelle. Ces deux op\Ž rations conjointes montrent que le morphisme  $\overline{\Cal M}_{g,n+1} \rightarrow \overline{\Cal M}_{g,n}$,  oubli de la section $n+1$, repr\'esente la courbe universelle $n$-piqu\Ž e.   Le r\'esultat suivant clarifie en partie la situation.

\proclaim  Proposition 6.13.  Le morphisme     (6.23)  est un isomorphisme au-dessus du sous-champ ouvert partout dense ${\Cal HM}_{g,G,\xi}$ sur son image;  il  induit un isomorphisme
d'espaces modulaires grossiers:
$\overline {HM}_{g,G,\xi}\buildrel\sim\over  \longrightarrow \overline {H} _{g,G,\xi}. $   

\dem Soit $\pi': C' \to S$ une $G$-courbe stable marqu\'ee par les
points de ramification.
Rappelons que cela impose en particulier que pour tout sous-groupe
cyclique $H$ de $G$, le diviseur de Cartier relatif   $Fix(H)_{hor}$ (Proposition 4.8) est
une somme de sections disjointes. Par oubli du marquage, suivi d'une contraction
\'eventuelle des composantes instables, on obtient une $G$-courbe stable $\pi: C \to S$.  Le
morphisme $\pi'$ factorise en $\pi' = \pi \gamma$, o\`u \/
$\gamma: C' \to C$    contracte les orbites de composantes lisses rationnelles \`a
 isotropie di\'edrale (Lemme 5.4). l'image d'une telle composante $E$ est un point double \`a  isotropie di\Ž drale $\Bbb D_m \,(m\geq 1)$,  si $m$ est le nombre de sections de points fixes qui rencontrent  $E$. Ces points sont les points fixes des reflexions  de $\Bbb D_m$. Notons  ainsi que $C' \cong C$ sauf si
$C$ poss\`ede des fibres avec des
points doubles \`a isotropie di\'edrale.  On peut  comme expliqu\'e ci-dessus reconstruire  $C'$ 
partant de   $C$ munie des images des  sections de points fixes,  en s\'eparant les
sections de points fixes qui coalescent en les points d'istropie di\'edrale.  Le morphisme $st$ est clairement un isomorphisme au-dessus du sous-champ ouvert form\'e des rev\^etement de courbes lisses.  Ce n'est cependant pas un isomorphisme en g\'en\'eral.
Si $C\to S = st \,(C'\to S)$, et si $Q \in C_s $ est un point double \`a isotropie di\'edrale $H = \Bbb D_m$, pour une r\'eflexion $\tau \in H$, on a not\'e que le diviseur des points fixes (horizontal) est de mani\`ere sp\'ecifique la r\'eunion de deux sections. 

Soit $\hat {\Cal O}_Q = \hat {\Cal O}_s [[x,y]]/(xy - a)$ l'anneau local compl\Ž t\Ž \/ de $C$ en $Q$. On suppose que $\tau (x) = y$ et $\tau (y) = x$. Si l'une des sections de points fixes est $x = b , y = c$, alors  $b = c$ et $a = b^2$. Si on se place sur $S = k[\epsilon], \, \epsilon^2 = 0$, on voit que si $C'$ et $C$ sont d\'efinies  sur un corps alg\Ž briquement clos $k$, alors l'application $\st$ au niveau des espaces tangents des d\'eformations universelles respectives n'est pas surjective. En fait les d\Ž formations qui sont dans l'image sont topologiquement triviales en les points doubles \`a isotropie di\'edrale.

Si pour tout sous-groupe cyclique  $H$ d'ordre deux le diviseur   $Fix(H)_{hor}$ est une somme de deux sections, alors on peut inverser $\st$.  L'action de $G$ sur les composantes exceptionnelles se r\'eduit \`a l'action d'un groupe di\'edral  $\Bbb D_m$.  On notera que $st$ est un monomorphisme car il est clair que 
 $\Aut_ G (C) \cong  \Aut_G (\st(C))$.   Finalement le fait qu'au niveau des espaces modulaires grossiers $st$ est un isomorphisme est clair, car il est propre, birationnel et bijectif entre vari\'et\'es normales.\hfill\break
\qed 

Dans la suite nous travaillerons   avec  le champ
   $\overline{\Cal HM}_{g,G,\xi}$, qui pour all\'eger les notations  sera
not\'e $\overline {\Cal H}_{g,G,\xi}$.   Rappelons que si $H \triangleleft G$,
on peut r\'ealiser
$\overline {\Cal H} _{g,G,\xi}$ comme une correspondance entre
$\overline {\Cal H}_{g,H,Res_H^G(\xi)}$ et  $\overline {\Cal H}_{g',
G/H, Cores_H^G(\xi)}$.   Si $H=1$ et $H=G$, on
obtient les
deux morphismes fondamentaux
$$\overline{\Cal M}_{g,r}\buildrel\imath\over \longleftarrow
\overline {\Cal H} _{g,G,\xi}\buildrel\delta\over \longrightarrow
\overline
{\Cal M}_{g',b} \tag (6.24)$$
    $r$ \'etant  le nombre de points de ramification, $b$ le nombre de
points de branchement, et
$g'$ le genre de la base. Le morphisme de gauche
est l'oubli de l'action de
$G$,  $\delta$ est le passage au quotient par
$G$.  Le morphisme $\delta$    sera appel\'e  le
morphisme discriminant. Son existence est justifi\'ee par  la
Proposition 5.2, qui assure que si $\pi: C \to S$ est un objet de
$\overline {\Cal H} _{g,G,\xi}$, alors la courbe quotient
$ C/G \to S$, marqu\'ee par le diviseur des points de branchement
(resp. piqu\'ee par les points de branchement)  est stable.  En r\'esum\'e:

  \proclaim Proposition 6.14.  1)  Le morphisme   $\overline
{\Cal H} _{g,G,\xi}\buildrel\imath\over \longrightarrow \overline
{\Cal M}_{g,(r)}$ est fini (repr\'esentable) non ramifi\'e; c'est une immersion r\'eguli\`ere locale dans le sens de Vistoli (\cite {66}, 1.20) .  \hfill\break
2)  Le morphisme
discriminant  $\overline {\Cal H} _{g,G,\xi}\buildrel\delta\over
\longrightarrow \overline {\Cal M}_{g',b}$ est propre quasi-fini  et plat \footnote {Non repr\'esentable si le centre $Z(G)$ est non trivial.}.\hfill\break
 3)  Au niveau des espaces  de modules grossiers, le
morphisme discriminant $\delta$   est fini,
surjectif.

\dem  1) Si $\Sigma \to S$ est une courbe stable marqu\'ee, la cat\'egorie fibre de $\imath$ en $C/S$ a pour objets au-dessus de $T\to S$ les plongements $G \hookrightarrow \Aut_T(\Sigma\times_S T)$.  Ces plongements correspondent aux $T$-points de $\Aut_S(\Sigma)\times_S \cdots \times_S \Aut_S(\Sigma)$ ($\vert G\vert$ facteurs), dont l'image est disjointe de la diagonale. Du fait que $\Aut_S(\Sigma)$
 est fini non ramifi\'e sur $S$, le premier point est clair. Comme les champs sont lisses, $\imath$ est une immersion locale r\'eguli\`ere dans le sens de Vistoli, i.e. repr\'esentable non ramifi\'e et localement d'intersection compl\`ete.
\hfill\break
2) La propret\'e de $\delta$ d\'ecoule 
du crit\`ere de r\'eduction stable pour les rev\^etements stables (Proposition 5.7).  La platitude 
d\'ecoule du th\'eor\`eme 5.5 donnant la structure locale   de $\pi$.
    Le caract\`ere \'etale de $\delta$ (au-dessus de ${\Cal M}_{g',b}$) d\'ecoule alors du
th\'eor\`eme 5.5. Le dernier point est clair.\hfill\break
\qed 

 \rema {6.5}    La
d\'efinition des rev\^etements stables donn\'ee ci-dessus permet de
retrouver d'une mani\`ere uniforme les objets \'etudi\'es par
  Abramovich, Corti et Vistoli \cite {1}, \cite {3}   sous le
terme {\it  balanced twisted stable maps}.  Les objets \'etudi\'es
par ces auteurs sont  les morphismes $\pi: D \longrightarrow
BG$ o\`u la source $D$, est la courbe  $D$ munie de sa structure
d'orbifold, i.e. donn\'ee par $D = C/G$. Dans cette structure, on
ignore les automorphismes induits par le centre de $G$. En
cons\'equence, on a (notation de \cite {3})  $$\coprod_{[\xi]} \overline
{\Cal H}_{g,G,\xi}// Z(G)  \buildrel\sim\over\longrightarrow
\overline {\Cal M}_{g',b} (BG) \tag (6.25)$$
$\lozenge$ 
\bigskip
  {\it 6.4.2. Le rev\^etement "universel"}  
\bigskip
  On continue d'explorer les aspects fonctoriels des champs de Hurwitz. Soit $[H,\chi]$ une classe  appara\^{\i}ssant 
dans  la donn\'ee de ramification $\xi$.  Soit $\Delta_{(H,\chi)}$ la composante de $\Delta_H$,
lieu des points fixes d'holonomie exacte $(H,\chi)$. Il y a une
action libre de $C_G(H)/H$ sur $\Delta_{(H,\chi)}$, le quotient
\'etant la composante $B_{[H,\chi]}$ de $B$, d\'efinissant  un
fibr\'e principal 
$  \Delta_{(H,\chi)} \longrightarrow
B_{[H,\chi]} $.
Supposons    que $B$ soit somme de sections  $Q_1,\dots, Q_b$  ($D$  est piqu\'ee)  
 contenant les points de branchement.   Par image r\'eciproque,
la  section $Q_i$, si  d'holonomie  $(H_i,\chi_i),$ d\'efinit un fibr\'e principal de
base $S$, et de groupe $C(H_i)/H_i$. Cela fournit   
un foncteur {\it d'\'evaluation}  en $Q_i$ \cite {45}: $$ev_{Q_i}:
\overline {\Cal H}_{g,G,\xi} \longrightarrow B(C_G(H_i)/H_i) \tag
(6.26)$$ 
 Soit $\Spec k \to   B(C(H_i)/H_i) $    l'atlas  correspondant au fibr\'e principal trivial, alors le 2-produit fibr\'e    $ {\Cal H}_{g,G,\xi} \times_{ B(C(H_i)/H_i) } \Spec k$
a pour objets les rev\^etements $\pi: C\to D$ \'equip\'es d'un point $P: S \to C$ qui est un point  de ramification d'holonomie $(H_i,\chi_i)$. 

Soit   $\imath: \overline {\Cal H}_{g,G,\xi}  \rightarrow \overline {\Cal
M}_{g,(r)}$ le morphisme  oubli\footnote{Rappelons que dans $\overline{\Cal M}_{g,(r)}$ le marquage est par parquets, d\'etermin\'es par la donn\'ee $\xi$.}de l'action de $G$, et consid\'erons  $\overline
{\Cal C}_{g,(r)} = \overline{\Cal M}_{g,(r),1} \rightarrow \overline {\Cal M}_{g,(r)}$ la courbe
universelle  (marqu\'ee) au-dessus de $\overline {\Cal M}_{g,(r)}$ \cite {48}.
 Une section de  $\overline {\Cal C}_{g,(r)}$ au-dessus
de $S$ est la donn\'ee  d'une courbe $r$-marqu\'ee $q: C\to S$  \'equip\'ee d'une section suppl\'ementaire  $P: S\to C$.  Il n'est pas n\'ecessaire d'exiger  
que cette section soit disjointe des points doubles ou des points marqu\'es. 
  Formons le carr\'e 2-cart\'esien $$ \commdiag{
\overline {\Cal C}_{g,(r)}  &\mapright^q& \overline {\Cal M}_{g,(r)}  \cr
 \mapup && \mapup\lft{\imath}   \cr
  \overline {\Cal C}_{g,G,\xi} & \mapright^q& \overline
{\Cal H}_{g,G,\xi} \cr}$$
Une  section  de  $\overline {\Cal C}_{g,G,\xi}$  au-dessus  de $S$
est  la  donn\'ee  d'un  $G$-rev\^etement 
 $q: C\buildrel\pi\over\rightarrow D\buildrel p
\over\rightarrow S$ \'equip\'e d'une section $P: S\to C$ de $q$. Le
groupe $G$ agit  (dans le sens de la section 6.1) sur le champ  $ \overline {\Cal C}_{g,G,\xi}$ par
$g (C\buildrel \pi \over\rightarrow D , P) =
(C\buildrel\pi\over\rightarrow D , g\circ P)$.
Soit  le morphisme $S \to \overline {\Cal H}_{g,G,\xi}$  d\'efini par le $G$-rev\^etement
$\pi: C\to D$ au-dessus de $S$.   Il est   clair qu'on a
une identification  de $G$-champs 
$  \overline {\Cal C}_{g,G,\xi}
\times_{\overline {\Cal H}_{g,G,\xi}} \, S \cong C,$
  justifiant le
fait que    $\overline {\Cal C}_{g,G,\xi}  \longrightarrow  \overline
{\Cal H}_{g,G,\xi}$ est la $G$-courbe universelle au-dessus de
$\overline {\Cal H}_{g,G,\xi}$. Le foncteur  
$$\Delta: (C\buildrel \pi
\over\rightarrow D , P) \mapsto (D , Q = \pi P), \quad  \overline
{\Cal C}_{g,G,\xi}  \longrightarrow \overline {\Cal C}_{g',b} $$
fait office de {\it rev\^etement universel}. 

Soit $A: U\to {\overline {\Cal H}_{g,G,\xi}}$ un atlas,   et notons $\pi:
C_U \to U$ le $G$-rev\^etement d\'efinissant $A$. La seconde
projection $p_2:  C_U \times_U C_U \rightarrow C_U$ d\'efinit un
$G$-rev\^etement, le groupe $G$ agissant sur le facteur de gauche
dans le produit  fibr\'e. La diagonale $\Delta: C_U \rightarrow C_U
\times_U C_U$ vue comme section de $p_2$, fait de  $p_2:  C_U
\times_U C_U \rightarrow C_U$ un objet de ${\overline {\Cal
C}_{g,G,\xi}}$ (on peut stabiliser cette section additionnelle si on veut), donc  conduit   \`a un morphisme $C_U \rightarrow
{\overline {\Cal H}_{g,G,\xi}}$ rendant le carr\'e suivant 
$$\commdiag{
\overline {\Cal C}_{g,G,\xi} & \mapright^q & \overline {\Cal H}_{g,G,\xi}
\cr
 \mapup &&  \mapup\lft{A} \cr
  C_U & \mapright^{ p_2} &    U } \tag (6.27)$$
  2-cart\'esien. En
particulier le morphisme $C_U \rightarrow  \overline {\Cal
C}_{g,G,\xi} $  d\'efinit un atlas de $ \overline {\Cal C}_{g,G,\xi}
$.  Pour faire le lien avec la base, notons le diagramme $2$-commutatif  
$$\commdiag{  \overline {\Cal C}_{g,G,\xi} & \mapright^\Delta & \overline {\Cal C}_{g',b}\cr
 \mapdown\lft{q} &&\mapdown\lft{p}\cr
  \overline {\Cal H}_{g,G,\xi}&  \mapright^\delta & \overline{\Cal M}_{g',b} }  $$
 Le morphisme  $\Delta$ n'est cependant pas le quotient de  $\overline {\Cal C}_{g,G,\xi}  $
par $G$.   Retenons cependant:

\proclaim  Proposition 6.15. Le morphisme   
$\psi: \overline
{\Cal C}_{g,G,\xi}  \longrightarrow\overline{\Cal H}_{g,G,\xi}\times_{\overline {\Cal M}_{g',b}} \, \overline{\Cal C}_{g',b} $
  est fini  de degr\'e $\vert G\vert $.

  \dem  En effet le champ  fibre de  $\psi$  en le rev\^etement $\pi: C\to D$ de base $S$, et muni de la section $Q: S\to D$, est  clairement repr\'esent\'e par le  sch\'ema $\pi^{-1}(Q)\subset C$.\hfill\break\qed 
  
  Noter que les sections universelles   $Q_\alpha \, (\alpha=1,\cdots,b)$  font de $\overline{\Cal H}_{g,G,\xi}\times_{\overline {\Cal M}_{g',b}} \, \overline{\Cal C}_{g',b}$ une courbe $b$-marqu\'ee de base  $\overline{\Cal H}_{g,G,\xi}$.

\bigskip

\beginsection  6.5. Champs de rev\^etements stables: cas non galoisien

\bigskip

   Dans cette section on construit  le champ de Hurwitz classifiant les rev\^etements  non galoisiens  stables \`a  monodromie fix\'ee.     Nous  le comparerons   au champ des rev\^etements admissibles \cite
{40}, \cite {42}, \cite {55}, \cite {68}, \cite{69}.  
\bigskip 
  {\it 6.5.1. Rev\^etements
admissibles}
  \bigskip

Le   point  cl\'e dans la construction    du champ de Hurwitz classifiant les rev\^etements galoisiens
 est  le fait de consid\'erer la base   comme
une courbe marqu\'ee par les points de branchement (lisses), qui en
cons\'equence ne sont pas autoris\'es \`a se rencontrer   
 sp\'ecialisation.  Les rev\^etements appel\'es
admissibles  par Harris-Mumford (loc.cit.p 57, et D\'efinition
6.15) satisfont \`a cette condition, mais cependant  ont   une mauvaise th\'eorie
des d\'eformations, en particulier la base de la d\'eformation verselle n'est en g\'en\'eral pas lisse, m\^eme pas normale. Ce d\'efaut a \'et\'e corrig\'e par Mochizuki \cite {55} et Wewers
\cite {68}, en proposant d'enrichir un tel rev\^etement d'une
log-structure, ce qui a pour effet  de r\'etablir la lissit\'e du
foncteur des d\'eformations.    

 Dans  la  construction qui suit  on proc\`ede en sens inverse.  Un rev\^etement stable   sera d\'ecrit, du moins \'etale-localement, comme
quotient d'une cl\^oture galoisienne, qui est un $G$-rev\^etement
stable marqu\'e par les points de branchement. 
 Fixer une collection de cl\^otures galoisiennes
locales  ajoute une structure suppl\'ementaire  qui cependant n'est pas \'equivalente \`a celle fournie par une log-structure. C'est une structure    plus faible.    Rappelons  la d\'efinition des rev\^etements admissibles
(\cite {40}, \cite {42} p 61, ou \cite {68}):
\proclaim D\'efinition 6.16. Soit une courbe stable marqu\'ee $D \to
S$, de genre $g'$, avec $S$
connexe,  le marquage \'etant d\'efini par un diviseur de Cartier
relatif $B/S$, \'etale de degr\'e $b$. Un rev\^etement admissible
$\pi: C\to D$ au-dessus de $S$,  de base $D/S$, est d\'efini
par:\hfill\break
\indent $\bullet$ i) Une courbe pr\'estable $C/S$, munie  de la
donn\'ee d'un morphisme fini surjectif
$\pi: C \to D$, et d'un diviseur de Cartier relatif $R \subset C$, \'etale sur $S$, et disjoint du lieu de non lissit\'e, tel que $\pi
(R) = B$. On fait l'hypoth\`ese que $\pi$ est  (mod\'erement) ramifi\'e exactement
le long de $R$. On suppose en outre que la monodromie  le long des fibres g\'eom\'etriques  est constante,  \'egale \`a celle fix\'ee initialement.\hfill\break
\indent $\bullet$ ii) (structure locale aux points doubles) Si $Q \in
D_s$ est un point double d'une fibre
g\'eom\'etrique, on suppose que localement pour la topologie \'etale
 la structure  du
morphisme $\pi$ en $P \in \pi^{-1}(Q)$ est d\'ecrite parÊ\footnote {\rm On peut exprimer la structure locale en disant que le syst\`eme de coordonn\'ees $(u,v)$ en $Q$ admet une racine $e$-i\`eme $(x,y)$ dans $ \hat{\Cal O}_{C,P}$ \cite {69}. On v\'erifie facilement en utilisant le lemme de Hensel que tout syst\`eme de coordonn\'ees en $Q$ admet   une racine $e$-i\`eme.}: 
$\pi_P^* : \hat{\Cal O}_{D,Q} \longrightarrow \hat{\Cal O}_{C,P}$
avec $\hat{\Cal O}_{C,P} =  {\Cal O}_s [[x,y]]/ (xy-t)$,\quad
$\hat{\Cal O}_{D,Q} =
  {\hat{\Cal O}}_s [[u,v]]/ (uv-\tau)$, o\`u $t\in \hat{\Cal M}_P, \tau
\in \hat{\Cal M}_Q$, et pour un certain $d\geq 1$,  
$u = x^d, v = y^d$ (donc $\tau = t^d$). 

Etale-localement on peut pr\'eciser  le sens de
l'hypoth\`ese $\pi$ est ramifi\'e le long de $R$. On peut en effet supposer
que $B = \sum_j Q_j$ est une somme de sections disjointes, et  de m\^eme 
$R = \sum_{i,j} P_{i,j}$, o\`u les $(P_{i,j})_i$ sont les points
d'image $Q_j$. On suppose alors que, l'indice de ramification le long
de $P_i$ \'etant $e_{i,j}$, on a $\pi^{-1} Q_j = \sum_i e_{i,j}
P_{i,j}$. Il  reste \`a expliquer la condition sur la monodromie
dans la d\'efinition 6.15,  condition qui n'est pas \`a priori d\'efinie le long des fibres singuli\`eres.  On passe par une cl\^oture galoisienne.

\bigskip

  {\it 6.5.2. Cl\^oture galoisienne: cas des courbes lisses} 

\bigskip

    L'op\'eration  fondamentale pour  la suite est l'op\'eration de cl\^oture galoisienne dans le contexte d'une famille de courbes  lisses
ou pr\'estables.   Pour une  famille de courbes,  m\^eme lisses, cette op\'eration  exige  quelques  pr\'ecisions.   Soit un rev\^etement $\pi: C \to D$, les courbes $C,D$ \'etant  lisses  et  d\'efinies sur un corps alg\'ebriquement clos $k$, de  points de branchement   $Q_1,\dots,Q_b$.  

 L'hypoth\`ese sur la monodromie, assertion (ii)
de la d\'efinition 6.15, revient \`a fixer  en premier le groupe de
Galois d'une  cl\^oture galoisienne de $\pi: C\to D$,  c'est \`a dire   le groupe de Galois  $G$ de l'extension $k(C)/k(D)$.  Soit   $\phi: Z \to D$  une  $G$-cl\^oture galoisienne de $\pi$.  Cela signifie  d'abord que le rev\^etement $\phi: Z\to D$ est galoisien, et que   le groupe
des automorphismes de $\pi: Z\to D$  est  identifi\'e \`a  $G$.  Il y a
d'autre part  une factorisation (non unique) de $\phi$ en
$\phi = \pi h : Z \buildrel h \over\longrightarrow  C\buildrel \pi
\over \longrightarrow  D $. La classe de conjugaison du sous-groupe
$H =\Aut (h) \subset G$ est  par contre bien d\'efinie.   Comme on a
une cl\^oture galoisienne 
$$ \bigcap _{s\in G} s H s^{-1} = 1 \tag (6.28) $$

 Si  $\phi': Z' \to D$ est une autre cl\^oture galoisienne, et si  $\phi'
= \pi h'$,  la th\'eorie de Galois, ou du groupe fondamental,  nous
assure qu'il y a un isomorphisme   $\psi:
Z\buildrel\sim\over\rightarrow Z'$,    pas n\'ecessairement  un $G$
isomorphisme, tel que $h'\psi = h$. 

Soit   $\xi\in R_+(G)$ la donn\'ee de ramification
du $G$-rev\^etement galoisien $\phi: Z \to D$.   Nous d\'efinissons la monodromie  de $\pi: C\to D$ de la mani\`ere suivante:
 
\proclaim D\'efinition 6.17.  Un type de monodromie est un triplet $m = (G,H,\xi)$ compos\'e de 
 i)  un groupe fini $G$,  ii)  un sous-groupe $H\subset G$ tel que $\bigcap_{s\in G} sHs^{-1} = 1$ (condition 6.28), et iii)  une donn\'ee de ramification $\xi \in R_+(G)$.   
  Le rev\^etement $\pi: C\to D$  d\'efini sur $k$ est  dit \`a  monodromie $m = (G , H , \xi),$  s'il existe une $G$-cl\^oture galoisienne $\phi = \pi h: Z \to C\to  D,$ telle     $\Aut (h) = H$,   de  donn\'ee de ramification   $\xi$.   On parlera alors de  $m$-cl\^oture galoisienne.

 Posons (voir \S\  2.2.2)
$$\Aut  (m) = \{ \theta \in \Aut (G) , \quad \theta (H) = H, \quad
\hbox{\rm et}\quad \theta
(\xi) = \xi \} \tag (6.29) $$
 Soit $Z(G)$ le centre de $G$.  Du fait  de  (6.28),   $H\cap Z(G) = 1$, de la sorte  $ H$ se r\'ealise par automorphismes int\'erieurs comme sous-groupe distingu\'e 
$H \vartriangleleft \Aut (m)$, on posera 
$$\Delta(m) = \Aut(m)/H \tag (6.30)$$
Noter   que  ${\overline N}_G(H)  =
N_G(H)/Z(G) $ est un sous-groupe distingu\'e de $\Aut (m)$ contenant $H$. 
Fixons un type de monodromie $m$,  et soit un   rev\^etement $\pi: C \to D$ d\'efini sur  $k$ (alg\'ebriquement clos).

 Deux $m$-cl\^otures galoisiennes  $(Z,h, \phi)$ et $(Z',h',\phi')$ sont isomorphes  s'il existe un $C$-isomorphisme $G$-\'equivariant $f: Z \buildrel\sim\over\rightarrow Z'$  qui est l'identit\'e sur $C$.   On sait que pour deux telles cl\^otures,  il existe  un isomorphisme $f: Z
\buildrel\sim\over\rightarrow Z'$ tel que $h' f = h$.  Cela d\'efinit
un automorphisme $\theta$ de $G$  unique  tel que   $f: Z\to {Z'}^{\theta}$ est $G$-\'equivariant. Notons que cela force \`a
avoir $\theta  \in \Aut(m)$.      

  Soit  $G(\pi)$ le groupo\"\i de des $m$-cl\^otures galoisiennes de $\pi: C\to D$. Le groupe $\Aut(m)$ agit \`a droite sur $G(\pi)$ par torsion de l'action.    Notons la remarque \'evidente:
 
\proclaim  Lemme 6.18. $G(\pi)$ est un groupo\"\i de \'equivalent \`a $\Aut (m)$.

\dem   
  Il suffit de noter que si $Z$ et $Z'$ sont deux cl\^otures galoisiennes, il existe un isomorphisme $f:Z\to Z'$ tel que $h'f = h$, et alors  un unique $\theta \in \Aut(m)$ tel que   $f: Z^\theta \buildrel\sim\over \to Z'$ soit $G$-\'equivariant.  Pour $\theta$ donn\'e, $f$ est unique car  $Z(G) \cap H = 1$.  \hfill\break
  \qed 
 
 Notons  que $\theta \in \Aut (m)$ fixe la classe d'\'equivalence du  $G$-rev\^etement galoisien   
$ \phi: Z\to D$ si et seulement si $\theta \in {\overline N}_G(H)$
(automorphisme int\'erieur), conduisant en accord avec la th\'eorie
de Galois au fait que $N_G(H) / H$ s'identifie au groupe des
automorphismes de $\pi: C\to D$.  

La d\'efinition d'une cl\^oture galoisienne s'\'etend  aux rev\^etements entre courbes lisses, de base quelconque. On conserve le triplet $m = (G,H,\xi)$.

\proclaim D\'efinition 6.19.
Soit  $\pi:C \to D$  un rev\^etement   entre $S-$ courbes lisses. Une $m$-cl\^oture
galoisienne de $\pi$ est un $G$-rev\^etement galoisien $\phi: Z\to D$
de base $S$, augment\'e d'une  factorisation  
 $\phi: \, Z\buildrel h\over\longrightarrow C\buildrel
\pi\over\longrightarrow D$.  
On demande   l'\'egalit\'e $\Aut (h) = H$, et que la  monodromie   le long des fibres g\'eom\'etriques  est  fixe \'egale \`a $m$.

Le $G$-rev\^etement
$\pi: C\to D$   est  dit \`a monodromie $m = (G,H,\xi)$  si la monodromie est constante \'egale \`a $m$ le long des fibres g\'eom\'etriques, par exemple s'il existe une $m$-cl\^oture galoisienne.
Les rev\^etements entre courbes lisses, de genres fix\'es, et  \`a
monodromie $m$ fix\'ee forment de mani\`ere \'evidente   un champ  de
Deligne-Mumford ${\Cal H}_{g,g',m}$,   except\'e si $g = g' = 0$ et $b = 2$.  C'est essentiellement une cons\'equence de la th\'eorie du groupe fondamental mod\'er\'e, comme not\'e dans \cite{42}, ou de mani\`ere plus d\'etaill\'ee dans \cite{69}.   Si $\phi: Z \to
D = Z/G$ est un $G$-rev\^etement galoisien  de donn\'ee de Hurwitz $\xi$,
le rev\^etement $C= Z/H \to  D$ est \`a monodromie $ m$,
et $\phi$ en est une $m$-cl\^oture galoisienne.  Cela d\'efinit un
morphisme de champs 
$$\imath: {\Cal H}_{g,G,\xi} \longrightarrow {\Cal
H}_{g',g,m} \tag (6.31)$$  
  Le groupe $\Aut(m)$ agit sur le champ ${\Cal H}_{g,G,\xi}$ (d\'efinition 6.1) par torsion de l'action.  Comme l'action du sous-groupe distingu\'e $H$ est triviale dans le sens de la section 6.1.1, il en r\'esulte une action de $\Delta(m)$, non stricte.    

 \proclaim  Th\'eor\`eme 6.20.  Le morphisme  (6.31) est  repr\'esentable. Il identifie  
   ${\Cal H}_{g',g,m} $   au champ  ''quotient''  
 $\overline \imath: {\Cal H}_{g',g,m}  \buildrel\sim\over  \rightarrow  [{\Cal H}_{g,G,\xi}/ \Delta(m)] \cong  [{\Cal H}_{g,G,\xi}/ \Aut(m)]//H.   $
  
 \dem  On supprime le pr\'efixe $m$.        
 Montrons d'abord que $\imath$ (6.31) est un epimorphisme, i.e. localement pour la topologie \'etale $\pi$ admet une cl\^oture galoisienne.  C'est en fait une cons\'equence de la th\'eorie du groupe fondamental (voir \cite {20} par exemple). On  peut pr\'ef\`ere utiliser un argument de d\'eformation, bien qu'essentiellement \'equivalent.  Soit $\pi:C_s\to D_s$ une fibre g\'eom\'etrique   d\'efinie sur une cl\^oture alg\'ebrique de $k(s)$. Soient   $(Q_j)_{1\leq
j\leq b}$ les points de branchement. Comme les courbes sont
lisses, on sait que le foncteur des d\'eformations formelles de $\pi_s:
C_s \to D_s$ est  isomorphe \`a celui de la courbe marqu\'ee $
(D_s,B=\sum_{j=1,\dots,b} Q_j)$.  Soit alors $h: Z\to C_s$ une cl\^oture galoisienne de $\pi: C_s\to D_s$. Consid\'erons la d\'eformation universelle $\Cal Z \to \Cal D$ du $G$-rev\^etement $\phi: Z \to D_s$, et soit $\Cal B\subset \Cal D$ le diviseur de branchement r\'eduit. On sait, voir par exemple le th\'eor\`eme 5.5, que $(\Cal D , \Cal B)$ est la d\'eformation universelle de la courbe marqu\'ee $(D_s,B_s)$.  On sait que la base de cette d\'eformation est formellement lisse.  On sait par ailleurs  que le foncteur des d\'eformations du rev\^etement $\pi_s$ peut \^etre identifi\'e \`a celui de la base marqu\'ee par les points de branchements \cite{68}. Donc  $\Cal Z/H \to \Cal D$ est la d\'eformation universelle de $\pi$, et $\Cal Z \to \Cal Z/H$ en est une cl\^oture galoisienne.  Des arguments  classiques d'alg\'ebrisation montrent que cette cl\^oture galoisienne existe sur un voisinage \'etale du point $s\in S$.  

Prouvons la repr\'esentabilit\'e de $\imath$. Il suffit de prouver que le morphisme diagonal    est un monomorphisme (\cite{51}, Prop. 4.4). En d'autres termes   si $\alpha \in \Aut(Z)$ pour un $G$-rev\^etement galoisien $\varphi: Z\to D$ de base $S$, et si $\imath(\alpha) = 1$, alors $\alpha = 1$.  Le groupe $\Aut_D(Z)$ des $G$-automorphismes de $Z$ qui induisent l'identit\'e sur $D$ est le groupe constant $G\times  S$. En effet ce groupe est fini non ramifi\'e sur $S$, contient $G\times S$. Comme sur les fibres g\'eom\'etriques il y  a \'egalit\'e, ces deux groupes sont \'egaux.  Donc par ce type d'argument $\alpha\in H\cap Z(G) = 1$.

 Passons \`a la seconde assertion.   
  Du fait de la propri\'et\'e universelle du quotient, le morphisme $\imath$ factorise en 
$ {\Cal H}_{g,G,\xi}  \to [{\Cal H}_{g,G,\xi}/\Aut (m)] \to {\Cal H}_{g',g,m}.$ Pour tout objet $Z\to D \in {\cal H}_{g,G,\xi}$, et tout $\sigma\in H$, on a $\sigma: Z\cong Z^\sigma$, on a $H\subset \Aut [Z]$, en notant $[Z]$ l'image de $Z$ dans le quotient. Il s'ensuit une factorisation de $\imath$
   $$\overline\imath:  [{\Cal H}_{g,G,\xi}/\Delta(m)] =  [{\Cal H}_{g,G,\xi}/\Aut(m)]//H  \to    {\Cal H}_{g',g,m}$$  
 qui est un epimorphisme. Pour conclure que $\overline\imath$ est un isomorphisme, il suffit de prouver que c'est un monomorphisme. Le champ de gauche Žtant un quotient, cela se r\'eduit ˆ prouver que  $\pi: C\to D$ et $\pi': C'\to D'$  sont deux rev\^etement de base connexe  $S$ et si $\phi:  Z  \to C\to D, \, \phi': Z'\to C'\to D'$  sont  deux  cl\^otures galoisiennes, alors on a un isomorphisme de  sch\'emas
 $$\coprod_{\theta\in \Delta(m)} \Isom_G(Z^\theta,Z')  \cong \Isom((C\to D),(C'\to D'))$$
 Dans la somme de gauche $\theta \in  \Delta (m)$ signifie un systme de repr\'esentants  de $\Aut(m)$ modulo $H$.  Le $S$-sch\'ema
$\Isom_{G,C}   (Z^\theta,Z')$  a pour points
 $ {\Isom}_{G,C} (Z^\theta,Z')(T) = \{ f: Z^\theta \times _ST
\buildrel\sim\over\longrightarrow
Z'\times_ST \, , \, h'f = h\}$.

Il est clair que $\Isom_{G}   (Z^\theta,Z')$ qui est un sous-sch\'ema ferm\'e de $\Isom (Z,Z')$, est   fini  et non ramifi\'e
sur $S$ \cite {20}, de sorte que les membres sont finis non ramifi\'es sur $S$. Le lemme 6.18 montre que sur les fibres g\'eom\'etriques  le morphisme est  bijectif. Pour conclure qu'on a bien un isomorphisme, il suffit de prouver que si $C = C'$, le sous-sch\'ema fibre au-dessus de l'identit\'e de $C/S$, soit $\coprod_{\theta\in \Delta(m)} \Isom_{G,C}(Z^\theta,Z') $ est isomorphe \`a $S$. Il est clair qu'il y a un $\theta$ unique tel que $\Isom_{G,C}(Z^\theta,Z')\to S$ est fini, non ramifi\'e  bijectif.   Il est en fait \'etale, donc un isomorphisme.  C'est local sur $S$, et devient clair si on passe \`a la d\'eformation universelle comme dans la premi\`ere partie, car elle est lisse, donc r\'eduite. \hfill\break\qed 

   \rema {6.6}
   Pour un rev\^etement  de base $S$, l'obstruction \`a l'existence   d'une $G$-cl\^oture galoisienne  dans le sens  (6.19), est dans $H^1(S,\Delta (m))$ (voir section 6.2).    $\lozenge$

   \bigskip

  {\it 6.5.3. Cl\^oture galoisienne: cas des courbes
(pr\'e)stables} 
\bigskip

La discussion pr\'ec\'edente sugg\`ere   qu'une d\'efinition alternative   des rev\^etements
admissibles,  les rev\^etements    stables,  est comme suit:

\proclaim D\'efinition 6.21.   Soit   $m  = (G , H ,\xi)$ un type de monodromie. Un rev\^etement  stable de base $D$ d\'efini  au dessus du sch\'ema
de   base $S$, de monodromie $m $, est un rev\^etement
obtenu par factorisation $ \pi: C = Z/H \longrightarrow D =  Z/G $  d'un $G$-rev\^etement galoisien  $\phi: Z\to D$ \`a donn\'ee de ramification $\xi$.

  Le  $G$-rev\^etement  
galoisien   $\phi: Z\to D$,   cl\^oture galoisienne de $\pi: C \to D$, doit \^etre  vu comme  un \'el\'ement  structurel du rev\^etement stable $\pi: C\to
D$.    
Noter que dans cette d\'efinition $Z$ \'etant stable marqu\'ee par
les points de ramification du
rev\^etement $Z \to  Z/G$, la courbe $D$ est en cons\'equence stable marqu\'ee
par les points de branchement de $\pi: C \to D$, et  $C$
est stable marqu\'ee par les  pr\'eimages des points de
branchement de $Z\to D$, ces points pouvant ne pas \^etre des points de ramification de $\pi$. D'autre part le rev\^etement $Z \to C$ est
un $H$-rev\^etement galoisien stable marqu\'e par un diviseur
$H$-invariant qui en g\'en\'eral contient strictement les seuls
points de ramification de $Z \to C$. 

Il est clair qu'un rev\^etement
stable est admissible.  Du fait de l'imposibilit\'e de recoller les uniformisations
galoisiennes locales
un rev\^etement stable  ne sera  donc  en d\'efinitive  que  localement pour
la topologie fppf (\'etale suffit)  de la forme 6.21.  

\proclaim D\'efinition 6.22.
Soit $m = (G,H,\xi)$ un type de monodromie. Un rev\^etement  stable au dessus de la base $S$, de monodromie  $m$, est   une   section  au-dessus de $S$ du  champ ''quotient''  
$$[\overline{\Cal H}_{g,G,\xi}/\Delta(m)] = [\overline{\Cal H}_{g,G,\xi} /\Aut (m)]//H \tag (6.32)$$  
Le champ de Hurwitz
compactifi\'e $\overline {\Cal H}_{g,g',m}$  param\'etrant les rev\^etements stables de monodromie $m$, entre courbes stables marqu\'ees de genres respectifs $g, \, g'$ est  par d\'efinition le champ quotient   
$\overline {\Cal H}_{g,g',m} = [\overline {\Cal H}_{g,G,\xi} / \Aut (m)]// H = [\overline{\Cal H}_{g,G,\xi}/\Delta(m)]  $.

  Soit un rev\^etement
stable   $\pi: C\to D$  d\'efini sur le corps
  $k$. On suppose   qu'il d\'erive
 d'un rev\^etement galoisien stable
$\tilde\pi : \tilde C \to D$.  Nous allons   comparer la d\'eformation formelle verselle du rev\^etement
admissible $\pi: C\to D$ \`a la d\'eformation verselle  de sa cl\^oture galoisienne $\phi: \tilde C \to D$, qui est aussi celle du rev\^etement stable  $\pi$, c'est \`a dire  enrichi  par
l'uniformisation galoisienne $$\phi: \tilde C \buildrel h \over
\rightarrow  C\buildrel \pi \over  \rightarrow  D$$
 Notons
$Q_1,\dots,Q_r$  les points doubles de $D$, et   $\{
Q_{\alpha,j} \}$ les points de $ C$ au dessus de $Q_\alpha$. Soit
 $\tilde Q_\alpha$ un quelconque point de $\tilde C$ au dessus
de $Q_\alpha$. On  notera que par d\'efinition les $Q_{\alpha,j}$
sont les points doubles de $C$.  Notons   $d_\alpha$ l'ordre du
stabilisateur de $\tilde Q_\alpha $, et enfin
$d_{\alpha,j}$ l'indice ("de ramification") d\'ecrivant la structure
locale (d\'efinition 6.15 (ii))
de $\pi$ en $Q_{\alpha,j}$.  On a le r\'esultat  suivant dont
la preuve est en tout point analogue au cas lisse:

\proclaim Lemme 6.23.  Sous les hypoth\`eses qui pr\'ec\`edent,  on
a $d_\alpha =
\ppcm _j{d_{\alpha,j}}$.

\dem Le r\'esultat  d\'ecoule de la cyclicit\'e de $I = G_{\tilde
Q_\alpha}$, et de la condition
(6.28). On a d'abord clairement pour tout $j$, $d_{\alpha,j } / q_\alpha$, donc
$\ppcm {\{d_{\alpha,j}\}_j}$ divise $d_\alpha$.  Par ailleurs la fibre
$\tilde\pi ^{-1} (Q_\alpha) = G.\tilde Q_\alpha \cong  G/I$
est une r\'eunion de $H$-orbites ${\Cal O}_1,\dots,{\Cal
O}_{r_\alpha}$, o\`u ${\Cal O}_j$ est la
$H$-orbite des points doubles de $\tilde C$ qui sont au dessus de
$Q_{\alpha,j}$. Soit
$\delta_j$ ''l'indice de ramification'' (l'ordre du stabilisateur)
relativement \`a $H$ des points
de ${\Cal O}_j$.; on a donc $\delta_j = {d_\alpha \over
d_{\alpha,j}}$. Le r\'esultat se ram\`ene
en fait \`a prouver que
$\pgcd {( \delta_1, \dots, \delta_{r_\alpha})} = 1$.
  Notant $\delta$ ce pgcd, si $K$ est le sous-groupe de $I$ d'ordre
$\delta$, on voit
de suite que $K \subset \bigcap_{s\in G} sHs^{-1} = 1$. D'o\`u la
conclusion. \hfill\break\qed 

D'apr\`es  Harris-Mumford la base de la d\'eformation universelle du rev\^etement
admissible $\pi: C \to D$ est
(\cite {42}, p 62):
$$R_\pi = W(k) [[ \tau_1,\dots,\tau_r,\dots,\tau_{3g'-3+b},
\{\tau_{\alpha,j}\} ]] / {\Cal I}
\tag (6.33) $$
o\`u $\Cal I$ d\'esigne l'id\'eal engendr\'e par les relations
$\tau_\alpha = \tau_{\alpha,j}^{d_{\alpha,j}}, \quad (\alpha =
1,\dots,r; \,\,
j=1,\dots,r_\alpha).   $
Dans cette description, $\tau_\alpha$ repr\'esente le param\`etre de
d\'eformation du point double
$Q_\alpha \in D$; on sait par ailleurs que la base de la d\'eformation universelle de la courbe
marqu\'ee $(D,\{ Q_\alpha\})$ est \quad $W(k) [[
\tau_1,\dots,\tau_{3g'-3+b} ]]$. L'anneau $R_\pi$  est r\'eduit, mais en g\'en\'eral pas int\'egralement clos dans son anneau total des fractions. La structure de $R_\pi$ peut \^etre pr\'ecis\'ee comme suit:

\proclaim Proposition 6.24. Posons $N = 3g' - 3 + b$.  Pour tout entier $d$ premier \`a $p$, soit $\mu_d \subset W$ le groupe des racines d-i\`eme de l'unit\'e. \hfill\break
i) L'anneau $R_\pi$ est r\'eduit, ses id\'eaux premiers minimaux sont les noyaux $\Cal P_\zeta$ des morphismes 
$$\varphi_\zeta= R_\pi \to W[[t_1,\cdots,t_N]], \,\,\tau_i \mapsto t_i^{d_i}, \,\, \tau_{i,j} \mapsto \zeta_{i,j}t_i^{r_{i,j}}$$
o\`u, $\zeta = \{\zeta_{i,j}\} \in \prod_{i,j} \mu_{d_{i,j}}, \,\, r_{i,j} = d_i/d_{i,j}$. On a $\Cal P_\zeta = \Cal P_{\zeta'}$ si et seulement si $\zeta'_{i,j} = \epsilon_i^{r_{i,j}} \zeta_{i,j}$ pour des racines de l'unit\'e $\epsilon_i \in \mu_{d_i}$. En particulier le nombre des id\'eaux premiers minimaux est 
$$N_\pi = \prod_i\left( {\prod_j  d_{i,j}\over d_i}\right)\tag (6.34)$$
ii) La fermeture int\'egrale de $R_\pi$ dans son anneau total des fractions est le produit de $N_\pi$ copies de $W[[t_1,\cdots,t_N]]$.

 \dem Posons $R = W[[t_1,\cdots,t_N]]$.  Si $\zeta'$ et $\zeta$ sont reli\'es comme indiqu\'e, il est clair que $\Cal P_\zeta = \Cal P_{\zeta'}$ car si $\psi$ est l'automorphisme de $R$ tel que $\psi (t_i) = \epsilon_i t_i$, alors $\varphi_{\zeta'} = \psi\varphi_\zeta$. Montrons que si $\varphi = \prod \varphi_\zeta$ alors $\ker \varphi = 0$. On peut se limiter \`a $N = 1$, ce qui permet de supprimer l'indice $i$, et de supposer que $1\leq j \leq m$ . Posons $\mu = \prod \mu_{d_j} $, et $\tau = \tau_1$. L'alg\`ebre $R_\pi$ a pour base sur $W[[\tau]]$ les mon\^omes $\prod_j \, \tau_j^{\alpha_j}, \, \alpha_j < d_j$. De la sorte un \'el\'ement du noyau est 
$$\xi = \sum_{k,\alpha = \alpha_1,\cdots,\alpha_m} \quad a_{\alpha,k} \prod_j \tau_j^{\alpha_j} $$
tel que pour  $l \geq 0$
$$   \sum_{k,\alpha = \alpha_1,\cdots,\alpha_m} \quad a_{\alpha,k} \prod_j \zeta_j^{\alpha_j}  = 0$$
la somme \'etant \'etendue aux couples $(k,\alpha)$ tels que $kd + \sum_j\alpha_jr_j = l$.  Si on regarde $ \prod_j \zeta_j^{\alpha_j} $ comme d\'efinissant un caract\`ere sur le groupe $\mu$, le th\'eor\`eme d'ind\'ependance des caract\`eres conduit \`a $a_{\alpha,k} = 0$. 

Observons   que $R$ est la fermeture int\'egrale de l'image de $\varphi_\zeta$.  Du fait
que  $\pgcd {(r_1,\cdots,r_m)} =
1$, il existe  avec  $ e_j \in {\Bbb Z}$, une relation de la forme:
$\sum_j e_j  r_j = 1$, qui entraine la relation
$ t = \prod_j \tau_j ^{e_j} $.  Le r\'esultat est alors clair.  De cela on tire le fait que si $\varphi_\zeta$ et $\varphi_{\zeta'}$ ont un noyau identique, alors il y a un $W$-automorphisme $\gamma$ de $R$ tel que $\varphi_{\zeta'} = \gamma\varphi_\zeta$. Alors n\'ecessairement $\gamma (t) = \epsilon t$ pour une certaine racine
d-i\`eme de l'unit\'e. Cela termine la preuve de la proposition\hfill\break
\qed 

La structure de l'anneau des d\'eformations $R_\pi$  explique   la non unicit\'e d'une cl\^oture galoisienne pour un rev\^etement  stable $\pi: C\to D$ d\'efini sur   $k$.

\proclaim Th\'eor\`eme 6.25.   Un rev\^etement admissible $\pi: C\to D$ au-dessus de $k$ admet une cl\^oture galoisienne, i.e. peut \^etre stabilis\'e.  
Le nombre de classes    de cl\^otures galoisiennes (sous l'action de $\Delta (m)$)  de  $\pi$ est $N_\pi = \prod_i\left( {\prod_j  d_{i,j}\over d_i}\right)$, c'est \`a dire le nombre de facteurs int\`egres de la fermeture int\'egrale de $R_\pi$. Dans cette correspondance, les facteurs (isomorphes \`a $W[[t_1,\cdots,t_N]]) $ correspondent bijectivement aux  anneaux de d\'eformation des cl\^otures galoisiennes respectives.

\dem   Admettons d'abord le premier point. Soit une $G$-cl\^oture galoisienne $Z \to C\to D$ de $\pi$, de  d\'eformation   universelle $\Cal Z\to \Cal D$. On sait (th\'eor\`eme 5.2) que la base de cette d\'eformation   s'identifie \`a $\Spec (R_Z)$, avec  
  $R_Z = W(k) [[t_1,\dots,t_r, \dots,t_{3g'-3+b}]]$. Posons $N = 3g' - 3+b$. Si $R = W[[\tau_1,\cdots,\tau_N]]$ est la base de la d\'eformation universelle $\Cal D$ de $D$, on peut supposer que le   morphisme discriminant $R \to R_Z$ est donn\'e par 
  $\tau_\alpha =
t_\alpha^{d\alpha}$ si $\alpha = 1,\dots,r$, et $\tau_\alpha =
t_\alpha$ si $\alpha > r$.

La consid\'eration de la d\'eformation $\Cal Z/H$ de $C$ conduit \`a une factorisation
$R \to R_\pi \buildrel\varphi\over \to R_Z$.  
On peut supposer que $\varphi$ est donn\'e par 
$$\tau_{\alpha,j} =
t_\alpha^{d_\alpha \over d_{\alpha,j}}, (\alpha = 1,\dots,r ;
j=1,\dots,r_\alpha)$$
Soit une autre cl\^oture galoisienne $Z' \to C\to D$, telle pour $\theta\in \Aut(m)$, on a un isomorphisme de cl\^otures 
$f: Z\buildrel\sim\over\rightarrow {Z'}^\theta$. Si $\imath': Z' \hookrightarrow {\Cal Z}'$ est le plongemenr canonique de $Z'$, la consid\'eration de $({\Cal Z}'\, , \imath' f)$ comme une autre d\'eformation universelle de $Z$, permet d'\'etendre $f$ en  un isomorphisme de diagrammes
$$  \commdiag{
  { \Cal Z }&  \mapright^f_\sim  &  {\Cal Z'}^\theta\cr
   \mapdown && \mapdown \cr
 \Spec R_Z & \mapright^f_\sim &  \Spec R_{Z'}
  \cr}$$
  Passant au quotient par $H$ il s'ensuit   un diagramme commutatif
   $$\commdiag{
   \Spec R_Z  &  \mapright^f_\sim &  \Spec R_{Z'} \cr
    \mapdown &&\mapdown \cr
     \Spec R_\pi & \mapright^=  &  \Spec R_\pi
  \cr} $$
  Cela montre que les deux cl\^otures galoisiennes conduisent \`a la m\^eme classe de morphisme $R_\pi \to W[[t_1,\cdots,t_N]]$, donc au m\^eme facteur de la fermeture int\'egrale. 

Prouvons la r\'eciproque. Cela revient \`a supposer qu'on a en mains les diagrammes commutatifs de dessus, exept\'e l'isomorphisme $f: \Cal Z \buildrel \sim\over\to \Cal {Z'}^\theta$. On peut  se ramener \`a $R_Z = R_{Z'}$.  De la sorte on a deux $G$-rev\^etements $\Cal Z $ et $\Cal Z'$ au-dessus d'une m\^eme base $R_Z$, avec $\Cal Z/H \cong \Cal C\otimes R \cong \Cal Z'$.  On sait que pour les cl\^otures galoisiennes induites au-dessus de la fibre g\'en\'erique $\nu$ (lisse) de $\Spec R_Z$, il y a $\theta \in \Aut (m)$ tel que  $\Cal Z_\nu \cong {{\Cal Z'_\nu}^\theta}$. Les courbes \'etant stables marqu\'ees, cet isomorphisme se prolonge en un isomorphisme global 
$\Cal Z \cong {\Cal Z'}^\theta$.
Par restriction \`a la fibre sp\'eciale, on a obtient  un isomorphisme $Z \cong {Z'}^\theta$ de cl\^otures galoisiennes de $\pi$. Pour terminer la preuve, reste \`a prouver que tout rev\^etement admissible $\pi: \overline C\to \overline D$ admet une cl\^oture galoisienne, i.e. peut \^etre stabilis\'e.

On peut dans ce but supposer qu'on a  une courbe stable marqu\'ee $D \to S = \Spec
(R)$, o\`u $R$ est un anneau de valuation discr\`ete complet, de
corps r\'esiduel $k$, le marquage \'etant d\'efini par un diviseur \'etale $B \to S$, et de fibre sp\'eciale $\overline D$.    On suppose que la fibre g\'en\'erique est lisse.
 Il y a des obstructions \`a
 trouver un rev\^etement admissible $C \rightarrow D$, relevant
$\overline C \rightarrow \overline D$. Par contre on peut
trouver un tel rev\^etement apr\`es
extension  finie de $R$. Quitte \`a agrandir $R$, faisons le choix
d'un tel rel\`evement $C \to D$.
Rappelons que si $B = \sum_j Q_j$, cela implique l'existence de
points $P_{i,j}$ au dessus de $Q_j$ en lesquels est concentr\'ee la
ramification. La courbe $C$ n'est pas forc\'ement stable marqu\'ee
par les $P_{i,j}$, mais elle l'est par l'ensemble de toutes les  pr\'eimages des $Q_j$. Soit
$\tilde C_K \rightarrow  C_K \rightarrow D_K$
une cl\^oture galoisienne de $C_K \to D_K$, les groupes de Galois
respectifs \'etant $G$ et $H$.  Quitte \`a agrandir si n\'ecessaire
$R$, soit $\tilde C \rightarrow D$ le mod\`ele stable marqu\'e.  Du
fait que $\tilde C /H$ est stable marqu\'ee par les pr\'eimages des
$Q_j$, on a $\tilde C /H = C$, il en d\'ecoule  une factorisation
$$ \tilde C \rightarrow C \rightarrow D\quad \hbox{ \rm et } \quad C =
\tilde C / H$$
Si on prend les fibres sp\'eciales, ceci montre que $\overline C$
poss\`ede une cl\^oture galoisienne, en d'autres termes, qu'on peut
enrichir le rev\^etement admissible $\pi: \overline C \to \overline
D$ en un rev\^etement stable.  
\hfill\break
  \qed

Il y a un  morphisme naturel de source le champ de Hurwitz $\overline {\Cal
H}_{g,g',m}$ de but le champ  $\overline {\Cal
H}^{\rm adm}_{g,g',m}$ des rev\^etements admissibles au sens
de Harris-Mumford \cite{42}.   On a le r\'esultat  (comparer avec l'approche de \cite{2}) 
 \proclaim  Proposition 6.26.
 $\overline {\Cal
H}_{g,g',m}$   est une  d\'esingularisation
(normalisation) de $\overline {\Cal
H}^{\rm adm}_{g,g',m}$.  
 \hfill\break 
\qed


 Le champ $\overline {\Cal H}_{g,g',m}$ (classifiant les rev\^etements stables) 
n'est  pas  en g\'en\'eral isomorphe au champ, plus ''fin'' des rev\^etements log-admissibles (comparer avec Wewers  \cite {69}).
   Si la base $S$ d'un rev\^etement stable est arbitraire, on ne peut
en  g\'en\'eral d\'ecrire un tel rev\^etement comme un quotient $\pi:
Z/H \to Z/G$, partant d'un $G$-rev\^etement  stable. On peut
cependant, partant de la d\'efinition abstraite (6.8),  donc d'un
objet de $\overline {\Cal H}_{g,g',m} (S)$, donner un contenu
plus visible \`a la structure d'un tel rev\^etement.  Les d\'efinitions  de la section 6.1 relatives aux quotients, disent qu'un
objet de ce champ  peut \^etre d\'ecrit
comme  d'une part la  donn\'ee d'un recouvrement $(S_i\to S)$ de $S$,
et ensuite pour tout $i$, d'un $G$-rev\^etement $Z_i \to Z_i/G$ de
base $S_i$, l'action de $G$  sur $Z_i$ \'etant fix\'ee \`a une
torsion  pr\`es par un \'el\'ement $\Aut(m)$.  L'objet inclut en plus
une donn\'ee de descente, qui en posant $S_{i,j} = S_i\times_S S_j$
se ram\`ene \`a  une collection d'isomorphismes:
$\sigma_{j,i}: {Z_i}_{\vert S_{i,j}} \buildrel\sim\over\rightarrow
{Z_j}_{\vert
S_{i,j}}^{\theta_{j,i}}$ 
qui v\'erifient  apr\`es restriction \`a $S_{i,j,k}$   la condition
de cocycle  
 $$ \sigma_{k,j}^{\theta_{j,i}}
\circ \sigma_{j,i} = \sigma_{k,i}, \quad \theta_{k,j}\theta_{j,i} =
\theta_{k,i}\tag (6.35)$$
 L'\'egalit\'e (6.35) est \`a prendre \`a un automorphisme interieur
pr\`es par un \'el\'ement  $h_{i,j,k} \in H$. En particulier,
passant aux quotients par $G$ et $H$, et faisant abstraction de la torsion de l'action
  on obtient un
diagramme
  $$\commdiag{
  {Z_i}_{\vert S_{i,j}}/H   &\mapright^{\sigma_{i,j}} &  {Z_j}_{\vert
S_{i,j}}/H \cr
   \mapdown &&\mapdown \cr
{Z_i}_{\vert S_{i,j}}/G   &\mapright^{\sigma_{i,j}}  & {Z_j}_{\vert S_{i,j}}/G
  \cr} $$

les isomorphismes horizontaux d\'efinissant des donn\'ees de descente
sur les   courbes d\'efinies localement par $(Z_i/H)$, et
$(Z_i/G)$. De la sorte on r\'ecup\`ere par descente \`a $S$ un
rev\^etement admissible $C\to D$, et aussi un rev\^etement $Z \to D$,
qui est localement un $G$-rev\^etement galoisien. Pour installer une
action globale de $G$, il y a lieu de passer \`a un recouvrement \'etale de $S$, alors dans ce cas $C = Z/G$ et $D = Z/G$.  

\bigskip
   {\it 6.5.4.  Le morphisme
discriminant: cas  g\'en\'eral} 
\bigskip

Il a \'et\'e observ\'e dans la section 3.1 que le
discriminant d'un  $G$-rev\^etement galoisien $\pi: C\to D$, entre
courbes lisses, de base $S$,  pouvait s'exprimer comme une somme avec
multiplicit\'es de diviseurs  de Cartier relatifs \'etales sur $S$.
Ces diviseurs  sont les images des diviseurs  (\'etales sur $S$):
$$\Delta_H = \hbox {\rm  diviseur des points fixes de stabilisateur
\'egal \`a } \, H$$de sorte que pour tout sous-groupe cyclique $K$
de $G$,   $Fix (K) = \sum_{K\subset C} \, \Delta_C$.  Ce fait
s'\'etend  \`a un $G$-rev\^etement  galoisien stable (voir remarque
section 4.1), en se limitant alors \`a la partie horizontale
$Fix_{hor} (K)$ (Proposition 4.8). Dans le cas g\'en\'eral, c'est \`a
 dire des rev\^etements non galoisiens, on a implicitement suppos\'e
 la validit\'e d'une telle description (D\'efinition 6.21).  Nous
revenons sur ce point. Soit un rev\^etement de base $S$, et degr\'e
d, qu'on peut sans perte de g\'en\'eralit\'e  supposer de la
forme
$\phi: Z\buildrel p\over\rightarrow C = Z/H
\buildrel\pi\over\rightarrow D = Z/G, $
donc de cl\^oture galoisienne au-dessus de $S$, 
$Z\buildrel\phi\over\rightarrow D = Z/G$.
 Soit  le diviseur de ramification (\S \ 4.6) d\'efini par
 $$R_\pi = \Div (\pi^*(\omega_{D/S})\to \omega_{C/S})\tag (6.36)$$    
Le diviseur de branchement est    $B_\pi =
\pi_* (R_\pi)$.   Notons  $[K]$ la $H$-classe de conjugaison du
sous-groupe cyclique $K$. 


\proclaim Proposition 6.27. i) Le
diviseur de ramification du rev\^etement $\pi: C\to D$ de base $S$, a
pour expression
$$R_\pi = \sum_{[K]} \, \left(   {{\vert K\vert }\over 
{\vert H\cap K\vert } }- 1\right) \, R_{[K]} \tag (6.37)$$
pour des
diviseurs de Cartier relatifs $R_{[K]}$ disjoints et \'etales sur $S$.
 Le terme  $  {{\vert K\vert }\over {\vert H\cap K\vert }}$
repr\'esente un indice de ramification.\hfill\break
ii) Le diviseur de branchement est $B_\pi = \sum_{[K] \nsubseteq H} \,\,  \left( d
\,-\, e_{H,K}\right) \, B_{[K]},$ expression dans laquelle
$e_{H,K} $ signifie  un indice de ramification, et $B_{[K]}$ un
diviseur de Cartier relatif \'etale sur $S$.  

\dem
{Preuve} i) On a clairement $R_\phi = p^* (R_\pi) + R_p$. Utilisant
l'expression de la section (3.2) pour le diviseur de ramification
dans le cas galoisien, c'est \`a dire $ R_\phi = \sum_K \left(
\vert K\vert - 1\right) \, \Delta_K$, et l'expression analogue pour
$R_p$, on obtient de suite 
$$p^*(R_\pi) = R_\phi - R_p =
\sum_{K,K\not\subset H} (\vert K\vert - \vert H\cap K\vert ) \Delta_K
\tag (6.38)$$
D'autre part les diviseurs $\Delta_K, \Delta_{K'}$ ont une m\^eme
image dans $C$ si et seulement si $K'$ est conjugu\'e \`a $K$ par
un \'el\'ement de $H$. Le nombre de ces conjugu\'es est $[H :
N_G(K)\cap H] $, et le degr\'e de $\Delta_K$ sur son
image $R_{[K]}$ est $[N_G(K)\cap H : K\cap H]$.
 Comme  $p_*p^* (R_\pi) = \vert H\vert R_\pi$, il vient 
 $$\vert H\vert R_\pi =  \sum_{ K, K\nsubseteq H} \,  (\vert K\vert  
 - \vert H\cap K\vert)  \,p_*(\Delta_K) = \vert H\vert \sum_{[K]} \, \left(  { {\vert K\vert }\over
{\vert H\cap K\vert }} - 1\right) \, R_{[K]} $$ 
\hfill\break ii) La preuve de 6.37 en tout point analogue est omise.\hfill\break
\qed

 \exa {6.7}  Soit le cas du champ  
de Hurwitz, classifiant les rev\^etements simples de degr\'e $n\geq 3$ de
$\Bbb P^1$, avec $b = 2g + 2n - 2$ points de branchement. Le groupe
de monodromie est $G = {\bf S}_n$, et $H =   {\bf S}_{n-1}$. Dans ce
cas $\Delta$ se r\'eduit \`a $H$ (automorphismes int\'erieurs). La
compactification lisse du champ de Hurwitz-Fulton est alors
$\overline {\Cal H}_{g,   {\bf S}_n, b(12)} /  {\bf  S}_{n-1}$. 

$\lozenge$

Soit $m = (G,H,\xi)$ un type de monodromie, et $b$ le degr\'e du diviseur de branchement.   Le champ de Hurwitz $\overline {\Cal
H}_{g,g', m}$ est  comme dans le cas galoisien la source  du morphisme  
discriminant: 
$$\delta : \, \overline {\Cal H}_{g,g', m}
\longrightarrow \overline {\Cal M}_{g',b} \tag (6.39)$$ explicitement
$\delta (\pi: C\to D)  = (D,B=\sum B_{[K]})$. Le morphisme   (6.39) est un morphisme plat, propre, quasi-fini et g\'en\'eriquement \'etale, i.e \'etale sur un sous-champ ouvert partout dense.  On peut voir $\delta$ comme un rev\^etement entre champs de Deligne-Mumford lisses. Sous ces conditions, il n'est pas difficile de d\'efinir le diviseur de ramification de $\delta$.
\proclaim Proposition 6.28. Soit un rev\^etement $\delta: {\cal H} \to {\cal M}$ entre champs de Deligne-Mumford lisses  sur $k$ (un morphisme quasi-fini, g\'en\'eriquement \'etale). Il existe un diviseur de Cartier effectif ${\cal R}\subset {\cal H}$ (le diviseur de ramification), tel que 
$$\omega_{\cal H} \cong \delta^*(\omega_{\cal M})\otimes {\cal O}({\cal R}) \tag (6.40)$$

\dem Le champ ${\cal H}$ (resp. ${\cal M}$) \'etant lisse, de dimension $N$ disons, le faisceau localement libre $\Omega^1_{{\cal H}/k}$ (de rang $N$) est d\'efini (resp.   $\Omega^1_{{\cal M}/k}$).  Soit $q: V\to {\cal M}$ un atlas. Formons le 2-produit fibr\'e 
$$\commdiag{ {\cal H} &\mapright^\delta & {\cal M}\cr
\mapup &&\mapup\lft{q}\cr
{\cal H}\times_{\cal M} U&\mapright&V\cr}$$
avec fl\^eches verticales \'etales surjectives. Si $U\to  {\cal H}\times_{\cal M} U$ est un atlas, de sorte que le morphisme compos\'e $p: U\to {\cal H}$ est un atlas de ${\cal H}$, on obtient le carr\'e 2-commutatif
$$\commdiag{{\cal H}&\mapright^\delta &{\cal M}\cr
\mapup\lft{p}&&\mapup\lft{q}\cr
U&\mapright^\sigma&V\cr}\tag (6.41)$$
Le morphisme $\sigma: U\to V$ est quasi-fini, g\'en\'eriquement \'etale, par cons\'equent  le conoyau $\Omega^1_{U/V}$ de 
$d\sigma: \sigma^*(\Omega^1_{V/k}) \to \Omega^1_{U/k}$
est de torsion. Noter que $\sigma^*$ une injection. On peut former  $R = \Div (\omega^1_{V/U})$, qui est un diviseur de Cartier effectif de $U$, d\'efinissant de mani\`ere \'evidente  un diviseur de cartier ${\cal R}$ sur ${\cal H}$.  Dans la situation de (6.39), si $C\to D$ est un rev\^etement, avec $r$ points doubles sur $D$, on d\'ecrit la struture locale de ${\cal R}$ comme suit. Soit $k[[t_1,\cdots,t_r,\cdots,t_N]] \hookrightarrow k[[\tau_1,\cdots,\tau_r,\cdots,\tau_N]]$ le morphisme induit entre les bases des d\'eformations universelles respectives, avec $t_i = \tau^{e_i}$ si $1\leq i\leq r$, et $t_j = \tau_j$ si $j>r$. Il est clair que l'\'equation locale de ${\cal R}$ est $\prod_{i,\, e_i>1} \, \tau_i^{e_i-1}$   .\hfill\break
\qed

\bigskip

\section7 { Graphes et rev\^etements} 

\bigskip

Dans cette section  le champ de Hurwitz $\overline {\Cal
H}_{g,G,\xi}$, ou l'espace  modulaire de Hurwitz $\overline
H_{g,G,\xi}$, est  \'etudi\'e d'un point de vue combinatoire.   
  Rappelons  (section 6) que
les points du bord  $\partial {\overline {\Cal H}}_{g,G,\xi} =
\overline {\Cal H}_{g,G,\xi} -  {\Cal H}_{g,G,\xi}$  (resp.
$\partial \overline H_{g,G,\xi}$)  sont repr\'esent\'es  par les
courbes stables munies d'une action stable de $G$ (resp. par les
classes \`a isomorphisme $G$-\'equivariant pr\`es). 
\bigskip

\beginsection 7.1. Graphes modulaires de Hurwitz

\bigskip

 La combinatoire d'une courbe pr\'estable se d\'ecrit   commod\'ement \`a l'aide du graphe dual,
appel\'e aussi graphe modulaire \cite {53}.  Fixons
les notations utilis\'ees dans cette section. Dans  la suite
 un graphe est  
un graphe fini.  
\proclaim D\'efinition 7.1. ( \cite{53}) i) Un graphe (fini, combinatoire)
$\Gamma$  est la donn\'ee de:   \hfill\break
i)  un ensemble  ${\Cal F}$
(de demi-ar\^etes)  muni d'une involution $\tau : {\Cal F} \to {\Cal
F}$\hfill\break
ii)  un ensemble    ${\Cal V}$ de sommets et d'une application
bord (ou incidence) : $\partial : {\Cal F} \to {\Cal V}$. Si $e \in
{\Cal F},\,\, \partial (e)$ est le sommet (initial) de $e$.

  Les orbites \`a deux \'el\'ements de $\tau$ constituent
les ar\^etes g\'eom\'etriques du graphe, tandis que les points fixes
de $\tau$ correspondent  aux feuilles (ou  pattes) du graphe. En modifiant l\'eg\`erement la terminologie, on a ainsi une partition ${\Cal F} = {\Cal E} \sqcup {\Cal L}$, \ ${\Cal E}$ (resp. ${\Cal L}$) \'etant l'ensemble des ar\^etes orient\'ees (resp. des pattes). La
valence d'un sommet $c$ est le nombre $v(c)$ des demi-ar\^etes
incidentes \`a ce sommet.  Si on ignore l'involution $\tau$, on obtient l'ensemble des drapeaux ${\Cal F}l (\Gamma)$, c'est \`a dire des couples $(v,e)\in {\Cal V}\times {\Cal F}, \, \partial (e) = v$.

   Un graphe {\it modulaire} est un graphe $\Gamma$, muni  en plus
de la donn\'ee d'une application  $g: {\Cal V} \to {\Bbb Z}_{\geq
0}$. Le graphe modulaire est dit {\it stable} si  pour tout sommet $c$ tel
que $g(c) = 0$ (resp. $g(c) = 1$), on $v(c) \geq  3$ (resp. $v(c)
\geq  1$).  Le graphe modulaire $\Gamma$ est
{\it marqu\'e} si les pattes sont index\'ees.  

    Un isomorphisme $\Gamma \buildrel\sim\over \to \Gamma'$ est la donn\'ee d'un couple de
bijections $f: {\Cal F }\to {\Cal F'}, \,\,h: {\Cal V} \to {\Cal V'}$, tel
que $h\partial = \partial f$, et $ f\tau = \tau' f$, pour les
graphes, avec en plus la condition $g' h = g$ pour les graphes
modulaires \footnote {Si pour toute demi-ar\^ete, f(e) $\ne \tau
(e)$, l'automorphisme $f$ est dit sans inversion. Nous ne consid\'ererons dans la suite que  des automorphismes sans inversion. Noter qu'il est possible de d\'efinir plus g\'en\'eralement un morphisme de
graphes $\Gamma \to \Delta$ \cite {6}. C'est encore la donn\'ee d'un couple
d'applications $f: {\Cal F}_\Gamma \to {\Cal F}_\Delta, \,\, h: { \Cal
V}_\Gamma \to {\Cal V}_\Delta$ avec les relations de commutation
analogues.}.  Ainsi, si  $\Gamma$ est un graphe (modulaire), on peut
d\'efinir le groupe $\Aut (\Gamma)$, et parler d'une action du groupe
fini $G$ sur le graphe  (resp. graphe modulaire) $\Gamma$.  Sous ces
conditions, il n'est pas difficile de d\'efinir le graphe (resp.
graphe modulaire) quotient $\Gamma /G$. Il est d\'efini par les
ensembles ${\Cal F} /G, \,{\Cal V} /G$, et les applications d\'eduites de
$\tau, \partial$, par passage au quotient. Il y a un morphisme de
graphes   $\Gamma \rightarrow \Gamma /G$, naturellement induit par le passage
au quotient dans les ensembles $\Cal F$ et $\Cal V$.
  Dans la suite  on notera par la m\^eme lettre, dans la mesure o\`u
cela n'induit pas de confusion, un graphe modulaire et le graphe
(ordinaire) sous-jacent. La d\'efinition  (7.1) conduit \`a la
d\'efinition du graphe modulaire $\Gamma (C)$  associ\'e \`a une
courbe pr\'estable marqu\'ee $(C , P)$  \cite {20}, \cite {53}.

Rappelons    cette d\'efinition.   L'ensemble $P$ est
dans cette construction l'ensemble des piq\^{u}res, ou points marqu\'es de $C$. Soit $\{C_i\}_{i\in I}$  l'ensemble  des composantes irr\'eductibles de $C$, ou ce qui est la m\^eme chose l'ensemble des
composantes $\{\tilde C_i\}$ de la normalisation $\tilde C$ de $C$.
Soit aussi  $\{b_\alpha\}_{\alpha \in \Lambda}$ l'ensemble  form\'e
des branches (points de $\tilde C$)  en les diff\'erents points
doubles de $C$, augment\'e des points marqu\'es.  De mani\`ere \'equivalente l'ensemble des points de $\tilde C$ qui se projettent sur
un point double, ou bien un point marqu\'e. Le graphe (dit dual)
$\Gamma (C)$ est  alors d\'efini par \hfill\break \indent i) ${\Cal F} =
\Lambda, \,\, {\Cal V} = I$,\hfill\break
\indent  ii)  $\tau (b_\alpha) = b_\beta$ ssi  $b_\alpha \ne b_\beta$ sont
les deux branches d'un m\^eme
point double, ou alors $\tau (b_\alpha) = b_\alpha$, et dans ce cas
$b_\alpha$ est un point marqu\'e,\hfill\break 
 \indent iii) $\partial b_\alpha = C_i$ ssi $b_\alpha$ est une branche
d'un point double, ou bien
un point marqu\'e,  situ\'e sur la composante $C_i$,\hfill\break
\indent  iv)  $ g(C_i)$ est le genre g\'eom\'etrique de $C_i$ (le genre de
$\tilde C_i$). 

De mani\`ere alternative on peut   voir les
demi-ar\^etes comme   des ar\^etes orient\'ees. Ce qui revient \`a
dire  que  les ar\^etes g\'eom\'etriques sont  d\'edoubl\'ees, en une
ar\^ete orient\'ee et  l'ar\^ete oppos\'ee (graphes \'epais\footnote
{ Dans un graphe \'epais, on impose un ordre cyclique sur les ar\^etes
incidentes en un sommet.}). Dans cette description  les
piq\^{u}res (les pattes du graphe)  sont alors  vues  comme des ar\^etes
orient\'ees d\'epourvues d'ar\^etes oppose\'es. L'involution
$\tau$, que nous noterons aussi $a \to \overline a$, renverse  de
fait l'orientation des ar\^etes, et fixe les pattes. Notons pour une composante $C_i$,
$l_i$  le nombre des  branches d'origine  $C_i$,  et $h_i$ le nombre
des piq\^{u}res d'origine $C_i$.  De sorte que la condition de
stabilit\'e s'\'ecrit de mani\`ere \'equivalente, si $g_i = g(C_i)$
est le genre :
$$2g_i - 2 + h_i + l_i > 0 \,\, ( i=1,\dots,s = \Card I) \tag (7.1) $$
Rappelons  aussi la relation bien connue qui donne le genre arithm\'etique de la courbe $C$, suppos\'ee connexe, en fonction des param\`etres contenus dans le graphe dual \cite {52},  c'est \`a dire le
nombre $d$  de points doubles (ar\^etes g\'eom\'etriques), et le
nombre de sommets $s$:
$$ g =  \sum_i g_i + \dim\, H_1( \Gamma) = \sum_i g_i + d - s + 1
\tag (7.2) $$

Revenons maintenant \`a la situation d'origine, donc  
$C$ est une $G$-courbe stable. Il est clair que $G$ agit   de
mani\`ere naturelle sur le graphe modulaire $\Gamma = \Gamma (C)$.  Du
fait de la stabilit\'e de  l'action, $G$ agit sans inversion, i.e.
pour toute branche $a$, et tout $\sigma \in G, \sigma \ne 1,
\,\,\sigma (a) \ne \overline a$. Il n'est pas difficile d'impl\'ementer sur $\Gamma$ la donn\'ee de Hurwitz.  Plus pr\'ecisemment  on attache \`a
 chaque \'el\'ement de $\Cal F$ une d\'ecoration mat\'erialis\'ee par  un couple  $(H,\chi)$, l'holonomie, $H$ \'etant le stabilisateur (cyclique) soit d'un
point non-singulier, soit d'une branche, et $\chi$  \'etant le
caract\`ere de l'action de $H$ sur l'espace tangent du point distingu\'e .  De plus 
si  l'ar\^ete $e$ est d\'ecor\'ee par $(H , \chi)$,  alors  l'ar\^ete oppos\'ee $\overline e \ne e$ est d\'ecor\'ee par $(H , \chi^{-1})$.
 Si on ne regarde que les classes de conjugaison des d\'ecorations qui de fait sont  constantes sur les
$G$-orbites on obtient  en particulier la  donn\'ee de Hurwitz.  
 En r\'esum\'e,  appelons graphe modulaire de Hurwitz, un graphe modulaire muni d'une action sans inversion d'un groupe fini $G$ qui satisfait aux conditions pr\'ec\'edentes.   Pour une ar\^ete, ou patte,  d\'ecor\'ee par  le  couple $(H,\chi)$, on parlera de ce couple comme \'etant l'holonomie. L'holonomie en $ge$ est $(gHg^{-1} , ^g\chi)$ (\S \ 2.2.1).
 Ce qu'on doit comprendre comme isomorphisme de graphes modulaires de Hurwitz est clair.

\proclaim D\'efinition 7.2.   On appelle graphe
modulaire de Hurwitz associ\'e au rev\^etement
galoisien stable $C\to D$, le graphe modulaire $\Gamma = \Gamma (C)$
muni   de l'action de $G$, et d\'ecor\'e  par les caract\`eres locaux attach\'es aux  orbites de points  de ramification (points sp\'eciaux)  et branches de $C$. Le graphe est marqu\'e si les orbites  de pattes sont index\'ees.\hfill\break
 \indent Soit $\Gamma$ un graphe modulaire de Hurwitz attach\'e \`a
$(g,G,\xi)$. Un $G$-rev\^etement $\pi: C \to D$ est dit de type
$\Gamma$, si $\Gamma (C) \cong \Gamma.$

On utilisera une notation identique pour le graphe modulaire de Hurwitz
du rev\^etement $\pi: C\to D$ et le graphe modulaire de la courbe
$C$.  Les graphes seront comme les rev\^etements, marqu\'es.

\bigskip

\beginsection 7.2. Graphes de groupes et rev\^etements de graphes

\bigskip

  Notre objectif   est de d\'ecrire
$\Gamma$, partant du graphe quotient,  qui est  est un graphe
ordinaire $\Delta = {\Gamma} /G$, et d'un morphisme de monodromie,
d\'efini   dans le but d'assurer la validit\'e
 du dictionnaire usuel.  La th\'eorie de Bass \cite {6} des rev\^etements
de graphes, est particuli\`erement  utile pour mener \`a
bien ce programme.

\bigskip

  {\it 7.2.1. Graphes modulaires quotients} 
\bigskip

Soit $C$ une courbe stable marqu\'ee d\'efinie sur $k$, munie d'une action stable de $G$, donc   $C$ est marqu\'ee par un diviseur $G$-invariant contenant le diviseur
des points de ramification (\S 4.3). Soit $\Gamma = \Gamma (C)$ le
graphe modulaire de Hurwitz attach\'e \`a $C$ (D\'efinition 7.2).
Notons pour d\'ebuter  la remarque \'el\'ementaire suivante:
 \proclaim Proposition 7.3.  Le graphe quotient $\Delta = \Gamma /G$
poss\`ede une structure canonique de graphe modulaire de genre $g'$,
\`a $b$ piq\^{u}res. Il s'identifie au graphe modulaire de la
courbe stable marqu\'ee $D = C/G$.

\dem Par d\'efinition du graphe combinatoire quotient \cite{6}   
$${\Cal V}(\Delta) =  {{\Cal V}(\Gamma)}/G  \quad \hbox{ \rm  et }
\quad  {\Cal F}(\Delta)  =   {\Cal F}(\Gamma)/ G \tag (7.3)$$   L'action de $G$
commute avec l'involution $\tau_\Gamma$, celle-ci induit donc une
involution sur
l'ensemble quotient, d\'efinissant  ainsi  la structure du graphe
quotient $\Delta$. Noter aussi
que les piq\^{u}res (points fixes de $\tau_\Delta $) du quotient,
sont les images des piq\^{u}res de $\Gamma$.   Pour finir de d\'ecrire $\Delta$ comme graphe modulaire, il  faut assigner une valeur
$g'_\alpha$  pond\'erant le sommet $v'_\alpha$ associ\'e \`a la
composante $D_\alpha$ de $D$. Par convention,  $g'_\alpha$ est le
genre (g\'eom\'etrique) de $D_\alpha$.  Il est clair que $g'_\alpha$
s'exprime en fonction des seules donn\'ees  de $\Gamma$. Soit une
composante $C_i$ de $C$, d'image $D_\alpha$, et soit $G_i$ le
stabilisateur de $C_i$, i.e. du sommet $v_i$ de $\Gamma$. Donc
$D_\alpha = C_i /{G_i}$.  L'expression cherch\'ee est la formule de
Riemann-Hurwitz, appliqu\'ee au rev\^etement ${\tilde C}_i \to
{\tilde D}_\alpha$: 
$$2g_i - 2 = \vert G_i\vert (2g'_\alpha - 2) +
\sum_\lambda B_\lambda  $$ 
o\`u  les $B_\lambda$ sont les
contribution des points de branchement, contributions d\'ependantes de la seule combinatoire du graphe de modulaire
de Hurwitz   $\Gamma$. Il est facile de prouver que l'on d\'efinit
ainsi une structure de graphe modulaire
(quotient) sur $\Delta$. Le dernier point est  clair. \hfill\break
\qed

Pour \'eviter toute confusion avec d'autres d\'efinitions  de rev\^etements
de graphes\footnote {On peut exiger, ce qui n'est pas le cas
dans notre d\'efinition, la propri\'et\'e de rel\`evement des
chemins (rev\^etements non ramifi\'es).}, insistons sur le fait que
l'action de $G$ sur les sommets et ar\^etes n'est  dans cette d\'efinition pas n\'ecessairement libre,  le quotient $\Gamma \to
\Gamma /G$ est donc \`a prendre dans le sens d'un rev\^etement ramifi\'e.  Avant de discuter en d\'etail le cas g\'en\'eral, il peut \^etre
utile de traiter deux cas particuliers importants; ce sont les
graphes modulaires de Hurwitz qui sont associ\'es aux composantes de
codimension un du bord de  $\overline {\Cal H}_{g,G,\xi}$, ou de
l'espace modulaire grossier de m\^eme nom. Ce sont les deux cas pour
lesquels $\Delta$ n'a qu'une seule ar\^ete g\'eom\'etrique,   
$\Delta$ est  donc soit un segment, soit une boucle.  
  \smallskip
$\bullet$\quad $\Delta$ {\sl est un arbre} 
\smallskip
Du fait de cette hypoth\`ese, il est  bien connu que la projection
$\Gamma \to \Delta$ admet une section \cite {6}; nous noterons  pour
simplifier par la m\^eme lettre $\Delta$ l'image d'une section
suppos\'ee choisie. Notant toujours ${\Cal V}$ l'ensemble des
sommets, et ${\Cal F}$ l'ensemble des ar\^etes orient\'ees, nous
noterons pour simplifier  $a\in \Gamma$ au lieu de $a \in {\Cal F}$,
m\^eme chose pour les sommets; la notation est analogue pour
$\Delta$. On a alors les identifications
$${\Cal F} =  \bigsqcup_{a\in \Delta} G /{H_a}\quad\hbox{\rm et}\quad {\Cal V} =
\bigsqcup_{v\in \Delta}  G /{G_v} \tag (7.4)$$
On note $G_v$ le stabilisateur du sommet $v$, et $H_a$ le
stabilisateur de l'ar\^ete (orient\'ee,
ou g\'eom\'etrique) $a$. D\'ecrivons la relation d'incidence. Si $e$
est une ar\^ete d'extr\'emit\'es
$g_1v_1$ et $g_2v_2$, il existe alors $a\in \Delta$ tel que $e = ga
\in Ga$; alors $g^{-1}g_1
v_1 = v_1$ et $g^{-1}g_2 v_2 = v_2$. Cela montre que l'ensemble des
ar\^etes  d'extr\'emit\'es
$g_1v_1$ et $g_2v_2$ est en bijection naturelle avec les classes
$gH_a$, $a$ \'etant l'unique
ar\^ete de $\Delta$ joignant $v_1$ et $v_2$, qui v\'erifient l'inclusion
$gH_a \subset g_1G_{v_1} \cap g_2G_{v_2} .$

 \exa{  7.1} (Le peigne) 

D\'etaillons  le cas particulier o\`u $\Gamma$ est un  "peigne", ce
qui correspond \`a choisir
$\Delta$ sous la forme d'une \'etoile \`a $k\geq 3$ branches. Les
notations pour $\Delta$ sont: \hfill\break $\bullet$ $v_0$  est la
racine de l'\'etoile, qui est  un sommet  de valence
$k$\hfill\break$\bullet$ $v_1,\dots,v_k$ sont les sommets terminaux
des $k$-branches.  

On doit installer sur $\Delta$ une
structure de graphe modulaire, et   impl\'ementer la donn\'ee de
Hurwitz. Pour la structure de graphe modulaire, on pose $g_0=
g_1=\dots=g_k = 0$. La donn\'ee de Hurwitz dans cet exemple est du
type "Harbater-Mumford" \cite {16}, \cite {30}, c'est \`a dire de la
forme
$$\xi = \sum_{i=1}^k ([H_i,\chi_i] + [H_i,\chi_i^{-1}]) \tag (7.5)$$ 

La forme sp\'ecifique de la donn\'ee de Hurwitz impose \`a chaque
sommet $v_1,\dots,v_k$ d'\^etre la source de deux piq\^ures. En effet si $C_i$ est une composante de $C$ au-dessus du brin d'indice $i$ de $D$, alors $C_i\to \Bbb P^1$ est ramifi\'e en trois points au plus. La forme de la donn\'ee de ramification impose que le nombre de ponts de branchement est deux, le point double \'etant d'isotropie triviale.

\vskip 10pt
\epsfxsize=2truein
\centerline{\epsfbox{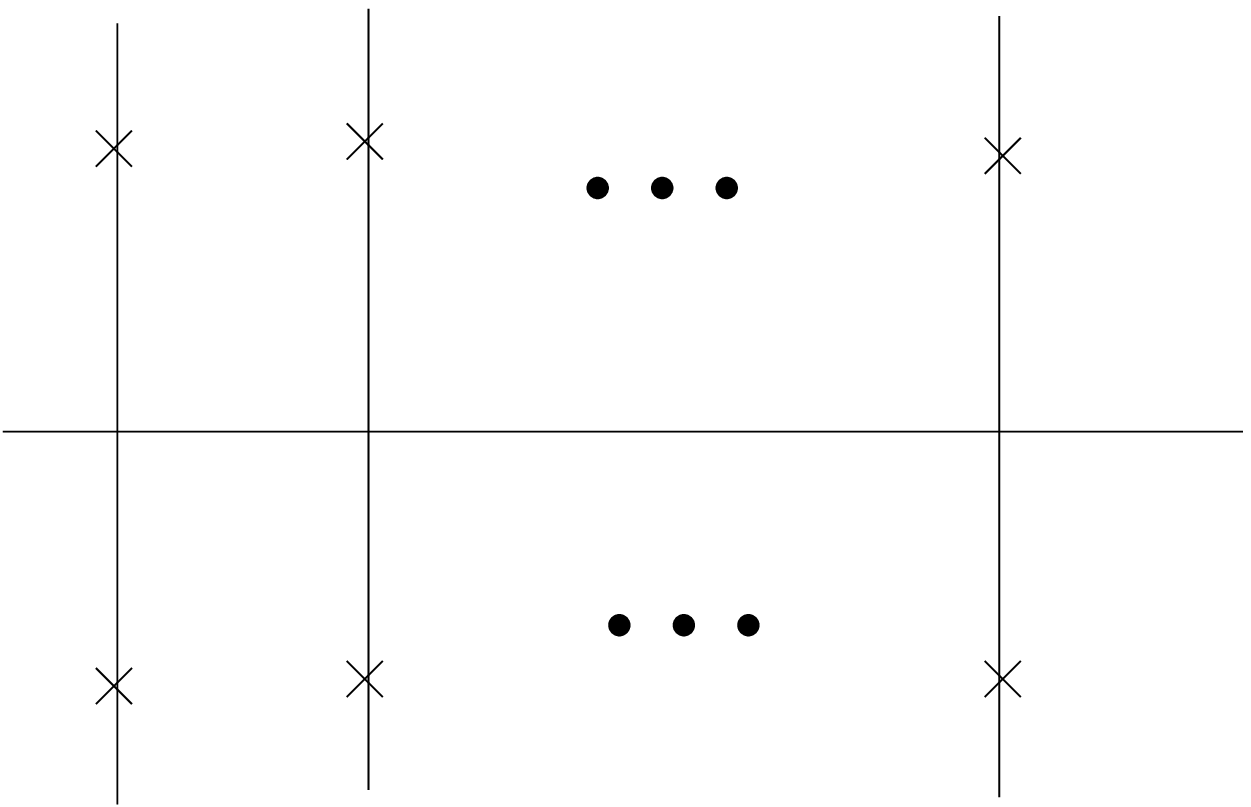}}
\vskip 10pt


Pour d\'ecrire la structure de graphe de Hurwitz sur $\Gamma$,
choisissons un sommet $s_i$ au
dessus de $v_i$ ($i=0,\dots,k$), et posons $H_i = G_{s_i}$; notons
que n\'ecessairement $H_0 =
1$. Observons aussi que $G$ est en cons\'equence engendr\'e par la
r\'eunion des $H_i$. On obtient finalement la description suivante  de 
   $\Gamma$:
$${\Cal V}_\Gamma = \big (\bigsqcup_{i=1}^k  G /{H_i}\big )
\,\,\bigsqcup G,\quad {\Cal
F}_\Gamma = \bigsqcup_1^{2k}G \tag (7.6)$$
Notons que le genre $g$ de $C$  est donn\'e  par
$$g = 1+ k\Card G - \sum_{i=1}^k [G : H_i] \tag (7.7) $$

\exa{ 7.2} (Le segment)  Dans cet exemple, $\Delta$ est    un segment de sommets
$v_1$ et $v_2$, et d'ar\^etes
orient\'ees $e$ et $\overline e$, $e$ pointant vers $v_2$.

\vskip 10pt
\epsfxsize=2truein
\centerline{\epsfbox{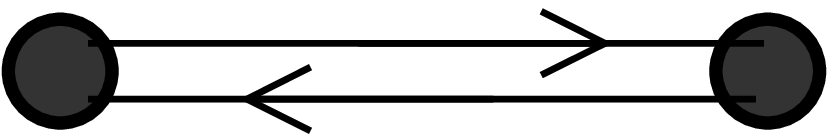}}
\vskip 10pt


On rel\`eve  $\Delta$ en un segment de $\Gamma$ de sommets $s_1$ et
$s_2$, de  stabilisateurs
respectifs $G_1$ et $G_2$. L'ar\^ete orient\'ee $a$ rel\`eve $e$:
soit $H$ le stabilisateur de $a$,
alors $H \subset G_1 \cap G_2$, et la connexit\'e impose que $G_1
\cup G_2$ engendre $G$. La
description du graphe (ordinaire)  $\Gamma$ est ais\'ee: on a
$$ {\Cal F}_\Gamma =  G/H \times \{\pm 1\}, \quad {\Cal V}_\Gamma =  G/{G_1}
\bigsqcup  G /{G_2} \tag (7.8) $$
Une ar\^ete g\'eom\'etrique est dans ce mod\`ele repr\'esent\'ee par
un \'el\'ement de $ G/H$; l'ar\^ete $gH$ joint les sommets $g_1G_1$
et $g_2G_2$  ssi $gH \subset g_1G_1 \cap g_2G_2$. L'involution est
donn\'ee par $\tau (gH,\epsilon) = (gH, -\epsilon)$. D'une  mani\`ere plus  g\'eom\'etrique, supposons que les composantes $C_1$ et
$C_2$ correspondent \`a $s_1$ et $s_2$.  Posons
$$\Ind_{G_1}^G (C_1) = C_1 \times_{G_1} G \tag (7.9) $$
Cela signifie qu'on effectue le quotient de $C_1\times G$ par
l'action de $G_1: \alpha (x,g) =
(\alpha x, g\alpha^{-1})$.   Le r\'esultat est une courbe avec $[G :
G_1]$ composantes irr\'eductibles
disjointes, toutes identiques \`a $C_1$, et munie de l'action de
$G$ donn\'ee par $ g[x,h] = [x,gh]$. La courbe $C$ est alors de la
forme
$$C = \Ind_{G_1}^G (C_1) \bigvee \Ind_{G_2}^G (C_2) \tag (7.10) $$
expression dans laquelle $\vee$ signifie que l'on  recolle les deux
facteurs en identifiant des
paires de points, selon la combinatoire des ar\^etes de $\Gamma$. Les
points doubles forment une
unique $G$-orbite; soit $p \in C_1 \cap C_2$ l'un d'eux, d'isotropie
$H$. Si $H$ agit sur la
branche $C_1$ par le caract\`ere $\chi$, et sur la branche $C_2$ par
$\chi^{-1}$ (stabilit\'e
oblige), alors la donn\'ee de Hurwitz $\xi$ de l'action de $G$ sur $C$ s'obtient par $\xi =
\xi_1 + \xi_2$, o\`u on suppose
que la donn\'ee d\'efinie par l'action de $G_1$ sur $C_1$ est
$\tilde\xi_1 = \xi_1 + [H,\chi]$, et
celle d\'efinie par l'action de $G_2$ sur $C_2$ est $\tilde\xi_2 =
\xi_2 + [H,\chi^{-1}]$. La
combinatoire est donc totalement contenue dans la donn\'ee
$(G_1,G_2,H ,\chi, \xi_1,\xi_2)$.  Noter qu'il n'est pas n\'ecessaire de s'assurer que le genre $g$ de $C$  est celui fix\'e initialement.  En effet le genre r\'esulte de la donn\'ee du  genre $g' = g'_1+g'_2$ de la  base $D = C_1/G_1\vee C_2/G_2$  joint \`a la donn\'ee de Hurwitz. On a donc  si $g_i$ est le genre de $C_i$ l'\'egalit\'e
$$g = [G:G_1]g_1+ [G:G_2]g_2 + [G:H] - [G:G_1]-[G:G_2]+1$$

\exa{ 7.3} (La boucle)  

  Si $\Delta$ est une boucle, alors  $G$ est une "extension " HNN. On
distingue deux cas selon que
$\Gamma$ a, ou n'a pas de boucle. Traitons d'abord le second cas: 
  
  $\bullet$   On suppose que  $\Gamma$ ne contient pas de boucle,
et que $\Delta =
\Gamma / G$ se r\'eduit \`a une boucle. Alors $\Gamma$ est d\'efini
par la donn\'ee (\`a conjugaison
pr\`es) par un triplet $(G_0,H,g_0)$ form\'e de deux sous-groupes $H$ et $G_0$ de $G$, $H$ \'etant
cyclique, et $H \subset G_0 \subset
G$, et d'un \'el\'ement $g_0 \in G, g_0 \not\in G_0$. On impose \`a
cette donn\'ee les conditions $   g_0 H g_0^{-1}\subset G_0$ et  $G_0 \cup \{g_0\}$ engendre $G$.
Il est facile de voir
qu'on reconstruit le graphe $\Gamma$   en posant:
$${\Cal F}_\Gamma = (G /H)\times \{\pm 1\}, \quad {\Cal V}_\Gamma =  G
/{G_0} \tag (7.11)$$
L'involution $\tau$ est comme dans l'exemple 7.2 la multiplication
par $-1$ sur le facteur de
droite. L'incidence est donn\'ee par
$$(\alpha H , +1) \to (\alpha G_0 , \alpha g_0 G_0)\quad,\quad
(\alpha H , -1) \to (\alpha g_0
G_0 , \alpha G_0) \tag (7.12)$$
en particulier $\dim H_1(\Gamma) = [G : H] - [G : G_0] + 1$. Par
exemple, si $G_0 = 1$, alors
$G$ est cyclique engendr\'e par $g_0$, et $\Gamma$ est un
$n$-circuit, avec $n = \Card  G$.
Par construction il y a  dans $\Gamma$ une seule orbite de sommets,
soit $v$ l'un d'entre eux
et soit  $G_0$  le stabilisateur de  
$v$. Si $C_0$ est la composante repr\'esent\'ee par $v$, la donn\'ee
de Hurwitz d\'efinie par l'action
de $G_0$ sur $C_0$ est de la forme
$$\tilde \xi = \xi + [H , \chi] + [H , \chi^{-1}]\in R_+(G_0) \tag (7.13)$$
L'hypoth\`ese de stabilit\'e de l'action qui impose aux caract\`eres locaux en $v$ et $g_0v$ d'\^etre
oppos\'es, ajoute une contrainte sur $\chi$, \`a savoir $\chi
(g_0^{-1} s g_0) = \chi^{-1} (s) ,
\,\,(s\in H)$.

$\bullet$ Reste le cas o\`u $\Gamma$ contient une boucle; il est
facile de voir que dans ce cas $\Gamma$ n'a
qu'un seul sommet, et que  toutes les ar\^etes sont des boucles, i.e.
$\Gamma$ est une fleur \`a $m\geq 1 $ p\'etales.

\bigskip
   {\it 7.2.2. Rev\^etements ramifi\'es de graphes
(Bass \cite {6})} \bigskip

 Pour reconstruire le graphe $\Gamma$
partant du quotient $\Delta$, et d'une donn\'ee additionnelle \`a
d\'efinir,  on compense le fait que l'action de $G$
n'est pas suppos\'ee libre, par une structure de graphe de groupes sur  le graphe quotient $\Delta$   (th\'eorie de Bass-Serre: \cite {6}, et
les r\'ef\'erences contenues dedans).   

Rappelons   qu'un graphe de
groupes $\Delta$ (tous les graphes sont finis), est la donn\'ee d'un
graphe $\Delta$, tel qu'\`a chaque sommet $v$ (resp. \`a chaque
ar\^ete orient\'ee $e$) est associ\'e un groupe (fini) $G_v$,
(resp. un groupe $H_e$); si $e$ pointe vers $v$, on a un morphisme
injectif $\partial _e^1: H_e \to G_v$, cette donn\'ee \'etant telle
que $H_e = H_{\overline e}$. On  pose $\partial _e^0 = \partial
_{\overline e}^1.$ On d\'efinit
alors le groupe fondamental d'un graphe de groupes $\Delta$ de la
mani\`ere suivante (loc.cit)

\proclaim D\'efinition 7.4.  Soit $F_\Delta$ le groupe libre de base
l'ensemble des ar\^etes orient\'ees ${\Cal F}_{\Delta}$; choisissons
aussi un arbre maximal $T \subset \Delta$. Le groupe fondamental
$\pi_1(\Delta;T)$ est  d\'efini comme \'etant  le quotient du produit
libre $F_{\Delta} \star
(\star_v G_v)$ par le sous-groupe distingu\'e engendr\'e par les relations
$$e = 1 \,\,\hbox{\rm  si} \,\,e\in T\,\,, \,\, e\overline e = 1\,\,(e\in
\Gamma) , \,\, \hbox{\rm  et}\,\,
e\partial_e^1 (h) \overline e \, = \, \partial_e^0 (h)\quad (e\in
\Delta, h\in H_e) \tag (7.14)$$

 Le groupe $\pi_1(\Delta;T)$  ne d\'epend pas, \`a isomorphisme pr\`es, du choix de $T$. On a la
description alternative suivante du groupe fondamental \cite {6};
appelons chemin de $\Delta$ une suite
$$ \gamma = (g_0,e_1,g_1,\dots, e_n,g_n) \tag (7.15)$$
o\`u  $(e_1,\dots,e_n)$ est une chemin "ordinaire",  $a_i$ d\'esignant le sommet terminal de
$e_i$ (donc aussi le sommet initial de $e_{i+1}$ si $i<n$), $a_0$ \'etant le sommet initial de
$e_1$, on a $g_i \in G_{a_i},\,\, (0\leq i\leq n)$. Si $a_0 = a_n,
\gamma$ est un lacet en $v=a_0$.
Ceci \'etant, on forme   le groupe des chemins
$$\pi (\Delta) = \left[ \left( \star_v G_v\right) \star
F_{\Delta}\right] \,/ \, R \tag (7.16)$$
  $R$ d\'esignant  le sous-groupe distingu\'e
engendr\'e  par les  relations
$$ e \overline e = 1 \,\, (e \in \Delta)\,\, , \,\,e\partial_e^1 (h)
\overline e \, = \,\partial_e^0 (h)\quad (e\in \Delta, h\in H_e) \tag
(7.17) $$
Si $\gamma$ est le chemin (7.15), on pose
$$\vert \gamma \vert = g_0e_1g_1 \dots e_ng_n \in \pi (\Delta) \tag (7.18) $$
Les \'el\'ements (7.18), lorsque $\gamma$ est un lacet en $v = a_0$,
forment un sous-groupe not\'e $\pi_1 (\Delta,v)$. On montre alors
\cite {6}  que la surjection $\pi (\Delta) \to \pi_1 (\Delta,T)$,
restreinte au sous-groupe $\pi_1(\Delta,v)$ est un isomorphisme
$\pi_1(\Delta,v) \cong \pi_1 (\Delta,T)$.

Formons $\Delta =  \Gamma /G$ le
graphe quotient, et soit $p:\Gamma \to \Delta$ le
morphisme quotient. On sait munir $\Delta$ d'une structure de graphe
de groupes, structure qui
capture les informations sur l'action de $G$, i.e. d\'efinissant sur
$\Delta$ une structure
"d'orbifold"  \cite {6}. Disons seulement, et cela sera suffisant
pour la suite, que si $s$ est un
sommet de $\Delta$, alors le groupe $G_s$ est essentiellement le
stabilisateur d'un sommet $v \in
\Gamma$ tel que $p(v) = s$; d\'efinition analogue pour $H_e$ si $e$
est une ar\^ete de $\Delta$. La construction n\'ecessite le choix de sections de $p$ au niveau des sommets et ar\^etes. On notera dans la suite $\Gamma // G$   le graphe $\Delta$ muni de la
structure de graphe de
groupes quotient. On a besoin du r\'esultat suivant de Bass-Serre (Bass \cite
{6}, Theorem  3.6):

\proclaim Th\'eor\`eme 7.5. Il existe une surjection naturelle $\psi
: \pi_1( \Gamma // G) \to G$,
appel\'ee morphisme de monodromie, telle que $\ker \psi \cong \pi_1
(\Gamma)$ (groupe fondamental du graphe ordinaire quotient). 
 \hfill\break\qed

La d\'efinition de $\psi$ mime la d\'efinition usuelle de l'action de
monodromie.  Noter que par simplicit\'e on ne fait pas r\'ef\'erence aux sommets de base. Le th\'eor\`eme 7.5 dit   qu'on
peut reconstruire le graphe $\Gamma$ muni de l'action de $G$, partant
 de la structure de graphe de groupes sur $\Delta = \Gamma
// G$, et du morphisme de monodromie, comme dans le cas classique. Le graphe $\Gamma$ appara\^{\i}t comme le quotient  du rev\^etement universel de $\Gamma//G$ par $\pi_1(\Gamma)$.

\proclaim D\'efinition 7.6.  Sous les hypoth\`eses pr\'ec\'edentes,
notons $D$ le sous-groupe distingu\'e de $G$ engendr\'e par les
stabilisateurs des sommets (sous-groupe de d\'ecomposition) , et notons $I$ le sous-groupe distingu\'e engendr\'e  par
les stabilisateurs des ar\^etes (sous-groupe d'inertie).

On a donc $I \triangleleft D \triangleleft G$. On peut pr\'eciser la
structure du groupe
fondamental ordinaire de $\Delta = \Gamma /G$, comme suit:
\proclaim Proposition 7.7.  Il y a une surjection naturelle $\phi:
\pi_1(\Gamma /G) \to
  G /D$, qui conduit \`a une suite exacte
$$\pi_1(\Gamma) \rightarrow \pi_1(\Gamma /G) \rightarrow
G/D\rightarrow 1 \tag (7.19) $$

\dem  D\'efinissons d'abord le morphisme $\phi$. Soit $a_0 \in
\Gamma$ un sommet, et posons $\overline a_0 =p(a_0)$. Si $\overline
\gamma = (\overline e_1,\dots,\overline e_n)$ repr\'esente un lacet
de $\Delta$ bas\'e en $\overline a_0$, on peut le relever en un
lacet $\gamma = (e_1,\dots,e_n)$ d'origine $a_0$. Soit $s$ le sommet
terminal de l'ar\^ete $e_n$; alors   $p(s) = p(a_0)$. Il y a donc $g
\in G$ tel que  $ s = ga_0$. On doit v\'erifier que la classe de $g$
modulo $D$ ne d\'epend pas du rel\`evement choisi. Soit en effet
$(e'_1,\dots,e'_n)$ un autre choix de rel\`evement; il y a $g_1 \in
G, e'_1 = g_1 e_1$. L'action de $G$ \'etant sans inversion, $g_1a_0 =
a_0$, donc $g_1 \in D$.
le rel\`evement $(e_1,g_1e'_2,\dots,g_1e'_n)$ diff\`ere de
$(e_1,\dots,e_n)$ par au plus $(n-1)$
ar\^etes. Par une r\'ecurrence imm\'ediate, on peut supposer que si
$s'$ est le sommet terminal de
$e'_n$, alors pour le sommet terminal de $g_1e'_n$, on a  $g_1s' \in
g D a_0$. On en tire imm\'ediatement que $s' \in gDa_0$, et donc $g$
et $g'$ sont dans la m\^eme classe modulo $D$. Le morphisme $\phi$
assigne \`a $\gamma$, la classe modulo $D$ de l'\'el\'ement $g$
ci-dessus. La v\'erification du fait que le morphisme $\phi$ est bien
d\'efini est \'evidente.

On sait par ailleurs que si $R$ est le sous-groupe distingu\'e de
$G$ engendr\'e par les
sous-groupes $G_s, s\in \Gamma // G$, alors
$$ {\pi_1(\Gamma // G)}/ R\,\cong \, \pi_1( \Gamma /G)$$
La construction montre qu'en fait $\phi$ co\"\i ncide avec le
morphisme d\'eduit de $\psi$ par
passage au quotient par les groupes $R$ et $D$ respectivement.  Cela
a pour traduction un diagramme commutatif:
 
  $$\harrowlength=45pt \varrowlength= 45pt\commdiag
{1  &\mapright &  \pi_1(\Gamma)  & \mapright & \pi_1(\Gamma // G)  &\mapright^\psi& G  &\mapright & 1  \cr
  &&\mapdown &&\mapdown &&\mapdown && \cr
  &&  \pi_1(\Gamma // G) /R   &\mapright^\sim&   \pi_1(\Gamma /G)   &\mapright^\phi&
  G /D  &&\cr }\tag (7.20)  $$  
De ce diagramme on tire un morphisme compos\'e  $\theta:
\pi_1(\Gamma) \to \pi_1(\Gamma
/G)$ qui n'est pas autre chose que le morphisme canonique d\'eduit de
$p$, et aussi l'\'egalit\'e $\psi (R) = D$. La conclusion vient alors
imm\'ediatement d'une chasse au diagramme. \hfill\break\qed

Passant \`a l'homologie, la suite exacte 7.20, donne imm\'ediatement comme corollaire:
 
\proclaim Corollaire 7.8.  Sous les hypoth\`eses pr\'ec\'edentes, on
a une suite exacte canonique:
$$ H_1(\Gamma , {\Bbb Z}) \longrightarrow H_1(\Gamma /G , {\Bbb Z})
\longrightarrow
\big( G/D \big)_{ab} \longrightarrow 1. \tag (7.21) $$
\qed

\bigskip

\beginsection  7.3. Groupe de Picard et rev\^etements

\bigskip

  {\it 7.3.1.  Sous-groupes de d\'ecomposition et d'inertie}

\bigskip

La proposition 7.2 admet une interpr\'etation g\'eom\'etrique simple
lorsque $\Gamma$ est le graphe dual (modulaire) d'une courbe pr\'estable, sur laquelle le groupe $G$ agit stablement. Soit C une telle
courbe (d\'efinie sur le corps alg\'ebriquement clos $k$), de
composantes irr\'eductibles
$C_1\dots,C_s$,  de sorte que $s$  repr\'esente le nombre d'ar\^etes
g\'eom\'etriques de $\Gamma$. Soit la normalisation  
 $\tilde C = \bigsqcup_i \tilde C_i \longrightarrow C$ de $C$. 
 On consid\`ere  $\tilde C$ comme munie de l'action induite de
$G$; m\^eme chose pour  $\Gamma$.  On note $G_i$ le stabilisateur de
$C_i$; si $p$ est un point double, on notera $G_p$ son stabilisateur
(cyclique). Soient comme
dans la section pr\'ec\'edente les sous-groupes distingu\'es $D$ et
$I$ (de d\'ecomposition et
d'inertie). On notera pour simplifier $H^1(\Gamma)$ le groupe
$H^1(\Gamma , {\Bbb Z})$. Rappelons  tout d'abord la description bien
connue du groupe de Picard de $C$; on le d\'ecrit au moyen de la
suite exacte
$$1\longrightarrow H^1(\Gamma) \otimes k^* \longrightarrow \Pic(C)
\longrightarrow \Pic (\tilde C)
= \prod_i \Pic (\tilde C_i) \longrightarrow 1 \tag  (7.21) $$
Consid\'erons  la courbe quotient $D = C/G$, qui est pr\'estable, et
soit $\pi: C \to  D = C/G$ le morphisme quotient.  Une suite exacte
analogue \`a celle utilis\'ee pour d\'ecrire $\Pic(C)$, appliqu\'ee
\`a $\Pic (C/G)$,  conduit \`a un diagramme commutatif
\ 
  $$\harrowlength=38pt \varrowlength=38pt \commdiag{
1  &\mapright&  H^1(\Gamma) \otimes k^*&  \mapright&  \Pic(C) &\mapright  & \prod_i
\Pic(\tilde C_i) & \mapright & 1 \cr
&& \mapup\lft{ \pi^* } &&   \mapup\lft{\pi^*}&&\mapup\lft{  \prod \pi_i^*} &&\cr
  1  &\mapright&
H^1(\Gamma/G) \otimes k^* &\mapright&  \Pic (C/G)&  \mapright&\prod_i \Pic({\tilde
C_i}/{G_i}) & \mapright & 1 \cr} \tag (7.22)$$ 
 dans lequel les fl\`eches verticales sont induites par les
morphismes quotients respectifs, not\'es $\pi$. On peut noter que la
fl\`eche verticale de gauche, c'est \`a dire $\pi^*$, est  celle
  qui d\'erive du morphisme de graphes
$\pi: \Gamma \to  \Gamma /G$. 

  On va pr\'eciser le contenu du  diagramme  (7.22) en \'etudiant la suite
exacte form\'ee par les noyaux des fl\`eches verticales. Si $A$ est
un groupe fini, on notera $\hat A$ le groupe des caract\`eres, soit
$ \hat A = \Hom (A , k^*)$, sous entendu que $\Card A$ est premier \`a
 la caract\'eristique de $k$.  Le r\'esultat suivant est important
pour la suite:

\proclaim Th\'eor\`eme 7.8. Le sous-groupe d'inertie de $(\tilde
C_i,G_i)$ \'etant not\'e  $I_i$, la suite exacte des noyaux du
diagramme pr\'ec\'edent s'identifie canoniquement \`a:
$$1 \longrightarrow  \widehat{G/D} \longrightarrow  \widehat{ G/I}
\longrightarrow \prod_i  \widehat {{G_i}/{I_i}}  \tag (7.23) $$

\dem  La preuve utilise la structure naturelle de $G$-faisceau
support\'ee  par $\pi_*({\Cal O}_C)$ sur $C/G$. Du fait de la pr\'esence de points doubles (\'eventuels) avec un groupe d'isotropie non
trivial, ce
faisceau n'est pas en g\'en\'eral localement libre, il est seulement
de type fini et sans torsion.
Cependant comme la ramification est mod\'er\'ee, on a toujours  une
d\'ecomposition en facteurs isotypiques (voir \S  3.2):
$  \pi_*({\Cal O}_C)  = \bigoplus_{\chi \in Irrep(G)}  {\Cal L}_\chi
\otimes V_\chi, $
  $V_\chi$ d\'esignant
l'espace  supportant  la repr\'esentation irr\'eductible
$\chi$, et ${\Cal L}_\chi$ \'etant un faisceau sans torsion donn\'e
par
$$ {\Cal L}_\chi = \pi_*^G \left( {\Cal O}_C \otimes \check {V_\chi}
\right) \tag (7.24) $$
Analysons  plus en d\'etail cette d\'ecomposition  au voisinage d'un
point double $q = \pi(p)$  de $D
=  C/G$, le point double $p$ \'etant de stabilisateur $H$.
Apr\`es localisation, et compl\'etion,  on est amen\'e \`a d\'ecomposer  $\pi_*(\widehat {\Cal
O}_C)_q$ en tant que $\big ( \widehat {\Cal O}_{D,q} , G\big )$-
module. La d\'ecomposition suivante est imm\'ediate:
$$\pi_*(\widehat {\Cal O}_C)_q \,\, \cong \,\, \Ind_H^G \big(
\widehat {\Cal O}_{C,p} \big)
\tag (7.25)$$
On est de la sorte ramen\'e \`a supposer que $\widehat {\Cal
O}_{C,p} =     {k[[x,y]]}/{(xy)}$, et que l'action sur les deux branches d'un g\'en\'erateur
$\sigma \in H$ de $H$, soit r\'eduite \`a la  forme usuelle,
$\sigma (x) = \zeta x\,,\, \sigma (y) = \zeta^{-1} y$, pour une
certaine racine
de l'unit\'e d'ordre $e = \Card H$. Ainsi, si on pose $u = x^{e}\,,
\,v = y^{e}$, on a
$\widehat {\Cal O}_{D,q} =    {k[[u,v]]}/{(uv)}$. La d\'ecomposition cherch\'ee est alors
totalement explicite
$$\widehat {\Cal O}_{C,p}\,\cong\,\widehat {\Cal O}_{D,q} \bigoplus
\big( \bigoplus_{j=1}^{e-1}
(x^j,y^{e-j}) \widetilde {\widehat {\Cal O}_{D,q}}\big ) \tag (7.26)  $$
o\`u le symbole tilde d\'esigne la normalisation. Notons $\mu_j$ la
repr\'esentation de degr\'e un de caract\`ere $\, \mu_j (\sigma) =
\zeta^j$, de sorte que la d\'ecomposition (7.26) \'equivaut \`a
$$\widehat {\Cal O}_{C,p}\,\cong\,\widehat {\Cal O}_{D,q} \bigoplus
\big( \bigoplus_{j=1}^{e-1}
\mu_j \otimes \widetilde {\widehat {\Cal O}_{D,q}}\big ) \tag (7.27) $$
Cela  donne finalement
$$\pi_*(\widehat {\Cal O}_C)_q \,\cong\, \widehat {\Cal
O}_{D,q}\otimes { \Ind}_H^G (1_H)
\bigoplus \big( \bigoplus_{j=1}^{e-1}\widetilde {\widehat {\Cal
O}_{D,q}}\otimes {
\Ind}_H^G (\mu_j) \big )  \tag (7.28) $$
De cette d\'ecomposition  on extrait le fait que  le facteur ${\Cal
L}_\chi$ est sans
torsion de rang \'egal au degr\'e de $\chi$. Par   r\'eciprocit\'e de Frobenius  un caract\`ere irr\'eductible de $G$
appara\^\i t dans $ \Ind_H^G (\mu_j)$ que si la restriction de $\chi$
\`a $H$
contient $\mu_j$. On voit ainsi que   ${\Cal L}_\chi$ est localement
libre en $q = \pi (p)$
si et seulement si la restriction de $\chi$ \`a $H$ est le caract\`ere trivial.\footnote {La d\'ecomposition 7.28 montre en fait que
$({\Cal L}_\chi )_q \cong  {\widehat {\Cal O}_{D,q}}^r \oplus
(\widetilde{ {\widehat {\Cal O}_{D,q}}})^s$, avec $r = \langle 1 ,
\chi\rangle_H$.} En particulier,
${\Cal L}_\chi$ est localement libre si et seulement si $\chi$ est
non ramifi\'e, signifiant par d\'efinition  $\chi (I) = 1$, i.e.
$\chi$ est trivial sur le sous-groupe d'inertie.

Pour terminer la preuve, on montre d'abord que le noyau de
$\pi^* : \Pic (C/G) \longrightarrow \Pic (C) $
s'identifie canoniquement \`a $\widehat { G/I}$. C'est un fait bien
connu, du moins dans
le cas lisse.  Adaptons la d\'emonstration  aux conditions  de la
pr\'esente situation. Soit ${\Cal L}$ un \'el\'ement du noyau, de
sorte que $\pi^*({\Cal L}) \cong {\Cal O}_C$  d\'efinit une $G$-lin\'earisation  de ${\Cal O}_C$.  La courbe   $C$ \'etant r\'eduite et
connexe, une telle lin\'earisation est d\'efinie par un caract\`ere
$\chi$ de $G$. Soit alors  ${\Cal O}_C (\chi)$, le faisceau
structural de $C$, ainsi $G$-lin\'earis\'e. On a sous ces
conditions   
$${\Cal L} \cong \pi_*^G \pi^* ({\Cal L}) \,\cong\,\pi_* \left( {\Cal
O}_C (\chi)\right) = {\Cal L}_\chi$$
Noter qu'alors $\chi$ est non ramifi\'e. Inversement, si $\chi$ est
un tel caract\`ere, on a un
morphisme canonique $\pi^* ({\Cal L}_\chi)  \to {\Cal O}_C$, entre
faisceaux  inversibles qui
est un isomorphisme sur un ouvert dense,  donc un isomorphisme.
Cela prouve que  le facteur isotypique ${\Cal L}_\chi$ de la d\'ecomposition (7.23),  dans le cas o\`u c'est un faisceau inversible,
c'est \`a dire si $\chi$ est non ramifi\'e,  est dans le noyau de
$\pi^*$.  L'identification cherch\'ee est sous une forme explicite
$\chi \in \widehat {  G/I} \rightsquigarrow {\Cal L}_\chi$.
Le m\^eme argument appliqu\'e aux composantes $\tilde C_i$ de la
normalisation donne
$$\ker {\pi^* : \Pic \big({\tilde C_i}/{G_i} \big) \longrightarrow \Pic ({\tilde
C_i})}\, \cong \, \widehat {{G_i} /{I_i}}\tag (7.29) $$
Le r\'esultat  d\'ecoule alors, apr\`es identification, de la suite
exacte des noyaux. \hfill\break
\qed 
  
\proclaim  Corollaire 7.9.  Le facteur isotypique $\Cal L_\chi$ de la
d\'ecomposition (7.23) est
localement libre (de rang \'egal au degr\'e de $\chi$) si et
seulement si $\chi$ est non ramifi\'e, soit $\chi (I) = 1$.
 \hfill\break
\qed

\rema{7.4}    Le corollaire pr\'ec\'edent
apporte une r\'eponse \'equivalente \`a l'alternative qu'on
rencontre usuellement lorsque dans une compactification, un faisceau
inversible d\'eg\'en\`ere  lors d'une sp\'ecialisation en un point du
bord. Pour comparaison avec l'\'etude du bord de l'espace modulaire
des courbes \`a spin, probl\`eme au parall\'elisme frappant, la
correspondance avec la terminologie employ\'ee par Jarvis \cite
{42},\cite {43} est
$$ {\chi} \quad\hbox{ \rm est}\quad \cases { \hbox{ \rm non  ramifi\'e   si} &$\, \, \chi\,$    est 
R-R \   (Ramond-Ramond) \cr
 \hbox { \rm  ramifi\'e  si
 }&  $\,\,\chi\,$    est  N-S 
(Neveu-Schwarz)\cr } $$ 
  $\lozenge $ 
 \bigskip
 
   {\it 7.3.2. Faisceaux sans torsion de rang un et rev\^etements
stables}  
\bigskip

Soit $p:D\to S$ une courbe pr\'estable.  
Un faisceau  coh\'erent $E$ sur $D$ est dit sans torsion de rang $n$,
relativement \`a $S$, si $E$ est plat sur $S$, et si sur chaque
fibre $D_s$, le faisceau induit $E_s$   est sans torsion de rang
$n$. En particulier $E$ est localement libre sur l'ouvert de lissit\'e\/ de $p$. Un th\'eor\`eme de Faltings (\cite{31} Theorem 3.5, voir
aussi \cite {43}) pr\'ecise la structure locale  de $E$ en un point
singulier, point en lequel $E$ n'est pas libre.    Pour d\'ecrire le r\'esultat   on peut choisir des coordonn\'ees locales le long des
branches, et supposer que $E$ est un module sans torsion de rang $n$
sur l'anneau local noetherien complet  $A = R[[x,y]]/(xy - \pi)$, o\`u
 $R = \hat{\Cal O}_s$ et $\pi \in {\Cal M}_R$. Si $\overline A =
A\otimes k$, alors on sait que  

$$E\otimes k \cong  {\overline A}^r
\oplus {\tilde{\overline A ^s}} \,\,(r+s=n)\tag (7.30)$$
 Comme  indiqu\'e dans \cite {31}, on peut se limiter au cas $n=s$.  La
construction d'une d\'eformation verselle de ${\tilde{\overline A
^n}} $ est comme suit. Noter que seul le cas $n= 1$ sera en fait
utilis\'e. Supposons que $P,Q \in \M_n(R)$ sont deux matrices
$n\times n$, telles que $P.Q = Q.P = \pi 1_n$. D\'efinissons deux
matrices  $2n\times 2n$, \`a coefficients dans $A$, par  
$$ \Phi
= \pmatrix {x.1 & P\cr Q & y.1 \cr}, \quad \Psi = \pmatrix{ y.1 &
-P\cr -Q & x.1 \cr} $$ 
Alors $\Phi\Psi = \Psi\Phi = 0$, et
$E(P,Q) = \ker \Psi = \Ima \Phi $ est sans torsion de rang $n$, de
r\'eduction $ {\tilde{\overline A ^n}} $. R\'eciproquement, tout
module sans torsion de rang $n$, de r\'eduction $ {\tilde{\overline A
^n}} $, est isomorphe \`a un $E(P,Q)$. On peut prouver  (loc.cit)
que si $P \equiv Q\equiv 0 \pmod {\Cal M}_R$, alors le couple $(P,Q)$
est d\'etermin\'e  \`a la transformation pr\`es 
$$ P' =
UPV^{-1}, \,\, Q' = VQU^{-1} \quad (U,V\in \GL_n(R))$$  Revenons \`a
 un $G$-rev\^etement galoisien stable $\pi: C \to D$ au-dessus de
la base $S = \Spec (R)$, l'anneau $R$ \'etant local noeth\'erien
complet. Soit $P \in C$ un point double de la fibre sp\'eciale, et
soit $Q = \pi(P)$. Il s'agit de d\'ecrire selon les termes de la
construction de Faltings, le module sans torsion $E = \widehat
{(\pi_* (\Cal O_C))}_Q$, de rang $n = \vert G\vert$.  Soient $x,y$
des coordonn\'ees locales le long des branches en $P$, telles qu'un
g\'en\'erateur $\tau$ du stabilisateur $H$ de $P$ agisse par  $\tau
(x) = \zeta x, \, \tau (y) = \zeta^{-1}y$, pour une certaine racine
$e$-i\`eme de l'unit\'e. On peut supposer que  $A = \widehat
{\Cal O_C}_P = R[[x,y]]/(xy - \pi) \quad (\pi\in {\Cal M}_R)$,   en
cons\'equence $B = A^H = \widehat {\Cal O_D}_Q = R[[u,v]]/(uv -
\pi^{e})$, avec $u=x^{e}, \, v=y^{e}$. Un \'el\'ement  $\xi\in A$
admet une \'ecriture unique de la forme
$$\xi = f(x) + g(y), \,\,
f(x)\in R[[x]], g(y) \in yR[[y]]$$ 
Il est facile de voir que le
facteur isotypique de caract\`ere $\zeta^\alpha$ de $A$, est le sous
$B$-module $\Cal L_\alpha$ form\'e par les \'el\'ements de la
forme 
$$ x^\alpha f(u) + y^{e-\alpha} g(v), \,\,( f\in R[[u]],
\,g\in R[[v]]) $$
 Il est ais\'e d'identifier ce module; on note
dans la suite $\beta = e - \alpha$. 

\proclaim Lemme 7.10. Le module
$\Cal L_\alpha$ est le conoyau de l'application $B^2\to B^2$ d\'efinie par la matrice  $\Phi = \pmatrix{ u & \pi^\alpha \cr \pi^\beta &
v\cr  }$; il est en particulier sans torsion de rang un, et
$\Cal L_\alpha \cong E(\pi^\alpha,\pi^\beta)$. 
\dem Soit
l'application $\Psi: B^2 \to A$ d\'efinie par la matrice $\pmatrix
{y^\beta &, -x^\alpha  } $. On a clairement $\Psi.\Phi = 0$.
Montrons que $\ker \psi = \Ima \Phi$. Soit $(f,g)\in B^2$ tel que
$y^\beta f - x^\alpha g = 0$. Ecrivons $f$ et $g$ sous la forme
$$\cases {f = f_+(u) + a_0 + f_-(v)\cr g = g_+(u) + b_0 + g_-(v)
\cr }$$
  les s\'eries $f_+,f_-,g_+,g_-$ \'etant sans terme
constant, et $a_0,b_0\in R$. Un calcul \'el\'ementaire
montre que l'\'equation $y^\beta f - x^\alpha g = 0$ implique les
deux \'egalit\'es: 
$$\cases {g_+(u) = {\pi^\beta\over u} f_+(u) -
b_0\cr
 f_-(v) = {\pi^\alpha \over v} g_-(v) - a_0\cr}  $$
  En
posant $\phi(u) = {f_+(u) \over u}, \, \psi(v) = {g_-(v)\over v}$, on
obtient bien que  $$ \pmatrix{ f \cr g  \cr}  = \pmatrix{ u &
\pi^\alpha \cr \pi^\beta & v\cr}  \pmatrix {\phi(u) \cr \psi(v)
 \cr} $$ \hfill\break\qed  

En utilisant les notations de dessus notons que
  $$\widehat {\Cal O}_{C,P} = \widehat {\Cal O}_{D,Q}
\bigoplus \left( \bigoplus_{\alpha = 1}^{e-1} E(\pi^\alpha,\pi^\beta)
\otimes W_\alpha\right)$$  o\`u  $W_\alpha$  d\'esigne la repr\'esentation
de degr\'e un de $H$ telle que $\tau\in H$ agisse par la
multiplication par $\zeta^\alpha$. D'o\`u par induction de $H$ \`a
 $G$ 
 $$ \widehat {\pi_*(\Cal O_C)}_Q = \widehat {\Cal O}_{D,Q}
\otimes \Ind_H^G(1) \bigoplus \left( \bigoplus_{\alpha=1}^{e-1}
E(\pi^\alpha,\pi^\beta) \otimes \Ind_H^G(W_\alpha)\right) \tag
(7.31)$$ 
Soit  pour tout caract\`ere irr\'eductible $\chi$ de $G$,
$n_{\alpha,\chi} $   la multiplicit\'e de $W_\alpha$ dans la
restriction de $\chi$ \`a $H$. Le facteur isotypique $\Cal L_\chi$
admet  finalement la structure suivante en le point double $Q \in
D$: 

\proclaim Proposition 7.11. La structure du  module sans
torsion $(\widehat {\Cal L}_\chi)_Q $ est donn\'ee par: $$ (\widehat
{\Cal L}_\chi)_Q  \cong \bigoplus_{\alpha=0}^{e-1}
E(\pi^\alpha,\pi^\beta)^{n_{\alpha,\chi}} \tag (7.32)$$
 \qed  
 
 On notera que la situation correspond  \`a une matrice $P$
diagonale. Le cas $G$ ab\'elien  sera d\'etaill\'e dans le
paragraphe 8.  Le fait que les faisceaux $\Cal  L_\chi$  peuvent \^etre sans torsion de rang un, non  inversibles, peut \^etre contourn\'e  par la construction suivante, inspir\'ee par la proc\'edure de stabilisation de Knudsen \cite {48} (voir aussi Jarvis \cite {44} \S 3.1.2). Soit une courbe pr\'estable $D \to S$, et soit $\Cal L$ un ${\Cal O}_D$-module sans torsion de rang un relativement \`a   $S$. Soit le $D$-sch\'ema
$$\rho: \tilde D = \Proj \,  (Sym^\bullet ({\Cal L}) \rightarrow D\tag (7.33)$$
Si $\Cal L$ est localement libre, alors $\tilde D = D$. Si $\Cal L$ n'est que sans torsion de rang un de lieu singulier
$Sing(\Cal L)$, alors le faisceau inversible ${\Cal O}_{\tilde D} (1)$ sur $\tilde D$ r\'esoud les singularit\'es de $\Cal L$ dans le sens suivant  (Jarvis \cite{44}):

\proclaim Proposition 7.12.   Le sch\'ema $\tilde D$ est une $S$-courbe pr\'estable. \hfill\break
i) On a $\rho_*({\Cal O}_{\tilde D}) = {\Cal O}_D$ et $\rho_* ({\Cal O}_{\tilde D} (1)) = \Cal L$. \hfill\break
ii) Si $j>0$ et $ n\geq 0,$ $R^j\rho_*({\Cal O}(n)) = 0$. \hfill\break
iii) La formation de $\tilde D$ ainsi que de $\rho_* ({\Cal O}_{\tilde D} (1))$ commute \`a tout changement de base.  De plus   $\omega_{\tilde D/S} = \rho^* (\omega_{D/S})$, de sorte que la restriction de $ \omega_{\tilde D/S}$ \`a toute composante exceptionnelle est triviale.
  
  \dem La preuve est   contenue dans les r\'ef\'erences pr\'ec\'edentes, particuli\`erement \cite {44}, Lemma 3.1.4, \cite {48}, Theorem 2.4. \hfill\break
\qed

Il est clair que dans cette construction le morphisme $\rho$ est un isomorphisme en dehors du lieu singulier de $\Cal L$. Il est utile d'avoir une description explicite de $\tilde D$ au voisinage d'un point singulier de $\Cal L$. Soit $Q \in D_s$ un tel point, et rempla\c cons $D$ par $\Spec \hat{\Cal O}_{D,Q}$ et $\Cal  L$ par $\hat {\Cal L}_D$. On suppose que  $\hat{\Cal O}_{D,Q} = \hat {\Cal O}_s [[x,y]] /(xy - \pi)$ et que $\hat {\Cal L}_D = E(p,q)$ avec $pq = \pi, \,\, p,q \in \hat {\Cal M}_s$. Du fait  de la pr\'esentation 
$${\hat{\Cal O}_{D,Q}}^2 \buildrel { \pmatrix{ y & -p\cr-q & x \cr }}\over\longrightarrow {\hat{\Cal O}_{D,Q}}^2 \longrightarrow E(p,q)\longrightarrow 0\tag (7.34)$$
on a localement en $Q$  
$$ {Sym^\bullet}_Q (E(p,q) = \hat {\Cal O}_Q[\xi , v]/ (y\xi - qv , -p\xi + xv) $$
En particulier  $\tilde D \subset \Bbb P^1\times D$ est le sous-sch\'ema d'\'equations 
$y\xi - qv = -p\xi + xv = 0.$
  La fibre au-dessus du point $Q$ est $E = \Bbb P^1$. Noter que les deux branches du point double $Q$ ont des transform\'ees strictes disjointes. Elles coupent $E$ en deux points qui sont des points doubles de la fibre de $\tilde D$, d'\'equations respectives  $ys = q$ et $xt = p$. Le morphisme $\rho$ est donc une contraction au sens de Knudsen\footnote{C'est une contraction particuli\`ere car les courbes exceptionnelles sont reli\'ees \`a deux autres composantes, et donc se contractent en un point double. La fibre au-dessus d'un point singulier de $\Cal L$ est donc de la forme $E_1\cup E\cup E_2$, $E$ \'etant la composante exceptionnelle, et $E_1\cap E = p, \, E_2\cap E = q$. Il est clair que le faisceau dualisant a une restriction \`a $E$ triviale.}, $D$ est la stabilisation de $\tilde D$. Il en d\'ecoule ais\'ement la relation $\omega_{\tilde D/S} = \rho^* (\omega_{D/S})$ (\cite {48}, Cor  1.5), en particulier la trivialit\'e de 
$\omega_{\tilde D/S}$ sur chaque composante exceptionnelle, ce qui est ailleurs clair directement. La d\'esingularisation de ${\cal L}$ ainsi construite est unique \`a un isomorphisme pr\`es \cite{44}.

\bigskip

\beginsection 7.4. Stratification canonique du bord

\bigskip
  
 Une propri\'et\'e appr\'eciable de  
   $\overline {\Cal M}_{g,n}$ est que le bord admet une stratification naturelle par le type combinatoire d'une courbe stable marqu\'ee (resp. piqu\'ee) (\cite{53} \ \S \ 2.7). Fixons un graphe modulaire $\Gamma$ de genre $g$, avec $n$
pattes (D\'efinition 7.1). Pour fixer les id\'ees on suppose le
graphe marqu\'e,   les pattes sont  alors num\'erot\'ees  de $1$ \`a
$n$.  Soit pour tout sommet $v\in V$, $h_v$ (resp.  $l_v$) le nombre
de pattes (resp. d'ar\^etes) incidentes en $v$. On d\'efinit un
morphisme  
$$\beta_\Gamma: \prod_{v\in V} \, \overline {\Cal
M}_{g_v,h_v+l_v} \longrightarrow  \overline {\Cal M}_{g,n}Ê\tag (7.35)$$ 
 de la mani\`ere suivante. On commence par num\'eroter l'ensemble $V$,  donc $V = \{v_1,\cdots,v_c\}$, et pour tout sommet $v =v_i$,  $e_1^v,
\dots, e_{l_v}^v$ les ar\^etes d'origine $v$. Noter que le graphe \'etant marqu\'e, les pattes incidentes \`a $v$ sont ordonn\'ees par l'ordre induit de $[1 , n]$, soit $p_1^v,\dots, p_{h_v}^v$ cette liste
ordonn\'ee. La donn\'ee d'un objet  $((C_i, x_1^{i},\dots
x_{h_i+l_i}^{i})_{v_i\in V})$ au-dessus de la base $S$ de $\prod_{v\in V}
\, \overline {\Cal M}_{g_v,h_v+l_v} $ peut s'interpreter comme un
"plongement", ou r\'ealisation de ${\Cal F}l(\Gamma)$:   
$$\imath: {\Cal F}l (\Gamma)
\hookrightarrow \bigsqcup_{v\in V} C_v = \bigsqcup_{v_i\in V} C_i  \tag (7.36)$$
Cela signifie  que $\imath$
identifie les ar\^etes de $\Gamma$ avec l'ensemble des points sp\'eciaux de  la somme disjointe $\bigsqcup_{v\in V} C_v$. Par convention les pattes d'origine $v$ ont pour images (ordonn\'ees) $x_1^v,\dots,
x_{h_v}^v$, et les ar\^etes d'origine $v$, ont pour images (ordonn\'ees) $x_{h_v+1}^v, \dots, x_{h_v+l_v}^v$.  R\'eciproquement un   plongement  (7.36)  conduit \`a un objet de  $\prod_{v_i\in V} \, \overline {\Cal
M}_{g_v,h_v+l_v}$.  Par identification des
sections $x_i^v, \,x_j^w$ qui correspondent \`a un couple d'ar\^etes
oppos\'ees (clutching morphism \footnote {Rappelons la d\'efinition de l'op\'eration de
recollement (clutching morphism) le long d'une paire de sections
(Knudsen \cite {48}, Manin \cite
{53} chapter 5). La construction se r\'esume par le r\'esultat
suivant:\hfill\break
  Soit $C'$ une courbe pr\'estable de base $S,$ dont on ne suppose pas
\`a priori les fibres
g\'eom\'etriques connexes. Soient deux sections $s_1,s_2 : S \to C'$,
dont les images sont form\'ees
de points non singuliers le long des fibres.  Il existe une
courbe pr\'estable $C'$, un
morphisme fini $p: C \to C'$ tel que $ps_1 = ps_2$, le couple
$(C',p)$ \'etant universel en un
sens \'evident. En particulier, il est d\'efini \`a  un
isomorphisme canonique pr\`es.  Une fibre g\'eom\'etrique \'etant
donn\'ee, si la fibre image de $C'$ est connexe, deux cas se
pr\'esentent:\hfill\break
$\bullet$ La fibre est irr\'eductible, dans ce cas on cr\'ee dans la
fibre correspondante de $C'$ une
boucle,    on a: $g' = g+1$.\hfill\break
$\bullet$ La fibre de $C$ est somme disjointe de deux courbes
connexes $C_1$ et $C_2$, avec
$s_i \in C_i (i=1,2)$. Dans ce cas la fibre (connexe) de $C'$ a pour
genre $g' = g_1+g_2$.}
  de Knudsen \cite {47}), on forme la courbe  $C = \bigsqcup_{v\in V}
\, C_v / \Gamma \in \overline {\Cal M}_{g,n}(S)$. Observer qu'il en
d\'ecoule une identification bien d\'efinie $\Gamma = \Gamma (C)$. On
note aussi qu'il y a une action \'evidente (\`a droite) de $\Aut
(\Gamma)$ (en fait de $\Aut ({\Cal F}l(\Gamma))$)  sur de tels plongements, i.e.  $\sigma. (\{ C_v\}_v ,
\imath) = (\{C_v\} ,   \imath.\sigma)$. Alors le morphisme
$\beta_\Gamma$ est un $\Aut (\Gamma)$-torseur  (dans le sens de la
section 6.1) sur son image \footnote { Le morphisme $\beta_\Gamma$
est fini, en particulier repr\'esentable  \cite {48}, de sorte que
son image, en tant que sous champ ferm\'e de $\overline {\Cal
M}_{g,n}$ est d\'efinie.}  $\overline {\Cal M}_{g,n} (\Gamma)$;
 qui est  en cons\'equence le quotient   
$$\prod_{v\in V} \,
\overline {\Cal M}_{g_v,h_v+l_v}  / \Aut (\Gamma) \tag (7.37)$$
Soit  maintenant $(\Gamma,G)$ un graphe modulaire de Hurwitz (D\'efinition 7.2). On consid\`ere le sous champ
localement ferm\'e  $ {\Cal H}_{g,G,\xi} (\Gamma)$ de
$\overline {\Cal H}_{g,G,\xi}$,   
dont les objets sont les $G$-courbes stables de type combinatoire
fix\'e  $\Gamma$, i.e. dont les
fibres g\'eom\'etriques sont de type $\Gamma$. Si $\Gamma
= \emptyset, {\Cal H}_{g,G,\xi}(\Gamma) =  {\Cal H}_{g,G,\xi}$.   Les
sous-champs (resp. sous-sch\'emas) localement ferm\'es ${\Cal
H}_{g,G,\xi} (\Gamma)$ (resp. les espaces de modules grossiers $ H_{g,G,\xi}(\Gamma)$) forment une stratification  de
$\overline {\Cal H}_{g,G,\xi}$  (resp. $\overline H_{g,G,\xi}$), analogue  \`a 
stratification    de $\overline {\Cal M}_{g,n}$ par le type
combinatoire  des courbes stables marqu\'ees \cite {53}. On a  
la d\'ecomposition
$$ \overline {\Cal H}_{g,G,\xi} = \bigsqcup_{\Gamma}  \,\, {\Cal
H}_{g,G,\xi}(\Gamma)  \tag (7.38)$$
 Le morphisme "discriminant",
associant \`a un rev\^etement sa
base marqu\'ee par les points de branchement
$ \delta : \overline {\Cal H}_{g,G,\xi} \longrightarrow \overline
{\Cal M}_{g',b}$
est  visiblement compatible aux stratifications,  ce qu'on peut traduire par le
fait que $\delta$ induit un morphisme  $ {\Cal H}_{g,G,\xi}(\Gamma)
\rightarrow {\Cal M}_{g',b} (\Delta), \,\, \Delta = \Gamma / G $. D'une
autre mani\`ere, si $r$ est le nombre de points de branchement, soit
le morphisme (oubli de l'action de $G$)  $\imath: \overline {\Cal
H}_{g,G,\xi} \rightarrow \overline {\Cal M}_{g,r}$. On a clairement
$  {\Cal H}_{g,G,\xi}(\Gamma) = \imath^{-1} \left(
  {\Cal M}_{g,r} (\Gamma)\right)$
Dans le terme de droite, $\Gamma$ est le graphe modulaire, all\'eg\'e de la donn\'ee de Hurwitz.

Avant de pr\'eciser la structure de la strate
${\Cal H}_{g,G,\xi}(\Gamma)$ des rev\^etements de type combinatoire $\Gamma$,  revenons  \`a un  graphe modulaire de Hurwitz $\Gamma$ (D\'efinition 7.2).
Soit $v \in \Gamma$ un sommet  et notons $st(v)$ l'\'etoile de $v$.  Les ar\^etes  \'el\'ements de $st(v)$, d\'ecor\'ees par  l'holonomie  $(H,\chi),$ sont      de deux sortes.
   D'une part les ar\^etes orient\'ees d'origine un sommet $v$ et
d'extr\'emit\'e $w$, $v=w$ \'etant possible (boucle), et d'autre
part les ar\^etes monovalentes correspondantes aux  pattes  issues de $v$. 
Faisons  tout d'abord le choix d'un syst\`eme
de repr\'esentants $v_1,\dots,v_k$ pour l'action de $G$ sur
 $V,$  ensemble des  sommets de $\Gamma$, i.e une section de $V\to V/G$.  Soit $G_{i}$ le stabilisateur de $v_i$, et notons $\st_i  = {\Cal E}_i \bigsqcup{\Cal L}_i$ l'\'etoile de $v_i$, union de  l'ensemble ${\Cal E}_i$ des ar\^etes (d\'ecor\'ees) issues de $v_i$ et de l'ensemble ${\Cal L}_i$   des pattes d'origine $v_i$. Du point de vue des $G$-ensembles, on a si ${\Cal A} = {\Cal L}\bigsqcup {\Cal E}$
 $$ {\Cal L} = \bigcup_i \Ind_{G_{i}}^G \, {\Cal L}_i, \,\,\, {\Cal E} = \bigcup_i G\times^{G_{i}} \, {\Cal E}_i$$
  Pour avoir l'analogue \'equivariant de (7.37), il nous faut   un plongement \'equivariant  ${\Cal F}l(\Gamma)$ dans ce qui remplace le second membre  de (7.36), c'est \`a dire  
  $$\tilde C =  \bigsqcup_{1\leq i \leq k}\Ind_{G_{i}}^G  C_i$$
    Comme seules les orbites de points sp\'eciaux sont index\'ees,  il est  nous faut  ˆ ce stade   choisir un marquage plus strict pour un rev\^etement galoisien. 
  
 Outre un syst\`eme de repr\'esentants des sommets, on va choisir pour chaque  classe de conjugaison  d'holonomie d\'ecorant une $G_i$-orbite de pattes ou ar\^etes  d'origine $v_i$, un repr\'esentant d'holonomie $(H_{i,\alpha},\chi_{i,\alpha})$ pour les ar\^etes, et $(H_{i,\lambda},\chi_{i,\lambda}) $  pour les orbites de pattes.  Ceci \'etant, on fait  ensuite le choix d'un point $e_{i,\alpha}$  (resp. $l_{i,\lambda}$) dans chaque orbite  d'ar\^etes qui a pour  holonomie  le repr\'esentant  choisi. D'une autre mani\`ere  on se donne   des sections de ${\Cal E}_i \to {\Cal E}_i/G_i, \,\, {\Cal L}_i \to {\Cal L}_i/G_i$ avec holonomie prescrite.

Soit maintenant un rev\^etement $\pi: C\to D$ de type combinatoire $\Gamma$, et 
$\widetilde C = \bigsqcup_v \,\, {\widetilde C}_v$ la normalisation
de $C$. La composante
${\widetilde C}_v$ munie de l'action du stabilisateur $G_v$ a pour
donn\'ee de Hurwitz
$$\tilde \xi_v = \xi_v + \eta_v \tag (7.39) $$
o\`u   $\xi_v$ repr\'esente la
contribution des  points de ramification de $C$, et $\eta_v$ est la
 contribution des ar\^etes
orient\'ees d'origine $v$, branches des points doubles.
Ceci nous conduit \`a  voir $\widetilde  C_v$ comme un
objet de $ {\Cal
H}_{g_v,G_v,{  \xi}_v,\eta_v}$  champ classifiant les courbes de
genre $g_v$, munies d'une
action de $G_v$ de donn\'ee de Hurwitz $ \tilde\xi_v =  \xi_v + \eta_v$. Dans cette d\'efinition  les courbes sont marqu\'ees   par un diviseur $G_v$-
invariant, ce qui signifie que seules les $G_v$-orbites de points sp\'eciaux sont num\'erot\'ees. 


   Avec la normalisation fix\'ee sur le graphe modulaire de Hurwitz  $\Gamma$,  soit un objet $C_i   \in {\Cal
H}_{g_{v_i},G_{v_i}, \xi_{v_i},\eta_{v_i}}$. Sur cette courbe, seules les orbites de points speciaux sont num\'erot\'ees. Pour tout $i$, et dans chaque orbite de tels points, on fait comme pour le graphe  le choix d'un point  $x_{i,\alpha}$ (resp.  $y_{i,\lambda}$)  dans chaque $G_i$-orbite, dont l'holonomie est celle s\'electionn\'ee.   On a donc  pour $1\leq \alpha \leq h_i$,  $x_{i,\alpha}$ d'holonomie $(H_{i,\alpha},\chi_{i,\alpha})$, et si $1\leq \lambda \leq l_i$ le point $x_{i,h_i+\lambda} = y_{i,\lambda}$ d'holonomie $(H_{i,\lambda}, \chi_{i,\lambda})$. Notons  alors pour tout  $1\leq i\leq k, \,  {\Cal H}^* _{g_{v_i},G_{v_i}, \xi_{v_i},\eta_{v_i}}$ 
le champ de Hurwitz modifi\'e, r\'esultat de ce marquage plus strict. 
Noter que l'oubli de ce choix  suppl\'ementaire conduit  \`a un morphisme 
 $${{\Cal H}^*}_{g_{v_i},G_{v_i}, \xi_{v_i},\eta_{v_i}}\to {\Cal H}_{g_{v_i},G_{v_i}, \xi_{v_i},\eta_{v_i}}$$
 de degr\'e $\prod_{\alpha,\lambda} \vert Z_{G_i}(H_{i,\alpha})/H_{i,\alpha}\vert \vert Z_{G_i}(H_{i,\lambda})/H_{i,\lambda}\vert$  ($Z_G(H)$ est le centralisateur de $H$ dans $G$). Ce morphisme est le quotient par un produit de groupes sym\'etriques   convenables. 

On d\'efinit  de m\^eme le champ ${\Cal H}^*_{g,G,\xi}(\Gamma) $, avec le morphisme ${\Cal H}^*_{g,G,\xi}(\Gamma)\to {\Cal H}^*_{g,G,\xi}(\Gamma)$ de degr\'e $\prod_{i,\lambda} \vert Z_{G_i}(H_{i,\lambda})/H_{i,\lambda}\vert$.  On peut formuler l'analogue de (7.37):

\proclaim Proposition 7.13.  Il y a  un    morphisme  naturel 
$ \xi_\Gamma :  \prod_{i=1}^k \,\,{\Cal
H}^*_{g_{v_i},G_{v_i}, \xi_{v_i},\eta_{v_i}}
 \longrightarrow   {\Cal H}^*_{g,G,\xi} (\Gamma)  $
qui  identifie le terme de droite  au quotient
de la source par le groupe
$\Aut_G(\Gamma)$.  Ce morphisme s'ins\`ere dans un diagramme
$$   \prod_{i=1}^k \,\,{\Cal
H}_{g_{v_i},G_{v_i}, \xi_{v_i},\eta_{v_i}}
  \longleftarrow   \prod_{i=1}^k \,\,{\Cal
H}^*_{g_{v_i},G_{v_i}, \xi_{v_i},\eta_{v_i}} \buildrel {\xi_\Gamma}\over \longrightarrow  {\Cal H}^*_{g,G,\xi} (\Gamma)\longrightarrow {\Cal H}_{g,G,\xi}(\Gamma) \tag (7.40)$$   

\dem Il faut tout d'abord d\'efinir le morphisme $\xi_\Gamma$. Cela revient \`a pr\'eciser la proc\'edure d'identification des branches par paires. Cette proc\'edure  doit \^etre   $G$-\'equivariante.  Elle  revient comme dans le cas $G = 1$, \`a plonger de  mani\`ere \'equivariante l'ensemble $ {\Cal F}l(\Gamma)$  sur l'ensemble correspondant des points sp\'eciaux de  $\bigsqcup_i {\Ind}_{G_{v_i}}^G (C_i) $,  puis \`a transporter l'involution canonique.      Il est clair que les  bijections entre points s\'electionn\'es  s'\'etendent en   un plongement $G$-\'equivariant 
$${\Cal F}l(\Gamma) \hookrightarrow  \tilde C = \bigsqcup_i {\Ind}_{G_{v_i}}^G (C_i) $$
typiquement $ge_{i,\alpha} \to gx_{i,\alpha}$ .  On peut alors transporter l'involution $\tau: a\to \overline a$  sur l'image dans le  terme de droite,  donnant de la sorte un
proc\'ed\'e d'identification par paires  des points distingu\'es de
$\tilde C$.  
Par recollement le long de ces paires de sections, on  obtient une
$G$-courbe stable marqu\'ee connexe de donn\'ee de ramification
$\xi$.   De mani\`ere plus concise  on peut \'ecrire:
$$C = \bigvee_i  \Ind_{G_{v_i}}^G (C_i)  \in {\Cal
H}_{g,G,\xi}(\Gamma) \tag (7.41) $$
On pose alors $\xi_\Gamma (\{ C_i\}) = C$. Cette construction   d\'efinit le foncteur $\xi_\Gamma$.  
Le reste est cons\'equence des remarques qui pr\'ec\`edent l'\'enonc\'e.\hfill\break \qed

 Notons que la codimension de  ${\Cal H}_{g,G,\xi} (\Gamma) $ est $ \Card \,{\Cal F}_{ \Gamma/ G} $, et que la construction est d\'efinie en fait au niveau des champs compactifi\'es.

  Les rev\^etements sont marqu\'es par les points de branchement r\'eels ou virtuels.   Supposons que  la donn\Ž e de ramification $\xi$ soit telle que l'un des sous-groupe $H_i$ soit trivial, par exemple $H_{b+1} = 1$, signifiant que le point  marqu\'e $Q_{b+1}$ n'est pas un point de branchement au sens strict. On a donc avec les notations  du \S 2, $\xi = \xi' + [1]$.  Par contraction (stabilisation) on peut effacer ce point.
\proclaim  Proposition 7.14. L'effacement du point de "branchement" d'indice $b+1$  d\'efinit  un morphisme de champs  $\rho:  \overline {\Cal H}_{g,G,\xi} \longrightarrow  \overline {\Cal H}_{g,G,\xi'}$.

\dem Soit un rev\^etement stable $\pi: C \to D$ de base $S$. L'oubli du point de branchement $Q_{b+1}$, suivi d'une stabilisation, donne un morphisme $\psi: D \to D'$. Soit d'autre part l'oubli  suivi d'une stabilisation des points de l'orbite r\'eguli\ re $\pi^{-1}(Q_{b+1})$, conduisant \`a un morphisme $\phi: C \to C'$. Par le caract\ re universel de cette op\Ž ration, l'action de $G$ se descend \`a  $C'$, et le morphisme $\psi \pi$ factorise par $C'$. Soit $\pi': C' \to D'$ cette factorisation. Le morphisme $\pi'$ factorise en $C' \to C'/G \to D'$ et $\psi$ en $D=C/G \to C'/G \to D'$. Comme la courbe $C'/G$ est stable marqu\Ž e par les images des points $Q_1,\cdots,Q_b$, il y a un morphisme $D'\to C'/G$ qui est clairement l'inverse de $C'/G \to D'$.  Donc $D' = C'/G$.  Le morphisme $\rho$ est  ainsi donn\'e par $\rho (C\to D) = (C'\to D')$. \hfill\break \qed

   On reprend les notations et hypoth\`eses de la section 7.2.1. Consid\Ž rons donc $G_i \subset G \,(i=1,2)$ deux sous-groupes, et un sous-groupe cyclique $H\subset G_1\cap G_2$. Soit des donn\'ees de ramification 
$$\xi_1^* = \xi_1 + [H,\chi]\in R_+(G_1), \quad \xi_2^* = \xi_2 + [H,\chi^{-1}]\in R_+(G_2)$$
Posons $\xi = \Ind_{G_1}^G \xi_1 + \Ind_{G_2}^G \xi_2\in R_+(G)$. Par recollement   le long de deux  sections, on obtient:
\proclaim Proposition 7.15. Il existe un morphisme de recollement le long de deux sections d'holonomie oppos\Ž es:
$$\rho: \overline {\Cal H}_{g_1,G_1,\xi_1^*} \times \overline {\Cal H}_{g_2,G_2,\xi_2^*} \rightarrow \overline {\Cal H}_{g_1+g_2,G,\xi} \tag (7.42)$$
Si $G = G_1=G_2$, et si $\xi^* = \xi + [H,\chi] + [H,\chi^{-1}]$ (Harbater-Mumford), par recollement des deux sections, on  a le morphisme analogue  $\rho: \overline {\Cal H}_{g-1,G,\xi^*} \rightarrow \overline {\Cal H}_{g,G,\xi}$.

 \dem  Soit $C_i\to D_i \,(i=1,2)$ un point du terme de gauche, et soient $Q_ i\in D_i$ les points de branchement d'holonomie   $[H,\chi_1 = \chi]$ et  $[H,\chi_2 = \chi^{-1}]$. Soit $P_i\in C_i$ un point d'holonomie exacte $(H,\chi_i)$. Par identification des paires de points $(gP_1,gP_2)$, on donne un sens  au rev\^etement  
$$C = \Ind_{G_1}^G C_1 \bigvee \Ind_{G_2}^G C_2\rightarrow D = D_1\bigvee D_2\tag (7.43)$$
Cette construction est bien d\'efinie \`a  isomorphisme unique pr\`es.  Elle d\'efinit le morphisme cherch\'e.  Une construction similaire    fonctionne dans le cas d'une donn\'ee de ramification du type {\it Harbater-Mumford} $\xi^* = \xi + [H,\chi] + [H,\chi^{-1}]$, et donne le second morphisme.\hfill\break\qed 
\medskip
Soit le morphisme discriminant $\delta: \overline {\Cal H}_{g,G,\xi}\to \overline{\cal M}_{g',b}$.  C'est un rev\^etement entre champs de Deligne-Mumford, dont a peut pr\'eciser le diviseur de ramification. Le r\'esultat est une simple traduction de la   proposition 6.28.
\proclaim Proposition 7.16.  Le diviseur de ramification du discriminant est 
$${\cal R} = \sum_{\Delta = NS} \,  (\vert H\vert - 1) \Delta\tag (7.44)$$
la somme \'etant \'etendue aux composantes NS du bord, donc celles telles que $H\ne 1$.\hfill\break
\qed

\beginsection 7.5. Type topologique d'un point du bord

\bigskip

Dans cette section $k = \Bbb C$. Soient   
  $H_1,\dots,H_{h(\xi)}$ les composantes irr\'eductibles de  $H_{g,G,\xi}$, $h(\xi)$ \'etant le nombre de Nielsen
(D\'efinition 2.4). Le   sch\'ema des modules
compactifi\'e  est  somme disjointe des "compactifications"
$\overline {  H}_{\alpha}, \,\, (\alpha = 1,\dots,h(\xi))$  (\S
6, Proposition 6.2). Si  $C\to D$ repr\'esente un point du bord, on peut
l\'egitimement chercher
\`a quelle composante  ce point  appartient,  c'est \`a  dire
trouver un proc\'ed\'e pour lire sur $C\to D$ le nom de la composante \`a
 laquelle le point est rattach\'e.  Si $C$ est
lisse, le type topologique qui  fixe la composante \`a laquelle ce
point appartient, est enti\`erement donn\'e par le morphisme de
monodromie $\psi:
\pi_1(D-\beta) \to G$ (Proposition 2.2). En fait une r\'eponse
similaire peut \^etre donn\'ee dans le cas singulier. Elle d\'ecoule
directement des r\'esultats de \cite {5}. Supposons donc  la courbe
$C$ singuli\`ere, et consid\'erons la d\'eformation universelle
\'equivariante, ici prise dans un sens analytique
$$\pi : {\Cal C} \longrightarrow D^{3g'-3+b},\quad (D = \{ z\in {\Bbb
C}, \vert z \vert < 1\}) \tag (7.45) $$
  avec par cons\'equent un isomorphisme (\'equivariant) $\pi^{-1}(0) =
C$. Si $t\in D$ est en dehors
du discriminant, la fibre ${\Cal C}_t$ est non-singuli\`ere, et
indique donc le type topologique
de la composante qui contient $C$.  En fait, et c'est le point que
nous allons retenir,  le type topologique peut \^etre lu directement
sur $C$. La raison est que le morphisme de monodromie associ\'e  \`a
 la fibre g\'en\'erale ${\Cal C}_t$, peut
\^etre construit directement \`a partir de donn\'ees lisibles sur
$C$. C'est essentiellement le
r\'esultat de (\cite {5}, thm 2.1), appel\'e  th\'eor\`eme de
Seifert-van Kampen.  Il est n\'ecessaire dans un tel \'enonc\'e de
remplacer le groupe fondamental ordinaire par le groupe fondamental
d'un graphe de groupes.

R\'esumons   la construction  de  (\cite {5},
\S  2.5). Soit   $\Gamma$ le graphe modulaire d\'efini par $C$, et
$\Delta = \Gamma/G$ le graphe quotient. Il y a sur $\Delta$  une structure de
graphe de groupes plus riche que celle utilis\'ee dans la section 7.2
 qui permet de reconstruire non seulement $\Gamma$ mais $C$  (loc.cit. D\'efinition  2.9). 
 
 Notons par $S_v,\, (v\in \Delta)$ les
composantes irr\'eductibles normalis\'ees de $D =   C/G$, et notons
$D_v \subset S_v$ les points qui ont pour images dans $D$ un point
double, donc les origines des branches. Notons aussi  $\beta_v$
l'ensemble des points marqu\'es port\'es par la composante $S_v$. Ces
points contiennent donc l'ensemble des points de ramification. On
note alors $\pi_v = \pi_1(S_v - (\beta_v \cup D_v) , *_v) $
le groupe fondamental de la courbe $S_v$ priv\'ee des points
exceptionnels, relativement \`a un
point de base $*_v$. La structure de graphe de groupes port\'ee par
$\Delta$ est d\'efinie ainsi:  tout d'abord, on associe \`a $v$ le
groupe $\pi_v$. Soit $(e,\overline e)$ une ar\^ete
  qui  pointe vers $v$, et $\overline e$ pointant
vers $w$. On pose alors $H_e =
H_{\overline e} = {\Bbb Z}$,  en fixant un g\'en\'erateur privil\'egi\'e $t_e$ tel que $t_e + t_{\overline e} = 0$.  On d\'efinit ensuite 
$$\partial_e^0(1) = \alpha (e)\quad,\quad \partial_e^1 (1) = \alpha
(\overline e)$$
expression dans laquelle  on note $\alpha (e)$ et $\alpha (\overline
e)$, des lacets de "pointe" autour des origines $p_v$ et $p_w$ des
branches d\'efinies par le point double $a$. La classe  de conjugaison  d'un tel lacet de pointe est bien d\'efinie, ce qui assure que le graphe de groupes est bien d\'efini \`a isomorphisme pr\`es.  Notons
$\Delta_\pi$ ce graphe de groupes  pour \'eviter toute confusion avec la d\'efinition de la section 7.2.
Rappelons    le r\'esultat   (\cite {5}, th\'eor\`eme  2.1),  de sp\'ecialisation du groupe fondamental:

\proclaim Proposition 7.17.  Le groupe fondamental $\pi_1
(\Delta_\pi)$ du graphe de groupes
$\Delta_\pi$ est isomorphe au groupe fondamental d'une fibre g\'en\'erique ${\Cal C}_t$ de la
d\'eformation universelle de $C$; en outre on peut choisir
l'isomorphisme de sorte qu'il pr\'eserve les classes de lacets autour
des piq\^ures. \hfill\break\qed

La preuve de la Proposition 7.17 pour laquelle on r\'ef\`ere  \`a  \cite {5},  fournit en outre un
algorithme pour expliciter un syst\`eme de g\'en\'erateurs canoniques
de "groupe de surface" pour le groupe fondamental. Pour lire
directement sur le rev\^etement $\pi: C \to D$, le nom de la
composante qui contient de point du bord, on voit qu'il suffit de
pouvoir reconstruire le morphisme de monodromie, au niveau de $D$,
donc sous la forme
$$\psi : \pi_1 (\Delta_\pi) \longrightarrow G \tag (7.45) $$
les images des lacets autour des piq\^ures \'etant soumises aux
contraintes impos\'ees par la donn\'ee de ramification (\S 2.4). En
r\'esum\'e:

\proclaim Proposition 7.18.  La correspondance   explicit\'ee
dans le lemme 2.1  reste
valable   en les points du bord de l'espace de Hurwitz. De
mani\`ere pr\'ecise,   le morphisme de monodromie (7.41)
permet de reconstruire le rev\^etement $\pi: C \to D$. La
classe double (voir 2.10)  d\'efinit    le type topologique, i.e.
le nom de la composante du sch\'ema de Hurwitz qui contient $\pi: C
\to D$.

\dem  La donn\'ee  de $\psi$  nous permet de  contruire  le rev\^etement $\pi: C\to D$.  En bref, la restriction de $\psi$ \`a $\pi_v$ donne un $G$-rev\^etement (non connexe peut \^etre) $C_v \to D$. Sur la courbe $C_v$ il y a une collection d'orbites exceptionnelles (les origines des branches) qui autorisent que l'on recolle les $C_v$ pour obtenir $C$. Noter que d'une autre mani\`ere  on r\'ecup\`ere le graphe de groupes  de la proposition 7.3, et donc $\Gamma$. On prend pour $G_v$ (resp. $H_a$)  l'image par $\psi$ de $\pi_v$ (resp. $H_a$).  utilise  alors le th\'eor\`eme 7.1
(th\'eor\`eme  de Bass \cite {6}).  Ensuite, c'est \`a quelques d\'etails pr\`es la r\'ep\'etition de la  construction utilis\'ee dans
la proposition 7.3.\hfill\break\qed

\bigskip

\section8 {Structures de niveau  sur les courbes stables }  

\bigskip

Dans cette section on revisite le champ compactifi\'e  des
courbes de genre $g$ (fix\'e) \'equip\'ees d'une structure de
niveau  $G$.    Le langage des champs de Hurwitz  permet donner
une d\'efinition  directe du champ    $\overline {\Cal M}_g(G)$ classifiant les courbes stables de genre $g$, munies d'une structure de niveau $G$.   En corollaire on obtient  une interpr\'etation modulaire claire de ses   points y compris au bord.    Le   champ $\overline {\Cal M}_g(G)$  bien que distinct du champ  $_G\overline {\Cal M}_g$ construit par  Deligne-Mumford \cite{18},\cite{19}  a  le m\^eme espace  des modules
grossiers.  Le premier qui est lisse est la d\'esingularisation du second.    Cette section   pr\'ecise   la section 7 de   \cite{61}.

 La d\'efinition classique d'une structure de niveau 
sur une courbe lisse, d\'efinie par le groupe fini $G$,  ne s'\'etend pas de
mani\`ere directe au cas singulier  (Brylinski \cite {14}, voir aussi \cite {20}, et
Oort-Van Geemen \cite {36}). Cela explique pourquoi la
compactification "naturelle" $_G\overline {\Cal M}_g$ de
Deligne-Mumford n'est pas d\'efinie comme l'espace des modules
grossiers attach\'e \`a un foncteur contravariant  naturel
d\'efini sur la cat\'egorie des courbes stables. Cela sugg\`ere
qu'il peut \^etre utile  d'attribuer  une interpr\'etation modulaire
 aux points du bord. C'est l'objectif de la pr\'esente section.  On
va observer que $\overline {\Cal M}_g(G)$ est  "essentiellement"
le champ des $G$-rev\^etements
principaux, \'eventuellement d\'eg\'en\'er\'es (stables), de base une
courbe stable. 

 Bien que r\'edig\'ee en des termes
assez diff\'erents notre approche  est  \'equivalente \`a celle sugg\'er\'ee
r\'ecemment par Abramovich, Corti et Vistoli \cite {1}, \cite {2}, \cite
{3}. Nous verrons aussi qu'elle permet de retrouver d'une mani\`ere
simple et directe, comme application des \S 7.1 et \S 7.2,  quelques
r\'esultats sur la structure du bord, essentiellement pour le niveau
ab\'elien $(n), \,n\geq 3$ (Boggi- Pikaart - De Jong \cite {13},
Oort-Van Geemen \cite {36}).

On fixe un groupe fini $G$ de cardinal $\vert G \vert$.  On note
pour abr\'eger ${\Cal H}_g(G)$ (resp. $\overline {\Cal H}_g(G))$, le
champ de Hurwitz classifiant les $G$-rev\^etements stables de base
une courbe lisse (resp. stable) de genre $g \geq 2$, de donn\'ee de
ramification  $\xi = \emptyset$, c'est \`a dire $\overline {\Cal H}_{g,G,\emptyset}$. Le cas $g=1$ peut \^etre inclus avec quelques modifications.

\bigskip

 \beginsection 8.1. Le champ des courbes stables avec structure de niveau $G$

\bigskip

  {\it 8.1.1. Structures de niveau sur les familles de courbes
lisses} 

\bigskip
Soit $C$ une courbe alg\'ebrique
projective et lisse  d\'efinie sur le corps alg\'ebriquement clos  $k$. On rappelle que $\vert G\vert \ne 0$ dans $k$. 
Une structure
de niveau $G$ sur $C$ est la donn\'ee d'une surjection ext\'erieure
$\phi: \pi_1(C) \longrightarrow G$ (\cite {14}, \cite {20}  D\'efinition 5.6).  On suppose que $G$  est quotient du groupe fondamental
$\pi_1(C)$ d'une telle courbe. On supposera en fait  un peu plus, \`a
savoir que $G$ est un quotient caract\'eristique (quotient par un
sous-groupe caract\'eristique). Le niveau $G'$ est dit dominer le
niveau $G$ si la surjection  $\phi: \pi_1(C) \longrightarrow G$
factorise par $G'$.  Le niveau ab\'elien $(n)$ correspond \`a $G
=  (  {\Bbb Z}/{n\Bbb Z})^{2g}$.  Une structure de
niveau   $(n)$ sur la courbe lisse $C$   revient  \`a sp\'ecifier  un
isomorphisme 
$$\alpha: \Pic (C) [n] \buildrel\sim\over\rightarrow
({\Bbb Z}/{n\Bbb Z})^{2g}\tag (8.1)$$
 ou si $k=\Bbb C, \,\,H^1(C, \Bbb
Z/n\Bbb Z) \buildrel\sim\over\rightarrow ({\Bbb Z}/{n\Bbb
Z})^{2g}$.     Un niveau  $G$ est {\it g\'eom\'etrique}  si notant $\pi_g$ le groupe
fondamental d'une surface de Riemann compacte de genre $g\geq 2$,
 le noyau d'une surjection $\phi: \pi_g \to G$ est ind\'ependant
de $\phi$. Les niveaux ab\'eliens sont g\'eom\'etriques, ainsi que
les niveaux di\'edraux \footnote {Le niveau di\'edral d'ordre $m \geq
2$, est d\'efini par $G =   {\pi_g}/{[\pi_g^{(2)},\pi_g^{(2)}]
}\pi_g^{(2m)}$, $[-,-]$ d\'esignant le sous-groupe des commutateurs,
et $\pi_g^{(k)}$ \'etant le sous-groupe engendr\'e par les
puissances d'ordre $k$ .}de Brylinski \cite {14} et Looijenga \cite
{52}.  On se limite dans la suite aux seuls niveaux
g\'eom\'etriques.

La d\'efinition d'une structure de niveau s'\'etend avec l'aide de la th\'eorie du groupe fondamental  de SGA 1, \cite
{34}) aux familles de courbes lisses (\cite {13}, \cite {14}, \cite
{20} \S 5). Soit $\pi: C\rightarrow S$ une courbe lisse au-dessus
de la base $S$, et de genre $g\geq 2$.  Le cas $g=1$ rentre dans ce
cadre, mais avec quelques pr\'ecautions suppl\'ementaires.  On se limite en cons\'equence \`a $g\geq 2$.    Soit une section 
  $s: S \longrightarrow C$   (une telle section
existe localement pour la topologie \'etale), alors on peut d\'efinir, $\Bbb L$ \'etant un ensemble de nombres premiers exclus,  un
groupe fondamental relatif  $\pi_1^{\Bbb L} (C/S,s)$  \cite{14}. Cela permet de d\'efinir   sur $S_{et}$ le faisceau localement
constant ${\Cal H}om^{ext} (\pi_1(C/S) , G)$ des
homomorphismes "ext\'erieurs". Alors (\cite {14} D\'efinition 2.3.1), \cite {20} D\'efinition 5.6), la d\'efinition pr\'ecise d'une structure de niveau
$G$ sur $C/S$ est la suivante:

\proclaim D\'efinition 8.1.  Une $G$-structure de niveau   sur une $S$-courbe lisse $C/S$, est la
donn\'ee d'une "surjection ext\'erieure" $\phi: \pi_1^{\Bbb L}(C/S)
\longrightarrow G$,  en des termes plus pr\'ecis, d'une section
globale   $\phi \in \Gamma (S , {\Cal
H}om^{ext} (\pi_1^{\Bbb L}(C/S) , G))$. m

Une structure de niveau\footnote {\rm 
Une structure de Teichm\"{u}ller de niveau $G$ dans \cite {14},\cite
{20}.} $\phi$ n'est donc que localement
(pour la topologie \'etale)  d\'efinie par un homomorphisme du
$\pi_1$ dans $G$.   La cat\'egorie des courbes lisses de genre $g\geq
2$, \'equip\'ees d'une structure de niveau $G$, d\'efinit un champ
alg\'ebrique $_G{\Cal M}_g$.    On a  le r\'esultat suivant (\cite {14}
Theorem 2.3.2, \cite {20} Lemma 5.7):

\medskip

\proclaim Th\'eor\`eme 8.2.  1)  Le champ $_G{\Cal M}_g$
est un champ alg\'ebrique de Deligne-Mumford.\hfill\break
  2)  Le morphisme $_G{\Cal M}_g \longrightarrow {\Cal M}_g$, oubli
de la structure de niveau, est fini, \'etale et  surjectif.\hfill\break
3)  Le nombre de composantes connexes de  $_G{\Cal M}_g$ est ind\'ependant
de $k$.\hfill\break
4)  Si le niveau $G$ domine le niveau $(n), n\geq 3$, alors
$_G{\Cal M}_g$ est repr\'esentable (un espace de modules fin).  
On a 
$${\Cal M}_g(G):= {\Cal H}_g(G)//Z(G) = _G{\Cal M}_g \tag (8.2)$$
 et le morphisme naturel ${\Cal M}_g (G) \to {\Cal M}_g$ est un $\Out
(G)$-torseur. \hfill\break
\qed

  La preuve de 3) n\'ecessite une
"compactification" de $_G{\Cal M}_g$(\cite {20}, cor 5.11).
 Faute d'une d\'efinition directe raisonable  du champ
compactifi\'e $  _G{\overline {\Cal M}}_g$, la proc\'edure
usuellement  retenue est de prendre pour d\'efinition du champ
compactifi\'e  $  _G{\overline {\Cal M}}_g$  la normalisation de $
\overline {\Cal M}_g [{1 \over{\vert G \vert}}]$ dans  $_G{\Cal M}_g$
(loc.cit. \S 5).  Seul le point 4) du th\'eor\`eme 8.2 demande une
v\'erification.  D\'ebutons par un
lemme \'el\'ementaire:

\proclaim Lemme 8.3.  1) Soient $\pi_i: \Sigma_i
\to C$ deux $G$-rev\^etements principaux de $C/S$, une $S$-courbe
lisse, avec $S$ connexe. Il existe $\theta \in \Aut (G)$, d'image
unique dans $\Out (G)$,  tel que le morphisme $Isom_{C,G}
(\Sigma_1^\theta , \Sigma_2) \rightarrow S$ soit un $Z(G)$-torseur.\hfill\break
2)  Si $\pi: \Sigma \to C$ est un $G$-torseur de base  $S$, $\Aut_C
(\Sigma) = G$, et $\Aut_{C,G} (\Sigma) = Z(G)$.

\dem Les notations, en particulier la d\'efinition du
foncteur $\Isom_{C,G} (\Sigma_1^\theta , \Sigma_2)$ sont celles  du th\'eor\`eme 6.20.  Des arguments identiques  conduisent au fait que pour un certain $\theta
\in \Aut (G)$, le morphisme 
$\Isom_{C,G} (\Sigma_1^\theta ,
\Sigma_2)\rightarrow  S  $ est fini, \'etale et surjectif.
Notons que  $Z(G)$ op\`ere naturellement sur $\Isom_{C,G}
(\Sigma_1^\theta , \Sigma_2)$. Comme cette action est simplement
transitive sur les fibres g\'eom\'etriques, on conclut que
$\Isom_{C,G} (\Sigma_1^\theta , \Sigma_2)/Z(G) = S$.\hfill\break
Le point 2) se traite de la m\^eme mani\`ere, en observant que la courbe
$\pi: \Sigma \to S$ \'etant lisse, $\Aut_C(\Sigma)$ est \'etale sur
$S$, et \'egal \`a $G$ sur les fibres g\'eom\'etriques, donc \'egal
\`a $G$.  Notons que cela prouve que dans 1), l'automorphisme
$\theta$ a une image bien d\'etermin\'ee dans $\Out(G)$,  d'o\`u le
lemme.\hfill\break\qed

Pour terminer la preuve du point 4) de (8.2),  on peut
invoquer la th\'eorie du groupe fondamental comme dans \cite
{13}, \cite {20}. On peut aussi, et de mani\`ere plus directe, noter
que si on forme le champ alg\'ebrique ${\Cal N}_{C/S} (G) = {\Cal H}_g(G)
\times_{{\Cal M}_g} S$ d\'efini par le morphisme $S\to {\Cal M}_g$, i.e.
$C/S$, le groupo\"\i de des  objets au dessus du $S$-sch\'ema $T$,
sont les $G$-fibr\'es principaux de base $C\times_S T$. Le Lemme 8.3,
dit que le champ ${\Cal N}_{C/S} (G)//Z(G)$ est repr\'esentable, et que
c'est un $\Out (G)$-torseur de base $S$. Noter qu'au sens strict, il
n'y a pas d'action de $\Out(G)$ sur  ${\Cal M}_g(G)$,   seulement
une action  faible (non stricte \S \ 6). Le lemme 8.3 dit que ce champ est en fait
isomorphe au champ form\'e des courbes de genre $g$, plus une
structure de niveau $G$, qui lui supporte une action  stricte de $\Out(G)$.\hfill\break
\qed 

 Le point de vue des sch\'emas de Hurwitz
(le point 4 de 8.2), conduit de  fait \`a renverser la D\'efinition
8.1, et  donc \`a consid\'erer une structure de niveau  comme d\'erivant d'un $G$-fibr\'e principal de base $C$. Cela sugg\`ere
qu'on doit  voir le champ  $  _G{\overline {\Cal M}}_g$ comme une
sorte de "compactification" de ${\bf B}G$, le champ classifiant du groupe
fini $G$ (on pourra comparer avec le point de vue voisin  mais plus sophistiqu\'e de
Abramovich-Corti-Vistoli \cite {2}).

\exa{ 8.1} 
Le niveau \'etant le ab\'elien $(n)$, donc $G =
({\Bbb Z} /{n\Bbb Z})^{2g}$, on notera le champ correspondant
${\Cal M}_g(n)$. La correspondance entre les deux d\'efinitions est
ais\'ee \`a expliciter. Si $\pi: \Sigma \longrightarrow C$ est un
$G$-rev\^etement principal, l'alg\`ebre $\pi_\star ({\Cal O}_\Sigma)$
se d\'ecompose canoniquement en une somme directe de faisceaux
inversibles (sous-faisceaux propres)
$ \pi_\star ({\Cal O}_\Sigma) = \bigoplus_{\chi \in \hat G} {\Cal
L}_\chi . $
Le choix d'une racine primitive $n$-i\`eme permet   d'identifier $G$
et $\hat G$;  par
suite  les faisceaux ${\Cal L}_\chi$ conduisent \`a une
identification de  $\Pic_{C/S} [n]$ et  $({\Bbb Z}/{n\Bbb
Z})^{2g}$. Un avantage de cette d\'efinition, est qu'elle garde un sens
comme on va le voir, si  on sp\'ecialise $C$ en un point du bord de
$\overline {\Cal M}_g$. Le nombre de composantes connexes de ${\Cal
M}_g(n)$ (le nombre de Nielsen), est  l'indice  
$[\GL_{2g} (\Bbb Z /n\Bbb Z) : \Sp_{g} ( \Bbb Z /n\Bbb Z)]$  \cite{20}. \ $\lozenge$

\bigskip

{\it 8.1.2. Structures de niveau sur les courbes stables }

\bigskip

Lorsque  $C$ est   singuli\`ere
(stable), le th\'eor\`eme 8.2 (4), sugg\`ere qu'une structure de
niveau $G$ sur $C$, doit pouvoir se d\'efinir en termes de $G$-fibr\'es "principaux"  de base $C$, mais  maintenant d\'eg\'en\'er\'es.  
\proclaim D\'efinition 8.4.  Soit $C$ une courbe stable de genre
$g\geq 2$ d\'efinie  sur $k$.
On d\'efinit un $G$- rev\^etement principal ou torseur,  stable (ou
d\'eg\'en\'er\'e) $\pi: \Sigma\longrightarrow C$, comme \'etant un
$G$-rev\^etement stable (dans le sens de la D\'efinition 6.7) , \`a
donn\'ee de Hurwitz  $\xi = \emptyset$.  Cela signifie   (D\'efinition 4.3)   que
$\Sigma$ est une courbe stable et   que l'action de $G$ sur
$\Sigma$  est stable, c'est \`a dire  libre en dehors  des
  points  doubles, et agissant "stablement" aux points doubles.
  Une $G$-pr\'estructure de niveau sur $C/S$,  est la donn\'ee d'une
classe d'\'equivalence (stricte) de $G$-rev\^etements principaux
stables de base $C$.

 Si d'une autre mani\`ere un $G$-fibr\'e principal de base $C$ lisse,
est interpr\'et\'e comme  un morphisme $C  \longrightarrow {\bf
B}G$, de $C$ dans le champ classifiant de $G$ \cite {1}, \cite {2},
on voit que cette d\'efinition diff\`ere de mani\`ere essentielle de
la D\'efinition 8.1 si $C$ d\'eg\'en\`ere.

  Soit $p: C \longrightarrow S$ une courbe stable de genre $g\geq
2$. Pour tout $S$-sch\'ema $T$,  soit $F_{C/S} (T) $
l'ensemble des $G$-{\it pr\'estructures de niveau} sur $C\times_ST$.
Le pr\'efaisceau $F_{C/S}$ n'est bien s\^ur pas en g\'en\'eral un
faisceau sur $S_{fppf}$,  m\^eme  si  le centre $Z(G)$ est non
trivial. Une des raisons est que le groupe des automorphismes d'un
$G$-torseur peut contenir strictement $Z(G)$. On d\'efinit  alors
$\Ni_{C/S,G}$, le faisceau  des $G$-structures de niveau  sur $p:
C\longrightarrow S$ comme \'etant le faisceau fppf associ\'e au
pr\'efaisceau $F_{C/S} $.  
 
 Na\"\i vement une
$G$-structure de niveau  sur $C/S$ est une section  $\alpha \in\Gamma (S ,
Ni_{C/S,G})$. Pour obtenir la bonne d\'efinition,  on  consid\`ere  
$$\overline {\Cal H}_g(G) := \overline {\Cal H}_{h,G,\emptyset} $$
    Si le niveau $G'$ domine le niveau $G$, signifiant que $G$
est un quotient de $G'$, il y a un morphisme naturel de champs
$ \overline {\Cal H}_g(G') \longrightarrow \overline {\Cal H}_g(G),$
et un morphisme analogue ${\overline M_g}(G') \longrightarrow
{\overline M_g}(G) $ entre sch\'emas modulaires. Le point  4) du Th\'eor\`eme 8.2 sugg\`ere que la compactification naturelle de $_G{\Cal
M}_g = {\Cal M}_g (G) = {\Cal H}_g(G)// Z(G)$ est 
$$ \overline{\Cal M}_g(G) = \overline {\Cal H}_g
(G)//Z(G) \tag (8.3)$$ 
Cela est justifi\'e par le r\'esultat
suivant:

\proclaim Th\'eor\`eme  8.5. Soit  ${\overline {_G{\Cal
M}}}_g$ la normalisation de $\overline {\Cal M}_g$ dans $_G{\Cal M}_g =
{\Cal M}_g (G) $ (compactification de Deligne-Mumford). Il existe un
morphisme   propre, birationnel 
 $$\overline{\Cal M}_g(G) := \overline {\Cal H}_g (G)//Z(G) \longrightarrow
_G\overline {\Cal M}_g \tag (8.4)$$
Les deux champs ont les m\^eme espaces grossiers de modules, en particulier  si le champ de droite est
repr\'esentable, par exemple si $G$ domine le niveau $(n), \,n\geq
3$,  c'est l'espace grossier de modules de $\overline{\Cal M}_g(G)$. 

\dem  Il est clair que le groupe des
automorphismes de tout objet de $\overline {\Cal H}_g (G)$ contient
$Z(G)$, ce qui permet de former le 2-quotient  $\overline {\Cal
H}_g (G)//Z(G) =\overline {\Cal M}_g(G)$  (voir \S \ 6.1). Ce
champ  est  de  Deligne-Mumford et  lisse (Proposition 6.5). Il
contient comme sous-champ ouvert dense   $_G{\Cal M}_g = {\Cal M}_g (G)
$. Il en r\'esulte, suite \`a la d\'efinition de la normalisation
(Deligne \cite {19}),  qu'il existe bien un morphisme naturel
$$\overline {\Cal H}_g (G)//Z(G) \longrightarrow _G\overline {\Cal
M}_g  \tag (8.5)$$
 \'etendant le morphisme $_G{\Cal M}_g \to {\Cal M}_g$. On peut pr\'ef\'erer un argument plus direct (\cite{61}, Theorem 7.2.3).
Ce morphisme est propre car les champs invoqu\'es le sont. 
Reste \`a prouver l'\'egalit\'e des espaces grossiers de modules.
On suppose que  $_G\overline {\Cal M}_g$ est repr\'esentable, par
exemple si le niveau $G$ domine $(n)$ avec $n\geq 3$ (Deligne \cite {19},
Proposition 3.5). Si $M$ (resp. $N$) sont les espaces grossiers de
modules respectifs,  alors  comme cons\'equence directe de la d\'efinition, on a  un morphisme naturel $M \to N =  _G\overline {\Cal M}_g$.
Ce morphisme est fini, et un isomorphisme sur un sous-sch\'ema
ouvert partout dense. Comme par construction $M$ est normal, ainsi
que $N$ (par d\'efinition), on a $M=N$. \hfill\break\qed

  \rema{8.2}   Le morphisme   $\overline
{\Cal M}_g(G) = {\overline {\Cal H}}_g(G) // Z(G) \longrightarrow
{\overline {_G{\Cal M}}}_g  $ est un isomorphisme sur les sous-champs ouverts form\'es par les structures de niveau sur les courbes
lisses de genre $g$. Noter que le champ de gauche est lisse, mais en
g\'en\'eral non repr\'esentable, et celui de droite est connu comme
repr\'esentable lorsque le niveau domine le niveau $(n), \,\, n\geq
3$ \cite {19}, mais non lisse en g\'en\'eral. Ces deux champs ont le
m\^eme espace modulaire grossier ${\overline M_g}(G)$ (fin pour
${\overline {_G{\Cal M}}}_g$ sous les conditions pr\'ec\'edentes).
  \hfill\break
$\lozenge$

Il est clair que $\overline
{\Cal M}_g(G) $ est d\'efini sur ${\Bbb Z}[{1\over \vert G\vert}]$.  Comme cons\'equence imm\'ediate de la construction, 
  chaque point poss\`ede  du champ  poss\`ede une interpr\'etation
modulaire pr\'ecise. Si $\pi: \Sigma \rightarrow C$ est un rev\^etement
principal (d\'eg\'en\'er\'e) de base la courbe
stable $C$, i.e. d\'efinissant une "structure de niveau" $G$ sur $C$,
le groupe des $G$-automorphismes du rev\^etement $\Aut_G (\Sigma)$
gouverne la structure locale de l'espace des modules. Sa d\'etermination explicite, en particulier savoir s'il contient
strictement ou non le centre de $G$, peut s'av\'erer difficile.  Le
fait que ce groupe puisse \^etre diff\'erent de $Z(G)$, est le fait
essentiel qui distingue le cas stable du cas lisse.  On peut \^etre
un peu plus pr\'ecis sur ce groupe, lorsque le niveau  $G$ domine le
niveau $(n), n\geq 3$,  comme le montre la proposition suivante,
cons\'equence facile du lemme de rigidit\'e  (par exemple:
Brylinski \cite {14},  Deligne \cite {19}, Oort \cite {36}) :
 
\proclaim Proposition 8.6. i)  Soit $\pi:\Sigma \longrightarrow C$ un
$G$-fibr\'e principal  de base la
$S$-courbe stable  $C$. Si $f \in \Aut_G (\Sigma)$, alors $f$ induit
l'automorphisme identit\'e de $C$, i.e. $\Aut_G (\Sigma) \subset
\Aut_C(\Sigma)$. \hfill\break
ii) Soit $C/S$ une courbe stable. Soit un niveau $G$ ($\vert G\vert \in {\Cal O}_S^*$). Alors $C/S$ a une $G$-struture de niveau (na\"\i ve) apr\`es extension fppf de $S$.

\dem On se ram\`ene de suite pour la premi\`ere assertion \`a $G =
({\Bbb Z}/{n{\Bbb Z}})^{2g}$.
 On consid\`ere  pour cela  une surjection $G \rightarrow   {\Bbb
Z}/{n{\Bbb Z}}$;  soit $H$ le
noyau. Le rev\^etement $\pi: \Sigma \rightarrow C$ factorise en
$\Sigma \rightarrow   \Sigma/H \rightarrow C.   $
Un $G$-automorphisme de $\Sigma$ induit un $  G/H$-automorphisme de la courbe
quotient $\Sigma \rightarrow   {\Sigma}/H$, ce qui justifie la r\'eduction annonc\'ee. Le
r\'esultat se teste sur les fibres g\'eom\'etriques, ce qui permet de
supposer aussi $S = \Spec \,(k)$,
$k$ corps alg\'ebriquement clos.
Le niveau \'etant maintenant le niveau $(n), \,n\geq 3,$ notons $h$
l'automorphisme de
$\pi: C\rightarrow S$ induit par le passage au quotient  de $f$; on a
donc $\pi f = h \pi$. En
particulier $f$ d\'efinit un automorphisme de ${\Cal O}_C$-alg\`ebre
$$ h^\star (\pi_\star ({\Cal O}_\Sigma)) \longrightarrow  \pi_\star
({\Cal O}_\Sigma) \tag (8.6)$$
commutant \`a l'action de $G$,  et donc un automorphisme qui fixe
chaque facteur
isotypique. En particulier, cela entra\^\i ne que $h$ induit
l'identit\'e sur $\Pic^0(C)[n]$, voir
Proposition 8.7 ci-dessous. Mais on sait  que sous cette condition le
lemme de rigidit\'e de Serre
entraine  $h = 1$ \cite {19}. Ainsi $ \Aut_G (\Sigma) \subset \Aut_C (\Sigma)$.

Prouvons le point ii). C'est local sur $S$. Soit un point g\'eom\'etrique $s\in S$. On montre d'abord que $C_s$ admet une $G$-structure de niveau. Soit une d\'eformation ${\Cal C}\to \Spec R$ de $C_s$ de base un anneau de valuation dicr\`ete complet  de corps r\'esiduel $k$, de corps des fractions $K$, \`a fibre g\'en\'erique lisse. On peut \'equiper apr\`es extension finie s\'eparable de $K$, et remplacement de $R$ par le normalis\'e,  la fibre g\'en\'erique $ C_\infty$ d'une $G$-structure de niveau $\Sigma_\infty \to  C_\infty$. Par r\'eduction stable on peut m\^eme supposer que $\Sigma_\infty$ se prolonge en une $G$-courbe stable $\Sigma \to \Spec R$. Alors il est clair que $\Sigma/G = {\Cal C}$ (voir \S 5.3).  Pour conclure, on reprend l'argument utilis\'e pour prouver l'existence locale d'une cl\^oture galoisienne (Th\'eor\`eme 6.21).  Soit une structure de niveau repr\'esent\'ee par un $G$-torseur stable $\Sigma_s\to C_s$.  La d\'eformation universelle de ce $G$-torseur existe sur une extension finie plate   $\Spec R^* \to \Spec R$  de la base de la d\'eformation universelle de $C_s$. Ce sont les points doubles de type NS qui sont responsables de la ramification de cette extension. Quitte \`a effectuer une extension \'etale de $S$ en $s$, on peut supposer que la d\'eformation universelle est d\'efinie sur $S$.  Il suffit alors d'effectuer le changement de base par $\Spec R^* \to \Spec R$ pour obtenir une structure de niveau (na\^\i ve).

Des arguments diff\'erents bas\'es sur le th\'eor\`eme de de Jong-Pickaart \cite{17},  permettent de prouver qu'une structure de niveau existe apr\`es une extension finie fid\`element plate (Romagny \cite{61}, Cor 7.2.4).\hfill\break\qed 

\bigskip

\beginsection 8.2 Le niveau ab\'elien ($n$) 

\bigskip

 {\it  8.2.1. Groupes de d\'ecomposition et d'inertie}

\bigskip

On fixe   pour niveau, le niveau ab\'elien $(n), n\geq 3$, 
donc $G =   ( {\Bbb Z}/{n{\Bbb Z}})^{2g}$. Les r\'esultats de la section 7.3 permettent de
d\'ecrire la structure combinatoire du bord. Pour cela, nous
reprenons les suites exactes (7.22) et (7.23) de la section 7.3, et
les traduisons dans le cas d'un $G$-fibr\'e principal. Soit  $\pi:
\Sigma \longrightarrow C =  \Sigma /G$
une structure de niveau $(n)$ sur la courbe stable $C$, d\'efinie sur
le corps alg\'ebriquement
clos $k$. La suite exacte (7.22), donne sous les pr\'esentes hypoth\`eses, et en
notant $\tilde C = \coprod_i \tilde C_i$ la normalisation de $C$, et
$\Gamma$ le graphe dual:
$$1\longrightarrow H^1(\Gamma)\otimes    {{\Bbb Z}\over {n\Bbb Z}}
\longrightarrow \Pic (C)[n]
\longrightarrow \prod_i \Pic (\tilde C_i)[n] \longrightarrow 1 \tag (8.7)$$
En particulier cela donne   
$\Pic (C)[n]
\buildrel\sim\over\longrightarrow  H^1_{et} (C ,  { {\Bbb Z}
\over {n\Bbb Z} })  \buildrel\sim\over\longrightarrow  ({{\Bbb Z}
\over{n\Bbb Z}})^{2g - h(\Gamma)}.  $ 
Notons que si $\Delta$ est le
graphe modulaire associ\'e \`a $\Sigma$, on a $\Gamma = \Delta
/G$. On peut peut rendre plus explicite les groupes de d\'ecomposition et d'inertie (D\'efinition 7.4); pour les notations
relatives \`a ces groupes, voir la section 7.2:
 
\proclaim Proposition 8.7.    Pour les groupes introduits ci-dessus, on  a:  \hfill\break
 i)  $  G/D \buildrel\sim\over\rightarrow ({\Bbb Z}/{n\Bbb Z}
)^{h(\Gamma)}, \quad   G/I  \buildrel\sim\over\rightarrow  (
{\Bbb Z}/{n\Bbb Z} )^{2g - h(\Gamma)}  \,
\buildrel\sim\over\rightarrow  \,\Pic\,(C)[n]$\hfill\break
ii) On a pour toute composante irr\'eductible $C_i$ de $C$,
  $I_i = I \cap G_i,\quad   {G_i}/{I_i} \buildrel\sim\over\rightarrow
({\Bbb Z} /{n\Bbb Z} )^{2g_i}$  et
  $D \buildrel\sim\over\rightarrow  ({\Bbb Z} /{n\Bbb
Z})^{2g - h(\Gamma)}$\hfill\break
iii) Si $C_i$ correspond \`a un sommet de valence $v_i$ de $\Gamma$, on a
$I_i \buildrel\sim\over\rightarrow  ({\Bbb Z}/{n\Bbb Z})^{v_i
- 1}, \quad G_i \buildrel\sim\over\rightarrow  ({\Bbb Z}/{n\Bbb
Z} )^{2g_i + v_i - 1}. $

\dem Consid\'erons le diagramme suivant \`a lignes exactes, et avec
les fl\`eches verticales
injectives
 $$\harrowlength=38pt  \varrowlength=38pt
 \commdiag {1 &\mapright&  H^1(\Gamma)\otimes  {\Bbb Z}/{n\Bbb
Z} &\mapright &\Pic (C)[n]  &\mapright & \prod_i \, \Pic (\tilde C_i)[n] &\mapright & 1\cr
 & & \mapdown &&  \mapdown &&\mapdown &&\cr
  1 &\mapright& \widehat {G/D}&  \mapright & \widehat {G/I}&
 \mapright &  \prod_I \,\widehat  {G_i/I_i} && \cr} \tag  (8.8)$$
  
L'exactitude de la suite horizontale du bas est le contenu du Th\'eor\`eme 7.2. On montre d'abord
que la fl\`eche verticale du milieu est une bijection. Cela se ram\`ene \`a prouver que si $L \in \Pic
(C) [n]$, alors $\pi^* (L) \cong {\Cal O}_\Sigma $ (voir la preuve du
Th\'eor\`eme 7.2). Supposons
$L$ d'ordre $d, \,\,d/n$, et consid\'erons le rev\^etement \'etale
connexe de degr\'e $d$, $\tau: C^*
\longrightarrow C$ d\'efini par $L$, c'est \`a dire tel que
$$ \tau_\star ({\Cal O}_{C^*}) = \bigoplus_{i=0}^{d-1} L^{i}\tag (8.9)$$
La courbe $C^*$ est connexe du fait que $L$ est d'ordre exact $d$, et
par ailleurs certainement stable. Notons que par construction $\tau^*
(L) \cong {\Cal O}_{C^*}$. Consid\'erons une d\'eformation du
$G$-rev\^etement $\pi: \Sigma \longrightarrow C$ de base $\Spec (R)$,
o\`u $R$ est un anneau de valuation discr\`ete  complet de corps
r\'esiduel $k$, \`a fibre g\'en\'erique lisse, donc d\'ecrite par
un diagramme
  
$$\commdiag{  Y &  \mapright &  X &  \mapright & \Spec \,(R) \cr
\cup & & \cup & & \cr 
\Sigma &  \mapright & C & &  \cr} $$
 
Du fait que $C^* \longrightarrow C$ est \'etale, on sait qu'il existe
une d\'eformation ${\Cal X}^*
\longrightarrow {\Cal X}$ de ce rev\^etement de base une quelconque
base prescrite, en particulier
${\Cal X}$;  on peut  se r\'ef\'erer par exemple au Th\'eor\`eme 5.2.
Notons $\eta$ le point g\'en\'erique de $\Spec\, (R)$, et $ Y_\eta$
(resp.
$ X_\eta$) les fibres g\'en\'eriques respectives.  Alors le rev\^etement
induit au dessus de $\eta$,
$ Y_\eta \longrightarrow  X_\eta$ factorise par $ X_\eta^*$. Cela
force $\pi:{\Cal
Y} \longrightarrow  {\Cal X}$ \`a factoriser par $ X^*$, et au bout
du compte $\pi: \Sigma
\longrightarrow C$ factorise par $C^*$, ce qui entraine bien $\pi^*
(L) \cong {\Cal O}_\Sigma$.

Prouvons maintenant les points 1) \`a 3).   Partant de
de l'\'egalit\'e  $\widehat{  G/I} \, = \, \Pic
(C)[n]$, et des identifications connues  
$$H^1(\Gamma)\otimes { 
{\Bbb Z}\over {n{\Bbb Z}}} \, \cong  (  {\Bbb Z}/{n{\Bbb
Z}})^{h(\Gamma)}, \,\,\hbox{ \rm et}\,\,\Pic (\tilde{C_i})[n] \cong
({\Bbb Z}/{n{\Bbb Z}})^{2g_i}$$
 le diagramme ci-dessus montre
que les applications verticales sont des bijections,  donc d'une part
$\widehat {  G/D} \,\cong \, (\frac {\Bbb Z}/{n{\Bbb
Z}})^{h(\Gamma)}$, et pour tout indice $i$, $\widehat{{G_i}/{I_i}} =
\Pic (\tilde {C_i})[n]$.   Par rapport \`a la suite exacte (7.23)
on a maintenant la suite exacte
$$1 \longrightarrow \widehat { G/D} \longrightarrow \widehat {  G/I}
\longrightarrow \prod_i \widehat {  {G_i}/{I_i}}\longrightarrow 1 \tag (8.10) $$
Par dualit\'e, on obtient une injection $\prod_i   {G_i}/{I_i}
\mapsto   G/I$, et
donc pour tout indice $i$, $G_i \cap I = I_i$. Comme $  {G_i}/{I_i}
\cong ({\Bbb Z}/{n{\Bbb Z}})^{2g-h(\Gamma)}$, on a  finalement $I \cong
({\Bbb Z}/{n{\Bbb
Z}})^{h(\Gamma)}$, et pour une raison identique, $D \cong  ({\Bbb Z}/{n{\Bbb
Z}})^{h(\Gamma)}$.

Notons maintenant $v_i$ la valence du sommet d'indice $i$ du graphe
$\Gamma$. Le groupe $I_i$ est engendr\'e par  $v_i$ \'el\'ements
$\{ \sigma_{i,\alpha} \}, (\alpha = 1,\dots,v_i)$, dont l'ordre
divise $n$, et soumis \`a une relation $\sum_\alpha
\sigma_{i,\alpha} = 0$. Dans le groupe  $I$ qui est le sous-groupe
engendr\'e par les $I_i$, un d\'ecompte naif des g\'en\'erateurs
conduit du fait que chaque $\sigma_{i,\alpha}$ appara\^\i t deux
fois, via une ar\^ete orient\'ee et l'ar\^ete oppos\'ee, \`a une
majoration
$$\sum_i v_i - A - (S-1) = h(\Gamma) \tag (8.11)$$
Les $S$ relations $\sum_\alpha \sigma_{i,\alpha}$ sont d\'ependantes,
car li\'ees par
$\sum_i ( \sum_\alpha \sigma_{i,\alpha} ) = 0,  $
ce qui explique la contribution $S-1$. Par ailleurs on sait que $I
\cong   ({\Bbb Z}/{n{\Bbb
Z}})^{h(\Gamma)}$, on voit donc ainsi qu'entre les g\'en\'erateurs
indiqu\'es, il ne peut y avoir
d'autres relations que les relations impos\'ees par le graphe (les
usuelles relations de courant).
Cela conduit aux \'egalit\'es
$I_{i,\alpha} = {{\Bbb Z}/{n{\Bbb Z}}}, \quad I_i = (
{\Bbb Z}/{n{\Bbb
Z}})^{v_i-1},  $ 
donc finalement $G_i \cong ({{\Bbb Z}/{n{\Bbb Z}}})^{2g_i+v_i-1}$
comme indiqu\'e.\hfill\break
\qed
\medskip
  {\it 8.2.2. Composantes irr\'eductibles du bord} 
  \bigskip
La description des composantes irr\'eductibles du bord de $\overline
M_g (n)$ est  
ais\'ee. La description g\'en\'erale donn\'ee dans les sections  7.2
et 7.4 montre que  les
composantes irr\'eductibles correspondent au choix d'un graphe
modulaire $\Gamma$, qui est soit
un segment, soit une boucle.

    {\bf segment:} \ \   
  dans ce cas, l'image de la composante dans $\overline M_g$ est
isomorphe au produit  $\overline M_{g_1} \times \overline M_{g_2}, \,\,(g_1 + g_2 = g).   $
 On a visiblement pour les groupes $I$
et $D$  
 $$I = 1, \,\,D = G =
\left ({\Bbb Z}/{n{\Bbb Z}} \right )^{2g} \tag (8.12)$$   De
mani\`ere plus pr\'ecise, le groupe $G$ se d\'ecompose en $G = G_1
\times G_2, \quad G_i
\cong ({\Bbb Z}/{n{\Bbb Z}})^{2g_i}, (i=1,2)$. Il en r\'esulte
en particulier que la composante
irr\'eductible correspondante s'identifie \`a $\overline M_{{g_1},1}
(n) \times \overline
{M_{{g_2},1}}(n)$, espaces modulaires de structure de niveau $(n)$
avec un point
marqu\'e. 

  {\bf boucle:}\ \  
 on a  (\S 7.2):  $$I =  {\Bbb Z}/{n{\Bbb Z}},\, \,
\hfill\break{ \rm et }\,\,D = ({\Bbb Z}/{n{\Bbb Z}})^{2g-1} \tag (8.13)$$
Avec les notations de la section 7.2 , l'\'el\'ement $g_0$ est
d'ordre $n$, et $G = D \times \langle g_0 \rangle $. La composante
correspondante du bord
s'identifie avec l'espace modulaire classifiant les structures de
niveau $(n)$ sur une courbe de
genre $g-1$ avec une  seule piq\^ure, c'est \`a dire
essentiellement \`a un espace de Hurwitz d\'efini par le groupe $G =
({\Bbb Z}/{n{\Bbb Z}})^{2g}$, et une donn\'ee de Hurwitz r\'eduite \`a un \'el\'ement d'ordre $n$. L'image d'une telle composante
dans $\overline M_g$ est la strate  not\'ee $\Delta_0 = \overline
{M}_{{g_1},2}$.

On peut \'etendre la description pr\'ec\'edente des composantes de
codimension un du bord \`a un groupe de niveau arbitraire, c'est \`a
 dire  au champ $\overline {\Cal M}_g (G)$.  Cela permet par
exemple de retrouver de mani\`ere naturelle  quelques r\'esultats de
Boggi-Pikaart (\cite {13}, \S 2).

Le niveau \'etant toujours le niveau ab\'elien $(n), n\geq 3$, on
obtient facilement en corollaire des m\'ethodes pr\'ec\'edentes, la
structure  connue des sous-groupes d'inertie pour les points doubles
\cite {36}.  On conserve dans le corollaire suivant les hypoth\`eses
de la Proposition 8.2:
 
\proclaim Corollaire 8.8.  Soit $Q_\alpha$ un point double de $C$,
et soit $H_\alpha$ le
stabilisateur d'un quelconque point double de $\Sigma$ au dessus de
$Q_\alpha$. Alors $H_\alpha = 1$ si et seulement si $Q_\alpha$
disconnecte  $C$, sinon  $H_\alpha =   {\Bbb Z}/{n{\Bbb
Z}}$.

\dem  Notons tout d'abord que le resultat est clair si $C$ n'a qu'un
seul point double, comme il r\'esulte de la preuve pr\'ec\'edente. On
va ramener le cas g\'en\'eral \`a ce cas particulier par un
argument de d\'eformation.    Supposons d'abord que le point double
$Q = Q_\alpha$ ne disconnecte pas $C$.  Dans ce cas la normalisation
partielle $C^*$ de $C$ en $Q$ est une courbe stable de genre $g-1$
marqu\'ee par deux points $Q'$ et $Q''$, et $C$ se d\'eduit de $C^*$
par l'identification  $Q' = Q''$. Cette construction, le "clutching"
morphism de Knudsen \cite
{48}, s'\'etend  \`a la d\'eformation universelle de
$(C^*,Q',Q'')$, et conduit \`a une d\'eformation  
${\Cal C}^*
\longrightarrow \Spec\, R^* $
 de $C$ dans laquelle le point double
$Q$ s'\'etend. Comme le clutching morphism est un morphisme repr\'esentable, fini et non ramifi\'e (loc.cit. cor 3.9), on voit que
cette d\'eformation n'est pas autre chose que la d\'eformation
universelle de $C$ qui pr\'eserve le point double $Q$. Elle peut se
d\'ecrire comme \'etant la restriction de la d\'eformation
universelle  ${\Cal C} \longrightarrow \Spec\, R$ de $C$ \`a
l'hypersurface $t_1 = 0$, si $t_1$ d\'esigne le param\`etre de d\'eformation de $Q$. On peut donc conclure que la fibre g\'en\'erique
de cette d\'eformation est une courbe avec un unique point double $Q$
qui reste du type "boucle".

Le niveau \'etant toujours le niveau ab\'elien $(n), n\geq 3$,
l'interpr\'etation modulaire  des points  du bord  permet de d\'ecrire  en des termes simples la ramification du morphisme $\overline {M_g} (n) \longrightarrow \overline M_g$ en d'autres
termes sa structure logarithmique \cite {55}. Fixons $C$ un point de
$\overline M_g $ sur le corps alg\'ebriquement clos $k$, $n$ \'etant
comme toujours inversible dans $k$. Soit $\Ni_G(C)$ l'ensemble des
structures de niveau "na\"\i ves" $G = ({\Bbb Z}/{n{\Bbb
Z}})^{2g}$ sur $C$ (D\'efinition 8.5).  Il y a une action \'evidente
de $\Out (G) = \GL_{2g} ({\Bbb Z}/{n{\Bbb Z}})$ sur   $\Ni_G
(C)$; alors:
\proclaim Lemme 8.9.  L'action de  $\Out (G) = GL_{2g} ({\Bbb
Z}/{n{\Bbb Z}})$ sur $\Ni_G
(C)$ est transitive.

\dem On suppose d'abord $S = \Spec\, k$ avec $k$ alg\'ebriquement
clos. Dans le cas
lisse  le r\'esultat est clair, du fait  qu'une structure de niveau
na\"\i ve est identique \`a une vraie structure de niveau.  Si
maintenant $C$ est singuli\`ere, le corollaire pr\'ec\'edant montre
que les points doubles de $C$ qui sont les  images des points doubles
de $\Sigma$ avec une isotropie non triviale, ne d\'ependent en fait
que de $C$, plus pr\'ecis\'ement du graphe  $\Gamma$, et pas de la
structure de niveau repr\'esent\'ee par $\pi: \Sigma \longrightarrow
C$.  En ces points l'indice de ramification est \'egal \`a $n$.
Soit maintenant $X$ la d\'eformation universelle \'equivariante  de $\Sigma$,
de base
 $$\Spec\, k[[t_1,\dots,t_r,\dots,t_{3g-3}]]$$  $t_\alpha$ \'etant le param\`etre de d\'eformation de l'orbite de points doubles
au dessus du point  $Q_\alpha$ de $C$. La courbe quotient $Y =   X/G$
se d\'eduit de la d\'eformation universelle ${\Cal C}$ de $C$,  ayant
pour  base   $\Spec\, k[[\tau_1,\dots,\tau_r,\dots,\tau_{3g-3}]]$,
par le changement de base donn\'e par $\tau_\alpha =
t_\alpha^{e_\alpha}$ (Th\'eor\`eme 5.2). Ce qui vient d'\^etre dit
montre que $Y$, \`a isomorphisme de d\'eformation pr\`es, ne d\'epend pas du choix de la structure de niveau.

Soient maintenant deux structures de niveau  $(n)$, $\pi_i: \Sigma_i
\longrightarrow C, \,(i=1,2)$
sur la courbe $C$;   soient aussi $X_i, (i=1,2),$ les d\'eformations
universelles \'equivariantes
respectives. On peut ainsi supposer que les d\'eformations $
{X_1}/G$ et $  {X_2}/G$ de $C$ sont
isomorphes,  donc que  ${X_1}/G\, \cong \,  {X_2}/G \,\cong Y $ est
une d\'eformation donn\'ee de $C$ de base  $ S = \Spec\,
k[[\tau_1,\dots,\tau_r,\dots,\tau_{3g-3}]]$.  On peut en outre, par
une sp\'ecialisation convenable,   la base \'etant ramen\'ee \`a $S
= \Spec (R)$, avec $R$ anneau de valuation discr\`ete complet de
corps r\'esiduel $k$,  supposer que les fibres g\'en\'eriques de $X_i
, (i=1,2)$ et de $Y$ sont  lisses et  d\'efinissent des structures de
niveau situ\'ees dans une m\^eme orbite de $\GL_{2g} ({\Bbb Z}/
{n{\Bbb Z}})$.  De mani\`ere \'equivalente, il existe $\theta \in
\GL_{2g} ({\Bbb Z}/{n{\Bbb Z}})$ tel que le sch\'ema $\Isom
(X_1^\theta , X_2)$ classifiant les isomorphismes \'equivariants  au
dessus de $Y$, soit fini, non ramifi\'e et surjectif.  Il y a donc
au moins une composante connexe de $\Isom_\pi (X_1^\theta,X_2)$ qui
domine $S$ et qui alors \'etale sur $S$, donc isomorphe \`a $S$.
On en d\'eduit un $Y$-isomorphisme $X_1^\theta \buildrel \sim \over
\longrightarrow X_2$, d'o\`u finalement un isomorphisme au dessus
de $C$,
$\Sigma_1^\theta\buildrel\sim\over \longrightarrow \Sigma_2$, d'o\`u
 la conclusion.\hfill\break \qed

Il a \'et\'e observ\'e que le morphisme ${\Cal M}_g(G)// Z(G)
\rightarrow {\Cal M}_g$ est un $\Out (G)$-torseur. Il n'en n'est plus
de m\^eme pour  le morphisme $\overline {\Cal M}_g(G)// Z(G) \rightarrow
\overline {\Cal M}_g$.  Au niveau des espaces grossiers de modules,
le morphisme $\overline M_g (n) \longrightarrow \overline M_g$ est un
rev\^etement galoisien de
groupe de Galois $\GL_{2g} ({\Bbb Z}/{n{\Bbb Z}})$. Pour pr\'eciser cela, soit un point  $\xi$ de  $\overline M_g (n)$ sur le corps
$k$, repr\'esent\'e par $\pi: \Sigma \longrightarrow C$. Le groupe
d'inertie  $I_\xi$ dans ce rev\^etement a la description suivante:
$$I_\xi = \{ \theta \in GL_{2g} ( {\Bbb Z}/{n{\Bbb Z}}) \, , \,
\Sigma^\theta \cong
\Sigma \} \tag (8.14)$$
Un \'el\'ement de $I_\xi$  d\'etermin\'e par $f: \Sigma
\buildrel\sim \over\longrightarrow
\Sigma^\theta$ d\'efinit donc un diagramme commutatif
 $$\commdiag{ \Sigma &\mapright^f& \Sigma^\theta \cr
   \mapdown\lft{\pi}&&  \mapdown\lft{\pi}\cr
    C &\mapright^h  &C \cr}   $$
dans lequel l'automorphisme horizontal du bas $h: C\buildrel
\sim\over \longrightarrow C$ est celui induit par $f$ par passage au
quotient. Notons $\widehat I_\xi$ le sous-groupe form\'e des couples
$(f,\theta)$ comme ci-dessus. On a la suite exacte
$$1 \longrightarrow \Aut_G (\Sigma) \longrightarrow \widehat I_\xi
\longrightarrow I_\xi
\longrightarrow 1 \tag (8.15) $$
Consid\'erons aussi le groupe   $\Aut_\pi (\Sigma)$ des
$C$-automorphismes (non \'equivariants \`a priori) de $\Sigma$ qui
commutent avec $\pi$, et qui pr\'eservent globalement $G$.  Du fait
que $n\geq 3$, le lemme de rigidit\'e implique que
$\Aut_G (\Sigma) \subset \Aut_\pi (\Sigma)$. L'application
$(f,\theta) \mapsto h$ est alors bien
d\'efinie; il en r\'esulte une  suite exacte
$$1\longrightarrow {\Aut_\pi (\Sigma)\over \Aut_G (\Sigma)}
\longrightarrow I_\xi \longrightarrow
\Aut (C) \longrightarrow 1 \tag (8.16) $$
Noter que la surjectivit\'e \`a droite d\'ecoule du Lemme 8.9. Si
$C$ est lisse, on a certainement
$\Aut_\pi (\Sigma) = \Aut_G (\Sigma) = G$, et $I_\xi = \Aut (C)$.
Dans le cas g\'en\'eral, on  a le
r\'esultat suivant:

\proclaim Proposition 8.10. Soit $m$ le nombre de points doubles de
$C$ (ar\^etes de $\Gamma$) qui ne  disconnectent pas $\Gamma$, alors
 $${\Aut_\pi (\Sigma)\over G} \buildrel \sim\over \rightarrow
\left({\Bbb Z}/{n{\Bbb Z}}
\right)^m \tag (8.17) $$

\dem Soit $f \in \Aut_\pi (\Sigma)$, il existe  $\theta \in \Aut (G)$
qui rende l'isomorphisme $f: \Sigma \buildrel\sim\over \longrightarrow \Sigma^\theta$ $G$-
\'equivariant. Si $X$ est une   
d\'eformation infinit\'esimale \'equivariante de $\Sigma$, et si
$\imath: C \hookrightarrow X$ est le
plongement \'equivariant correspondant, le couple $(X , \imath
f^{-1})$ d\'efinit donc une
d\'eformation not\'ee $X_f$ de $\Sigma^\theta$. Si ${\Cal X}
\longrightarrow  \Spec\,(S)$ avec $S =
k[[t_1,\dots,t_{3g-3}]]$ est la d\'eformation universelle  \'equivariante de $\Sigma$, alors
${\Cal X}^\theta \longrightarrow \Spec\,(S)$ est celle de
$\Sigma^\theta$. Il vient alors de $f$, un
automorphisme $\hat f$  de $\Spec\,(S)$, qui conduit au diagramme cart\'esien
  $$\commdiag{ {\Cal X}_f  & \mapright^f&  {\Cal X}^\theta  \cr
   \mapdown&&\mapdown\cr
  \Spec\,(S) & \mapright^{\hat f}&    \Spec\,(S)\cr}  \tag (8.18)$$
o\`u $h$ induit $f$ sur les fibres au dessus de l'origine $0\in S$.
Comme $\pi f = \pi$, les
d\'eformations de $C$ obtenues par passage au quotient par $g$ sont
\'equivalentes, ce qui force
$\hat f$ a \^etre l'identit\'e sur $\Spec \,(R)$ base de la d\'eformation universelle de $C$. Ainsi $\hat
f$ appartient au groupe de galois du rev\^etement $\Spec\,(S)
\longrightarrow \Spec\,(R)$, groupe  identifi\'e \`a $\left({\Bbb Z}/n{\Bbb Z}\right) ^m$. Noter que si $\hat f = 1$, alors par
examen au point g\'en\'erique, on obtient en cons\'equence que $f$
est une translation par un \'el\'ement de $G$.

R\'eciproquement partons d'un \'el\'ement $\hat f$ du groupe de
Galois  $\left({\Bbb Z}/n{\Bbb
Z}\right) ^m$, et notons ${\Cal X}_{\hat f}$ la d\'eformation d\'eduite du carr\'e cart\'esien
  $$\commdiag{ {\Cal X}_f & \mapright^h&  {\Cal X}  \cr
   \mapdown &&\mapdown\cr
  \Spec\,(S) &  \mapright^{\hat f}  & \Spec\,(S) \cr} $$
Par passage au quotient par $G$, on obtient un isomorphisme
$  {{\Cal X}_{\hat f} }/G \buildrel\sim\over\longrightarrow   {\Cal
X}/G.  $
  Par le Lemme 8.9 on conclut qu'il y a au dessus de $\Spec\,(S)$ un
isomorphisme \'equivariant
${\Cal X}_{\hat f} \buildrel\sim\over\longrightarrow {\Cal X}^\theta$
  pour un $\theta$ convenable. Le diagramme cart\'esien qui en d\'ecoule, analogue \`a (8.18),
  assure que tout \'el\'ement du groupe de galois est dans l'image de
$\Aut_\pi\,(\Sigma)$.
\hfill\break
\qed 

En g\'en\'eral pour un rev\^etement principal d\'eg\'en\'er\'e, et
pour un niveau arbitraire (dominant $(n), \,\,n\geq 3$) l'une des
deux  inclusions
$$G \subset \Aut_G (\Sigma) \subset \Aut_\pi \,(\Sigma) \tag (8.19)$$
peut  \^etre stricte, eventuellement les deux. Notons que l'anneau
local compl\'et\'e de ${\overline M_g}(G)$ au point d\'efini par
$\Sigma$ est
$ \hat{\Cal O}_\Sigma \cong k[[\tau_1,\dots,\tau_{3g-3}]]^{\Aut_G
(\Sigma)}.  $
La singularit\'e eventuelle de cet anneau local provient donc du
groupe $ {\Aut_G (\Sigma)}/G$ et de son action sur l'espace tangent
$H^1_G (\Sigma , \Theta_\Sigma)$.
  Comme application  donnons une
preuve diff\'erente du r\'esultat connu suivant (Oort-Van Geemen
\cite {36}):
\proclaim Proposition 8.11.  Soit un rev\^etement principal d\'eg\'en\'er\'e $\pi: \Sigma \rightarrow C$
d\'efinissant une structure de niveau $(n), \,n\geq 3$. On suppose
que   $C$ a deux composantes irr\'eductibles se coupant
transversalement en deux points, alors on a
$$ {\Aut_\pi\, (\Sigma)\over G}\,\, \buildrel\sim\over\rightarrow
\,\,{\Bbb Z}/{n{\Bbb
Z}} \tag (8.20)$$

\dem  On suppose donc $C = C' \cup C''$, avec $C'$ (resp. $C''$)
lisse de genre $g'\geq 1$ (resp $g''\geq 1$), et $g = g'+g''+1$. Soit
$C' \cap C'' = \{Q'\,,\,Q''\}$. Comme $h(\Gamma) = 1$, les
stabilisateurs des points doubles de $\Sigma$ , donc au dessus
de $Q'$ ou $Q''$, sont \'egaux  disons \`a $H \cong {{\Bbb Z}
/{n{\Bbb Z}}} $. Soit $\Sigma'$,
resp. $\Sigma''$ une composante de $\Gamma$ qui rel\`eve $C'$, resp.
$C''$, et supposons que $P' \in \Sigma' \cap \Sigma''$ soit au dessus
de $Q'$. Si $f \in \Aut_G (\Sigma)$, on a $f(P') \in GP'$, on
peut donc se ramener ne modifiant $f$ par un \'el\'ement de $g\in G$
\`a $f(P') = P'$. Alors $f(\Sigma') =
\Sigma'$ et $f(\Sigma'') = \Sigma''$. La restriction de $f$ \`a
$\Sigma'$ commute \`a l'action de $G'$,
donc est d\'efinie par un \'el\'ement de $G'$, du fait que $\Sigma'$
est lisse.  Par une nouvelle
r\'eduction on peut supposer de plus  que $f\vert \Sigma' = id$.
Alors on doit avoir
$f\vert \Sigma'' \in H \cong { {\Bbb Z}/{n{\Bbb Z}}}. $
Il est facile de voir qu'inversement un automorphisme  $f$ de
$\Sigma$  d\'efini par ces deux conditions centralise l'action de
$G$. Le r\'esultat en d\'ecoule.

Si maintenant $P'' \in \Sigma''$ est un point au dessus de $Q''$, il
est clair que les actions de
$f$ en les espaces tangents en ces deux points sont  d\'efinies par
des racines $n$-i\`eme de l'unit\'e
 oppos\'ees. Donc l'action de  $ {\Aut_G (\Sigma)}/G$ sur l'espace
de la d\'eformation \'equivariante universelle de $\Sigma$ est
$\tau_1 \mapsto \epsilon\tau_1, \,\tau_2 \mapsto \epsilon^{-1}
\tau_2, \,\tau_\alpha \mapsto
\tau_\alpha \quad (\alpha \geq 3)$
D'o\`u vient le fait qu'en ce point ${\overline M_g} (n)$ a une
singularit\'e  $A_{n-1} \times $(vari\'et\'e lisse de codim 2).
\hfill\break
\qed

Le calcul explicite du groupe $\Aut_G (\Sigma)$, g\'en\'eralisant le
calcul de la proposition 8.11,  est
possible avec un niveau g\'en\'eral, et  sous une condition de g\'en\'ericit\'e de $C$. La condition requise est celle qui pour $f\in
\Aut_G(\Sigma)$ assure que l'automorphisme induit $h$ de la courbe
$C$, point\'ee par les points de branchement, est \'egal \`a
l'identit\'e. De mani\`ere g\'en\'erale, apr\`es avoir fait le
choix d'une orientation du graphe  $\Gamma$, on a le r\'esultat:
\proclaim Proposition 8.12.  Soit $\Aut_{G,\pi} (\Sigma) = \{f \in
\Aut_G(\Sigma), \,\, \pi f =
\pi\}$.  On a une suite  exacte:
$$1\longrightarrow  \Aut_{G,\pi} (\Sigma) \buildrel
{j}\over\longrightarrow \prod_v G
\buildrel\partial\over\longrightarrow \prod_e
G/{(\bigcap_{\omega\mapsto e}} \,G_\omega)
\tag (8.21)$$
o\`u \quad  $\partial ((\sigma_v)_v) \,= \, (\sigma_{e(1)}^{-1}
\sigma_{e(0)})_e$,  et    $e(0)$, $e(1)$ \'etant  respectivement
l'origine et l'extr\'emit\'e de $e$.

\dem D\'efinissons l'application $j$. Par hypoth\`ese l'automorphisme
$h$ de $C$ induit par un \'el\'ement $f\in \Aut_{G,\pi} (\Sigma)$ est
suppos\'e \^etre l'identit\'e. D\`es lors pour tout sommet $v$ de
$\Gamma$ la courbe  $\Sigma_v = \pi^{-1} (C_v)$ est fix\'ee par $f$,
ce qui entra\^\i ne que l'automorphisme d\'efini par la  restriction
de $f$ \`a cette courbe co\"\i ncide avec la restriction d'un \'el\'ement $\sigma_v \in G$.  En effet, quitte \`a modifier $f$ par
un \'el\'ement de $G$, on peut supposer que $f$ fixe une composante
irr\'eductible $\Delta_v$ de $\Sigma_v$. Sur cette composante $f$ est
alors donn\'e par un \'el\'ement  de $G$, en fait du centre du
stabilisateur de cette composante, donc est \'egal \`a cet \'el\'ement partout. On  d\'efinit alors l'injection $j$ par la collection
$j(f) = (\sigma_v)_v$.   

La condition de coh\'erence qu'on doit
exiger d'une collection $(\sigma_v)$ pour qu'elle d\'efinisse un
automorphisme de $\Sigma$ par recollement des $\sigma_v$, est que
pour tout point
double $Q$ de $C$, donc toute ar\^ete  $e$ de $\Gamma$,  et si $a$ et $b$
sont les extr\'emit\'es de $e$,  alors $ \sigma_a (P) = \sigma_b (P)
 $ pour tout point double  $P \in \Sigma$ d'image
$Q$.  Cette condition qui \'equivaut   \`a  ce que $\sigma_b^{-1}
\sigma_a$ fixe tout point double au dessus de $Q$, assure que les
automorphismes $\sigma_v$ des $\Sigma_v$ se recollent en un
$G$-automorphisme de $\Sigma$. Cette condition s'exprime bien par
$\partial ((\sigma_v)) = 1$. \hfill\break\qed 

\bigskip

\section 9  {  Rev\^etements cycliques}  

\bigskip

   Le champ des
rev\^etements cycliques  admet une description particuli\`ere.\footnote{ 
Il est instructif de   comparer ce champ  avec le champ des
courbes \`a spin \'etudi\'e en d\'etail par Jarvis \cite {43},
\cite {44}.}.       Si la base est $\Bbb P^1$ on donnera
une description tr\`es  concr\`ete du sous-champ ouvert classifiant
les rev\^etements non d\'eg\'en\'er\'es, comme champ quotient, g\'en\'eralisant le r\'esultat de Arsie-Vistoli \cite {4}. La construction est parall\`ele \`a celle de Jarvis \cite {43}, \cite {44}, qui traite des courbes \`a spin.
\bigskip
\beginsection 9.1.   Rev\^etements cycliques versus racines d'un faisceau inversible

\bigskip

 {\it  9.1.1 Racines et quasi-racines d'un faisceau inversible}

\bigskip
 
On fixe un entier $r\geq 2$.
 Le concept de racine $r$-i\`eme  se d\'efinit en deux temps \cite {43}, \cite{44}:  
 \proclaim  D\'efinition 9.1. 
 Soit $D$
  une courbe pr\'estable.
  Une quasi-racine $r$-i\`eme d'un faisceau inversible $\kappa$ sur $
D$ est la donn\'ee d'un couple
$({\Cal L} , \Phi)$, o\`u ${\Cal L}$ est un faisceau sans torsion
de rang un, et $\Phi$ un morphisme
$$\Phi: {\Cal L}^{\otimes r} \longrightarrow \kappa \tag (9.1)$$
On demande  que   $\Phi$  soit un isomorphisme en tout point o\`u ${\Cal L}$
est libre.  \hfill\break
   Une quasi-racine est une racine $r$-i\`eme si on a en plus on a:\hfill\break
i)   En un point singulier de ${\Cal L}$, i.e. un point en lequel
${\Cal L}$ n'est pas libre, le conoyau de $\Phi$ est de rang $r-1$.\hfill\break
 ii)  $r  \deg \, {\Cal L} = \deg {\kappa } \quad (\deg\,  {\Cal L} = \chi({\Cal L}) - \chi({\Cal O}_C))$.

   Dans le contexte des rev\^etements cycliques  on est amen\'e \`a consid\'erer les racines d'ordre $r$ d'un faisceau inversible  ${\Cal O}_D (- \sum_i m_i Q_i)$.
Dans ce cas les conditions ci-dessus \'etant impos\'ees, en particulier la
condition iii)
$$  \deg ({\Cal L}) = 2g - 2 - \sum_i m_i \tag (9.2)$$
Jarvis \cite {43}, \cite {44}  a montr\'e que le champ  classifiant
les donn\'ees  $(D , {\Cal L}, \{Q_i\})$
, o\`u ${\Cal L}$ est une racine
$r$-i\`eme comme  cela vient d'\^etre d\'efini, est un champ de
Deligne-Mumford s\'epar\'e et
lisse au dessus de $\Spec \,{\Bbb Z}[{1\over r}]$.  Un tel r\'esultat
exige  qu'une d\'efinition convenable des familles de racines $r$-i\`eme de faisceaux inversibles soit propos\'ee. Cette d\'efinition
comme on va le voir  est en fait tr\`es naturelle dans le cadre des
champs de Hurwitz.  Elle utilise au pr\'ealable la description  par
Faltings ( \S 7.3.3, voir \cite {31}, \cite {44})  de la d\'eformation verselle des modules   sans torsion sur l'anneau local d'un
point double.   

Comme il est expliqu\'e dans \cite {44}, la donn\'ee  seule d'une racine $r$-i\`eme de   
$\omega_D (- \sum_{i=1}^n b_i Q_i)$, c'est \`a dire la version
relative de la d\'efinition 9.1,  d\'efinit certes un champ alg\'ebrique,
mais qui sauf si $r$ est premier n'est pas lisse. Pour  contourner ce d\'efaut  on doit lorsque
la courbe stable $D$ est singuli\`ere  rendre plus pr\'ecise  la structure de
$Spin$  en fixant la structure locale du faisceau ${\Cal  L}$ au
voisinage d'un point double $P$.  Soit $S$ la base, et $s\in S$
l'image de $P$. Posons $R = \hat {\Cal O}_{S,s}$, et soit $x,y: \, xy
= \pi\in {\Cal M}_R$ un syst\`eme de coordonn\'ees en $P$. On exige de
mani\`ere pr\'ecise que la structure locale de ${\Cal L}$, c'est \`a
 dire apr\`es passage \`a l'anneau local compl\'et\'e (ou un
voisinage \'etale),   soit de la forme\footnote { Ce sont les pure
$n$ roots de \cite {44} Definition 5.1.3. Cette restriction est n\'ecessaire pour que le foncteur des d\'eformations infinit\'esimales
du triplet $(D,{\Cal L},\Phi)$ soit formellement lisse (loc.cit.
Proposition 5.4.4). Cette condition est clarifi\'ee dans \cite {3}.} (voir \S\  7.3.3 pour les notations)
$${\hat {\Cal L}}_Q \cong E(\tau^{a} , \tau^b),\;\;a+b = r, \quad \pi
= \tau^r \quad (\tau \in {\Cal M}_R) \tag (9.3)$$

\bigskip

 {\it 9.1.2. Description   des rev\^etements cycliques } 

\bigskip

  Fixons   une donn\'ee de
ramification $R = \sum_{i=1}^s b_i [H_i,\chi_i] $ (D\'efinition 2.1) relative au groupe  $G =  {\Bbb
Z}/{n{\Bbb Z}}$, et soit $b = \sum_i b_i$. 
  Si on fixe un g\'en\'erateur $\sigma$ de $G$,   on sait que cela  \'equivaut  \`a  sp\'ecifier l'ordre
$e_i$ de $H_i$, ainsi que  les entiers $k_i, \,1\leq k_i<e_i$
(l'holonomie locale)  d\'efinissant les classes $[H_i,\chi_i]$ par
$\chi_i (\sigma_n ^{n/{e_i}}) =  \zeta_n^{k_in/{e_i}}$. Fixons   
 une racine  primitive $n$-i\`eme de l'unit\'e.

Soit un $G$-rev\^etement   $\pi: C \longrightarrow   D \cong C/G$, $C$ et $D$ \'etant lisses.
 Il n'y a pas lieu de pr\'eciser
pour le moment si les points de branchement sont piqu\'es ou marqu\'es. Un tel rev\^etement, \`a groupe de Galois cyclique, est
totalement d\'etermin\'e par la $({\Cal O}_D , G)$ alg\`ebre coh\'erente (\S\  7.3.3)
$$\pi_\star ({\Cal O}_C ) \, = \, \bigoplus_{j=0}^{n-1} \,\,{\Cal
L}_j\tag (9.4) $$
o\`u la d\'ecomposition du membre de droite est la d\'ecomposition
 facteurs isotypiques. On notera ${\Cal
L}_j$  le sous-faisceau propre de l'op\'erateur $\sigma$ de valeur
propre  $\zeta_n^j$.   Dans la
suite , on posera ${\Cal L}_1 = {\Cal L}$.
La multiplication dans $\pi_\star ({\Cal O}_C )$ induit des morphismes
${\Cal L}_j \otimes {\Cal L}_k \longrightarrow {\Cal L}_{j+k},  $
les indices \'etant pris modulo $n$. Il en r\'esulte en particulier
pour tout $j$, et  $d'\,\vert \;d$ des morphismes
$$\Phi_j: {\Cal L}^{\otimes j} \rightarrow {\Cal L}_j, \quad
\Phi_{d',d}: {\Cal L}_{d'}^{\otimes d/d'}
\rightarrow {\Cal L}_d \tag (9.5)$$

Ces morphismes d\'ecrivent la structure d'alg\`ebre sur  $\pi_\star
({\Cal O}_C )$. Le morphisme $\Phi = \Phi_n: {\Cal L}^n \rightarrow
{\Cal O}_D$  permet d'identifier la donn\'ee de ramification  du rev\^etement $\pi: C\rightarrow D$.   Limitons
nous pour d\'ebuter au cas lisse.  Soit  $D$  une courbe lisse d\'efinie sur  $k$; rappelons en premier la description classique des rev\^etements cycliques de base $D$, en
termes du triplet $(D,\Cal L,\Phi)$:
\proclaim Proposition 9.2.     Si au point $Q_i \in D$ l'indice de
ramification est $e_i$, l'action de $\sigma^{n/ e_i}$  au dessus de
$Q_i$ (l'holonomie locale)   est d\'ecrite par $k_i,\,\,1\leq
k_i<e_i$, o\`u  $\nu_i$  est tel que $1\leq \nu_i<e_i$,
$\nu_i k_i\, \equiv\, 1 \,\,(e_i)$ et $m_i\, = \,{n\over
e_i}$, alors
  $\Div (\Phi) = B$, soit
$ \Phi ({\Cal L}^n ) \,= \, {\Cal O}_D \,( - \sum_{i=1}^b m_i
\,\nu_i \,Q_i). $  
  Il y a un unique isomorphisme
d'alg\`ebre  qui est l'identit\'e sur $\Cal L$:
$$\pi_* ({\Cal O}_C) = \bigoplus_{i=0}^{n-1} {\Cal L}_i
\buildrel\sim\over \longrightarrow \bigoplus_{i=0}^{n-1}
  {\Cal L}^{\otimes i} \left( \left[ { { iB}\over n}\right] \right) \tag
(9.6)$$

\dem  Le  premier point   est    classique  (voir par exemple
\cite {56}). Rappelons   bri\'evement les arguments.  Soient $Q$
un point de branchement,  $\Cal O = \Cal O_Q$ l'anneau local de $D$
en $Q$,   $v$ la valuation associ\'ee. Soit  $\overline {\Cal O} =
\bigcap_{P\to Q} \Cal O_P$ la cl\^oture int\'egrale de $\Cal O$ dans le
corps des fonctions de $C$.  
 Soit la d\'ecomposition en sous-espaces propres ${\Cal O} =
\oplus_{i=0}^{n-1}  {\Cal O} \xi_i$; posons  $\xi := \xi_1$.  Notons que
$\xi^{i} = h_i \xi_i$ avec $h_i\in {\Cal O}$. Soit par ailleurs
$\xi_i^n = f_i\in {\Cal O}$. Du fait de la normalit\'e de $\overline
{\Cal O}$, il est clair que  $v(f_i) < n$. Posons $v(f_1) = d =  
{n \over e} \nu$, o\`u $e = \pgcd (n,d), \,1\leq \nu < e$. Notons que $d >
0$ car $Q$ est un point de branchement. Si $\pi\in {\Cal O}$ est une
uniformisante, et si $P$ est un point de $C$ au-dessus de $Q$, l'\'egalit\'e $e v_P(\xi) = \nu v_P(\pi)$ dit que $e / v_P(\pi)$. Comme
$\tau^{n/e} $ est une unit\'e de $\Cal O$, o\`u $\tau  =  
{{\xi^{e}}\over {\pi^\nu}}$, on voit que le nombre des points de $C$
au-dessus de $Q$ est au moins $  n \over e$. Il en d\'ecoule que ce
nombre est $ { n \over e}$, et que l'indice de ramification  est $e$. On
a par ailleurs $\xi^{in} = {h_i}^{n} {\xi_i}^n$ donc $nv(h_i) +
v(\xi_i^n) = nd$, et comme $v(\xi_i^n)  < n$,  on en tire $v(h_i) =
[{ {id}\over n}]$ et $v(\xi_i^n) = n\langle { {id} \over n}\rangle $.  
Il en r\'esulte  la description:
$$ {\Cal L}_i = {\Cal L}^{\otimes {i}} \,( \sum_{j=1}^b \, \,[{i\nu_j
\over e_j}]\,Q_j ) \;\buildrel\sim\over\rightarrow \;{\Cal
L}^{\otimes {i}} ([{iB\over n}])  \tag (9.7)$$
 De mani\`ere plus pr\'ecise il y
a un unique isomorphisme qui est l'identit\'e sur ${\Cal
L}^{i}$. 

On peut conclure que le rev\^etement $\pi:
C\rightarrow D$ est de fait totalement d\'etermin\'e par le couple
$({\Cal L},\Phi)$.  Notons d'abord que le terme de droite dans
(9.7) est une ${\Cal O}_D$-alg\`ebre, la multiplication \'etant
donn\'ee par les applications naturelles$${\Cal L}^{i} [ {{ iB}\over n}]
\otimes {\Cal L}^j[ { {jB}\over n}] \longrightarrow {\Cal L}^{i+j} [ 
{{(i+j)B}\over n}] \tag (9.8)$$
l'indice $i+j$  \'etant pris modulo $n$. Il est clair
que l'unique isomorphisme qui est l'identit\'e sur
$\bigoplus_{i=0}^{n-1} {\Cal L}^{i}$ est un isomorphisme d'alg\`ebre.\hfill\break
\qed

La description   (9.2) d'un rev\^etement
cyclique  entre courbes lisses s'\'etend facilement \`a une famille
de courbes lisses.  Il suffit  de donner un sens au terme
de droite dans (9.6). Supposons que dans la donn\'ee de
ramification les valeurs distinctes prises par les $n_i = { n
\over {e_i}}\nu_i$ sont $p_1 < \dots < p_s$.    Il est pr\'ef\'erable
pour  la suite  de mettre le diviseur $B$ sous la forme  $B =
\sum_{i=1}^s p_i B_i$,  les   $B_i$ \'etant  maintenant
disjoints deux \`a deux. Posons $\deg (B_i)  = k_i$ qui repr\'esente  le nombre d'occurence de l'holonomie $p_i$, de sorte que
$\sum_i k_ip_i = \sum_j n_j = nm$, l'entier $m$ \'etant d\'efini par
cette \'equation.

 La base $S$ du rev\^etement est maintenant
arbitraire.  Pour tout $1 < p =  { n \over e }\nu, \,\,
1\leq \nu < e, \, \pgcd (e,\nu) = 1$, soit  $\Delta_p \subset C$ le
sous-sch\'ema  des points de $C$ d'holonomie $p$.   On sait (section
4)  que si $p=p_i$,  $\Delta_p$ est un diviseur de Cartier relatif
\'etale de degr\'e $k_i  { n \over {e_i}}$ sur $S$.  Si $p = p_i$, on notera $\Delta_p
= \Delta_i$. Le groupe quotient de  $  {\Bbb Z}/
{n\Bbb Z}$ par le sous-groupe d'ordre $e_i$, agit librement sur
$\Delta_i$. Le quotient  $B_i \subset D$ de $\Delta_i$  est donc un
diviseur de cartier relatif   \'etale de degr\'e $k_i$ sur $S$. On a alors:

\proclaim Proposition 9.3. Soit $\pi: C\to D$ un rev\^etement
cyclique entre courbes lisses de base  $S$. Avec les notations pr\'ec\'edentes,
en particulier ${\Cal L} := {\Cal L}_1$, on a    $B =
\sum_{i=1}^s p_i B_i$,  et: \hfill\break   
1)   Il y a un unique isomorphisme ${\Cal
L}_i \buildrel\sim\over\rightarrow {\Cal L}^{\otimes i} \left[ 
{{iB}\over  n}\right]$ qui est l'identit\'e sur ${\Cal L}^{\otimes i}$. 
Mieux il  y a un isomorphisme unique de l'alg\`ebre  "cyclique" $E  =
\pi_* ({\Cal O}_C)$  sur  
$$\bigoplus_{i=0}^{n-1} {\Cal L}^{\otimes i}
\left[  {{iB}\over  n}\right] \tag (9.9)$$ 
2) le champ de Hurwitz
${\Cal H}_{g,n,\xi}$ dont les objets sont les rev\^etements cycliques
de degr\'e $n$, de base une courbe lisse de genre $g$, et de donn\'ee de ramification $\xi$ est isomorphe au champ dont les objets sont
les triplets $(D,{\Cal L},\Phi)$, o\`u $D$ est lisse de genre $g$, et
$({\Cal L},\Phi)$ est une racine d'ordre $n$ de ${\Cal O}_D(-B)$, le
diviseur $B = \sum_{j=1}^s p_j B_j$ \'etant comme  prescrit par la
donn\'ee $\xi$.

\dem   Du fait que
$B = \sum_{j=1}^s p_j B_j$, avec $B_j$ \'etale de degr\'e $k_j$ sur
la base $S$, on peut poser 
$$\left [ { { iB}\over n}\right] = \sum_j
\,\,[ { {ip_j}B_j\over  n}] \tag (9.10)$$
 L'assertion 1) ainsi que  la
premi\`ere partie de 2) se v\'erifient  par un argument local
analoque \`a celui utilis\'e dans la preuve 9.2 (voir note de
bas de page). On peut d'une autre mani\`ere tout d\'eduire de la Proposition 9.2. En
effet si la base $S$ est int\`egre, l'\'egalit\'e de diviseurs de
Cartier relatifs $\Div (\Phi_i) = \left[ { {iB}\over  n}\right] $ d\'ecoule du fait qu'il y a \'egalit\'e sur toute fibre g\'eom\'etrique. Comme toutes ces constructions sont compatibles aux
changements de base, le cas g\'en\'eral se d\'eduit de ce cas
particulier par changement de base \`a partir de la d\'eformation
universelle d'une quelconque fibre. Noter que l'isomorphisme invoqu\'e
 est unique.
 Pour finir la preuve de 1),   notons que le
terme de gauche  dans (9.10) poss\`ede une structure naturelle de ${\Cal O}_D$-alg\`ebre, l'identification est alors claire.

 Pour 2) on observera
d'abord que la cat\'egorie form\'ee des triplets $(D,{\Cal L},\Phi)$
est fibr\'ee en  groupo\"{\i}des, et en fait est un champ.
L'assertion se r\'esume  donc \`a v\'erifier que le foncteur
indiqu\'e est \`a la fois un monomorphisme et un \'epimorphisme
\cite {51}. Comme monomorphisme signifie que pour deux rev\'etements
$\pi: C\to D$ et $\pi': C'\to D'$ de base $S$, le morphisme 
$$\Hom
(\pi,\pi') \rightarrow \Hom \left( (D,{\Cal L},\Phi),(D',{\Cal
L'},\Phi')\right)$$
est bijectif, il est clair que cela d\'ecoule de
2). Pour v\'erifier que ce foncteur est un \'epimorphisme, il suffit
de noter que si on part d'un objet $(D/S,{\Cal L},\Phi)$ comme il est
indiqu\'e, alors 
${\Cal A} = \bigoplus_{i=0}^{n-1}\,  {\Cal L}^{\otimes
i} \left[ {{iB}\over n}\right] $
 a  une structure naturelle de
$O_D$-alg\`ebre, il suffit alors de poser $C = \Spec_{{\Cal O}_D} ({\Cal
A})$.  \hfill\break\qed

 \bigskip
 \beginsection 9.2.  Rev\^etements cycliques stables
  
  \bigskip
   {\it 9.2.1 Rev\^etements cycliques
stables et quasi-racines d'un faisceau inversible}  \bigskip

  Si $\pi: C
\to D$ est un rev\^etement stable d\'efini sur $k$,  avec $D$   singuli\`ere,   on peut former de nouveau la d\'ecomposition
 
   $$\pi_\star ({\Cal O}_C ) \, = \, \bigoplus_{j=0}^{n-1} \,\,{\Cal
L}_j \tag (9.11)$$
Le ${\Cal O}_D$ module $\pi_\star ({\Cal O}_C )$ n'est
plus en g\'en\'eral localement libre de rang $n$, mais seulement sans
torsion de rang $n$.  Rappelons dans le contexte des rev\^etements
cycliques la description   pr\'ecise des singularit\'es des ${\Cal
L}_j$ (\S \ 7.3.3). Consid\'erons  $P\in C$ un point double de stabilisateur
d'ordre $e \vert n$, et soit $Q$ le point double image dans $D$. Si
$\sigma$ est un g\'en\'erateur fix\'e de $G$, on notera  $\tau =
\sigma^{n\over e}$ le g\'en\'erateur distingu\'e du stabilisateur
$H = G_P$. La description locale au point  $Q$ de la d\'ecomposition
9.11  est  comme suit (\S 7.3.3).  

Choisissons  des "coordonn\'ees
locales"  $x,y$ le long des branches telles que  
$$\tau (x) = \zeta
_e^k x \,,\, \tau (y) = \zeta_e^{e-k} y \tag (9.12)$$
    conditions exprimant
la stabilit\'e de l'action.  Notons  par un tilde la normalisation
$k[[x]]\times k[[y]]$ de ${\hat{\Cal O}_Q}$.  Si un caract\`ere de
$H$ est identifi\'e \`a sa valeur sur le g\'en\'erateur $\tau$,
 rappelons la  d\'ecomposition  (7.26):

\proclaim Lemme 9.4. La d\'ecomposition en $H$-sous faisceaux
propres de ${\hat{\Cal O}_P}$ est:   $$ {\hat{\Cal O}_P} \;=\;
{\hat{\Cal O}_Q} \;\bigoplus \; \left ( \bigoplus_{\alpha = 1}^{e-1}
{\widetilde{\hat{\Cal O}_Q}} (x^\alpha,y^{e-\alpha}) \right ) \tag (9.13)$$
\qed  

   En application des r\'esultats g\'en\'eraux  de
la section 7.3, on peut d\'ecrire le lieu singulier  des faisceaux ${\Cal
L}_j$ (voir en particulier le corollaire 7.9)\footnote{
 Noter que  ${\Cal L}_j$ et ${\Cal L}_k$ ont m\^eme lieu
singulier ssi $(j,n) = (k,n)$. }:
\proclaim Proposition 9.5. Le faisceau ${\Cal L}_j$ est localement
libre en un point double $Q \in
D$ si et seulement si $e$ \'etant l'ordre du sous-groupe d'isotropie
en  un point $P$ au dessus
de $Q$,  $j \equiv 0 \pmod e$.

  Posons ${\Cal L} = {\Cal L}_1$,  et analysons maintenant les
morphismes (9.5) induits par la
multiplication, et plus particuli\`erement  les morphismes $\Phi_j:
{\Cal L}^{\otimes j} \longrightarrow {\Cal L}_j\;\; (1\leq j\leq n).$

\proclaim Proposition 9.6.
  Si $Q\in  D$ est du type NS ($e > 1$), il y a deux cas:  \hfill\break
i) Le point $Q$ est une singularit\'e de ${\cal L}_j$,  alors
la dimension de la fibre en $Q$ du conoyau de $\Phi_j$ est $j-1$, et
l'id\'eal de ${\cal O}_Q$ image de  $(\Phi_j)_Q$ est non principal.\hfill\break
ii) Dans le cas contraire, si $e \vert j$, la dimension du conoyau
de $\Phi_j$ en $Q$ est encore $j-1$, mais l'image de $\Phi_j$ est
alors un id\'eal principal, i.e. libre de rang un. \hfill\break
 Si $j=n$, le morphisme $\Phi: {\Cal L}^{\otimes n}
\longrightarrow {\Cal O}_D ( - B) $ induit par la multiplication fait de ${\Cal L}$ une
racine $n$-i\`eme de ${\Cal O}_D ( - B)$.

  \dem Cela r\'esulte d'un calcul \'el\'ementaire (voir la section
7.3.2). Traitons le premier cas, donc $e\;\vert\; j$. Il r\'esulte de
la proposition 9.5 , et avec des notations \'evidentes, qu'on peut
prendre  la "r\'esolvante"
$X_j = \sum_{i=1}^{n/e} \; \zeta_e^{-ij} \sigma^{ij} (x^{j\nu} ,
y^{e-j\nu}) $
comme g\'en\'erateur local du ${\widetilde{\hat{\Cal O}_Q}}$ module
principal $(\hat{{\Cal L}_j})_Q$. Le m\^eme r\'esultat vaut pour $j=1$.
On voit alors imm\'ediatement que l'image par $\Phi_j$ de
$X_1^{\otimes j}$ est
$$ (\Phi_j)_Q \;(X_1^{\otimes j}) \;=\; (u^{[{j\nu\over
e}]}\;,\;v^{j-1-[{j\nu \over e}]}) X_j
 $$
Du fait que la fibre en $Q$ de ${\Cal L}^{\otimes j}$ modulo le
sous-module de torsion
n'est pas libre, car le point est NS, les deux membres de l'\'egalit\'e (9.14) sont sans torsion de rang
un, non libres, d'o\`u le conoyau de $\Phi_j$ en $Q$  s'identifie \`a
$$ {{k[[u]] \times k[[v]]} \over {(u^{[{j\nu\over e}]},
v^{j-1-[{j\nu\over e}]}})}$$ la dimension de ce conoyau  est donc
$j-1$. Pour la derni\'ere assertion, la seule chose \`a v\'erifier est l'\'egalit\'e $n\,\deg ({\Cal L}) = - \deg (B)$. Comme cette \'egalit\'e est invariante par d\'eformation, de m\^eme la formation de $\Cal L$, il suffit de la v\'erifier dans le cas lisse, ce qui est clair.\hfill\break
\qed

Il est  naturel de   demander  si le morphisme $\Phi$ provient par $\rho_*$ d'un morphisme analogue $\tilde\Phi: {\Cal O}_{\tilde D} (n) \rightarrow  {\Cal O}_{\tilde D}$ (voir Proposition 7.12 et la construction qui pr\'ec\`ede).  La r\'eponse est donn\'ee essentiellement par (\cite {44}, Proposition 3.1.5).
\proclaim Proposition 9.7. Sous les conditions pr\'ec\'edentes il existe un unique morphisme
$$\tilde \Phi: {\Cal O}_{\tilde D} (n) \rightarrow  {\Cal O}_{\tilde D}( - B) \tag (9.15)$$ tel que $\rho_* (\tilde \Phi) = \Phi$.

\dem Soit $U$  l'ouvert image r\'eciproque du compl\'ementaire du lieu singulier de $\Cal L$. Le morphisme $\rho$ est un isomorphisme sur  $U$. Comme $U$ est sch\'ematiquement dense,  le morphisme $\tilde\Phi$ s'il existe est unique.  Il suffit donc de le construire au-dessus d'un voisinage d'un point singulier $Q$, appel\'e $D$ pour simplifier. De la description de $\tilde D$ comme sous-sch\'ema de $\Bbb P^1\times D$ donn\'epar les \'equations 
$$\xi y - qv = \xi p - xv = 0$$
 on note que $\tilde D = U_1 \cap U_2$, avec $U_1 = \Spec (\hat {\Cal O}_Q [s] / (sy - q , sp - x)), \,\, U_2 = \Spec ( \hat {\Cal O}_Q [t] / (y - qt , p - xt))$. Le cocycle qui d\'efinit  ${\Cal O}_{\tilde D} (n)$ est $s_{1,2} = t^n\in \Gamma (U_1\cap U_2 , {\Cal O}_{\tilde D})$. Il suffit donc de trouver des fonctions r\'eguli\`eres $\phi_i$ sur $U_i \, (i=1,2)$ avec $\phi_1 = t^n \phi_2$. Mais on a $p = \pi^\nu , q = \pi^{e - \nu}$. Si $a = m\nu , b = m(e - \nu)$, alors $p^b = q^{a}$. Il suffit de prendre $\phi_1 = y^b$ et $\phi_2 = x^{a}$. La v\'erification est imm\'ediate. \hfill\break\qed 
 
  Avec les notations pr\'ec\'edentes,  \`a un point de type NS, on
attache une paire d'entiers (non ordonn\'es) {\it le symbole }   $(a , b)$ :
$$a = m\nu,\;b=n - m\nu,\; m = {n\over e},  \quad  a + b
= n\tag (9.16)$$
On a $1\leq \nu < e$, et $\pgcd (\nu ,  e) = 1$. En particulier
$\pgcd (a , n) = m$. Le symbole  d\'ecrit    l'action  du
stabilisateur d'un point  $P \in C$ au dessus de $Q$, sur les
branches en $P$. Notons que l'image de $\Phi$ en $Q$ est l'id\'eal
(non principal) $(u^{a} , v^b)
\widetilde {\hat {\Cal O}_Q}$ de colongeur $n-1$.


  Soit maintenant un object $\pi: C \rightarrow  D = C/G$  du champ de
Hurwitz $\overline{\Cal H}_{g,n,\xi}$, au dessus du sch\'ema de base $S$. Dans
la notation du champ on omet   le genre
$g$ (fix\'e) de $C$, ainsi que la donn\'ee de ramification $\xi$. Notons $q$ et $p$ les morphismes
structuraux des $S$-sch\'emas $C$ et $ D$. Soit  de nouveau la d\'ecomposition en facteurs isotypiques:
$$\pi_\star ({\Cal O}_C ) \, = \, \bigoplus_{j=0}^{n-1} \,\,{\Cal
L}_j\tag (9.17) $$
les faisceaux ${\Cal L}_j$ \'etant des faisceaux sans torsion  de
rang un (relatifs) sur $S$, i.e.
plats sur $S$ et les fibres au-dessus d'un point g\'eom\'etrique \'etant sans torsion de rang un. La structure locale de ${\Cal L} =
{\Cal L}_1$ peut \^etre pr\'ecis\'ee de la mani\`ere suivante:
\proclaim  Proposition 9.8. Soit $Q \in  D_s$ un point de la fibre
au dessus du point
g\'eom\'etrique $s\in S$. On suppose que $Q$ est un point singulier
de ${\Cal L}$,  de symbole $(a , b)$; alors  $\hat {\Cal L}_Q $ est
un module sans torsion de rang un au dessus de $\hat {\Cal O}_s$ de
type $E( \tau^{a} , \tau^b)$.

  \dem  Il suffit d'appliquer la Proposition 7.11.\hfill\break \qed 

 On  prouve  que  le champ des rev\^etements
cycliques de degr\'e $n$ et de donn\'ee de ramification
fix\'ee $\xi$ (not\'e  $ \overline{\Cal H}_{g,n,\xi}$), est isomorphe au champ  dont les objets sont les
courbes stables marqu\'ees (ou piqu\'ees)\footnote {la courbe $D$ est
marqu\'ee par le diviseur $B = B_1+\dots+B_s$, le diviseur relatif
$B_i$ \'etant \'etale de degr\'e $k_i$ sur la base $S$, et  les
$B_i$ \'etant disjoints deux \`a deux. Rappelons que $B_i$ est le
lieu des points de branchement d'holonomie fix\'ee $p_i$.} $(D ,  B =
\sum_{j=1}^s p_j B_j)$  \'equip\'ees d'une racine $n$-i\`eme ${\Cal
L}$ du faisceau ${\Cal O}_D (-  B), \, p_i = m_i \nu_i$. 
  \proclaim  Th\'eor\`eme 9.9.   Le foncteur $(\pi:C \to D) \mapsto (D,{\Cal
L},\Phi)$ \'etablit un isomorphisme de $ \overline{\Cal H}_{g,n,\xi}$
sur le champ dont les objets sont les triplets $(D,{\Cal L},\Phi)$, o\`u
 $D$ est une courbe stable marqu\'ee par $B = \sum_{i=1}^s B_i$
(voir note 28), et $({\Cal L},\Phi)$ est une racine, dans le sens (9.1)n d'ordre $n$ de
${\Cal O}_D(-\sum_{i=1}^s p_iB_i)$.

\dem Le fait que la
cat\'egorie fibr\'ee en groupo\"{\i}des d'objets les triplets
$(D,{\Cal L},\Phi)$ soit un champ alg\'ebrique  s\'epar\'e et lisse
est une cons\'equence directe des r\'esultats de Jarvis (\cite {44},
Theorem 5.3.1). Essentiel est le r\'esultat qui d\'ecrit la d\'eformation universelle d'une racine d'ordre $n$ (loc.cit, et \cite {43} Theorem 2.3.2). Ce r\'esultat  dit que si $(D,{\Cal L},\Phi)$ est
une courbe stable marqu\'ee par le diviseur $B_1+\dots+B_s$, munie
d'une racine d'ordre $n$ de ${\cal O}(-B), \, (B = \sum_i p_iB_i)$, la base de la d\'eformation universelle est une alg\`ebre de s\'eries formelles
$k[[\tau_1,\dots,\tau_b,t_{b+1},\dots,t_{3g-3+b}]]$, et le
morphisme qui a une d\'eformation de ce  triplet associe la  d\'eformation correspondante de la courbe marqu\'ee $D$, s'identifie \`a
 $$\psi: \Spec
k[[\tau_1,\dots,\tau_b,\tau_{b+1},\dots,\tau_{3g-3+b}]]
\longrightarrow \Spec k[[t_1,\dots,t_b,\dots,t_{3g-3+b}]]\tag
(9.18)$$
o\`u $\psi^* (t_j) = \tau_j^{e_i} $ si $j\leq b$, et $\psi^*
(t_j) = \tau_j$ si $b<j\leq 3g-3+b$. Si $Q$ est un point singulier de
$\Cal L$ de signature $(a,b), \, a+b=n$, on pose $e =  { n \over{\pgcd
(a,b)}}$. On notera que ce r\'esultat dit que si le triplet $(D,{\Cal
L},\Phi)$ de base $\Spec k$  d\'erive du rev\^etement  $\pi: C\to D$,
alors $\pi$ et $(D,{\Cal L},\Phi)$ ont "m\^eme" th\'eorie des d\'eformations (comparer avec le Th\'eor\`eme 5.5). Cela signifie que si
$\pi: {\Cal C}\to {\Cal D}$ est la d\'eformation universelle du rev\^etement
$\pi: C\to D$, alors  le triplet image $({\Cal D},{\Cal L},\Phi)$
est la d\'eformation universelle de $(D,{\Cal L},\Phi)$. Ceci \'etant,
arm\'e de cet argument, la preuve est analogue \`a celle de la
proposition 9.4 (3). 

Montrons que le foncteur est un monomorphisme.
Cela revient d'abord \`a prouver que si un $S$-automorphisme $f$ de
$C$ induit l'identit\'e sur $(D,{\Cal L},\Phi)$, alors $f=1$. Comme le
groupe $\Aut (C/S)$ est fini non ramifi\'e, il suffit de le voir
sur les fibres, donc si $S = \Spec \, k$. Il est clair que $f$ est
l'identit\'e en dehors des points doubles de $C$, donc $f=1$. On
suppose maintenant avoir deux objets $\pi: C\to D$ et $\pi: C'\to D$
avec une m\^eme image $(D, {\Cal L},\Phi)$, on prouve que les deux
rev\^etements sont isomorphes. Du fait de l'unicit\'e de
l'isomorphisme, le probl\`eme est local sur $S$, on peut donc supposer
que  tout est d\'efini au-dessus de la base de la d\'eformation
universelle, commune aux trois objets, disons de la fibre $(D_s,{\Cal
L}_s,\Phi_s)$. Cela permet de supposer que $S$ est lisse (connexe), et
que les fibres g\'en\'eriques des courbes en jeu sont lisses. 

Soit le $S$-sch\'ema $\nu: \Isom (\pi , \pi') \longrightarrow S$ 
form\'e des isomorphismes induisant l'identit\'e sur $(D,{\Cal L},\Phi)$.  On
sait que $\nu$ est  fini non ramifi\'e sur $S$. On sait par
ailleurs   que $\nu$ est un monomorphisme dominant, c'est donc un
isomorphisme.  

Reste \`a voir que localement  sur $S$ tout triplet
$(D,{\Cal L},\Phi)$ de base $S$ provient d'un rev\^etement.   Il suffit
de relever un triplet d\'efini sur la base $S = \Spec\  k$, ceci du fait que   l'identification
des foncteurs de d\'eformations permet alors de relever localement
une famille arbitraire. Dans le cas ponctuel, on d\'eforme le triplet
$(D,{\Cal L},\Phi)$ \`a une base $\Spec\,  R$, o\`u $R$ est un anneau
de valuation discr\`ete complet de corps r\'esiduel $k$, de fibre g\'en\'erique $D_K$ lisse. Quitte \`a \'epaissir $R$, on peut supposer
que la fibre g\'en\'erique du  triplet provient du rev\^etement
$C_K\to D_K$.  Le mod\`ele stable de ce rev\^etement  (\S \ 5.3) r\'epond \`a la question.
\hfill\break\qed

\medskip

  {\it 9.2.3.  Composantes irr\'eductibles du bord} 

\bigskip

La description des composantes irr\'eductibles 
du bord de $ \overline{\Cal H}_{g,n,\xi}$  d\'ecoule 
directement des r\'esultats de la section 7.2.  On doit consid\'erer une
base   qui  correspond soit \`a un segment, soit une boucle.

\indent  {\bf Le segment:} \  Il s'agit de d\'ecrire les
diviseurs du champ $  \overline{\Cal H}_{g,n,\xi}$ d\'etermin\'es par les graphes modulaires de Hurwitz $\Gamma$  qui recouvrent
le segment. D\'ecrire un rev\^etement stable $C\rightarrow  \Sigma$
de groupe $G  =  {\Bbb Z}/{n{\Bbb Z}}$ de base $\Sigma =
\Sigma_1 \cup \Sigma_2$, courbe avec deux composantes lisses de
genres respectifs $1\leq g_1\leq g_2, \;\;g_1+g_2 = g$, un seul point
double, et de donn\'ee de ramification $\xi$, est \'equivalent \`a
la donn\'ee d'une racine $n$-i\`eme du faisceau  ${\Cal O}_\Sigma (-
B)$, $B$ \'etant le diviseur de branchement.   Il s'av\`ere cependant
beaucoup plus commode dans les descriptions qui suivent de conserver
le point de vue des rev\^etements. Dans la suite la base sera not\'ee
indiff\'eremment  $D$ ou $\Sigma$. Le diviseur de branchement sera
not\'e de la mani\`ere suivante 
$ B = \sum_{\alpha = 1}^b
m_\alpha \nu_\alpha Q_\alpha,  $
 avec  $m_\alpha = { n\over
{e_\alpha}}$,  $1 \leq \nu_\alpha < e_\alpha $ et $\pgcd (e_\alpha ,
\nu_\alpha) = 1$. On notera que  
$ \sum_{\alpha = 1}^b m_\alpha
\nu_\alpha \equiv 0 \pmod {n}.  $

Le rev\^etement $C\rightarrow \Sigma = \Sigma_1\cup\Sigma_2$ comme
sp\'ecifi\'e ci-dessus impose \`a $C$ a avoir la topologie suivante
(comparer avec la discussion g\'en\'erale  de la section 7.4):
$$C  = \Ind_{G_1}^G \;C_1\;\vee\;\Ind_{G_2}^G\,C_2 \tag (9.19)$$
pour deux courbes lisses $C_1$ et $C_2$ induisant des rev\^etements
$C_1\rightarrow \Sigma_1$ et $C_2\rightarrow \Sigma_2$ de groupes
respectifs $G_1$ et $G_2$, sous-groupes de $G$.
Le stabilisateur d'un point double $P$ (on suppose que $P \in C_1\cap
C_2$) de $C$ \'etant
$H\subset G_1\cap G_2$. Soit $n_i = \# G_i \;(i=1,2),\;e = \# H$.
Donc $e \;\vert \;  n_1,n_2$, et
du fait de la connexit\'e de $C$, $\ppcm\, (n_1,n_2) = n$. Il est
ais\'e de voir que si le diviseur de
branchement $B_i^\star$ de $C_i \rightarrow \Sigma_i$ est  \'ecrit
sous la forme
$$B^\star_i =  \sum_{\alpha =1}^{b_i} m_\alpha^{i}\nu_\alpha^{i}
Q_\alpha^{i}\; + \; { {n_i}\over
{e} }\nu^{i}Q^{i}\;=\; B_i +   {{n_i}\over{e}} \nu^{i}Q^{i} \;\;(i=1,2) \tag (9.20)$$
avec $m_\alpha^{i} = { { n_i}\over {e_\alpha^{i}}}$, les  $Q_\alpha^{i} ({\rm resp.} Q^{i}) \;(i=1,2)$
\'etant  les points qui apr\`es
identification donnent le point double $Q_\alpha ({ \rm resp.} Q)$. Le diviseur de
branchement du rev\^etement  $\pi: C \rightarrow \Sigma$ est alors
$$ B =  { {n}\over {n_1}}\left(\sum_{\alpha =1}^{b_1}
m_\alpha^{1}\nu_\alpha^1  Q_\alpha^{1}\right)\; +
 { {n} \over{n_2}} \left(\sum_{\alpha =1}^{b_2} m_\alpha^{2}
\nu_\alpha^2 Q_\alpha^{2}\right)
\;=\;{ {n}\over {n_1}} B_1 +  { {n} \over {n_2}} B_2 \tag (9.21)$$
Notons  comme cons\'equence de (7.20) les congruences
$\sum_{\alpha =1}^{b_i} m_\alpha^{i}\nu_\alpha^{i} \;+\;  
{{n_i}\over {e} }\nu^{i}\; \equiv \;0\;\;
\hbox{ \rm mod}\; n_i\;\;(i=1,2) $
qui si $B_i$ est connu, d\'eterminent de mani\`ere unique $e$ et les
$\nu^{i}$. Par exemple $  { {n_i}\over {e}} \nu^{i}$ est le
reste modulo $n_i$ de  $\sum_{\alpha =1}^{b_i}
m_\alpha^{i}\nu_\alpha^{i} $; notons que $\nu^1 + \nu^2 = e$.   En d'autres termes  la signature du
point  double $Q$, obtenu par identification de $Q^1$
et $Q^2$ est  donn\'ee par (9.15)  
$$  {a =  { {n}\over {e}} \nu^1, b =
 { {n} {e}\over \nu^2}}\tag (9.22)$$  
 On a $a + b = n$  et  $ { n \over e} = \pgcd (a , b)$.   La   partition not\'ee dans
la suite $\pi$, du diviseur de branchement (de mani\`ere \'equivalente de la donn\'ee de ramification) en
  $$  {(\pi) \,\,\, B =   { {n}\over {n_1}} B_1 +  { {n} \over{n_2}} B_2 }
\tag (9.23) $$
d\'etermine tous les autres param\`etres i.e. $m^1\nu^1 = {
{n_1}\over {e^1}}\nu^1$ et $m^2 \nu^2 =  {{n_2}\over  {e^2} }\nu^2$, restes
respectifs modulo $n_i$
de $\sum_{\alpha =1}^{b_i} m_\alpha^{i}\nu_\alpha^{i} \;\;(i=1,2)$.
Rappelons que dans l'\'ecriture
pr\'ec\'edente, $1 \leq  \nu^{i} < e^{i}$ et $(\nu^{i} , e^{i}) = 1$.
La congruence initiale
$$\deg (B) \;=\;   { {n}\over {n_1}}\left(\sum_{\alpha =1}^{b_1}
m_\alpha^{1}\nu_\alpha^{1}\right)
\;+\; { {n} \over{n_2}}\left(\sum_{\alpha =1}^{b_2}
m_\alpha^{2}\nu_\alpha^{2}\right)\;\equiv
\;0\;\; \pmod  n\tag (9.24)$$ entra\^{\i}ne imm\'ediatement   $e^1 = e^2$.
Notons par ailleurs que la
stabilit\'e impose $2g_i - 2 + b_i \geq  0 \;\;(i=1,2)$.    La signature $(a,b)$ du point double est en fait d\'etermin\'ee par les seuls supports de la partition $\pi$  (9.23), i.e. la partition des points de
branchements  $[1 , b] = I_1 \coprod I_2$.
Nous noterons dans la suite  $\delta_{g_1,g_2,\pi} $ ou $\delta_{g_1,g_2,B_1,B_2} $  la composante irr\'eductible du bord correspondante \`a  ces choix. Son image dans $\overline {\Cal M}_{g',b}$ est la
composante $\delta_{g_1,g_2,I_1,I_2}$. Noter que $\delta_{g_1,g_2,B_1,B_2} = \delta_{g_2,g_1,B_2,B_1} $.

\medskip

  \indent{\bf La boucle: \  ($g'\geq 1$) }  
On suppose que la base du rev\^etement est maintenant une courbe
irr\'eductible $\Sigma$ de genre  g\'eom\'etrique $g' - 1$, avec un
unique point double $Q$.  Notons $\tilde \Sigma$ la normalisation de
$\Sigma$ et soient $Q^1, Q^2$,
les deux points de $\tilde \Sigma$ pr\'e-images de $Q$. On sait (\S
7, 7.13) que pour l'espace total $C$ du rev\^etement, il y a deux
possibilit\'es:\hfill\break
\indent   i) La courbe $C$ est irr\'eductible, disons avec $d =  {
n\over e}$ points doubles. Alors la normalisation $\tilde C$ est irr\'eductible, et le rev\^etement cyclique $\tilde C \rightarrow \tilde
\Sigma$ est ramifi\'e avec $b + 2$  points de branchement, ceux
provenant des points de branchement de $C \rightarrow \Sigma$ et
$Q_1, Q_2$. Le diviseur de  branchement est donn\'e  par
$$\tilde B =  B + m_1\nu_1 Q^1 + m_2\nu_2 Q^2 \tag (9.25)$$
$B$ \'etant le diviseur de branchement  de $\Sigma \rightarrow C$.
Notons  que $m_1\nu_1 + m_2\nu_2 = n  $, en fait 
${{n}\over  {e}} = m_1 = m_2$  et  
$a = m_1\nu_1, b = m_2\nu_2$ est la signature du point double $Q$.   Nous noterons  $\delta_{0,a,b}$,   la composante d\'ecrite de cette mani\`ere.\hfill\break
\indent   ii) La courbe $C$ est r\'eductible. Alors il y a dans $C$
une seule $G$-orbite de composantes irr\'eductibles, fixons $C_0$ l'une
d'entre elles. Il y a aussi une seule orbite de points doubles, soit
$P \in C_0$ l'un de ces points. Notons $G_0$ le stabilisateur de
$C_0$ et $H$ celui de $P$, d'ordres respectifs $n_0 < n$ et $e$, et
$e\vert n_0$. Soit $g_0\in G\setminus G_0$ un \'el\'ement tel que $P
\in C_0 \cap g_0(C_0)$.  On sait que $G_0\cup g_0$
engendre $G$ du fait de la connexit\'e de $C$ (\S 7.2), et que par ailleurs
$H \subset G_0 \cap g_0 G_0 g_0^{-1} = G_0$ puisque le groupe est
ab\'elien. Notons que l'\'el\'ement $g_0$ est d\'etermin\'e
modulo $G_0$, et au changement $g_0 \mapsto g_0^{-1}$ pr\`es. La
courbe $C$ est alors obtenue
par identification des couples de points $(gP,gg_0^{-1}P) \in
Ind_{G_0}^G \;(C_0)$.  Le diviseur de
branchement de $C_0\rightarrow \Sigma_0$ est visiblement  de la forme
$$B_0 +  {{n_0} \over {e} }\nu^1 Q_1 +  { {n_0} \over{e}} \nu^2 Q_2,\;
\quad   \nu^1  + \nu^2  =  e
\tag (9.26) $$
avec n\'ecessairement la  factorisation de $B$ d\'esign\'ee (de
nouveau)  $\pi$:
 $$  (\pi): \,\, B = {  n \over {n_0}} B_0   \tag (9.27)$$
La signature du point double $Q$ est  $a = {  n \over {n_0}} a_0, \;b =
{  n \over {n_0}} b_0$, avec
 $ a_0 =  {{n_0}\over {e}} \nu^1 \;, \;b_0 =  { {n_0} \over {e} }\nu^2.  $
Noter que ce choix n'est soumis \`a aucune contrainte, i.e.
$e\vert n_0$ \'etant fix\'e, on choisit
$\nu_1 < e$,  premier \`a $e$, alors $\nu_2$ s'en d\'eduit du
fait de l'\'egalit\'e  (9.26). En sens
inverse, on a $e =  { {n_0} \over {\pgcd (n_0 , a_0)}}$.

 La notation retenue pour cette composante
irr\'eductible est $\delta_{0,\pi,a,b}$. Le cas i) correspond  \`a  $n_0 = n$. La notation regroupe  donc les deux cas. \hfill\break
\qed 
   
 \medskip

\beginsection 9.3.  Rev\^etements
cycliques de $\Bbb P^1$

\bigskip

 Lorsque la base est $\Bbb P^1$  on
 identifie  le champ (ouvert) de Hurwitz \`a un champ quotient
  (comparer avec \cite {4} qui se limite aux rev\^etements
uniformes).

\bigskip

  {\it 9.3.1. Quelques calculs de torseurs}  

\bigskip

 Il est
  connu que la classification des objets au dessus de la cat\'egorie des $k$-sch\'emas  ($k$ corps, ou anneau de base),  \'etale-localement  (ou fpqc) d'une "marque" donn\'ee (terminologie de
Demazure-Gabriel \cite {21}),  est \'equivalente \`a  la donn\'ee
du champ classifiant $\B G$, avec  pour groupe structural  $G$,  le
groupe des automorphismes de l'objet marque.  Une d\'efinition
plus pr\'ecise   est comme suit.  Soit $\Cal M$ un champ (alg\'ebrique)
au dessus de $k$. Soit $X$ un $k$-sch\'ema, et soit un objet
$P\in {\Cal M}(X)$, l'objet "marque". Disons que   $Q \in {\Cal
M}(U)$ est de marque $P$, s'il y a un morphisme $\alpha: U\to X$, un
recouvrement (\'etale) $U'\to U$, tels que $Q\times_U U' \cong
P\times_X U'$, donc si  $\Isom (Q , P\times_X U) (U') \ne
\emptyset$.

On d\'efinit un champ ${\Cal M} (P)$, le champ des objets
de marque $P$, dont la cat\'egorie  des objets au-dessus de $U$ sont
les couples $(Q,\alpha: U\to X)$, avec la propri\'et\'e de dessus,
traduisant le fait que $Q$ est de marque $P$, les morphismes \'etant
les morphismes \'evidents. Ce champ n'est qu'une  pr\'esentation
diff\'erente du champ classifiant $\B (G/X)$, si $G = \Aut (P)$ est
le groupe alg\'ebrique des automorphismes de $P$ (\cite {51},
2.4.2).
\proclaim Proposition 9.10.  Sous les hypoth\`eses qui pr\'ec\`edent,  on a un isomorphisme de champs ${\Cal M} (P)
\buildrel\sim\over \rightarrow  \B (G/X),  $
$G = \Aut_k(P)$ \'etant   le  $X$-groupe alg\'ebrique des automorphismes de
$P$.  En particulier, si $G$ est lisse de type fini, le champ $\Cal
M(P)$ est alg\'ebrique.

\dem  Rappelons bri\`evement la d\'efinition de l'isomorphisme \cite {21}. Si $(Q,\alpha)
\in {\Cal (P)}$ est un objet de marque $P$  au-dessus de $S$, il est
clair que le $S$-sch\'ema $\Isom_S (P\times S , Q\times_X S)$  est  de
mani\`ere naturelle un  $ G$ fibr\'e principal de base $S$, d\'efinissant le foncteur 9.33.  En sens inverse  si $E \to S$ est un
$\Aut _X(P)\times_X  S$ torseur,  le produit contract\'e 
$$E  \times^{\Aut (P)\times_X S}  \,(P\times_X S)$$
 fournit un objet de
marque $P$ de base $S$. Ces deux constructions sont inverses l'une de
l'autre. D'une autre mani\`ere le champ ${\Cal M}(P)$ est visiblement
une gerbe  sur $X$, qui  est neutre car $P$ en est une section sur
$X$. Le r\'esultat d\'ecoule alors de (\cite {51}, Lemme
(3.31))\hfill\break\qed  

 Par exemple   on  peut  prendre la marque $k^n$, avec ${\Cal M}$ le champ des
faisceaux localement libres de rang $n$ sur les $k$-sch\'emas (\cite
{51}, Th\'eor\`eme 4.6.2.1). Alors le principe qui vient d'\^etre
rappel\'e   dit  que ce champ est  isomorphe \`a $\B \GL (n)$.
Plus significatif pour la suite est le cas du groupe projectif lin\'eaire $\PGL (n) = \Aut \Bbb P^n$. Le champ classifiant $\B \PGL (n)$
est   isomorphe au champ dont les objets sont les fibrations en $P^n$
au dessus d'un $k$-sch\'ema. Le cas
utile pour la suite est un m\'elange de ces  exemples.  Fixons un
faisceau localement libre  $V$ de rang $r \geq 1$  sur $P^1$, donc en
vertu d'un th\'eor\`eme de Grothendieck de la forme 
$V = {\Cal
O}(n_1)\oplus \dots \oplus {\Cal O}(n_r)$,  
avec $n_1\leq \dots \leq
n_r$. Soit  la cat\'egorie fibr\'ee en groupo\"\i des  ${\Cal P}$, dont
les objets    au-dessus de  $S$ les couples $(D\to S, E)$, avec $D\to
S$ une fibration en $\Bbb P^1$, et $E$ un fibr\'e vectoriel sur $D$
de type $V$ le long des fibres, c'est \`a dire, localement sur $S$,
de la forme $\oplus_{i=1}^r {\Cal L}_i$, o\`u ${\Cal L}_i$ est
inversible de degr\'e relatif $n_i$ .  Noter que cela a un sens
car, localement pour la topologie \'etale, $D = \Bbb P_S^1$.
L'identification des objets  est un exercice ais\'e:
\proclaim
Lemme 9.11.    La cat\'egorie fibr\'ee en groupo\"\i des
$\Cal P$ est un champ alg\'ebrique.  La marque   des objets  de
$\Cal P$ est $(\Bbb P^1,V)$. \hfill\break
Le groupe   des automorphismes  $G$ du
mod\`ele, est   le groupe qui lin\'earise universellement le fibr\'e
$V$.  Ce groupe  s'ins\`ere donc dans une suite exacte$$ 1
\rightarrow H = \Aut (V) \rightarrow G \rightarrow \PGL (1)
\rightarrow 1 \tag (9.28)$$

\dem Le
premier point se r\'esume essentiellement \`a montrer que si
$(D,E)$ et $(D',E')$ sont deux objets de $\Cal P$ au-dessus du sch\'ema $S$, alors le faisceau sur $S_{et}$ 
$$ (T\to S) \mapsto \Isom_T
\left ((D_T,E_T) , (D'_T,E'_T)\right )$$ est repr\'esentable  par un
$S$-sch\'ema de type fini.  On peut proc\'eder en deux \'etapes. Soit
d'abord le foncteur $\Isom_S (D,D')$ qui est  repr\'esent\'e un
torseur $\Cal I$ sous le groupe $\PGL (1)$. Si $\Phi: D\times _S {\Cal
I} \buildrel\sim\over\rightarrow D'\times_S {\Cal I}$ est l'isomorphisme
universel, on est ramen\'e visiblement  \`a prouver que le
foncteur  d\'efini sur la cat\'egorie des $\Cal I$-sch\'emas 
$
(U\to {\Cal I}) \mapsto \Isom_U (E_U, \Phi^*(E')_U)$
est repr\'esentable, en fait  par un $\Cal I$ sch\'ema affine sur $\Cal I$.
C'est  un fait connu.  Identifions  maintenant  la marque des
objets de $\Cal P$. Il est clair que toute fibration en $\Bbb P^1$ de
base $S$ est localement (pour la topologie \'etale) de la forme $\Bbb
P^1_S$. Supposons donc $D = \Bbb P^1_S$. Ceci \'etant,  le  module
localement libre  $E$ \'etant  de type $V$,    quitte \`a localiser
si n\'ecessaire,  on peut le supposer de la forme   
$E =  {\Cal
O}_{\Bbb P^1_S}(n_1)\oplus \dots \oplus {\Cal O}_{\Bbb P^1_S} (n_r) =
p^* (V)$,   
o\`u  $p: \Bbb P^1_S \to \Bbb P^1_k$ est le morphisme de
changement de base.  Ce qui  montre que  la marque est bien $(\Bbb
P^1,V)$. Pour d\'ecrire le groupe des automorphismes de l'objet
marque,  rappelons que si un groupe alg\'ebrique $G$ agit sur une
vari\'et\'e $X$, l'action \'etant   $\mu: G\times X \to X$, une
$G$-lin\'earisation d'un faisceau coh\'erent $F$, est un isomorphisme
$\Phi:  p_2^* (F) \buildrel\sim\over\longrightarrow \mu^* (F)$
v\'erifiant une relation de cocycle bien connue (par exemple \cite {57},
voir aussi \S \ 6.1). En g\'en\'eral un faisceau  $F$ localement libre
de rang $n$  sur $X$ n'est pas $G$-lin\'earisable, cependant il le
devient si le groupe $G$ est agrandi convenablement.  De mani\`ere
pr\'ecise, l'argument rappel\'e en d\'ebut de preuve montre qu'il
existe un groupe alg\'ebrique $\tilde G$, et un morphisme de
groupes $\tilde G \to G$, tel que les $G$-lin\'earisations de $F$
correspondent de mani\`ere bijective aux sections de $\tilde G \to
G$.  En fait on applique la construction  du d\'ebut au couple de
faisceaux $(\mu^*(F), p_2^*(F))$, $\mu$ et $p_2$ \'etant
respectivement l'action et la projection $G\times X \to X$ , de
sorte que 
$$\tilde G = \Isom (p_2^* (F) , \mu^* (F))$$
 D'une mani\`ere
simplifi\'ee, les points de $\tilde G$ sont les couples $(g,\phi)$,
avec $g\in G$, et $\phi: F\buildrel\sim\over\rightarrow g^*(F)$. La
loi de composition  \'etant  
$$(g,\phi)(h,\psi) =
(gh,h^*(\phi)\psi) \tag (9.29)$$ 
Si on a pour tout $g\in G, \;g^*(F) \cong F$,
alors  ce groupe s'ins\`ere dans une extension $1\rightarrow \Aut (F) \rightarrow \tilde G \rightarrow
G\rightarrow 1.$ Appliqu\'e \`a la situation pr\'esente, cela donne la seconde partie du lemme.\hfill\break\qed 

Noter que l'extension (9.28)  n'est  en g\'en\'eral pas scind\'ee
du fait que ${\Cal O}(n)$ est $\PGL (1)$-lin\'earisable  si  et seulement
si $n$ est pair; il n'est  que $\GL (2)$-lin\'earisable.  Soit
$\alpha_g: {\Cal O}(1) \buildrel\sim\over\rightarrow g^* ({\Cal O} (1))$
la lin\'earisation tautologique de ${\Cal O} (1)$ qui induit une
lin\'earisation   de $V$, 
$ \oplus_{i=1}^r \,
\alpha_g^{\otimes n_i}:  V \buildrel\sim\over\rightarrow g^* (V)
\quad (g\in \GL (2)). $
Noter que $\lambda 1_2$ agit
diagonalement  sur $V$ par la matrice diagonale $   {\rm diag} (\lambda^{n_1},
\dots ,\lambda ^{n_r})$ On notera cependant  que la lin\'earisation
tautologique  de $V$ sous $\GL (2)$ se descend dans tous les cas  en
une  lin\'earisation sous $\GL (2) / \mu_d$, si $d $ est le pgcd des
$n_i$.   Le groupe  $H$ est ais\'e \`a d\'ecrire.  Supposons la
partition $(n_i)$ de $r=\sum_{i=1}^r n_i$,  prenant $p$ valeurs
distinctes $k_1 < \dots < k_p$,  $k_i$ apparaissant $\alpha_i$ fois.
Alors $V$ a une filtration canonique 
 $$0
\ne F^p \subset \dots \subset F^1 = E,\quad  F^j/F^{j+1} = {\Cal
O}_{\Bbb P^1_S}(k_j)^{\oplus ( \alpha_j)}\tag (9.30) $$ 
En fait si $V =
\oplus_{j=1}^p {\Cal O}_{\Bbb P^1}(k_j)^{\oplus (\alpha_j)} $, on a
$F_j =  \oplus_{m=j}^p {\Cal O}_{\Bbb P^1}(k_m) ^{\oplus (\alpha_m)} $.
En particulier un automorphisme de $E$ respecte filtration $F_.$.
Il est  alors imm\'ediat  que $H = \Aut (E)$ est  le produit
semi-direct $$ H = U \rtimes \prod_{j=1}^p \GL (\alpha_i) $$
d'un groupe unipotent $U$ par un produit de groupes lin\'eaires. Dans certains cas, $G$ est en fait un produit semi-direct.
Supposons   $\alpha_1=\dots \alpha_p = 1$. Soit une
relation de Bezout $1 = \sum_j m_j  { { k_j} \over d}$. Consid\'erons
le caract\`ere  de $H$ donn\'e par   $\chi :  H \rightarrow
\prod_{j=1}^p \GL (\alpha_i)\buildrel\psi\over\rightarrow G_m $
avec $\psi ((g_j)) = \prod_j  \det (g_j)^{m_j}$, et soit
$K $ le noyau de $\chi$.  On a la remarque \'evidente:

\proclaim Lemme 9.12.  Sous les conditions indiqu\'ees, le groupe $H$   est
le  produit semi-direct  $H =  K\rtimes \GL (2) /\mu_d$.  Si de plus
$d$ est pair, la suite (9.28) est scind\'ee.  

\dem
  Notons $\alpha: \GL (2) \to G$ la section induite par la
lin\'earisation canonique. Alors  les \'el\'ements de $K$ qui sont
aussi dans l'image de $\alpha$, sont les $  {\rm diag}(\lambda^{n_1},
\dots ,\lambda^{n_r})$, avec $\lambda^{\sum m_in_i} = \lambda^d =
1$.  Cette intersection est donc r\'eduite \`a l'identit\'e. Le
reste est clair.  Pour le point deux, notons que si les $n_i$ sont
tous pairs, alors on obtient une section de $G \to \PGL (1)$ par factorisation
de $g\in \GL (2) \mapsto \oplus_i\,  \alpha_g^{\otimes n_i} \det
(g)^{-n_i/2}$.\hfill\break\qed   
\bigskip  

  {\it 9.3.2.
Interlude: G\'eom\'etrie du discriminant  des formes binaires  }
  \bigskip
  
   L'espace affine des formes binaires admet une stratification naturelle, dont la
g\'eom\'etrie  est intimement li\'ee \`a
celle des sch\'emas de Hurwitz  classifiant les rev\^etements
cycliques de  $\Bbb P^1$.   Soit $N$ un entier fix\'e
dans toute cette section,  et soit $\mu = n_1\leq n_2\leq  \dots
\leq n_r$ une partition (fix\'ee) de poids $N$.  Une forme binaire
(non nulle) 
$$ a_0X^N + a_1X^{N-1}Y + \dots + a_NY^N \tag (9.31)$$
est repr\'esent\'ee par le point $(a_0, \cdots ,a_N) \in \Bbb
A_*^{N+1}$. On notera    $X_\mu = X_{(n_1,\dots,n_r)}$ la sous
vari\'et\'e   (localement ferm\'ee) de l'espace affine $ \Bbb
A^{N+1}$ des formes  binaires de degr\'e $N$,  form\'ee de celles
qui sont  le produit de  $r$ facteurs lin\'eaires  non
proportionnels,   de multiplicit\'es  respectives $n_1\leq n_2\leq
\dots \leq n_r$.  On parlera de $\mu$ comme du type combinatoire de
la forme  binaire.    Il est facile de  d\'ecrire $X_\mu$ comme un c\^one au dessus d'un produit d'espaces projectifs.   On \'ecrira  pr\'ef\'erablement  la partition $\mu$ sous la forme duale
$\mu = m_1^{k_1} \dots m_s^{k_s}$, signifiant que la fonction  $i\to
n_i$ prend les valeurs distinctes $m_1 < \dots < m_s$ respectivement
$k_1, \dots ,k_s$ fois. Une forme binaire de type $\mu$ peut alors
\^etre \'ecrite  sous la forme plus rigide 
$F = \prod_{i=1}^s
f_i^{m_i}, \,\,(f_i,f_j) = 1 \;{\rm si} \; i\ne j , \; \deg (f_i )=
k_i $ 
,  $f_i$ \'etant sans facteur multiple. La d\'ecomposition   est g\'en\'eriquement bien  d\'efinie \`a l'action du groupe diagonalisable  pr\`es 
 $$T_\mu = \{
\lambda = (\lambda_1, \dots ,\lambda_s), \;\prod_{i=1}^s \lambda_i^{m_i} =
1\}\tag (9.32)$$ 
Pour d\'ecrire  avec plus de  pr\'ecision  la
stratification  de $\Bbb A_*^{N+1},$ de strates les  sous-sch\'emas
$X_\mu$ ,    rappelons (voir par
exemple \cite {47})   quelques    propri\'et\'es du morphisme
de Vi\`ete.  Soit  $ Z_N =  (\Bbb A_*^2)^N$  la vari\'et\'e  des
$2N$-uples
  $$(u_1,v_1,\dots ,u_N,v_N)  , \quad  (u_i,v_i) \ne (0,0) , \;1\leq
i \leq N   \tag (9.33)$$
  Il y a  sur $Z_N$ une action \'evidente d'une part du groupe  $G_m^N$, et d'autre part du groupe
sym\'etrique  $S_N$. Il en  r\'esulte une action du produit
semi-direct   $ G_m^N \rtimes S_N$ sur  $  (\Bbb A_*^2)^N$,
telle que  $   Z_N / G_m^N \rtimes S_N = (\Bbb P^1)^N / S_N = \Bbb
P^N$. 
  Soit le sous-groupe $T_N  = \{ (t_1,\dots , t_N) \in G_m^N /
\prod_i t_i = 1\}$.  Les choses s'expriment mieux en termes de
c\^ones au dessus d'une vari\'et\'e. Rappelons que si $X$ est une
vari\'et\'e  normale, et si $\Cal L$ est un faisceau inversible
ample sur $X$, alors l'alg\`ebre
  $$R = \bigoplus_{d \geq 0} \; \Gamma (X , {\Cal L}^{\otimes d} ) \tag (9.34)$$
  est de type fini, et normale.  La vari\'et\'e affine  $\Spec (R)$
est le c\^one de base $X$ d\'efini par $\Cal L$. Si $X$ est lisse,  le sommet est le seul point
singulier, point singulier conique, et  il est imm\'ediat que le
fibr\'e  vectoriel $V({\Cal L}) \rightarrow X$,   induit une r\'esolution de la
singularit\'e conique, via le morphisme naturel $V({\Cal L}) \to
\Spec (R)$. Le diviseur
exceptionnel est la section nulle de la fibration, et son compl\'ementaire est   isomorphe au c\^one \'epoint\'e  $\Spec (R) - 0$.
Notons par ailleurs  que $\Spec (R) - 0 / G_m = X$, et $\Pic (\Spec R
- 0) =  {\Pic (X)}/{[\Cal L]}$.  Ainsi la vari\'et\'e not\'ee
$\Bbb A_*^{N+1}$ des formes binaires de degr\'e $N$ peut \^etre
identifi\'ee avec  le c\^one \'epoint\'e de base $\Bbb P^N$,
relatif \`a ${\Cal O}(1)$.
Rappelons que le morphisme  de Viete (relations coefficients-racines)  $v:  Z_N = (\Bbb A^2_*)^N \rightarrow \Bbb
A^{N+1}_*$ est d\'efini comme suit: 
$$ (u_1,v_1,\dots ,u_N,v_N)
\mapsto \prod_{i=1}^N (u_iX - v_iY) = a_0X^N + \dots +a_NY^N  \tag
(9.35)$$
Ce morphisme  permet de justifier l'identification
 classique suivante:
\proclaim Lemme 9.13.     Le quotient  $ X_N = Z_N / T_N$
est le c\^{o}ne \'epoint\'e de base $(\Bbb P^1)^N$ d\'efini par le
faisceau inversible ${\Cal O} (1,\dots,1)$, de plus    $Z_N / T_N \rtimes
S_N =  \Bbb A^{N+1}_*$.    

\qed

 Pour d\'ecrire  $X_\mu$  on forme 
d'abord le plongement "diagonal" de $(\Bbb A^2_*)^r$ dans $(\Bbb
A^2_*)^N$ i.e. 
$$(u_1,v_1,\dots ,u_r,v_r) \mapsto (\overbrace
   {u_1,v_1,\dots ,u_1,v_1} , u_2,v_2,\dots
,\overbrace    { u_r,v_r,\dots, u_r,v_r}) \tag
(9.36)$$
le premier couple r\'ep\'et\'e  $n_1$ fois, etc. On note
que l'image  de (9.36) est la sous vari\'et\'e des points fixes de
$S_{n_1} \times \dots \times S_{n_r}$.  Le sous groupe de $T_N\rtimes
S_N$ qui agit sur  cette image,  est le normalisateur de $S_{n_1}
\times \dots \times S_{n_r}$, donc est le groupe $T_\mu \rtimes
Norm_{S_N} (S_{n_1} \times \dots \times S_{n_r})$, avec  $T_\mu$
donn\'e par
$$T_\mu = T_N^{S_{n_1} \times \dots \times S_{n_r}} =
\{({t_1,\dots t_1} , t_2,\dots
,\dots,t_r)\}, \;t_1^{n_1}\dots t_r^{n_r} = 1\tag (9.37)$$
Noter que
cette action est libre g\'en\'eriquement,  mais pas libre partout. On
a utilis\'e  en passant   le lemme \'evident suivant:
\proclaim  Lemme 9.14. Soit un produit semi direct $G = H \rtimes K$, et soit
$\Delta \subset K$ un sous-groupe de $K$.  Le normalisateur de
$\Delta$ dans $G$ est $H^\Delta \rtimes N_K(\Delta)$.
 \hfill\break\qed

Le groupe diagonalisable $T_\mu$
n'est connexe que si $d = \pgcd (n_1,\dots,n_r) = 1$, sinon  $T_\mu
/ T_\mu^\circ =   {\Bbb Z} /{d\Bbb Z}$ . Noter qu'\`a cela
correspond   une action libre du groupe diagonalisable  $T_\mu \cong \ker
(G_m^r \rightarrow G_m) , \, (t_1,\dots,t_r) \mapsto t_1^{n_1}\dots
t_r^{n_r}$ sur $(\Bbb A^2_*)^r$ qui s'exprime par 
$$(t_1,\dots,t_r).(u_1,v_1,\dots,u_r,v_r) = (t_1u_1,t_1v_1,\dots
,t_ru_r,t_rv_r)  $$
et qui rend le plongement  (9.36) \'equivariant.  Il en d\'ecoule  une action r\'esiduelle  du
normalisateur  $N_\mu := Norm_{S_N} (S_{n_1} \times \dots \times
S_{n_r})$ sur   $Z_\mu = (\Bbb A^2_*)^r) / T_\mu$ .   Rappelons que
les valeurs distinctes prises par les $n_i$   sont  $m_i =
n_{k_1+\dots +k_i}, \;\;m_1 < \dots  < m_s$. On trouve pour
le normalisateur $N_\mu$   le wreath
product   
$N_\mu = (S_{m_1} \wr S_{k_1} )\times \dots \times
(S_{m_s} \wr S_{k_s}).  $ 
 En r\'esum\'e on a le fait
connu suivant:  

\proclaim Proposition 9.15. 1)  Le
quotient  $Z_\mu = (\Bbb A^2_*)^r) / T_\mu$ s'identifie au c\^one  \'epoint\'e au dessus de $(\Bbb P^1)^r$ relatif \`a ${\Cal L} = {\Cal O} (n_1,\dots, n_r)$.  \hfill\break  
2)  Le quotient $C_\mu = Z_\mu / N_\mu$
s'identifie au  c\^one \'epoint\'e au dessus de $\Bbb P^{k_1}
\times \dots \times \Bbb P^{k_s}$ d\'efini par le faisceau
inversible ${\Cal O} (m_1,\dots ,m_s)$.  Le morphisme $ C_\mu  \to
\Bbb A_*^{N+1}$ est fini birationnel sur son image, il d\'efinit  la
normalisation de la strate ferm\'ee $\overline {X_\mu}$. \hfill\break
 3)  La  vari\'et\'e $X_\mu$ est isomorphe au compl\'ementaire dans   $
C_\mu = Z_\mu / N_\mu$ de l'hypersurface  $\Delta (f_1. \dots .f_s)
= 0$ ($\Delta (-)$ \'etant le discriminant usuel)

 \dem  Pour 1), on note que les fonctions
invariantes sous $T_\mu$ sont les fonctions polynomiales
$f(u_1,v_1,\dots,u_r,v_r)$ telles que  pour $t\in T_\mu$:
$$f(t_1u_1,t_1v_1,\dots,t_ru_r,t_rv_r)  =
f(u_1,v_1,\dots,u_r,v_r)$$ 
  si $f$ est multihomog\`ene de
multidegr\'e $(\alpha_1,\dots,\alpha_r)$, le  poids  de $f$ est $t
\mapsto t_1^{\alpha_1}\dots t_r^{\alpha_r}$, donc on doit avoir
$(\alpha_1,\dots ,\alpha_r)  = e(n_1,\dots,n_r)$ pour un $e\geq 1$.
L'assertion 1) en d\'ecoule. D'une autre mani\`ere, on peut noter que
la projectivisation du plongement (9.36), est le plongement
polydiagonal    $(\Bbb P^1)^r \mapsto (\Bbb P^1)^N$   et que la
restriction de ${\Cal O}(1,\dots,1)$ est ${\Cal O}(n_1,\dots,n_r)$.
L'assertion d\'ecoule alors du lemme 3. \hfill\break    Pour 2), on
note que  l'action du groupe fini $N_\mu$,  se r\'eduit \`a
l'action de  $S_{k_1}\times \dots \times S_{k_s}$, ce qui conduit
au fait que les points de $Z_\mu / N_\mu$  ont pour coordonn\'ees
(multihomog\`enes), les coefficients des produits partiels
$$f_1(X,Y) = \prod_1^{k_1} (u_iX - v_iY) = a_0^1X^{k_1} + \dots
a_{k_1}^1 Y^{k_1} \;,  \dots     \tag (9.38)$$  
 On peut identifier
un point de ce quotient \`a un $s$-uple de formes  $[f_1,\dots
,f_s]$, de degr\'es respectifs $k_1,\dots,k_s$.  L'identification \'etant    
$$[f_1,\dots,f_s] = [g_1,\dots,g_s] \Longleftrightarrow
g_i = \lambda_i f_i \;\;(\lambda_i \in k^*), \prod_{i=1}^s
\lambda_i^{m_i} = 1 \tag (9.39)$$   
Comme dans la description  pr\'ec\'edente, on retrouve  bien la description des points du c\^one de
base   $\Bbb P^{k_1} \times \dots \times \Bbb P^{k_s}$  relatif au
faisceau ${\Cal O }(m_1,\dots,m_s)$.  D'une autre mani\`ere le
quotient $Z_\mu/N_\mu$ est isomorphe \`a $\Bbb A_*^{k_1+1}\times \dots
\times \Bbb A_*^{k_s+1} /  T'$, le tore $T'$ \'etant le noyau de
$G_m^s \rightarrow G_m , \;(\lambda_1,\dots,\lambda_s) \mapsto
\prod_{i=1}^s \lambda_i^{m_i}$; l'identification s'obtient  en
effectuant  le quotient partiel de $Z_\mu$ par le  sous tore  de
$T_\mu$ donn\'e par les \'equations $\prod_{i=1}^{k_1} t_i = 1,
\dots,\prod_{i=1}^{k_s} t_i = 1$.    Le morphisme  $Z_\mu / N_\mu
\rightarrow \Bbb A_*^{N+1}$ n'\'etant autre que $(f_1,\dots , f_s)
\mapsto f_1^{m_1}\dots f_s^{m_s}$.  Tout est plus clair sur le
diagramme   
$$\harrowlength=40pt    \commdiag{    Z_\mu =
(\Bbb  A_*^2)^r / T_\mu &  \mapright  & Z_\mu/N_\mu & \mapright^\Phi & \Bbb A_*^{N+1} \cr     \mapdown & &
 \mapdown & &  \mapdown \cr   
  (\Bbb P^1)^r &  \mapright &
\Bbb P^{k_1}\times \dots \Bbb P^{k_s} & \mapright^\phi & \Bbb P^N \cr} $$   
 Les fl\`eches  $\Phi, \phi$
sont celles explicit\'ees  lors de la preuve,  c'est \`a dire
$$\phi(D_1,\dots,D_s) = \sum_{i=1}^s m_iD_i, \;\; \Phi
([f_1,\dots,f_s])  = [f_1^{m_1}\dots f_s^{m_s}] \tag (9.40)$$
\qed

 \rema{9.2} 
     Les relations d'incidence entre les strates $X_\mu \subset {\Bbb A}_*^{N+1}$ se d\'ecrivent  facilement
au moyen d'un ordre  entre partitions de poids  donn\'e $N$ \cite {47}:
$\eta = (\eta_1\leq \dots \leq \eta_l)  \leq \mu = (n_1\leq \dots
\leq n_r) \Longleftrightarrow \eta_j = \sum_{i\in I_j} n_i $,     les
$I_j$ formant une partition de $[1,r]$. L'adh\'erence de la
strate $X_\mu$ est     $\overline {X_\mu} = \bigcup_{\eta \leq \mu}
\; X_\eta$.    Le morphisme $C_\mu = Z_\mu /N_\mu \rightarrow \Bbb
A_*^{N=1}$ a pour image $\overline {X_\mu} $,   est  la d\'esingularisation  de $\overline {X_\mu} $.
Apr\`es projectivisation, ce morphisme se r\'eduit \`a     $\Bbb
P^{k_1} \times \dots \times \Bbb P^{k_s} \longrightarrow \Bbb P^N$
qui est le compos\'e d'un plongement de Veronese, suivi d'une
projection lin\'eaire; ce morphisme est la normalisation de son
image.   
 \exa{9.3}   Le cas $\mu =
(1,\dots,1,2)$,  soit $k_1 = N-2, k_2 = 1, m_1 = 1, m_2 = 2$,
correspond a l'hypersurface discriminant (priv\'ee de l'origine),
qui a pour normalisation le c\^one \'epoint\'e de base $\Bbb
P^{N-2}\times \Bbb P^1$, relativement \`a ${\Cal O}(1,2)$.  Le
morphisme de normalisation est $[f_1,f_2] \mapsto F = f_1f_2^2$. Si
$f_1$ a un facteur lin\'eaire multiple, on peut l'\'echanger avec
$f_2$, ce qui montre que ce morphisme n'est pas bijectif sur l'image.
     Regardons comme exemple  le cas des formes quartiques  (9.31) $(N =
4).$        Le discriminant $\Delta$   caract\'erise les formes
avec un facteur lin\'eaire double,   donc $\mu = (1,1,2)$. Pour
  $\mu = (1,3)$,   $\overline X_{1,3}$ est    le c\^one  \'epoint\'e au dessus de $\Bbb P^1 \times \Bbb P^1$ d\'efini par
${\Cal L} = {\Cal O}(1,3)$.  La th\'eorie classique  donne des \'equations pour $\overline X_{1,3}:     \,P(F) = Q(F) = 0$,     si on
note $P$ et $Q$ les covariants      $$P =   {1 \over 6}
(a_2^2-3a_1a_3+12a_0a_4), \;\;Q = a_0a_2a_4- { 3 \over 8} a_0a_3^2 -
 { 3 \over 8} a_1^2 a_4+   1 8 a_1a_2a_3-  { 1\over  36} a_2^3$$     Dans
le cas  $\mu = (2,2)$, donc $m_1 = 2, \;s=1, \; k_1 =  2$, $X_{2,2}$
repr\'esente les formes  quartiques carr\'e d'une forme
quadratique; dans ce cas  aussi on obtient  des \'equations en \'ecrivant que la forme hessienne $H(F)$ est proportionnelle \`a $F$.
Enfin si $\mu = (4)$ qui est le cas  des formes avec un seul facteur
lin\'eaire, les \'equations  de $X_4$  qui est la strate ferm\'ee, se
r\'esument  comme on le voit facilement \`a $ H(F) = 0$.

\bigskip    
  {\it  9.3.3.  Le champ des rev\^etements
cycliques de la droite projective}   
  
   \bigskip
On fixe la ramification d'un rev\^etement cyclique de degr\'e $n$
(\S \ 9.1)   sous la forme  d'une suite de diviseurs de $n$, $\;1 < e_i
\vert n \;(1\leq i\leq r)$, et pour tout $i$, un entier $1\leq \eta_i
< e_i$, tel que $\pgcd (e_i,\eta_i) = 1$; on posera   $n_i =  { n\over 
{e_i}} \eta_i$.  On supposera que $n_1 \leq \dots \leq n_r$, d\'efinissant ainsi une partition $\mu$ de poids $N = nm = \sum_i n_i$.
          
          Soit   $\pi: C\to D $  un rev\^etement cyclique \`a
ramification fix\'ee de type $\mu$, de base $S$,  les fibres g\'eom\'etriques de $D\to S$ \'etant  isomorphes \`a  $\Bbb P^1$.   On sait
(Proposition 9.4) que le rev\^etement $\pi:C \to D$ peut \^etre
reconstruit en  partant de la base $D$, \'equip\'ee de la structure
additionnelle $({\Cal L},\Phi)$.     Il sera pr\'ef\'erable dans la suite de  noter la partition $\mu$ comme une suite
strictement croissante  avec r\'ep\'etitions, donc  de consid\'erer
les entiers $m_j \;(1\leq j\leq s)$ (distincts)  d\'efinis par la
partition $n_1\leq \dots \leq n_r$ comme  cela a \'et\'e fait dans
la section 9.  Cela permet d'\'ecrire   le diviseur $B$ tel
que $\Div (\Phi) = B$  sous la forme $B = \sum_{i=1}^s m_i B_i$, les
diviseurs $B_i$ \'etant  maintenant disjoints deux \`a deux. Posons
$\deg (B_i) = k_i$, de sorte que $\sum_i k_im_i = \sum_j n_j = nm$,
l'entier $m$ \'etant d\'efini par cette \'equation.  

Revenons \`a  la vari\'et\'e lisse $X_\mu$ partie ouverte du c\^one
$Z_\mu/N_\mu$. Le groupe $\GL (2)$ agit de mani\`ere naturelle sur
l'espace vectoriel des formes binaires de degr\'e donn\'e, i.e.
$g.F(X,Y) = F(g^{-1}X,g^{-1}Y)$, mais cette action ne se descend pas
en g\'en\'eral  en une action de $\PGL (1)$. L'entier $m$ \'etant
comme introduit au dessus,  on notera  cependant qu'il y  a une
action du groupe $\GL(2)/\mu_m$ sur le c\^one \'epoint\'e
$Z_\mu/N_\mu$, de mani\`ere explicite
 $$g[f_1,\dots,f_s] =
[g.f_1,\dots,g.f_s] \tag (9.41)$$
Cette d\'efinition est coh\'erente
du fait que  $\sum_{i=1}^s k_im_i = nm$. De mani\`ere \'equivalente,
le faisceau ${\Cal O}(m_1,\dots,m_s)$ sur $\Bbb P^{k_1}\times \dots
\times \Bbb P^{k_s}$ est canoniquement $\GL (2)/\mu_m$-lin\'earis\'e. La description qu'on a en vue du champ ${\Cal H}_{n,\mu}$ des rev\^etements
cycliques de degr\'e $n$ de $\Bbb P^1$, \`a ramification
de type $\mu$, et qui dans le cas des courbes g\'en\'eralise \cite
{4},  est:

\proclaim Th\'eor\`eme 9.16.  Les points de branchement \'etant marqu\'es (non index\'es), on a:
$${\Cal H}_{n, \mu}  = [X_\mu /(\GL(2)/\mu_m)] \tag (9.42)$$
Si les points de branchement sont piqu\'es (ordonn\'es) on a, pour une action de $G_m$ qui sera pr\'ecis\'ee
dans la preuve, et $Y_\mu$ \'etant le compl\'ementaire dans $Z_\mu$de
la grosse diagonale  $ {\Cal H}_{n, \mu}  = [Y_\mu / G_m].  $

\dem  On commence par d\'efinir un
morphisme  ${\Cal H}_{n,\mu} \to  [X_\mu / (\GL(2)/\mu_m)]$, qui apr\`es v\'erification, sera un isomorphisme. Soit un rev\^etement $p: C\to
D$, de base $S$. La courbe $q: D\to S$ est un fibr\'e en coniques,
et le sous-faisceau propre ${\Cal L} \in \Pic (D)$,  de $\pi_* ({\Cal
O}_C)$,  relatif \`a valeur propre $\zeta_n$,  est de degr\'e $-m$
le long des fibres. Par ailleurs la multiplication d\'efinit un
morphisme injectif  $\Phi:  {\Cal L}^{\otimes n} \to {\Cal O}_D$,  le
diviseur $B$  de $\Phi$ ayant par d\'efinition  un type combinatoire
fix\'e le long des fibres.  Ce qu'on traduit par   
$$B = \Div
(\Phi) = \sum_{1<e\vert n, \; 1\leq \eta <e} \; {n\over e} \eta
B_{e,\eta} \tag (9.43)$$
 On r\'ecup\`ere  d'une
part un torseur    $P =  \Isom ((\Bbb
P^1_S , V_S) , (D\to S , L))$ sous $G = \GL(2)/\mu_m$, avec  $V = {\Cal O}(-m)$.  Du fait de la
forme du diviseur de la section $\Phi$, cette section d\'efinit
d'autre part un morphisme  $P \to X_\mu \subset \Bbb A^N$, qui est
par  construction  clairement $G$-\'equivariant.  Cette construction
donne donc un foncteur   ${\Cal H}_{n,\mu} \rightarrow [X_\mu/G] $
Pour prouver que c'est un isomorphisme,  il suffit de voir que
c'est simultan\'ement un monomorphisme et un \'epimorphisme \cite
{51}. Il est possible en fait d'exhiber un foncteur quasi-inverse,
comme  la preuve de dessous le sugg\`ere. Pour tout sch\'ema $S$,  on
doit s'assurer que la restriction de ce foncteur \`a la  cat\'egorie  des sections au-dessus  de $S$ est pleinement fid\`ele.  Cela
se ram\`ene \`a voir que tout  $S$-isomorphisme de rev\^etements
$f: C\buildrel\sim\over\rightarrow C'$ correspond de mani\`ere
bijective \`a un isomorphisme des objets images par le foncteur
consid\'er\'e. Notons $h: D\buildrel\sim\over\rightarrow D'$ le
$S$-isomorphisme entre les courbes quotients. Il d\'efinit d'une part
un isomorphisme entre les  $G$-torseurs associ\'es
$P\buildrel\sim\over\rightarrow P'$, et un isomorphisme $\  {\Bbb
Z}/{n\Bbb Z}$-\'equivariant 
$$\psi: {\Cal  E} = \oplus_{i=0}^{n-1}
{\Cal L}_i \buildrel\sim\over\longrightarrow h^* ({\Cal E}'  =
\oplus_{i=0}^{n-1} {\Cal L}'_i )\tag (9.44)$$
Notons $\psi_i: {\Cal L}_i
\buildrel\sim\over\rightarrow h^*({\Cal L}'_i)$ sa restriction \`a
${\Cal L}_i$.  Rappelons la relation $\Phi' . \psi_1^{\otimes n} = \Phi
$; il s'agit donc de voir que l'isomorphisme initial $f$ peut \^etre
reconstitu\'e demani\`ere unique en partant de la donn\'ee
$(h,\psi_1)$.  C'est exactement le contenu de la Proposition 9.4.

Dans le cas des points de branchement ordonn\'es la donn\'ee qui se substitue \`a ($\Bbb P^1, {\Cal O}(-m)$)  
est $(\Bbb P^1, (Q_1,\dots,Q_b),{\Cal O}(-m))$, le groupe $G$ est donc $G =
G_m = \Aut ({\Cal O}(-m))$. Il est \'evident que cela revient dans la
preuve pr\'ec\'edente \`a remplacer $X_\mu$ par $Y_\mu$, dont les
points correspondent \`a la donn\'ee de $b$ points distints
$Q_1,\dots,Q_b$ de $\Bbb P^1$, et d'une forme binaire $f$ de
diviseur des z\'eros $\Div (f) = \sum_j n_jQ_j$. L'action de $G_m$
d\'erive de "l'action"  $\lambda [(u_1,v_1),\dots,(u_b,v_b)] =
[\lambda^{1/m}(u_1,v_1),\dots,\lambda^{1/m}(u_b,v_b)]$. Le reste est
clair.\hfill\break \qed 

 \exa{9.4}  Si $n=3$ (voir \cite {4}),
rev\^etements  galoisiens triples  de la droite projective,
 le diviseur $B$ est de la forme  $B = B_1+2B_2$,  avec $B_1$ et
$B_2$ disjoints et sans multiplicit\'es; on a  entre les degr\'es
respectifs la relation $k_1+2k_2 = 3m$. On a en fonction de la
prescription g\'en\'erale:
$\Div (\Phi_2) = \left[ {2B\over
3}\right] =  2B_1+B_2$. On voit que $X_\mu$ est dans ce cas un ouvert
du c\^one de base $\Bbb P^{k_1}\times \Bbb P^{k_2}$ d\'efini par
${\Cal O}(1,2)$. Dans \cite {4}  Arsie et Vistoli   ont trait\'e le
cas g\'en\'erique  $\mu = (1^N)$; les rev\^etements  correspondants
sont appel\'es par eux uniformes.  Dans  cette situation, on a
$s=1$, et $X_{(1^N)}$ est  simplement le c\^one des formes de degr\'e
$N$ sans racine multiple.

 \bigskip

\section{10}  {Groupe de Picard et classes tautologiques} 

\bigskip

La position du champ de Hurwitz  comme \'etablissant une
correspondance entre champs de courbes marqu\'ees (ou piqu\'ees)
$$\commdiag{ \overline { \Cal M}_{g,r} 
  \mapleft^{\imath  } 
   \overline
{\Cal H}_{g,G,\xi}  \mapright^{ \delta}   \overline{\Cal M}_{g',b} \cr} \tag (10.1)$$ 
 sugg\`ere  l'existence sur ce champ    de
fibr\'es tautologiques li\'es \`a la sp\'ecificit\'e des rev\^etements galoisiens. 

On notera pour tout groupe fini   
  $G$, $(- , -)_{\vert G}$ le produit scalaire de deux caract\`eres   de $G$,  et   $R(G)$ l'anneau
des repr\'esentations (ou des caract\`eres) de $G$.  Si $v\in \hat G$,
on notera $V_v$ l'espace vectoriel qui supporte la repr\'esentation
de caract\`ere (irr\'eductible) $v$. Pour tout r\'eel $x$,  on posera
$\langle x \rangle = x - [x]$.

\bigskip

\beginsection 10.1. Fibr\'e de Hodge

\bigskip

  {\it 10.1.1.  $G$-fibr\'es vectoriels sur le champ de
Hurwitz} 

\bigskip

  Soit le  morphisme canonique
(repr\'esentable) oubli de l'action de $G$   
$\imath: \overline
{\Cal H}_{g,G,\xi}\longrightarrow \overline {\Cal M}_{g,(r)} $ 
o\`u comme  pr\'ecedemment, $r$ d\'esigne  le nombre de points de ramification, et $b$ le nombre de points de branchement.
  Si $\pi: C\to D$ est un rev\^etement on supposera    la
base  piqu\'ee par $b$-points contenant les points de branchement.   
 La courbe $C$ est  marqu\'ee par
 les pr\'eimages  des points de
branchement. Le marquage de ces points par paquets se d\'eduisant de la num\'erotation des points
de branchement.    Rappelons qu'un  
faisceau  sur  $\overline{\Cal H}_{g,G,\xi}$ est   l'assignation pour tout rev\^etement $q:
C\buildrel\pi\over\rightarrow D\buildrel p \over\rightarrow S $ d'un faisceau coh\'erent $E (\pi)$ sur la base $S$, de mani\`ere compatible avec les morphismes. Soit  la remarque: 
\proclaim  Lemme 10.1. Soit $F$ un faisceau coh\'erent
(resp. un fibr\'e vectoriel) sur $\overline {\Cal M}_{g,r}$. Le
faisceau $\imath^* (F)$ sur $\overline {\Cal H}_{g,G,\xi}$ est un
$G$-faisceau coh\'erent (resp. un $G$-fibr\'e vectoriel),
relativement \`a l'action triviale de $G$. Dans le second cas  il
se d\'ecompose en la somme directe de ses composants isotypiques
$$\imath^* (F) = \bigoplus_{\chi\in \widehat G}  F_\chi \otimes
V_\chi \tag (10.2)$$
les $F_\chi$ \'etant des faisceaux  localement
libres sur $\overline {\Cal H}_{g,G,\xi}$, et la  somme directe \'etant index\'ee par les repr\'esentations irr\'eductibles de $G$.

 \dem  Si $F$ est un faisceau sur $\overline {\Cal
M}_{g,r}$, donc d\'efini par une collection de faisceaux coh\'erents
$F(q) = F(C\buildrel q\over \rightarrow S)$, on a  pour tout
morphisme de $\overline{\Cal M}_{g,r}$, i.e.  un carr\'e cart\'esien 
$$\commdiag{ C & \mapright^f &C'\cr \mapdown\lft{q}&&\mapdown\lft{q'}\cr   S  &\mapright^h& S' \cr} $$
un
isomorphisme $\phi_{f,h}: F(q) \buildrel \sim\over \longrightarrow
h^* (F(q')),  $
les $\phi_{f,h}$ \'etant assujettis \`a
satisfaire \`a une relation de cocycle  $\phi_{f'f,h'h} =
h^*(\phi_{f',h'}). \phi_{f,h}. $    Si maintenant $q: C\to S$
provient  d'un $G$-rev\^etement $\pi: C\to D$, alors tout \'el\'ement
$g\in G$ d\'efinit un automorphisme de $C\to S$ dans $\overline {\Cal
M}_{g,r}$, mais pas dans le champ de Hurwitz, sauf si $g$ est dans le
centre. Il en r\'esulte un isomorphisme 
$\phi (g): F(q)\buildrel \sim\over\longrightarrow F(q)$ qui du fait  de la relation de
cocycle repr\'esente une $G$-lin\'earisation sur le faisceau image
r\'eciproque $\imath^*(F)$ d\'efini par $\imath^* (F) (\pi) =
F(q)$. Si $F$ est  de plus localement libre,  il se d\'ecompose  en facteurs
isotypiques
$$F(q) = \bigoplus_{\chi\in \widehat G}  F_\chi (q)
\otimes V_\chi \tag (10.3)$$
avec $F_{\chi}(q) = {\Cal H}om_G (\Cal O_S
\otimes V_\chi , F)$. Comme les morphismes de $\overline {\Cal
H}_{g,G,\xi}$  sont $G$-\'equivariants, il en d\'ecoule une action de
$G$ sur $\imath^* (F)$, en d'autres termes les facteurs $F_\chi (q)$
d\'efinissent un fibr\'e vectoriel sur le champ de Hurwitz. De la
sorte la d\'ecomposition qui pr\'ec\`ede induit donc la d\'ecomposition en sous-faisceaux propres (10.1).  \hfill\break\qed  

 Un $G$-fibr\'e vectoriel sur le champ
de Hurwitz  aura pour sens   un fibr\'e donn\'e par une d\'ecomposition
(10.2). Un $G$-morphisme entre deux $G$-fibr\'es vectoriels est un
morphisme qui commute aux actions respectives de $G$. Soient $F$ et $H$
deux $G$-fibr\'es vectoriels sur le champ de Hurwitz. Nous dirons qu'ils
sont disjoints si une m\^eme repr\'esentation irr\'eductible $V_\chi$
n'appara\^\i t pas simultan\'ement dans $F$ et dans $H$. Notons le
fait \'evident suivant:

\proclaim  Lemme 10.2. Supposons avoir une
suite exacte   $0\rightarrow E \rightarrow F \rightarrow H
\rightarrow 0$, les faisceaux \'etant des $G$-fibr\'es vectoriels, et
les morphismes  des $G$-morphismes. Si $E$ et $H$ sont disjoints,
alors $F = E \oplus H$.  

\dem  L'hypoth\`ese signifie que
pour une repr\'esentation irr\'eductible $V_\chi$ donn\'ee, on a soit
$E_\chi = F_\chi$, soit $H_\chi = F_\chi$. Le r\'esultat en d\'ecoule. \hfill\break\qed 

\bigskip

  {\it 10.1.2. D\'ecomposition du fibr\'e de Hodge}
 
\bigskip

Rappelons \cite {40}
que le fibr\'e de Hodge  $\Bbb E_{g',b}$ (resp.$ \Bbb E_{g,r}$) sur
$\overline{\Cal M}_{g',b}$ (resp.  $\overline{\Cal M}_{g,r}$) est
le faisceau dont les sections sur l'objet $p: D \rightarrow S$
(resp. $q: C\to S$), courbe stable marqu\'ee de genre $g'$ (resp.
$g$), est le faisceau localement libre  de rang $g'$ (resp. $g$) sur
$S$  
$$\Bbb E (p) = p_*(\omega^1_{D/S}) \tag (10.4) $$  
Notons que
par dualit\'e  de Serre ce faisceau est dual de $\R^1q_* ({\Cal
O}_D)$, et que le morphisme trace permet une identification $R^1q_*
(\omega_{D/S}) = \Cal O_S$. La classe de Hodge $\lambda \in \Pic
(\overline {\Cal M}_{g',b})$ est par d\'efinition $\lambda = \det
(\Bbb E) = c_1(\Bbb E)$; on d\'efinit plus g\'en\'eralement
$\lambda_i = c_i (\Bbb E)$. Pour \'eviter toute confusion, on notera
dans la suite $\Bbb E_{bs}$ (resp. $\lambda_{bs}$) le fibr\'e vivant sur la base.

Soit un rev\^etement galoisien stable de groupe $G$, $\pi: C
\rightarrow D$  de base $S$,  $p: D \rightarrow  S$ et $q: C
\rightarrow S$ les morphismes structuraux.   Sous
ces conditions  les faisceaux   $q_* \omega_{C/S}$  et
$\R^1q_* ({\Cal O}_C)$ sont des $G$-faisceaux  pour l'action triviale de $G$,  localement libres, et en dualitŽ.  Le faisceau localement
libre de rang $g = g_C$ (genre de  $C$) $q_* \omega_{C/S}$ d\'efinit un
$G$-faisceau sur le champ  $\overline{\Cal H}_{g,n}^\xi$ que nous
noterons  dans la suite ${\Bbb E}_{tot}$.  Ce fibr\'e
vectoriel est exactement l'image r\'eciproque du fibr\'e de Hodge
sur $\overline {\Cal M}_{g,r}$, comme explicit\'e en (10.1.1). Le
fibr\'e de Hodge sur $\overline{\Cal M}_{g',b}$ \'etant not\'e
$\Bbb E_{bs}$, on veut comparer $\Bbb E_{tot}$ et $\delta^* (\Bbb
E_{bs})$. Soit  d\'ecomposition en facteurs isotypiques  (10.1):
$${\Bbb E}_{tot} \;= \; \bigoplus_{v\in \hat G} \; {\Bbb E}_v \otimes
V_v \tag (10.5)$$
La fibre de cette d\'ecomposition en  $C\to D$  est simplement la d\'ecomposition en facteurs irr\'eductibles de $H^0(C,\omega_C)$. La formule de Chevalley-Weil  (\S \ 3.2) permet de pr\'eciser le rang de chaque facteur, la
donn\'ee de ramification   $\xi = \sum_{i=1}^b [ H_i , \chi_i]$ \'etant fix\'ee.
  \proclaim Proposition 10.3. Soit la d\'ecomposition du  fibr\'e de Hodge ${\Bbb E}_{tot}$
  en facteurs isotypiques
${\Bbb E}_{tot} \;=\; \bigoplus_{v \in \hat G} \;{\Bbb E}_v \otimes
V_v $.  On a  en particulier  ${\Bbb E}_1 = \delta^* ({\Bbb E}_{bs})$. Le
rang   de ${\Bbb E}_v$ est donn\'e  si $v \ne 1$ par la formule de
Chevalley-Weil 
$$rg ({\Bbb E}_v) \;= dim\;(v)\left( g'-1 + \sum_i (1 -  { 1\over
{e_i}}) \right) \;-
\;\sum_i\sum_{l=1}^{e_i-1} ( {{e_i-l} \over {e_i}}) (v^\vee ,
\chi_i^{l-1})_{\vert H_i}  \tag (10.6)$$
Si $G$ est cyclique\footnote {  Par exemple si $v =  n-1$, on trouve
pour le rang de la composante correspondante du fibr\'e de Hodge
$g'-1+\sum_{i=1}^b  { {\nu_i} \over{e_i} }= g' - 1 +  { {\sum_i m_i\nu_i} \over n}$.} d'ordre $n$, (10.6) se sp\'ecialise en ($0 \leq \langle - \rangle < 1$ d\'esignant la partie
fractionnaire):
$$rg ({\Bbb E}_v) \;=\; g' - 1  +\sum_{i=1}^b \! \langle  {
{-v\nu_i} \over {e_i} }\rangle \tag (10.7)$$

  \dem  Il suffit de se placer sur $S = \Spec k$. Si $\pi:
C\rightarrow\Sigma$ est un rev\^etement du type indiqu\'e d\'efini
sur    $k$ ($\vert G\vert \ne 0 \in k$), il suffit, pour avoir
l'expression (10.6), d'expliciter pour tout $v$ le rang du facteur
isotypique $H^0(C , {\Cal O}_C)_v$;  de mani\`ere \'equivalente,
trouver la trace de Lefschetz  (voir \S 3.2)
$$L_G (\Omega_C) = [H^0(C , \Omega_C)] \;-\;[H^1(C , \Omega_C)] \in R(G)$$
Si $\pi_* (\Omega_C) = \bigoplus_{v\in \hat G} \;E_v \otimes V_v$,
cette derni\`ere s'exprime par
$L_G (\Omega_C) \;=\; \sum_v \;\chi (E_v) [V_v] \in R(G).$
Cette expression  montre qu'il suffit de trouver le degr\'e de
$E_v$ et d'appliquer alors le
th\'eor\`eme de Riemann-Roch (voir \S \ 3.2). Le degr\'e s'obtient
facilement au moyen de la suite exacte
$$0\rightarrow \pi^* (E_v) \rightarrow \Omega_C \otimes V_v^\vee
\rightarrow {\Cal Q}_v \rightarrow
0$$
 $ {\Cal Q}_v$ \'etant concentr\'e sur les orbites de  points
de ramification,   plus
pr\'ecis\'ement    ${\Cal Q}_v = \bigoplus_i
\;\Ind_{H_i}^G \;\left(({\Cal Q}_v)_{P_i}\right)$    o\`u $P_i$
est un quelconque point d'une l'orbite de type $[H_i , \chi_i]$.  De
cette observation on tire facilement (voir \S \ 3.2 pour un calcul
 analogue)
$$ \deg\;(E_v)\;= \; { {dim\;(v)\;(2g_C - 2)}\over  {\# G}}\; - \;
\sum_{i=1}^b \sum_{l=1}^{e_i-1}
{  {e_i - l}\over {e_i}} \;(v^\vee \;,\; \chi_i^{l-1})_{\vert H_i}$$
  La formule de Riemann-Hurwitz donne $  {{2g - 2}\over {\# G} } = 2g' -
2 + \sum_{i=1}^b (1 -  { 1\over
{e_i}})$,  d'o\`u par substitution dans la formule de Riemann-Roch
$$\chi (E_v) = \dim \;(v)\;\left( g'-1 + \sum_i (1 -  { 1\over {e_i}})
\right) \;- \;\sum_i
\sum_{l=1}^{e_i-1} ( { {e_i-l}\over {e_i}}) (v^\vee , \chi_i^{l-1})_{\vert H_i} $$
Si maintenant $G$ est cyclique  d'ordre $n$ engendr\'e par
$\sigma$, alors pour tout entier $0 \leq v < n$,  $v$ d\'esignant
la repr\'esentation de caract\`ere $\chi_v (\sigma) = \zeta_n^v$, on
a  $({\chi_v}^\vee , \chi_i^{l-1})_{\vert H_i} = 0$ sauf si $l-1$ est le
reste modulo $e_i$ de $ - v\nu_i$. D'o\`u  le
r\'esultat. 
 Pour confirmer la validit\'e de  
l'expression  (10.7), v\'erifions que $\sum_v \dim (v) rg (\Bbb E_v) = g_C =
g$. La somme sur les caract\`eres irr\'eductibles de la partie qui
suit le signe moins, donne $$\sum_{i,l} { {e_i-l}\over {e_i}} (\sum_v
(\dim v) v^\vee , \chi_i^{l-1})_{H_i}= \sum_{i,l} (\chi_{\Bbb
C[G]} , \chi_i^{l-1})_{H_i}=\sum_{i,l} ( { {e_i-l}\over  {e_i}})
 { {\# G}\over {e_i}} = \# G \sum_i  {e_i(e_i-1) \over 2}$$La premi\`ere somme se simplifie du fait de $\sum_v \dim (v)^2 = \# G$; les
choses \'etant mises ensemble, la formule suit.\hfill\break \qed

\exa   {10.1}  Soit le cas $n=3$, i.e les rev\^etements triples galoisiens de genre $g$ de $\Bbb P^1$. Dans ce cas une donn\'ee de ramification \'equivaut \`a une partition des points de branchement $[1,g+2] = \Lambda_1\sqcup \Lambda_2$, avec $\vert \Lambda_1\vert \equiv \vert\Lambda_2\vert \pmod 3$  (voir \S\  2.3).  Le fibr\'e de Hodge est une somme directe
$\Bbb  E = \Bbb E_1\oplus \Bbb E_2$
o\`u 
$$rg (\Bbb E_1) ={\vert \Lambda_1\vert + 2\vert\Lambda_2\vert \over 3} - 1,\, \,\, rg (\Bbb E_2) ={2\vert \Lambda_1\vert + \vert\Lambda_2\vert \over 3} - 1\tag (10.8)$$
$\lozenge$
 
Il y a une  construction
duale.  Soit un rev\^etement galoisien
stable $\pi: C \rightarrow D$ de groupe $G$. Le faisceau $\pi_* ({\Cal
O}_C)$ est un $({\Cal O}_D , G)$ module sans torsion de rang $\# G$,
de formation
compatible \`a tout changement de base. Il se d\'ecompose  en facteurs isotypiques (section 4)
$$\pi_* ({\Cal O}_C) \;=\; \bigoplus_{v\in\hat G} \;{\Cal L}_v
\otimes V_v \tag (10.9)$$
Il est imm\'ediat  que ${\Cal L}_1 = \pi_* ({\Cal O}_C)^G = {\Cal
O}_D$, et que le rang de ${\Cal
L}_v $ est \'egal \`a $\dim (v)$.
\proclaim Proposition 10.4. On a  $p_* ({\Cal O_C}) = {\Cal O}_S$,
et  si $v\ne 1$,  on a $\;p_* ({\Cal
L}_v ) = 0$ et
$${\Bbb E}_{v^\vee} \;\cong \R^1p_* \;({\Cal L}_v )^\vee \tag (10.10)$$

\dem  Tout est clair, sauf peut-\^etre le dernier point. Notons que
par dualit\'e de Grothendieck-Serre on a un isomorphisme canonique
$q_* (\omega_{C/S} ) \;\cong \; \R^1 q_* ({\Cal O}_C) ^\vee$
qui est   de mani\`ere \'evidente   un $G$-isomorphisme. Par ailleurs
$$\R^1 q_* ({\Cal O}_C) = \R^1 p_* \pi_* ({\Cal O}_C) =
\bigoplus_{v\in \hat G} \R^1 p_* ({\Cal L}_v \otimes V_v) =
\bigoplus_{v\in \hat G} \R^1 p_* ({\Cal L}_v) \otimes V_v$$
  Par identification  des facteurs isotypiques, on obtient le r\'esultat.\hfill\break\qed

 On peut moduler cette construction en se reportant  \`a la section 6.4.2,  en particulier en consid\'erant le carr\'e 2-cart\'esien
 $$\commdiag {\overline{\Cal C}_{g',b}&\mapright^\pi & \overline{\Cal M}_{g',b}\cr
 \mapdown &&\mapdown\cr
 \overline{\Cal C} &\mapright^\pi &\overline{\Cal H}_{g,G,\xi}\cr}\tag (10.11)$$
 Une section de base $S$ du champ $\overline{\Cal C}$ est un rev\^etement $\pi: C\to D$ de base $S$, \'equip\'e d'une section $Q: S\to D$ additionnelle. On d\'efinit un  module sans torsion   $\Bbb L_v \, (v\in \hat G)$ sur $\overline{\Cal C}$, de rang $\deg(v)$,  par la prescription 
 $$\Bbb L_v (\pi) = Q^*({\Cal L}_v) \tag (10.12)$$
Invoquant le morphisme  $\psi:  \overline{\Cal C}_{g,r}\to   \overline{\Cal C} = \overline{\Cal H}_{g,G,\xi}\times_{  \overline{\Cal M}_{g',b }} \overline{\Cal C}_{g',b}$ (proposition 6.14), on voit que $\psi_* ({\Cal O}_{ \overline{\Cal C}_{g,r}}) = \bigoplus_{v\in \Irrep(G)} \, \Bbb L_v \otimes V_v  $.
Avec cette d\'efinition il est que 
 ${ {\Bbb E}^\vee}_{v^\vee}  =   \R^1\pi_* ({\Bbb L}_v)  = -\pi_!(\Bbb L_{v }).   $
    Sous les hypoth\`eses de la proposition
10.1 nous poserons pour tout $v \in \hat
G$
$$ \lambda_v \,=\, \det\,({\Bbb E}_v) \in Pic(\overline{\cal H}_{g,G,\xi}) \tag (10.13)$$  
On notera par la m\^eme lettre la premi\`ere classe de Chern de $\Bbb E_v$.
Avec ces notations   on a   $\lambda_1 = \delta^* (\lambda_{bs})$, et $ \lambda = \lambda_{tot} =  \bigotimes_{v\in
\hat G} \,   \lambda_v^{\otimes \deg (v)}.  $

    La dualit\'e de Serre-Grothendieck  donne une autre construction des ${\Bbb E}_v$.  Soit un rev\^etement galoisien
stable $\pi: C \rightarrow D$ de groupe $G$. Le faisceau $\pi_* ({\Cal
O}_C)$ est un $({\Cal O}_D , G)$ module sans torsion de rang $\# G$,
de formation
compatible \`a tout changement de base. Il se d\'ecompose  en facteurs isotypiques (\S \  4)
$\pi_* ({\Cal O}_C) \;=\; \bigoplus_{v\in\hat G} \;{\Cal L}_v
\otimes V_v . $
Il est imm\'ediat  que ${\Cal L}_1 = \pi_* ({\Cal O}_C)^G = {\Cal
O}_\Sigma$, et que le rang de ${\Cal
L}_v $ est \'egal \`a $\dim (v)$.

\proclaim Proposition 10.5.  On a  $p_* ({\Cal O_C}) = {\Cal O}_S$,
et  si $v\ne 1$,  on a $\;p_* ({\Cal
L}_v ) = 0$ et
$${\Bbb E}_{v^\vee} \;\cong \R^1p_* \;({\cal L}_v )^\vee \tag (10.14)$$

\dem  Tout est clair, sauf peut-\^etre le dernier point. Notons que
par dualit\'e de Grothendieck-Serre on a un isomorphisme canonique
$q_* (\omega_{C/S} ) \;\cong \; \R^1 q_* ({\Cal O}_C) ^\vee$
qui est    un $G$-isomorphisme. Par ailleurs
$$\R^1 q_* ({\Cal O}_C) = \R^1 p_* \pi_* ({\Cal O}_C) =
\bigoplus_{v\in \hat G} \R^1 p_* ({\Cal L}_v \otimes V_v) =
\bigoplus_{v\in \hat G} \R^1 p_* ({\Cal L}_v) \otimes V_v$$
  Par identification  des facteurs isotypiques, on obtient le r\'esultat.   \hfill\break \qed

  Rappelons que  classiquement on d\'efinit un   faisceau inversible  not\'e  ${\cal L}_i \in \Pic (\overline {\Cal
M}_{g,n})$     par 
$${\cal L}_i (q: D\to S) = \sigma_i^* (\omega_{D/S}) \tag (10.15)$$
en notant   $\sigma_i$ la i-\`eme section. Pour  \'eviter toute  confusion  avec les faisceaux introduits ci-dessus,  on notera par la m\^eme lettre   $\psi_i$ le faisceau inversible ${\cal L}_i$ et sa premi\`ere classe de Chern $c_1({\cal L}_i)$.   La m\^eme lettre  d\'esignera  aussi pour simplifier \`a la fois $\psi_i\in \Pic (\overline {\Cal M}_{g',b})$, ainsi  que son image r\'eciproque dans
$\Pic (\overline {\Cal H}_{g,G,\xi})$.     La structure du
groupe $\Pic (\overline {\Cal M}_{g,n})$ (resp. $\Pic({\Cal M}_{g,n})$) est
connue \cite{40}, \cite{53}.  Si
$g\ne 1,2$, ces  groupes sont   libres de rang fini\footnote { Si
$g\geq 3$, les faisceaux inversibles $\lambda$, et $\psi_i \,(1\leq
i\leq n)$ forment une base de $\Pic ({\Cal M}_{g,n}$); si on ajoute \`a
 ces faisceaux inversibles ceux associ\'es aux composantes du bord
(les $\delta$), on obtient une base de $\Pic (\overline {\Cal
M}_{g,n})$,  loc.cit. Le groupe $\Pic (\overline {\Cal M}_{0,n})$ est
libre de rang  $2^{n-1} - { n \choose 2} - 1$ (Keel, voir \cite{53})}.     
Soit le morphisme de groupes 
$$\delta^*:  \Pic ( \overline {\Cal
H}_{g,G,\xi}) \longrightarrow  \Pic (\overline {\Cal M}_{g',b}) \tag
(10.16)$$
Si le groupe de droite est libre (de rang fini), par exemple si $g' = 0$,    
 il n'est pas difficile de montrer que ce morphisme est injectif (voir par exemple \cite
{43}). Naturellement   le groupe $\Pic
(\overline {\Cal H}_{g,G,\xi})$ est discret que dans
des cas exceptionnels, essentiellement  $g' = 0$ et $G$ ab\'elien.
Si $g'=0$ et $G$ ab\'elien, un fait g\'en\'eral
qui s'applique aux champs de Deligne-Mumford donne   $\Pic  (\overline {\Cal
H}_{g,G,\xi})\otimes \Bbb Q \cong  \Pic (\overline
{H}_{g,G,\xi})\otimes \Bbb Q$,  alors  $\Pic  (\overline {\Cal
H}_{g,G,\xi})\otimes \Bbb Q \cong  \Pic  (\overline {\Cal M}_{0,b})
\otimes \Bbb Q$. Le groupe  $\Pic  (\overline {\Cal H}_{g,G,\xi})$
ne diff\`ere     de $\Pic  (\overline {\Cal M}_{0,b})$ que
par la torsion, qui est non triviale en g\'en\'eral. Ils ont un m\^eme rang  $2^{b-1} -  { b \choose 2} - 1$ \cite{53}.
\bigskip

 \beginsection 10.2.  Les fibr\'es en droites $\psi_{i,\chi}$ et $\mu_{i,v}$

\bigskip

 Examinons   les fibr\'es tautologiques  qui d\'erivent  des fibres au-dessus de
$Q_1,\dots,Q_b$,  les points marqu\'es de la base, et portant l'information sur la ramification. Rappelons que   les points de branchement sont inclus dans les $Q_i$.  Notons $\sigma_i: S \to D$ le
morphisme d\'efinissant $Q_i$. Il est clair que  le $({\Cal
O}_S,G)$-faisceau ${\Cal O}_{\pi^{-1} (Q_i)}$ est localement libre de
rang $\vert G\vert$.  On a  en tenant
compte de (10.9) la d\'ecomposition ${\Cal O}_{\pi^{-1} (Q_i)}
\cong \bigoplus_{v\in \widehat G} \sigma_i^* ({\Cal L}_v) \otimes V_v $. 
On pose  
  $$ {\Cal L}_{i,v} (\pi) = \sigma_i^* ({\Cal L}_v) \tag (10.17)$$ 
 
  Cette construction d\'efinit  un   fibr\'e vectoriel  ${\Cal L}_{i,v}$ sur  $\overline {\Cal
H}_{g,G,\xi}$ de rang $\deg (\chi)$.   Si $Q_i $ est un point de
branchement  d'holonomie $[H_i,\chi_i]$,   il est possible de d\'evisser les fibr\'es ${\Cal L}_{i,v}$ en une somme directe de
sous-fibr\'es portant une certaine part de l'information sur la
ramification.  Fixons un repr\'esentant $(H_i,\chi_i)$ de l'holonomie
en le point de branchement d'indice $i$. Soit alors le diviseur de
Cartier relatif $\Delta_i$, fini et  \'etale sur $S$, lieu des points
de la fibre $\pi^{-1}(Q_i)$ exactement d'holonomie $(H_i,\chi_i)$.
Alors $\Delta_i \to S$ est un torseur de groupe $C_G(H_i)/H_i$ (voir
3.1.2),    le diviseur $\pi^{-1}(Q_i)$ ayant  pour
expression
$$\pi^{-1}(Q_i)  = e_i \left( \Ind_{C_G(H_i)}^G  \Delta_i
\right)  \quad (e_i = \vert H_i\vert)  $$
En particulier
${\Cal O}_{\pi^{-1}(Q_i)} = \Ind_{C_G(H_i)}^G \left({ \Cal
O}_{e_i\Delta_i} \right)$. On est  ainsi ramen\'e    \`a  pr\'eciser la
structure de $({\Cal O}_S,G)$-module de $\Ind_{C_G(H_i)}^G \left({\Cal
O}_{e_i\Delta_i}\right) $, c'est \`a dire  la structure du $({\Cal
O}_S,C_G(H_i))$-module ${\Cal O}_{e_i\Delta_i}$. Pour $0\leq k <e_i$
soit la suite exacte de $({\Cal O}_S,C_G(H_i))$-modules

$$0\rightarrow
{\Cal O}_{\Delta_i}(-k\Delta_i) \rightarrow {\Cal O}_{(k+1)\Delta_i}
\rightarrow {\Cal O}_{k\Delta_i} \rightarrow 0\tag (10.18)$$

Soit  $\chi$ un caract\`ere irr\'eductible   de $C_G(H_i)$ avec $\rho_\chi:
C_G(H_i) \to \GL (W_\chi)$ la repr\'esentation associ\'ee. Du fait
que $H_i$ est dans le centre de $C_G(H_i)$,    si $h\in H_i,$ la restriction de
$\rho_\chi (h)$ \`a $H_i$ est une homoth\'etie $\rho_\chi (h) =
\chi_i^j(h)$,  pour un certain $j \in [0,e_i-1[$. Nous
parlerons de $\chi_i^j$ comme du poids de la restriction de $\chi$ \`a
 $H_i$.  Le groupe  $H_i$ agit sur le faisceau conormal ${\Cal
O}_{\Delta_i} (-\Delta_i)$  au moyen du caract\`ere $\chi_i$, il est
alors imm\'ediat de d\'eduire de la suite exacte (10.18), que les deux
faisceaux localement libres $ {\Cal O}_{\Delta_i}(-k\Delta_i)$ et
${\Cal O}_{k\Delta_i}$ sont "disjoints", et donc  vus comme $({\Cal
O}_S,C_G(H_i))$-modules on a $  {\Cal O}_{(k+1)\Delta_i}  \cong  {\Cal
O}_{\Delta_i}(-k\Delta_i) \oplus {\Cal O}_{k\Delta_i} $.

L'action de $H_i$ sur ${\Cal O}_{\Delta_i} (-k\Delta_i)$ est
donn\'ee par la caract\`ere $\chi_i^k$, il en r\'esulte  que toute
fibre de ${\Cal O}_{\Delta_i} (-k\Delta_i)$ s'identifie \`a  la repr\'esentation
$\Ind_{H_i}^{C_G(H_i)}(\chi_i^k) $, il s'ensuit  la d\'ecomposition 
$${\Cal O}_{\Delta_i} (-k\Delta_i) = \bigoplus_{\chi\in \widehat {C_G(H_i)},\chi_{\vert H_i}={\chi_i}^k}   E_{i,\chi} \otimes W_\chi \tag (10.19)$$  
la somme portant sur les caract\`eres irr\'eductibles de
$C_G(H_i)$ de restriction $\chi_i^k$   \`a $H_i$. On en d\'eduit  ais\'ement, en mettant bout \`a bout les d\'ecompositions (10.18), la d\'ecomposition 
$${\Cal O}_{k\Delta_i} = \bigoplus
_{\chi\in \widehat {C_G(H_i)}} E_{i,\chi} \otimes W_\chi \tag
(10.20)$$
o\`u   la somme directe porte sur les caract\`eres irr\'eductibles de $C_G(H_i)$ dont la restriction \`a $H_i$ est  le
caract\`ere ${\chi_i}^\alpha$, avec  $0 \leq \alpha \leq k-1$.  Finalement on a la d\'ecomposition
$${\Cal O}_{e_i \Delta_i} = \bigoplus_{k=0}^{e_i-1} \,
{\Cal O}_{\Delta_i} (-k\Delta_i)  = \bigoplus_{\chi\in \widehat
{C_G(H_i)} } E_{i,\chi} \otimes W_\chi 
\tag (10.21) $$
 
 Cette construction    d\'efinit  pour tout  $\chi\in \widehat {(C_G(H_i))}$   un fibr\'e   
  $  E_{i,\chi}$ de rang $\deg(\chi)$.  Explicitons les relations qui relient les  fibr\'es ${\Cal L}_{i,v}$ et les $E_{j,\chi}$.
  
\proclaim Proposition 10.6. Avec les
notations  pr\'ec\'edentes,  $E_{i,\chi}$  est un
fibr\'e vectoriel de rang $\deg (\chi)\,\,(\chi\in \widehat
{(C_G(H_i))}$.  Le fibr\'e vectoriel  ${\Cal   L}_{i,v}  $ est de
rang $\deg v$, il se  d\'ecompose en  
$${\Cal L}_{i,v} =
\bigoplus_{\chi \in \widehat {C_G(H_i)}} \,  E_{i,\chi}
^{\,n_{v,\chi}} \quad  \hbox{\rm o\`u}  \quad n_{v,\chi} =
(\chi,v\vert_{C_G(H_i)}) \tag (10.22)$$

\dem  
Dans l'anneau des repr\'esentations de $G$   soit la  d\'ecomposition en facteurs irr\'eductibles
$$\Ind_{C_G(H_i)}^G (W_\chi) \, = \, \bigoplus_{v\in \hat G}
\,\,V_v^{n_{v,\chi}} $$
 avec  $n_{v,\chi} = \langle
\Ind_{C_G(H_i)}^G (W_\chi) , V_v \rangle_{\vert G}\, =\,  \langle
W_\chi , {V_v}_{\vert C_G(H_i)}\rangle_{\vert C_G(H_i)} $.
   Par insertion de cette relation dans (10.21), et induction, on obtient en regroupant les termes  exactement (10.22).
 \qed
 
  Les d\'eterminants respectifs $\det({\Cal L}_{i,v})$ et $\det(E_{i,\chi})$
seront  not\'es $\mu_{i,v}$ et $\psi_{i,\chi}$, ceci pour $v\in \hat G, \,\, \chi \in \widehat{C_G(H_i)}$.  

Dans la suite, on notera par la m\^eme lettre   un faisceau inversible, et sa premi\`ere classe de Chern, par exemple $\mu_{i,v}$ et $\psi_{i,\chi}$. On remarquera que les fibr\'es vectoriels donn\'es par la d\'ecomposition ${\Cal O}_{\Delta_i} (-\Delta_i) = \bigoplus_{\chi\in \widehat {C_G(H_i)},\chi_{\vert H_i}={\chi_i}}   E_{i,\chi} \otimes W_\chi $ se substituent  aux  images r\'eciproques des ${\cal L}_P\in \Pic(\overline{\cal M}_{g,r})$, qui ne sont pas d\'efinis car les points de ramification ne sont pas num\'erot\'es.  Cependant la relation
$$\bigoplus_{P\in \Delta_i} {\cal L}_P =   \bigoplus_{\chi\in \widehat {C_G(H_i)},\chi_{\vert H_i}={\chi_i}}   E_{i,\chi} \otimes W_\chi $$
a un sens. Lorsque $G$ est cyclique (\S \ 10.4), tout devient plus simple.

\bigskip
  \beginsection 10.3. Relations  de Riemann-Hurwitz d'ordre sup\'erieur   

\bigskip
  {\it 10.3.1. Calculs dans l'anneau de Chow } 
  
  \bigskip

     Fixons quelques conventions. Si   $\Cal M$ est un champ de Deligne-Mumford de dimension pure $ n = \dim {\Cal M}$, on notera  $A^\bullet (\Cal M)$ l'anneau de Chow de $\Cal M$ \`a coefficients rationnels. On retient les notations de Mumford \cite{58}, voir aussi  Vistoli \cite{68}, en particulier l'identification $A^k({\Cal M}) = A_{n-k}({\Cal M})$ (loc. cit. \S\  3). Si $Z$ est un sous-champ ferm\'e int\`egre de $\Cal M$ de dimension $k$, on notera   $[Z]\in A_k({\Cal M})$ le cycle correspondant, et  $[Z]_Q$ la classe fondamentale de $Z$. Si  $e(Z)$ est  l'ordre du groupe des automorphismes d'un point g\'en\'erique de $Z$, alors \cite{58}, 
$[Z]_Q = e(Z)^{-1} [Z]$.  L'identit\'e de $A^\bullet (\Cal M)$ est donc $[{\Cal M}]_Q$. Si $M$ est un espace des modules grossiers le morphisme ${\Cal M }\to M$ induit un isomorphisme $A^\bullet(\Cal M)\buildrel\sim\over \to A^\bullet(  M)$ \cite{68}.   Dans la suite  ${\Cal M }= \overline{\Cal M}_{g,n}$.  

Soit la courbe universelle $q:\overline{\Cal C}_{g,n}= \overline{\Cal M}_{g,n+1} \to \overline{\Cal M}_{g,n}$; le morphisme $q$ est  le morphisme not\'e  usuellement  $\pi_{n+1}$, oubli du $n+1$-i\`eme point. Posons   $K = c_1(\omega_q)$ et   soit  $D_i$ (ou $D_{i,n+1}$) l'image de la  section universelle d'indice $i, \,\, (1\leq i\leq n)$.  Il est commode d'utiliser, outre les classes kappa  de Mumford-Morita-Miller \cite {58},   les classes kappa de Mumford   (Arbarello-Cornalba \cite {4})   
$ \tilde\kappa_l = q_\star (K ^{l+1}) \in
A^l(\overline {\Cal M}_{g,n}), \quad \hbox{\rm resp.} \quad\kappa_l = q_\star(c_1(\omega_q(\sum D_i))^{l+1}). $
 Elles sont reli\'ees par   \cite{4}
$$\kappa_a = \tilde \kappa_a+ \sum_{i=1}^n \psi_i^{a}\tag (10.23)$$
 Les classes
$\kappa_l$ satisfont  relativement au morphisme $\pi_n: \overline{\Cal M}_{g,n} \to  \overline{\Cal M}_{g,n-1}$,   \`a la  relation simple
$$\kappa_l = \pi_n^*(\kappa_l) + \psi_n^l \tag (10.24)$$
   Rappelons que $D_{i,n+1}$   est le diviseur  dont les points sont obtenus en collant une ''bulle''   au point d'indice $i$.  Les classes $\psi_i$ satisfont  relativement \`a $\pi_{n+1}$ \`a la relation importante (Lemme de comparaison)
 $$\psi_i = \pi_{n+1}^*(\psi_i) + D_{i,n+1} \, (1\leq i\leq n) \tag (10.25)$$
Soit le champ ${\Cal H}_{g,G,\xi}$  source du   rev\^etement universel $q: \overline{\Cal C}_{g,G,\xi} \to \overline{\Cal H}_{g,G,\xi},$ et  $\omega = \omega_q$. On   d\'efinit des classes 
$\kappa_l \, (\hbox{\rm resp.} \tilde\kappa_l) \in A^l(\overline{\Cal H}_{g,G,\xi})\otimes \Bbb Q$.
  Pour les distinguer des pr\'ec\'edentes, on notera  $  {\kappa'}_l  \,\,( \hbox{\rm resp.} \,\, \tilde\kappa'_l) \in    A^l(\overline {\Cal M}_{g',b})\otimes \Bbb Q$ celles de m\^eme nom  qui vivent sur la base, c'est \`a dire dans $A^l(\overline{\Cal M}_{g',b})$.   Du fait du carr\'e cart\'esien (6.24), on voit que les classes $\kappa_l$ et $ \tilde\kappa_l$ proviennent des classes de m\^eme nom de  $\overline{\Cal M}_{g,(r)}$.  
 Le morphisme $\imath: \overline {\Cal H}_{g,G,\xi} \to \overline{\Cal M}_{g,(r)}$ \'etant une immersion locale r\'eguli\`ere (prop. 6.14) on notera que  le morphisme de Gysin $\imath^!$ est d\'efini au niveau des cycles.  Au niveau de la premi\`ere classe de Chern, et pour prendre en compte les \'el\'ements de torsion  du groupe de Picard,  il   est utile d'utiliser le formalisme du produit d'intersection de  Deligne  (voir  \cite {9}, Jarvis \cite
{43},\ \S \ 3).  

Rappelons que cette construction associe \`a toute courbe prestable $q: C\to S$, et  deux
faisceaux inversibles ${\Cal L}, {\Cal M}$ sur $C$,  un faisceau
inversible $\langle {\Cal L},{\Cal M}\rangle \in \Pic (S)$  par 
$$\langle {\Cal L},{\Cal M}\rangle   = \det Rq_* ({\Cal L}\otimes {\Cal M}) \otimes (\det Rq_* {\Cal
L})^{-1} \otimes (\det  Rq_* {\Cal M})^{-1} \otimes \det Rq_* {\Cal
O}_C) \tag (10.26)$$
 Le d\'eterminant de la cohomologie $\det Rq_* ( - )$ se r\'eduit  dans (10.27) \`a 
$\det R^0 q_* (-) \otimes
(\det R^1q_* (-))^{-1}$.
 On  montre que l'op\'eration
$\langle {\Cal L},{\Cal M}\rangle $ est bilin\'eaire (bimultiplicative) 
sym\'etrique, et fonctorielle relativement aux isomorphismes (\cite
{9},\cite {18}). Si $\omega = \omega_{C/S} $ est le faisceau
dualisant relatif,  l'identit\'e suivante  (en notation additive)
est satisfaite (Th\'eor\`eme de Deligne-Riemann-Roch  (\cite
{18}, \cite {43}): 
$$2 \det Rq_* {\Cal L} = \langle {\Cal L},{\Cal L}
\rangle - \langle {\Cal L},\omega \rangle + 2 \det Rq_* \omega \tag (10.27)$$
 Si  $D$ est un diviseur de Cartier relatif, et ${\Cal M} = {\Cal O}_C(D)$,  il d\'ecoule imm\'ediatement de la d\'efinition
jointe \`a la suite exacte $O\to {\Cal O}(-D) \to {\Cal O} \to {\Cal
O}_D\to 0$ la relation \cite {9}  
$$\langle {\Cal L},{\Cal M}\rangle  =
{\det} _{{\Cal O}_S} (q_* ({\Cal L}\otimes  {\Cal O}_D)) \otimes
({\det}_{{\Cal O}_S} q_* ({\Cal O}_D))^{-1}\tag (10.28)$$
    Dans ce contexte la classe  $\kappa = \kappa_1 $ est repr\'esent\'ee  par le
faisceau inversible d\'efini  par  le produit d'intersection  $\kappa_1 = \langle
\omega , \omega \rangle $.
 D'abord une remarque pr\'elimiminaire.  
\proclaim Lemme 10.7. Soit  $q:
C\buildrel\pi\over\rightarrow  D\buildrel p \over \rightarrow S$ un
rev\^etement entre $S$-courbes lisses, de degr\'e $n$. Pour tout faisceau inversible $\Cal L \in
\Pic (D)$, on  
$$\langle \pi^*({\Cal L}),\pi^*({\Cal L})\rangle =  n \langle {\Cal L},{\Cal L} \rangle \tag (10.29)$$  
 
\dem  Du fait que $\pi$ est fini, on a $Rq_* (-) = Rp_* \pi_* (-)$,
et donc pour tout faisceau inversible  $\Cal M$ sur $D$, $Rq_* (\pi^*
({\Cal M}) = Rp_* ({\Cal M} \otimes \pi_* ({\Cal O}_C))$.   L'hypoth\`ese dit que $E = \pi_*({\Cal O}_S)$ est localement libre de rang $n$. Le  prob\`eme se  r\'eduit  \`a v\'erifier  qu'il y a  pour tout fibr\'e
 vectoriel $E$ de rang $n$, et tous faisceaux inversibles $\Cal L,
\Cal M$ sur $D$, un isomorphisme canonique (unique) 
$$n\langle  {\Cal
L}, {\Cal M }\rangle  =    \det Rp_* ({\Cal L}\otimes {\Cal M}\otimes  E) -
\det Rp_* ({\Cal L}\otimes E )   - \det Rp_* ({\Cal M} \otimes E) + \det
Rp_* (E)  $$ 
En suivant les arguments de (\cite {19}, \S \ 9.5), on
peut se ramener \'etale-localement au cas o\`u $E$ admet une
filtration  par des sous-fibr\'es, \`a quotients successifs de rang
un.    Supposons que $E$ admette une filtration 
$0 = E_0 \subset E_1\subset \cdots \subset E_n = E,$
$E_j$ \'etant localement libre de rang $j$, et ${\Cal L}_j = E_j/E_{j-1}$  localement libre de rang un.   On obtient, la  notation \'etant additive
$$\eqalign{\langle \pi^*({\Cal L}),\pi^*({\Cal M})\rangle  = & \sum_j \det Rp_*
({\Cal L}\otimes {\Cal M}\otimes {\Cal L}_j) - \sum_j  \det Rp_* ({\Cal
L}\otimes {\Cal L}_j)  \cr
& - \sum_j \det Rp_* ({\Cal M} \otimes
{\Cal L}_j) + \sum_j \det Rp_* ({\Cal L}_j)\cr}$$
ce qui tenant compte du
fait que les ${\Cal L}_j$ sont de rang un, peut s'\'ecrire $\sum_j
\langle {\Cal L}\otimes {\Cal L}_j,{\Cal M}\rangle - \sum_j  \langle {\Cal M}
, {\Cal L}_j\rangle .$
La conclusion d\'ecoule dans ce cas de la
bilin\'earit\'e du produit d'intersection.  \hfill\break  \qed

 \bigskip
 {\it 10.3.2. Relations de Riemann-Hurwitz d'ordre sup\'erieur}
 
 \bigskip
  Soit   de nouveau  la  $G$-courbe universelle $\pi : \overline {\Cal
C}_{g,G}\longrightarrow \overline {\Cal H}_{g,G,\xi}$ d\'efinie dans la section 6.4.2.  Elle
  s'ins\`ere dans un carr\'e $2$-commutatif    
$$\commdiag{ \overline {\Cal C}_{g',b}&\mapright^p& \overline {\Cal
M}_{g',b}\cr
 \mapup\lft{\Delta} &&\mapup\lft{\delta}\cr
   \overline {\Cal C}_{g,G,\xi} &\mapright^q &
\overline {\Cal H}_{g,G,\xi}  \cr}\tag (10.30)  $$   

 Notons  $\omega  = \omega_{\overline {\Cal C}_{g,G,\xi}/\overline
{\Cal H}_{g,G,\xi} }\in \Pic (\overline {\Cal C}_{g,G,\xi})$ le
faisceau inversible dualisant relatif du morphisme $ q$. Ce  faisceau est visible sur un atlas. Soit $\pi_U: C_U \to D_U$  un $G$-rev\^etement de base $U$ d\'efinissant   un
atlas de  $\overline {\Cal H}_{g,G,\xi} $. Alors le rev\^etement obtenu par changement de base 
$p_2: \tilde {C}_U = C_U\times_U C_U\to C_U,$
muni de la section diagonale (il n'est pas n\'ecessaire de stabiliser) fait de $C_U$ un atlas de  ${\Cal C}_{g,G,\xi}$.          L'incarnation de $\omega$ sur l'atlas $C_U \to \overline
{\Cal C}_{g,G,\xi}$  est $\omega_{C_U/U}$.   En notant 
  $\omega'$ le faisceau $  \omega_{\overline {\Cal
C}_{g',b}/\overline {\Cal M}_{g',b}}$,   on a  $\Delta^*(\omega')
(C_U) = \pi_U^*(\omega_{D_U/U})$.  

Pour tout indice $i$ avec $ 1\leq i \leq b$,
soit $R_i$ la partie du diviseur de ramification au dessus de la
section $Q_i$.  Pour d\'efinir $R_i,$  il suffit de le  r\'ealiser au niveau de l'atlas $C_U\to {\Cal C}_{g,G,\xi}$.   Soit le rev\^etement $p_2: \tilde{C}_U \to C_U$. 
 Le diviseur relatif $\tilde{R}_i = \Ind_{H_i}^G \Delta_i\subset \tilde{C}_U$ est  alors bien d\'efini. Rappelons  que ce diviseur  est \'etale de degr\'e $[G:H_i]$ sur $C_U$. Alors   $R_i = \Delta^*({\tilde R}_i)$, et ${\tilde R}_i = R_i\times_u C_U$. On peut donc consid\'erer $R_i$ comme un diviseur de 
$\overline {\Cal C}_{g,G,\xi}$, \'etale de degr\'e $[G:H_i]$ sur $\overline {\Cal H}_{g,G,\xi}$. Le diviseur de ramification est alors 
$$R = \sum_i (e_i-1)R_i \subset \overline{\Cal C}_{g,G,\xi} \tag (10.31)$$

Notons $B_i\subset \overline {\Cal M}_{g',b}$ le diviseur image de la section $Q_i$. La formule de ramification transcrite au  niveau des champs  a pour
forme:
\proclaim  Proposition 10.8.  i)  Pour tout  $1\leq i\leq b$, on a  l'\'egalit\'e de diviseurs de Cartier  $  \Delta^*(B_i) = e_iR_i$.\hfill\break
ii) Dans le groupe de Picard de
$\overline {\Cal C}_{g,G,\xi} $ on a  la relation (formule de ramification)  
 $$\omega = \Delta^* (\omega') \otimes \left(
\otimes_{i=1}^b {\Cal O} (R_i)^{\otimes {e_i-1} }\right)\tag
(10.32)$$

\dem Pour le premier point, notons le diagramme commutatif
$$\commdiag{ \overline {\Cal C}_{g,G} &\mapright^{\Delta}&   \overline {\Cal
M}_{g',b}\cr
 \mapup &&\mapup\cr
   C_U &\mapright^{\pi_U }&
 D_U\cr}   $$ 
 dans lequel les fl\^eches verticales sont repr\'esent\'ees respectivement par $p_2: \tilde C_U\to C_U $ et $D_U\times_U D_U \to D_U$ muni de la section diagonale.  Pour identifier $\delta^*(B_i)$, il suffit d'examiner l'image r\'eciproque dans $C_U$, qui est $\pi_U^*(B_i) = e_iR_i$.\hfill\break
ii)  Le faisceau $\Delta^*(\omega')$ a pour description  relativement \`a l'atlas $C_U$ l'image r\'eciproque $\pi_U^*(\omega_{D_U/U})$. De la sorte la relation (10.32) n'est que  la  simple traduction de la formule
de ramification appliqu\'ee  \`a $\pi: C_U \to D_U$,
jointe \`a la description du diviseur de ramification  (voir remarque
4.1, et \S \ 6.4.2). \hfill\break\qed  
 
 Soit $\tau_i: R_i\hookrightarrow \overline{\Cal C}_{g,G,\xi}$ l'inclusion.  On peut \'enoncer le r\'esultat principal de cette section, un  analogue  sup\'erieur  de la formule de Riemann-Hurwitz, relation qui exhibe une propri\'et\'e remarquable des classes kappa\footnote{Le r\'esultat est \'enonc\'e pour les rev\^etements galoisiens, mais il est  facile de  valider  la preuve dans le cas non galoisien.}:
\proclaim Th\'eor\`eme 10.9. Dans $A^l( \overline{\Cal H}_{g,G,\xi}) $ on a  pour tout $l\geq 0$ (on rappelle que les coefficients sont rationnels)  la relation de Riemann-Hurwitz:
$$\tilde\kappa_l = \vert G\vert \delta^* (\tilde\kappa'_l) + \sum_i \, (-1)^{l+1} (1 - e_i^{l+1}) {\tau_i}_* (c_1(N_{R_i})^l) \tag (10.33)$$
d'une autre  mani\`ere $ \kappa_l = \vert G\vert \delta^*(\kappa'_l).  $

 \dem  De la relation de Riemann-Hurwitz (10.32), on tire en prenant la  premi\`ere classe de Chern des deux membres
 $$K = \Delta^* (K') + \sum_i (e_i - 1) c_1({\Cal O}(R_i))$$
 En particulier en  \'elevant \`a  la puissance $l+1$  des deux membres de cette \'egalit\'e, et tenant compte aussi du fait que $R_i$ et $R_j$ sont ''disjoints'' si $i\ne j$
$$ K^{l+1} = \Delta^* ({K'}^{l+1}) + \sum_i \sum_{j=1}^{l+1}  { {l+1}\choose j } (e_i - 1)^j \Delta^* ({K'}^{l+1-j}) c_1({\Cal O}(R_i))^j $$
D'une autre mani\`ere, en ins\'erant l'\'egalit\'e $\Delta(K')^{l+1-j} = (K - \sum_i (e_i - 1) c_1({\Cal O}(R_i)))^{l+1-j}$ si $1 \leq j \leq l+1$   on obtient  
$$\Delta^*({K'}^{l+1-j}) c_1({\Cal O}(R_i))^j \,= \,
\sum_{k=0}^{l+1-j}  { {l+1-j}\choose k} (1 - e_i)^k K^{l+1-j-k} c_1({\Cal O}(R_i))^{k+j}$$
On a $ K^{l+1-j-k} c_1({\Cal O}(R_i))^{k+j} = {\tau_i}_*(c_1({\omega_\pi}_{\vert R_i})^{l+1-j-k} c_1({\Cal O}_{R_i}(R_i))^{k+j-1}$, ce qui compte tenu  du fait que ${\omega_\pi}_{\vert R_i} \cong {\Cal O}_{R_i}(R_i)$   donne finalement 
$$\Delta^*({K'}^{l+1-j}) c_1({\Cal O}(R_i))^j \, = \sum_{k=0}^{l+1-j}  {{l+1-j}\choose  k} (1 - e_i)^k (-1)^{l+1-j-k}   {\tau_i}_* c_1({\Cal O}_{R_i} (R_i))^{l}$$
  Le coefficient  devant  ${\tau_i}_* \, c_1({\Cal O}_{R_i}(R_i))^{l}$ qui en r\Ž sulte est
 $$ (-1)^{l+1-j} \sum_{k=0}^{l+1-j}   { {l+1-j}\choose  k} (e_i - 1)^k =  (-1)^{l+1-j} e_i^{l+1-j}$$
 En reportant cela dans l'expression 10.26, on obtient   que le coefficient  devant l'expression  $ \tau_! \,c_1({\Cal O}_{R_i} (R_i))^{l}$ est 
 $$\sum_{j=1}^{l+1}  { {l+1}\choose j }  (e_i - 1)^j (-e_i)^{l+1-j}  = (-1)^{l+1} - (-e_i)^{l+1}$$
 Appliquons    $\pi_*$ aux deux membres, on voit alors qu'il suffit  alors de prouver que  
 $\pi_* (\Delta^*({K'}^{l+1})) = \vert G\vert \delta^* (\tilde\kappa_l)$.
  Soit le morphisme fini  $\psi: \overline {\Cal C}_{g,G,\xi} \rightarrow \overline {\Cal H}_{g,G,\xi} \times_{\overline {\Cal M}_{g',b}} \overline {\Cal C}_{g',b}$. Si  $p_1,p_2$ sont  les projections du produit   cart\Ž sien sur les deux facteurs, on  note que 
 $$\pi_* \psi^* p_2^* ({K'}^{l+1}) = {p_1}_* \psi_* \psi^* {p_2}^* ({K'}^{l+1}) $$
 Mais par la formule de projection, jointe  au fait que  $\psi$ est fini  de degr\'e $\vert G \vert$, on a $\psi_* \psi^* (-) = \vert G\vert (-)$ (voir \cite{58} Prop 3.8). Par ailleurs,   du fait de la platitude des morphismes, on a la propri\'et\'e de  changement de base (\cite {66}, Lemma 3.9) 
 $${p_1}_* p_2^* ({K'}^{l+1}) = \delta^* ({\tilde\kappa'}_l) $$
 la  relation (10.33)  en d\'ecoule. On passe facilement de cette relation \`a la forme compacte  qui relie les classes $\kappa$. Notons d'abord que par d\'efinition
 $$\kappa_l = \tilde\kappa_l + \sum_i (q_{\vert R_i})_* \left(c_1(N^\vee_{\vert R_i})^l\right)$$
   Comme  $q^*(N^\vee_{B_i}) = (N^\vee_{R_i})^{e_i}$,  on a
 $$e_i^l(q_{\vert R_i})_*\left(N^\vee_{\vert R_i})^l\right) = 
   q_{\vert R_i})_*\left(N^\vee_{\vert B_i})^l\right) 
     = m_i\psi^l \,\, (m_i = {\vert G\vert\over e_i})$$
 d'o\`u finalement  \ \ $\kappa_l = \vert G\vert \delta^*(\tilde\kappa'_l) + \sum_i e_i^{l+1} \, {m_i\over e_i^l} \, \delta^*(\psi_i)^l = \vert G\vert \,\delta^*(\kappa'_l).$
 \hfill\break
 \qed 
 
 Si $l = 0$, la relation se r\'eduit \`a   la formule de Riemann-Hurwitz  usuelle 
 $2g - 2 = \vert G\vert (2g' - 2) + \sum_i (1 - {1\over e_i})$,
on encore 
$2g-2+r = \vert G\vert (2g'-2+b).$
   En effet on sait que $\tilde \kappa_0 = (2g - 2)[\overline{\Cal H}_{g,G,\xi}], \, \tilde\kappa'_0 = [\overline{\Cal M}_{g',b}]$ et comme $R_i\to \overline{\Cal M}_{g',b}$ est un torseur sous $G/H_i$, la somme se r\'eduit \`a  $\sum_i (e_i - 1){\vert G \vert \over e_i} [\overline{\Cal M}_{g',b}]$.
   
\rema {10.2}    Si $l=1$, on peut   donner une preuve diff\'erente de   (10.34)  exploitant l'\'egalit\'e  $\tilde \kappa_1 = \langle \omega , \omega \rangle$,  avec  la bimutiplicativit\'e   du produit d'intersection de Deligne.
   On \'evalue   $\langle \omega , \omega \rangle $ de deux mani\`eres. La premi\`ere  est   simplement
  la formule de Mumford pour la courbe $C_U\to U$, et la
seconde  revient \`a  substituer \`a $\omega$ le second membre de la relation de Riemann-Hurwitz
(10.32).   On se limite au cas   $G$ cyclique,  le cas d'un groupe $G$ arbitraire pouvant se traiter de mani\`ere analogue,   avec quelques modifications.  Un  probl\`eme vient cependant  du fait qu'on ne peut appliquer directement le lemme 10.7 \`a $\Pic(\overline{\Cal H}_{g,G,\xi})$ du fait de la non platitude de $\Delta$, i.e  de $\pi_*({\Cal O}_C)$ si  $C$  est singuli\`ere. Pour valider le lemme, il faut d\'efinir le produit d'intersection $\langle \Cal L , \Cal M\rangle$ si  $\Cal L$ ou $ \Cal M$ d\'eg\'en\`ere en un faisceau   sans torsion de rang un sur la $S$-courbe $C$. 
   
   Gr\^ace \`a la proposition 7.12, si   $\Cal L$ est sans torsion de rang un, on sait qu'il y a une d\'esingularisation $\rho: \tilde C\to C$ de ${\cal L}$ avec un faisceau inversible  canonique ${\Cal O}(1)$ tel que ${\Cal L} = \rho_*({\Cal O}(1))$.  On a  $R\rho_*({\Cal O}(1)) = {\Cal L}$  (Proposition 7.12).   Si $\tilde q = q\rho: \tilde C \to S$ est  le morphisme compos\'e, on note   que 
   $Rq_* ({\Cal L}) = R\tilde q_*({\Cal O}(1))$ est un complexe parfait, de sorte que $\det (Rq_*({\Cal O}(1))$ a un sens, ce qui permet de  d\'efinir $\langle \Cal L , \Cal M\rangle$  comme \'etant 
   $$\langle {\Cal L} , {\Cal M}\rangle = \langle {\Cal O}(1) , \rho^*({\Cal M})\rangle \tag (10.34)$$
   Avec ces conventions, si $\Cal P\in \Pic(C)$, pour calculer le produit d'intersection $\langle \Cal L\otimes \Cal P , \Cal M\rangle$, il suffit de noter que ${\Cal L}\otimes {\Cal P} = \rho_*({\Cal O}(1)\otimes \rho^*({\Cal P})$.   D\`es lors la relation d'additivit\'e 
   $\langle \Cal L\otimes \Cal P , \Cal M\rangle = \langle \Cal L  , \Cal M\rangle + \langle \Cal P    , \Cal M\rangle$
   est clairement v\'erifi\'ee. De cette remarque d\'ecoule la validit\'e du lemme 10.7  pour un rev\^etement stable, non n\'ecessairement lisse.   Si $G$ est cyclique d'ordre $n$ on peut  alors  justifier comme pr\'ec\'edemment  la relation 
$$\langle \omega' , \omega'\rangle = \det Rq_*({\omega'}^2\otimes {\Cal L}_j) - 2\det Rq_*(\omega'\otimes {\Cal L}_j) + \det Rq_*({\Cal L}_j) \tag (10.35)$$

\bigskip
\beginsection  10.4.   Rev\^etements cycliques 

\bigskip
  {\it 10.4.1.  Relations entre les $\psi$ et les $\mu$} 
\bigskip
Dans le cas d'un groupe de Galois cyclique
les choses se simplifient de mani\`ere notable.  Soit  $G =   {\Bbb Z}/ {n\Bbb Z}$, la classe de 1 correspondant
\`a l'automorphisme $\sigma$,  on a donc   $\widehat G = G$. Pour
tout $j\in \Bbb Z$  soit $\chi_j$ le caract\`ere tel que
$\chi_j(\sigma) = \zeta_n^j$.    L'holonomie en un
point de branchement $Q_i$ est d\'efinie par l'entier $k_i, \, 1\leq
k_i < e_i, \, \pgcd (k_i,e_i) = 1$;  soit   $1\leq \nu_i < e_i$ tel que $k_i\nu_i \equiv 1 \pmod {e_i}$.

 Rappelons que    $\mu_{i,\alpha} = \psi_{i,\alpha}$ (proposition 10.6).
   Les faisceaux inversibles $\psi_{i,\alpha}$ sont reli\'es par des relations simples: 
\proclaim  Proposition  10.10.  Soit
$\alpha =  le_i + k < n$, avec $0\leq k  < e_i$. Alors
$$\psi_{i,\alpha} = \psi_{i,k} \otimes \psi_{i,e_i}^{\otimes
l},\quad  \quad \psi_{i,e_i}^{\otimes m_i} = 0 \quad
(\psi_0 =  0) \quad \hbox { \rm  et} \,\,
 \psi_{i,\alpha} =
\psi_{i,1}^\alpha \psi_i^{-[ {\alpha \nu_i \over   e_i}]}\tag (10.36)$$   
Le morphisme $ev_i^*: \Pic(\B(C_G(H_i)/H_i))\to  \Pic (\overline {\Cal H}_{g,n,\xi})$  est  injectif. 
 
\dem   Le morphisme $\pi: \Delta_i \rightarrow S$
est un $G/H_i$-torseur, il en d\'ecoule imm\'ediatement que  si
$\alpha = le_i, \,\,\psi_{i,\alpha} = \psi_{i,e_i}^{\otimes l}$, et
$\psi_{i,e_i}^{\otimes m_i} =  0$. En consid\'erant la d\'ecomposition du  faisceau $\pi_* ({\Cal O}_{\Delta_i} (-k\Delta_i) )$ (voir 10.19), on voit facilement que les facteurs, index\'es par les
$j \equiv kk_i \pmod {e_i}$, sont $\psi_{j + le_i} = \psi_{i,j}
\otimes \psi_{i,e_i}^{\otimes l} $. On peut justifier  ces relations
d'une autre mani\`ere, en
observant qu'au voisinage de la section  $Q_i$ on a ${\Cal L}_j =
{\Cal L}^{\otimes j} [{ j\nu_i\over e_i} Q_i]$. On en tire imm\'ediatement par image
r\'eciproque et pour tout indice $\alpha$  la relation $\psi_{i,\alpha} =
\psi_{i,1}^\alpha \psi_i^{-[ { \alpha  \nu_i\over e_i}]}$.  
\hfill\break
Prouvons  le derni\`ere  point.  Rappelons que $ev_i$ est l'\'evaluation (6.26).   On
 se limite \`a  $g' = 0$, qui est le seul cas consid\'er\'e dans la suite. Il suffit de prouver que $\psi_{i,e_i}$ est d'ordre $m_i$.  On est ramen\'e \`a exhiber un rev\^etement $\pi :C \to D$ entre courbes lisses, 
de base  $\Spec \, k(\Delta)$, o\`u $\Delta$ une courbe projective lisse,  tel que le rev\^etement
galoisien \'etale $\Delta_i \to \Spec \,\Delta$ soit connexe,  avec la ramification indiqu\'ee.      Choisissons la courbe $\Delta $
de genre $h $ tel que $ 2h \geq b$. Comme les $e_i$ sont premiers \`a la caract\'eristique $p$ si $p>0$, cela ne pose pas de probl\`eme si $\Delta$ est g\'en\'erique. On prouve l'assertion pour
l'indice $i=1$. On peut certainement  trouver ${\Cal L}_1, \dots, {\Cal
L}_b \in \Pic (\Delta)$ d'ordres respectifs $e_1,\dots,e_b$, et tels
que la somme $\sum_{i=1}^b \Bbb Z{\Cal L}_i$ soit directe. Du fait que
$g'=0$, $\ppcm (e_1,\dots,e_b) = n$.

Soit un diviseur $D_i$ tel que ${\Cal L}_i = {\Cal O}(D_i)$. Il y a une fonction rationnelle
$\phi_i \in k(\Delta)^*$ telle que $e_i D_i = \Div (\phi_i)$. Si $D
= \sum_{i=1}^b \nu_i D_i$, on a 
$nD = \Div (\prod_{i=1}^b
\phi_i^{m_i\nu_i}),$
 et la classe de $D$, c'est \`a dire $\otimes_i {\cal L}_i^{\nu_i}$,  est d'ordre $n$. 
  On observe  aussi que  la fonction
rationnelle  $\psi =  \prod_{i=2}^b \phi_i^{m_i\nu_i}$ est d'ordre
$m_1$ dans $k(\Delta^*)/k(\Delta^*)^{m_1}$. En effet si pour un diviseur
strict $d$ de $m_1$, on a $\psi^{m_1\over d} = \gamma^{m_1}$, alors on peut  supposer que  
$ \psi = \gamma^d$. D\`es lors $\prod_{i=1}^b
\phi_i^{m_i\nu_i}$ serait une puissance  $d$-i\`eme avec $d/m$, ce qui \`a
exclure  car $D$ est d'ordre $n$. 
Consid\'erons alors le corps de fonctions
d'une variable $E = k(\Delta) (t,\xi)$ o\`u 
$$\xi^n =
(t-\phi_1)^{m_1\nu_1}\prod_{i=2}^b (t - \phi_1+\phi_i)^{m_i\nu_i} $$
On va voir qu'il est de
degr\'e $n$ sur $k(\Delta)(t)$, et qu'il y a une seule place $\wp$
au-dessus de $t=\phi_1$, d'indice de ramification $e=e_1$ et de
degr\'e r\'esiduel  $m_1$.   Soit $v_\wp$ la valuation normalis\'ee
en $\wp$.  On a $nv_\wp (\xi) = m_1\nu_1 e$ donc $e$ est un
multiple de $e_1$.  Il suffit pour conclure de voir
que le degr\'e r\'esiduel est au moins \'egal \`a $m_1$.  Comme
 $\left( \xi^{e_1} / (t-\phi_1)^{\nu_1} \right)^{m_1} =
\prod_{i=2}^b (t - \phi_1+\phi_i)^{m_i\nu_i}, $  dans le corps r\'esiduel
$k(\wp)$ la classe de $u =  \xi^{e_1} / (t-\phi_1)^{\nu_1}     $ v\'erifie   ${\overline u}^{m_1} = \prod_{i=2}^b \phi_i^{m_i\nu_i} = \psi$.   
Il en d\'ecoule que son degr\'e est $\geq m_1$, donc \'egal \`a $m_1$.   Le r\'esultat en d\'ecoule.\hfill\break\qed  

Dans la suite, le groupe \'etant toujours cyclique, on posera
$\mu_i = \psi_{i,1}$, qui est   le faisceau  tel que  
$\mu_i (\pi) = \sigma_i^* (\Cal L). $
On remarque  que $p_* ({\Cal L}) = 0$, plus g\'en\'eralement $p_* ({\Cal L}_j) = 0$ si $1\leq j<n$, de sorte que  $p_!({\Cal L}_j) = - [\R^1p_*({\Cal L}_j)]$. Les classes $\mu_i$, $\psi_i$ et
$\psi_{i,j}$ sont reli\'ees par des relations universelles simples:

\proclaim Proposition 10.11.    i)   On a dans $\Pic (\overline {\Cal
H}_{g,n,\xi})$    les relations
$$   \psi_{i,j} = j\mu_i - [  {j\nu_i \over  e_i}] \psi_i,\, \,\,  n
\mu_i \,= \, m_i\nu_i \psi_i,\, \,\,  n\mu_{i,k_i} = m_i\psi_i, \,\,  \hbox {\rm et}\,\, \, m_i\psi_{i,e_i} = 0 \tag
(10.37) $$
 
  \dem Notons $Q_i$ l'image de la section $\sigma_i$, de sorte que le
faisceau conormal \`a $Q_i$ est  ${\Cal O}_{Q_i}(-Q_i) \cong
\Omega^1_{D /S} \otimes {\Cal O}_{Q_i}$. Comme pour $i\ne j$, les
sections $\sigma_i$ et $\sigma_j$ sont disjointes, il vient
$$ \sigma_i^* \left({ \Cal O} (-Q_j)\right)  = \cases {\psi_i   
& si $ i=j$\cr
0  &  si $i\ne j$\cr}  \tag (10.38)$$
Du fait que $\Cal L$ est une racine $n$-i\`eme de ${\Cal O}_D ( -
\sum_{j=1}^b \, m_j\nu_j Q_j)$
(D\'efinition 9.1), il vient par image r\'eciproque, le groupe de
Picard \'etant not\'e additivement
$ n \sigma_i^* ({\Cal L}) \cong m_i\nu_i \psi_i.  $
En utilisant (9.8), on obtient imm\'ediatement la relation
$\sigma_i^* ({\Cal L}_j) = j\sigma_i^* ({\Cal L}) - [ { {j\nu_i}
\over {e_i}}]\psi_i$, c'est \`a dire  $\psi_{i,j} = j\mu_i - [ {
{j\nu_i} \over{e_i}}] \psi_i$
comme indiqu\'e.  On a  alors 
$\nu_i\mu_{i,k_i} = \nu_ik_i - \nu_i[{k_i\nu_i\over e_i}] \psi_i = \mu_i$.
 En multipliant  les deux membres par $m_i$,  la relation $n\mu_{i,k_i} = m_i\psi_i$  suit.   La derni\`ere relation d\'ecoule de la proposition 10.10.

 \qed 

Le th\'eor\`eme de Deligne-Riemann-Roch   permet   de
mettre en \'evidence  des  relations simples dans   $\Pic ({\Cal H}_{g,G,\xi})$ entre les  classes
$\lambda_j \,\,(1\leq j\leq n-1)$ et les classes $\psi_{i,\alpha}$,
et finalement les classes $\lambda_i$ et $\psi_j$ seules.  Le
raisonnement s'inspire de (\cite {43}, Theorem 3.3.4).   Un raffinement tenant compte du bord sera explicit\'e dans la section suivante. 
Notons  que  pour $1\leq j\leq n-1$ (\S \ 10.1.2)  
$\lambda_{n-j} = \det (\Bbb
E_{n-j}) = (\det R^1p_*({\Cal L}_j)^{-1} = \det Rp_* ({\Cal L}_j)$.   Le
th\'eor\`eme  de Deligne-Riemann-Roch permet d'\'evaluer le dernier
terme.  
 \proclaim Proposition 10.12. Pour tout $1\leq j \leq
n-1$, on a dans  $\Pic ({\Cal H}_{g,G,\xi})$ les relations: 
$$ 2n\lambda_{n-j} - 2n\lambda'  =  - \sum_{\alpha=1}^b jn  \langle
{jm_\alpha \nu_\alpha \over n}\rangle \mu_\alpha +  \sum_{\alpha=1}^b
n \langle {jm_\alpha \nu_\alpha \over n}\rangle\left( 1+[{
{jm_\alpha\nu_\alpha}\over  n}]\right) \psi_\alpha \tag (10.39)$$

\dem     Partons de l'identit\'e valable
pour tout $j$,   ${\Cal L}_j = {\Cal L}^{\otimes j} [ { {jB}\over  n}]$.
Supposons d'abord $j = n$,  donc  ${\Cal L}^{\otimes n} [ B] = \Cal O$.
En utilisant la notation additive, on a 
$$0 = \langle n{\Cal L} + {\Cal
O}(B) , \omega'\rangle = n\langle {\Cal L} , \omega'\rangle +
\sum_{\alpha=1}^b m_\alpha\nu_\alpha \langle Q_\alpha ,
\omega'\rangle  =   n\langle {\Cal L} , \omega'\rangle  +
\sum_\alpha m_\alpha\nu_\alpha \psi_\alpha $$
De la m\^eme
mani\`ere on  obtient imm\'ediatement 
$$0 =  \langle n{\Cal L} + B ,
{\Cal  L}\rangle 
  =  n\langle {\Cal L} , {\Cal L}\rangle + n  \rangle {\Cal
O}(B),{\Cal L}\rangle = n\langle {\Cal L},{\Cal L}\rangle +\sum_\alpha
m_\alpha\nu_\alpha \mu_\alpha$$ 
Le th\'eor\`eme de
Deligne-Riemann-Roch donne par ailleurs  $$
2n \lambda_{n-j} =  
\, 2n \det Rp_* ({\Cal L}_j) =  n\langle {\Cal L}_j , {\Cal L}_j\rangle - n\langle {\Cal
L}_j , \omega' \rangle + 2n \lambda'$$
ce qui  par substitution
de ${\Cal L}_j = j{\Cal L} + [{jB\over n} ]$ conduit \`a
$$2n\lambda_{n-j} -
2n\lambda' =  n \{ j^2\langle {\Cal L},{\Cal L}\rangle +2j\langle {\Cal
L}, [{jB\over n}]\rangle + \langle [{jB\over n}] , [{jB\over n}]
\rangle\}  - nj \langle {\Cal L},\omega' \rangle - n^2\langle  [{jB\over
n}],\omega' \rangle$$ 
donc en tenant compte des expressions du d\'ebut 
$ 2n\lambda_{n-j} - 2n\lambda'  =  $
$$-j^2 \sum_\alpha m_\alpha
\nu_\alpha\mu_\alpha    + 2jn \sum_\alpha [{jm_\alpha \nu_\alpha \over
n}] \mu_\alpha  +\left(  - n  \sum_\alpha [{jm_\alpha \nu_\alpha \over
n}] ^2      + j \sum_\alpha m_\alpha \nu_\alpha
  -  \sum_\alpha n[{jm_\alpha \nu_\alpha \over n}]
\right)\psi_\alpha   $$
soit finalement, $\langle - \rangle$ d\'esignant la
 partie fractionnaire, cette  somme se r\'eduit  \`a 
 $$
2n\lambda_{n-j} - 2n\lambda' =    - \sum_{\alpha=1}^b jn  \langle
{jm_\alpha \nu_\alpha \over n}\rangle \mu_\alpha +  \sum_{\alpha=1}^b
n \langle {jm_\alpha \nu_\alpha \over n}\rangle\left( 1+[ {
{jm_\alpha\nu_\alpha}\over  n}]\right) \psi_\alpha $$ 
qui est l'\'egalit\'e
annonc\'ee. \hfill\break
 \qed  

On peut retrouver la relation de Riemann-Hurwitz \`a l'ordre  $l = 1$ (comparer avec le th\'eor\`eme 10.9)
  \proclaim Proposition 10.13. Dans $\Pic(\overline{\Cal H}_{g,G,\xi})$ on a la relation    $ \kappa_1 = n \kappa'_1$.

 \dem   Partons de la relation de Riemann-Hurwitz $\omega = \Delta^*(\omega') \otimes \Cal O(R)$ (10.32). On a  en utilisant   la notation additive
$$\tilde \kappa_1 = \langle \omega , \omega\rangle = \langle \Delta^*(\omega') , \Delta^*(\omega')\rangle + 2\langle \Delta^*(\omega') ,{ \Cal O}(R)\rangle + \langle {\Cal O}(R) , {\Cal O}(R)\rangle $$
puis en tenant compte de la remarque qui pr\'ec\`ede, et du fait que $R_i$ et $R_j$ sont disjoints si $i\ne j$
$$\tilde\kappa_1 = n\delta^*\langle \omega' , \omega'\rangle +2\sum_j \langle\,\delta^*(\omega'), {\Cal O}(R_j)\rangle - \sum_j\,(e_j-1)^2\langle {\Cal O}(-R_j) , {\Cal O}(R_j\rangle$$
 Un d\'eterminant \'etant pris sur $\overline{\Cal H}_{g,G,\xi}$, on a par d\'efinition et utilisation des relations de la section 10.2, 
   $\langle\delta^*(\omega'), {\Cal O}(R_j)\rangle 
 = \det \left({\Cal O}_{R_j} \otimes \omega'\right) - \det {\Cal O}_{R_j} = \sum_j m_j\psi_j$.
 Pour l'autre terme on a $\langle {\Cal O}(-R_j) , {\Cal O}(R_j\rangle = \det \left({\Cal O}_{R_j}(-R_j)  \right) - \det {\Cal O}_{R_j} = \sum_{l=0}^{m_j-1} \psi_{j,k_j+le_j}  - \det {\Cal O}_{R_j}. $
 Du fait de la proposition 10.10, cette  derni\`ere expression se r\'eduit \`a
 $$= m_j\psi_{j,k_j} +  {m_j(m_j-1)\over 2} \psi_{j,e_j} - \det {\Cal O}_{R_j} = m_j\psi_{j,k_j}$$
  En sommant sur $j$ on obtient \quad $\kappa_1 = \langle \omega,\omega\rangle + \sum_j m_j\psi_{j,k_j} = $
  $$2\sum_j (e_j-1)m_j\psi_j - \sum_j (e_j-1)^2 m_j\psi_{j,k_j}  + \sum_j m_j\psi_{j,k_j}$$
et  en ins\'erant la relation $m_j \psi_j = n\psi_{j,k_j}$, on obtient finalement
  $\kappa_1 =   n\delta^*\langle \omega',\omega'\rangle + n\sum_j \psi_j = n\delta^*(\kappa'_1)$\hfill\break
\qed 
\bigskip

{\it 10.4.2. Faisceaux inversibles associ\'es aux composantes du
bord}  
\bigskip
Fixons les notations en ce qui concerne le bord de $\overline{ \Cal M}_{g,n}$ \cite {40}. Si $0\leq i\leq [g/2]$ et   $A\subset [1,n]$, on note 
$\Delta_{i,A}$ la composante irr\'eductible de $\partial \overline{\Cal M}_{g,n}$ correspondante  \`a un segment, l'un des sommets \'etant pond\'er\'e par $i$, et marqu\'e  par les
indices appartenant \`a $A$. L'autre sommet est donc pond\'er\'e
par $g-i$, et marqu\'e par le compl\'ementaire $B$ de $A$. Si $i = { g \over 2}$, alors $\Delta_{{  g \over2},A} = \Delta_{ { g \over 2},B}$, et si $i = 0$, la stabilit\'e   impose $\vert A\vert \geq 2$.  La composante associ\'ee \`a la boucle est not\'ee $\Delta_0$.    Ces composantes sont des diviseurs de Cartier, du fait de la lissit\'e  de $\overline{ \Cal M}_{g,n}$.   Elles d\'efinissent des faisceaux inversibles  qui seront   seront not\'es  ${\Cal O}(\Delta{i,A})$, ou simplement par la m\^eme lettre.  Pour \'eviter toute confusion, si $\Delta$ est un diviseur de Cartier effectif dans l'un des champs consid\'er\'es, on notera $[\Delta]\in A^1$ le cycle (ou classe) associ\'e.

 Rappelons que le bord $\overline {\Cal H}_{g,G,\xi} - {\Cal
H}_{g,G,\xi}$  admet une stratification (\S 7.2) compatible avec le morphisme $\delta: \overline {\Cal H}_{g,G,\xi}\to \overline{\Cal M}_{g',b}$ .  
 Les  strates qui ne sont pas en g\'en\'eral les composantes  irr\'eductibles, sont index\'ees
par les graphes modulaires de Hurwitz de type $(g,n)$ \`a une seule orbite d'ar\^etes
(\S 7.2). 

A un tel graphe  $\Gamma$  est associ\'e un diviseur de Cartier  $\Delta_\Gamma$
contenu dans le bord du champ de Hurwitz, et comme le champ de Hurwitz est lisse,   un faisceau
inversible ${\Cal O}(\Delta_\Gamma)$.    Si   $G$ est  cyclique, ce qui simplifie la combinatoire du bord, on va mettre 
en \'evidence  quelques  relations universelles simples, analogues \`a
celles observ\'ees par Jarvis pour les courbes \`a spin, qui
relient les faisceaux inversibles ${\Cal O}(\Delta)$ (les notations sont celles
de la section 9.2.3).  L'argument   g\'en\'eral
  s'applique  \`a un groupe de Galois arbitraire,  mais avec des notations plus compliqu\'ees.  Pour un groupe cyclique  rappelons  que  les graphes de Hurwitz sont index\'es par des partitions $\pi$ du diviseur de branchement, et qu'une telle partition induit une partition $(A,B)$ de l'ensemble des points de branchement. Si $n=p$ est premier, il y a une unique partition $\pi$ subordonn\'ee \`a $(A,B)$. La notation retenue pour  le diviseur  correspondant est $\Delta_\pi$   et ${\Cal O}(\Delta_\pi)$ d\'esigne le faisceau inversible associ\'e.   On note aussi  $e_\pi$ l'indice d'inertie  d'un point  double d'une courbe  g\'en\'erique dans $\Delta_\pi$, et on pose $m_\pi = {\vert G\vert \over e_\pi}$. Dans la proposition suivante  le groupe de Galois est arbitraire.
 
  \proclaim Proposition 10.14.  Dans le
groupe  $\Pic(\overline {\Cal H}_{g,G,\xi})$ on a
les relations o\`u le produit tensoriel porte sur les partitions subordonnŽ\'ees \`a $(A,B)$ 
$$\cases {i) \;\;\delta^* ({\Cal O} (\Delta_{i,A,B}))\; = &   $\,\bigotimes_\pi  {\Cal O}(\Delta_{ \pi})^{\otimes e_\pi}$\cr
ii)\;\; \delta^* ({\Cal O}(\Delta_0))\; = &$\bigotimes_{\pi }  
\;{\Cal O}(\Delta_{0,a,b,\pi})^{\otimes e_\pi}$ } \tag (10.40)$$

 \dem Soit un point $\pi: C\rightarrow D$ (g\'en\'erique)  dans une
composante o\`u  la base $ D$ est suppos\'ee avoir un
seul point double, dans le cas i)  deux composantes irr\'eductibles,
une seule dans le cas ii). La d\'eformation \'equivariante
universelle de $C$ a pour base une alg\`ebre de s\'eries formelles
$\Spec (k[[\tau_1,\dots,\tau_N]])$, o\`u $N = 3g-3+b$, et $\tau =
\tau_1$  repr\'esente  le param\`etre de l'orbite des points doubles
(Th\'eor\`eme 5.5). Soit la  base $\Spec (k[[t_1,\dots,t_N]])$  de la d\'eformation universelle de  $(D,Q)$, $t=
t_1$ \'etant le param\`etre de $Q$, avec les relations
$t = \tau^{e}, \,\,t_j = \tau_j\quad   {\rm si} \quad j\geq 2$.
Du fait que les \'equations locales des diviseurs irr\'eductibles
mentionn\'es dans i) sont dans ces
coordonn\'ees locales, respectivement $t = 0$ et $\tau = 0$, la
relation en d\'ecoule. La preuve de
ii) est analogue.\hfill\break
 \qed 

On peut  aussi  \'etudier la situation   du bord de $\overline {\Cal H}_{g,n,\xi}$ relativement \`a l'immersion locale non ramifi\'ee $\imath: \overline {\Cal H}_{g,n,\xi}\to \overline{\Cal M}_{g,r}$.  Rappelons que  seules   les orbites de points de ramification   sont index\'ees,  et non les points eux-m\^emes. 
\proclaim Proposition 10.15. Soit $\Delta$ le diviseur  du bord de $\overline{\Cal M}_{g,r}$. On a
$\imath^*(\Delta) = \sum_\pi m_\pi \Delta_\pi$. 

\dem Soit $\pi: C\to D$ un rev\^etement  correspondant \`a un point g\'en\'eral de $\delta_\pi$, donc le graphe modulaire de $D$ est un segment, ou une boucle. Il y a donc dans $C$ une unique orbite de points doubles, de cardinal $m = m_\pi$. Dans la d\'eformation universelle non \'equivariante de $C$, soient $t_1,\cdots,t_{m}\in R$ les param\`etres de d\'eformation de ces points. Si $t_i$ est le param\`etre de $P_i$, on  sait que le stabilisateur $H$ de $P_1$ laisse invariant $t_1$, et  qu'on a $t_i = g^*(t_1)$ si $g(P_1) = P_i$. L'\'equation de $\Delta$   dans la carte $\Spec \, R$ est  $t_1\cdots t_m = 0$. Mais la carte locale correspondante de $C\to D$ dans $\overline {\Cal H}_{g,n,\xi})$ est $\Spec \, R_G$, o\`u $R_G = R/J$, $J$ \'etant l'id\'eal engendr\'e par les $g(a)-a),\ \ a\in R, g\in G$. De la sorte l'\'equation locale de $\imath^*(\Delta)$ en $C\to D$ est 
$\prod_{i=1}^m t_i   = u^m \epsilon$
o\`u  $u $ est l'image de $t_1$, et $\epsilon$ est inversible. Comme $u = 0$ est l'\'equation locale au point consid\'er\'e du diviseur $\delta_\pi$, la conclusion suit.\hfill\break\qed 

   Faber et Pandharipande \cite {30} ont
fait observer que l'utilisation du th\'eor\`eme de Riemann-Roch   par Mumford \cite{58},   marche  dans un cadre plus g\'en\'eral,  en particulier      s'applique  au rev\^etement  universel  $\pi: {\Cal C} = \overline{\Cal C}_{g,G,\xi} \to \overline{\Cal H}_{g,G,\xi} = \Cal H$. 

Rappelons la formulation de (\cite {30}, \S 1.1). Les notations sont celles de loc.cit, en particulier  les  $B_k$  sont les nombres de Bernoulli.   Le r\'esultat qui donne le caract\`ere de Chern de $E$ s'\'enonce dans le pr\'esent contexte:
$${\rm ch} (E) = g + \sum_{l=0}^\infty {B_{2l}\over (2l)!} \left( \tilde\kappa_{2l-1} + {1\over 2}\imath_* \sum_{i=0}^{2l-2} (-1)^{i} \psi^{i} \overline {\psi}^{2l-2-i}\right) \in A^*({\Cal H}) \tag (10.41)$$
En particulier si $l = 1$   on obtient la relation importante (voir 10.4.2), dans laquelle $\delta =  j_* ([{\Cal S}]_Q)$  est la classe fondamentale du bord  $\tilde\kappa_1+ \delta  = 12\lambda.   $

 On peut aussi obtenir ces relations  par transfert,  simplement en appliquant $\imath^*$ \`a la relation correspondante dans  le champ $ \overline{\Cal M}_{g,r}$.  Soit le diviseur $\delta = \sum_\pi m_\pi \Delta_\pi \subset  \overline{\Cal H}_{g,G,\xi}$, alors de la proposition 10.15 vient
\proclaim  Proposition 10.16. On a dans   $\Pic(\overline{\Cal H}_{g,G,\xi})$ la relation de Mumford
$$\tilde\kappa_1 = \langle \omega,\omega\rangle - \delta\tag (10.42)$$
\qed 
 
Supposant $G$ cyclique (\S \ 9).  On peut en suivant le raisonnement de (\cite {43}, Theorem 4.3.8)   mettre en \'evidence une collection de  relations naturelles qui relient    les classes tautologiques avec certaines des composantes du bord, donc  qui vivent dans $\Pic (\overline {\Cal H}_{g,n,\xi})$.  Ces relations  traduisent  l'anomalie (localis\'ee sur le bord du champ de Hurwitz) qui   mesure sur un rev\^etement $\pi: C \to D$ de base $S$, la non trivialit\'e de $\Theta = {\Cal L}^{\otimes n} (B)$.

Reprenons   la construction de la section \S 7.3.3. Elle produit une courbe pr\'estable $r: \tilde D \to S$ et une factorisation par   $\rho: \tilde D \to D$. Le morphisme $\rho$ est une contraction au sens de Knudsen, il remplace tout point singulier de $\Cal L$ par une courbe exceptionnelle $E = \Bbb P^1$.  Soit ${\Cal O}_{\tilde D} (1)$ le faisceau inversible tautologique sur $\tilde D$; rappelons que $\rho_* ({\Cal O}_{\tilde D} (1)) = {\Cal L}$.

On notera   $B$ l'image inverse $\rho^* (B)$; noter que ce diviseur est disjoint du lieu exceptionnel, la notation est donc sans cons\Ž quence. Formons le faisceau inversible sur $S$
$$\langle {\Cal O}_{\tilde D} (1) , {\Cal O}_{\tilde D} (- n - B)\rangle \tag (10.43)$$
Il est clair que cette construction \'etant compatible aux changements de bases, fournit un faisceau inversible sur $\overline {\Cal H}_{g,n,\xi}$. Si $C$ est lisse  il est trivial,  d\`es lors ce  faisceau inversible  doit s'exprimer comme combinaison lin\'eaire  des classes des composantes irr\'eductibles du bord. Le r\'esultat est qui compl\`etement analogue \`a   (\cite {48}, Proposition 4.3.8)  reste valide  en fait avec des modifications \'evidentes pour ${\Cal L}_j$  si $j, \,\, (1\leq j < n)$ est premier \`a $n$. Posons $\Xi =  {\Cal O}_{\tilde D} (- n - B) $. On a  le r\'esultat important (comparer avec Jarvis \cite{43}:
\proclaim Proposition 10.17. On a 
$$\Xi = {\Cal O} \left( \sum_{\pi=NS}  {  ab\over  m} \,\, \Delta_\pi \right) \quad {\rm et }\quad  \langle \Xi , \tilde \omega \rangle = 0 \tag (10.44)$$

 \dem Noter que le terme de gauche ne fait intervenir que les composantes NS du bord.
  La preuve  reprend exactement les calculs de  loc.cit,  il n'y a pas lieu de les r\'ep\'eter. La seconde relation d\'ecoule du fait que $\tilde \omega = \rho^*(\omega)$, donc $\tilde \omega$ est trivial sur chaque diviseur exceptionnel.\hfill\break\qed 

En \'evaluant d'une autre mani\`ere  le produit d'intersection $\langle {\Cal O}_{\tilde D}(1) , \Xi \rangle$  on obtient des relations importantes entre les faisceaux inversibles tautologiques introduits pr\'ec\'edemment.    Dans le cas du champ des courbes hyperelliptiques de genre $g$ ces relations  se r\'eduisent essentiellement comme on va le voir  \`a la relation de Cornalba-Harris (\cite {40}, \S \ 6).
\proclaim  Th\'eor\`eme 10.18.  On a dans $\Pic (\overline {\Cal H}_{g,n,\xi})$, si $1\leq j<n, \,  (j,n) = 1,$ la relation
 $$\sum_{\pi=NS}    {a(j)b(j)\over  m} \,\, {\Cal O}(\Delta_\pi) = \tag (10.45)$$
$$     2n(\lambda' - \lambda_{n-j}) - \sum_{\alpha=1}^b jn  \langle
{jm_\alpha \nu_\alpha \over n}\rangle \mu_\alpha +  \sum_{\alpha=1}^b
n \langle {jm_\alpha \nu_\alpha \over n}\rangle\left( 1+[
{jm_\alpha\nu_\alpha\over  n}]\right) \psi_\alpha  $$

\dem Il suffit de prouver le r\'esultat pour $j = 1$; le cas $j$ premier \`a $n$ s'en d\'eduit par des modifications \'evidentes. Evaluons d'abord le produit d'intersection  (10.43). On a 
$$\langle    {\Cal O}_{\tilde D} (1) , {\Cal O}_{\tilde D} (- n - B)\rangle \rangle =   
  - n \langle  {\Cal O}_{\tilde D} (1)  ,  {\Cal O}_{\tilde D} (1) \rangle - \langle  {\Cal O}_{\tilde D} (1) , {\Cal O}_{\tilde D} (B)\rangle  =  - n \langle  {\Cal O}_{\tilde D} (1)  ,  {\Cal O}_{\tilde D} (1) \rangle  - \sum_\alpha m_\alpha \nu_\alpha \psi_\alpha$$
Par le th\'eor\`eme de Deligne-Riemann-Roch  \quad $ n \langle   {\Cal O}_{\tilde D} (1) ,  {\Cal O}_{\tilde D} (1) \rangle  =$
$$ 2n \det Rr_* ({\Cal L}) + \langle  {\Cal O}_{\tilde D} (n + B) , \tilde \omega\rangle - \langle  {\Cal O}_{\tilde D} (B)  , \tilde \omega \rangle  - 2n \det Rr_* \tilde \omega $$
Comme  $ \langle  {\Cal O}_{\tilde D} (n + B) , \tilde \omega\rangle  = {\Cal O} $ on obtient $\langle  {\Cal O}_{\tilde D} (B)  , \tilde \omega \rangle  = \sum_\alpha  m_\alpha \nu_\alpha \psi_\alpha$, 
 et finalement
 $$\det Rr_* (\tilde \omega) = \det Rp_* (\omega) = \lambda'$$    \qed 

Tenant compte de la relation  $n\mu_i = m_i\nu_i \psi_i$,  et apr\`es multiplication par $n$, on r\'ecup\`ere bien par restriction \`a   ${\Cal H}_{g,n,\xi}$ les relations (10.38).  Soit toujours $n=p$ premier, et de plus $g' = 0$ (rev\^etements de $\Bbb P^1$).  Par sommation des relations (10.45), on peut exprimer  $\lambda$ comme combinaison lin\'eaire des classes $\psi_\alpha$ et $\delta_\pi$. La relation obtenue \'etend la relation de Cornalba-Harris \cite{40} qui correspond \`a $p=2$.
\proclaim Corollaire 10.19.  On a dans $\Pic (\overline{\Cal H}_{g,p})$ la relation 
$$ p^2\left({p^2 - 1\over 6}\right)\sum_{\pi = NS} \, \delta_\pi = - 2p^2\lambda + p\left({p^2 - 1\over 6}\right) \sum_\alpha \psi_\alpha \tag (10.46)$$  
 
 \dem  On effectue la somme sur $j$ des relations (10.45), et on multiplie par $p$ les deux membres. On a $\lambda' = 0$, de sorte que la somme des termes de gauche dans (10.45) se r\'eduit \`a \'evaluer la somme $\sum_{j=1}^{p-1} a(j)b(j)$. Rappelons que $a(j) = p\langle {ja\over p}\rangle $, et $b(j) = p - a(j)$. La somme est  visiblement  
 $$p(\sum_{r=1}^{p-1} r - \sum_{r=1}^{p-1} r^2) = {p(p^2-1)\over 6 }$$
  Pour le terme de droite, notons que $\psi_\alpha = \nu_\alpha \mu_\alpha$. Donc apr\`es multiplication par $p$,  la somme des termes  de droite se ram\`ene facilement \`a
  $$\sum_\alpha \psi_\alpha\left(\sum_j  p^2\langle {j\nu_\alpha \over p}\rangle  -  ( p\langle {j\nu_\alpha\over p}\rangle)^2\right) = {p(p^2-1)\over 6} \sum_\alpha \psi_\alpha$$ 
  \hfill\break
  \qed

   Les relations (10.45) et (10.46) deviennent particuli\`erement simples  dans le cas hyperelliptique.  Soit $\overline{\Cal H}_g$ le champ des courbes hyperelliptiques de genre $g \geq 1$, les points de Weierstrass \'etant num\'erot\'es de $1$ \`a $2g+2$. La position remarquable de $\overline{\Cal H}_{g}$  est r\'esum\'ee par diagramme
 $$\overline{\Cal M}_{g,2g+2}\buildrel\imath\over \longleftarrow \overline{\Cal H}_g\buildrel\delta\over \longrightarrow \overline{\Cal M}_{0,2g+2}\tag (10.47)$$
 Un principe g\'en\'eral valable  pour un champ de Hurwitz  avec $g'=0$,  est que  les relations de Riemann-Hurwitz sup\'erieures jointes aux relations  (10.45) et (10.46), ram\`enent en principe un calcul dans le champ $\overline{\Cal H}_g$ \`a un calcul dans $\overline{\Cal M}_{0,2g+2}$.   

Illustrons ce principe en   montrant comment retrouver  la relation de Cornalba-Harris (\cite {40},\ \S \ 6) qui exprime dans le lieu hyperelliptique la classe $\lambda$ comme combinaison lin\'eaire des composantes du bord. Cette relation  d\'erive naturellement de la relation (10.46) avec $p=2$.  
\proclaim Proposition  10.20.
Dans le groupe $\Pic (\overline {\Cal H}_g)$ on a la relation (en notation additive)
$$8(2g+1)\lambda = 4\sum_{ \pi = R}  \, \alpha(g+1 - \alpha)     {\Delta_\pi}^R \, + \, 8\sum_{ \pi = NS} \,  \beta(g - \beta)\    {\Delta^{NS}_\pi} \tag (10.48)$$
 
\dem Pour $p=2$   (10.46) devient
$2\sum_{\pi = NS} \, \Delta_\pi = -8\lambda + \sum_{\alpha=1}^{2g+2} \psi_\alpha$.
En multipliant les deux membres par 2g+1, on peut invoquer la relation connue entre les classes $\psi$ et les composantes du bord dans $\overline{M}_{0,n}$
$$2(n - 1) \sum_{\alpha=1}^n \psi_\alpha = \sum_{(I,J), \vert I\vert = j} \, j(n - j) [{\Delta'}_{I,J}]\tag (10.49)$$
 Dans cette relation, les partitions sont ordonn\'ees. On ram\`ene cette relation dans $\Pic (\overline{\Cal H}_g)$ par $\delta^*$, en tenant compte de la proposition (10.14), qui dit que $\delta^*(\Delta'_\pi ) = \Delta_\pi$ dans le cas R, et $2\Delta_\pi$ dans le cas NS. Cela conduit  apr\`es arrangement \`a
 $$8(2g+1)\lambda = 4\sum_{ \pi = R}  \, \alpha(g+1 - \alpha)     {\Delta_\pi}^R \, + \, 2\sum_{ \pi = NS} \, 2\left[ (2\beta+1)(2g+1 - 2\beta)-(2g+1)\right]\    {\Delta^{NS}_\pi} $$
expression qui apr\`es simplification du coefficient entre crochets est  exactement le r\'esultat annonc\'e.\hfill\break\qed 
 
  Soit  $ \overline{H}_g $ le lieu des courbes hyperelliptiques, i.e le champ quotient $[\overline{\Cal H}_g/S_{2g+2}]$, avec le morphisme d'oubli de l'involution $\imath: \overline{H}_g \to \overline{\Cal M}_g$. Noter que  $\imath$ est seulement un plongement au-dessus de $H_g$, et aux points g\'en\'eriques des composantes du bord.   La relation (10.48) pouss\'ee dans $A^1( \overline{H}_g)$ est exactement la relation de Cornalba-Harris (\cite {40} , \S 6).
  \proclaim Proposition 10.21. Dans  $A^1( \overline{H}_g)$, on a la relation
 $$ (4g+2)\lambda = 
    {g\over 2} [\Delta^R_1]_Q +  \sum_{\alpha=2}^{[{g+1\over2}]}
  {\alpha(g+1 - \alpha) }   [\Delta^R_\alpha]_Q+ 2\sum_{\beta=1}^{[{g\over 2}]}  {\beta(g - \beta)}    [\Delta^{NS}_\beta]_Q   \tag (10.50) $$

   \dem Soit $\imath: \overline{\Cal H}_g \to  \overline H_g$ le morphisme quotient par $
  S_{2g+2}$. Noter que le fibr\'e de Hodge sur $\overline{\Cal H}_g$ est l'image r\'eciproque du fibr\'e de m\^eme nom sur  $  \overline{H}_g$. Prenons la premi\`ere classe de Chern des deux membres de la relation (10.47), on obtient
  $$  8(2g+1)\lambda = 4\sum_{ \pi = R}  \, \alpha(g+1 - \alpha)     {[\Delta_\pi]^R}_Q \, + \, 8\sum_{ \pi = NS} \,  \beta(g - \beta)\    {[\Delta_\pi]^{NS}}_Q$$
  $$= 4\sum_{ \pi = R}  \, {\alpha(g+1 - \alpha)   \over 2}  {[\Delta_\pi]^R}  \, + \, 8\sum_{ \pi = NS} \,  {\beta(g - \beta)\over 4}  {[\Delta_\pi]^{NS}}$$
Car dans le cas $R$, le groupe des automorphismes d'un point g\'en\'eral d'une composante est  d'ordre deux, par contre il est d'ordre quatre dans le cas $NS$.    Si  $\pi$ (partition non ordonn\'ee) est subordonn\'ee \`a la partition $(I,J)$ de $[1,2g+2]$, posons $j = \vert I\vert \leq \vert J\vert, \, \, j=2\alpha$ dans le cas pair, et $j = 2\beta+1$ si impair,  il vient par application de  $\imath^*$,  
   $$ (2g+2)! \, (16g+8) \lambda =  
     4\sum_{ \pi = R}  \, {\alpha(g+1 - \alpha)\over 2}   {2g+2\choose  j}i_*{[\Delta_\pi]^R}  \,  +   \,8\sum_{ \pi = NS} \,  {\beta(g - \beta)\over 4} {2g+2\choose j}     i_*{[\Delta_\pi]^{NS}} $$
   Il est clair que le degr\'e de $\Delta^R_\pi$ (resp.  $\Delta^{NS}_\pi$) sur son image $\Delta^R_j$ (resp. $\Delta^{NS}_j$) est $j!(2g+2-j)!$ sauf si $j=2$, car dans ce cas une courbe hyperelliptique $C\in \Delta^R_\pi$ a une composante rationnelle, munie d'une involution. Les deux points fixes de l'involution sont les points marqu\'es. Les points d'intersection $P',\, P''$ avec l'autre composante, de genre $g-1$ sont \'echang\'es par l'involution. Il y a une seconde involution qui a pour point fixe $P', \, P''$ et qui commute avec la premi\`ere. Donc l'\'echange des deux points marqu\'es $P', \, P''$ conduit \`a une courbe isomorphe. Il vient finalement en rappelant que le cas $R$ correspond \`a $j$ pair
$$ (2g+2)! \, (16g+8) \lambda =  
     4\sum_{ \pi = R}  \, {\alpha(g+1 - \alpha)\over 2}   {2g+2\choose  j}i_*{[\Delta_\pi]^R}  \,  +   \,8\sum_{ \pi = NS} \,  {\beta(g - \beta)\over 4} {2g+2\choose j}     i_*{[\Delta_\pi]^{NS}}  $$
   Il est clair que le degr\'e de $\Delta^R_\pi$ (resp.  $\Delta^{NS}_\pi$) sur son image $\Delta^R_j$ (resp. $\Delta^{NS}_j$) est $j!(2g+2-j)!$ sauf si $j=2$, car dans ce cas une courbe hyperelliptique $C\in \Delta^R_\pi$ a une composante rationnelle, munie d'une involution. Les deux points fixes de l'involution sont les points marqu\'es. Les points d'intersection $P',\, P''$ avec l'autre composante, de genre $g-1$ sont \'echang\'es par l'involution. Il y a une seconde involution qui a pour point fixe $P', \, P''$ et qui commute avec la premi\`ere. Donc l'\'echange des deux points marqu\'es $P', \, P''$ conduit \`a une courbe isomorphe. Il vient finalement 
      $$\eqalign{(2g &+2)! \,(16g+8)\lambda = \cr
   & 4g {(2g+2)!\over 4} [\Delta^R_1] + 4\sum_{\alpha=2}^{[{g+1\over2}]}
  {\alpha(g+1 - \alpha)\over 2}   (2g+2)! [\Delta^R_\alpha]+ 8\sum_{\beta=1}^{[{g\over 2}]}  {\beta(g - \beta)\over 4} (2g+2)!  [\Delta^{NS}_\beta] \cr}   $$
   En revenant aux classes fondamentales des composantes 
  $$= 4g {(2g+2)!\over 2} [{\Delta_1^R}_Q] + 4\sum_{\alpha=2}^{[{g+1\over2}]}
  {\alpha(g+1 - \alpha) }   (2g+2)! [{\Delta_\alpha^R}]_Q+ 8\sum_{\beta=1}^{[{g\over 2}]}  {\beta(g - \beta)} (2g+2)!  [{\Delta_\beta^{NS}} ]_Q    $$
  et finalement en simplifiant par  $8(2g+2)!$,  on obtient la relation de Cornalba -Harris 
  $$(8g+4)\lambda = 
     g [{\Delta_1^R}]_Q +  \sum_{\alpha=2}^{[{g+1\over2}]}
  2{\alpha(g+1 - \alpha) }   [{\Delta_\alpha^R}]_Q+ 4\sum_{\beta=1}^{[{g\over 2}]}  {\beta(g - \beta)}    [{\Delta_\beta^{NS}}]_Q    $$
 \qed

\medskip
 {\it 10.4.3.  Int\'egrales de Hodge - Hurwitz} 
\bigskip
 
 On se limite en premier  au champ des courbes hyperelliptiques.
Une  autre illustration du principe  de correspondance  utilis\'e  ci-dessus revient \`a observer que   des calculs d'intersection avec le lieu hyperelliptique $\overline{H}_g\subset \overline {\Cal M}_g$ une fois transport\'es sur le champ hyperelliptique se ram\`enent   \`a des calculs dans $\overline{\Cal M}_{0,2g+2}$. Cela montre que l'objet naturel pour les calculs est $\overline{\Cal H}_g$ plut\^ot que $\overline {H}_g$.  
La m\'ethode utilis\'ee est une amplification du calcul de Faber-Pandharipande (\cite {30}). Elle conduit    simplement \`a l'\'evaluation d'int\'egrales    
$$  \int_{\overline{{\Cal H}^1}_g} \! \kappa_1\mu_1^{2g-2}  = {(2g -1)^2\over 2^{2g}(2g+1)!  }\quad { \rm  et }\quad  \int_{\overline{{\Cal H}^1}_g} \!\! \mu_1^{2g-1}  = {1\over 2^{2g}(2g+1)! }  \tag (10.51)$$
en notant  $\overline{{\Cal H}^1_g}$ le champ des courbes hyperelliptiques avec un point de Weierstrass marqu\'e.  Rappelons que le faisceau inversible $\mu_i \in \Pic(\overline{\cal H}_g)$ (sa classe de Chern)  est l'image r\'eciproque du faisceau ${\cal L}_i\in \Pic(\overline{\cal M}_{g,2g+2})$ (\S \ 10.2). On conserve la notation $\mu_i$ pour \'eviter toute confusion avec les $\psi_i$ qui proviennent du bas.  A.J. Bene a obtenu des r\'esultats  analogues, mais par une approche combinatoire  bas\'ee sur   cellulation de $\overline{\Cal H}_g$ au moyen des graphes \'epais,  et le sch\'ema d'int\'egration de Penner   \cite {10}. A la diff\'erence de  notre m\'ethode, il travaille avec les classes combinatoires $W_a$ de Witten, et \'etudie l'intersection de ces classes avec le lieu hyperelliptique.

 Rappelons que    dans le champ $\overline{\Cal H}_g$ les points de Weierstrass sont  marqu\'es de $1$ \`a $2g+2$. Il est commode de noter $\overline{{\Cal H}^1}_g$ le champ dont les objets sont les courbes hyperelliptiques  stables marqu\'ees par un seul point de Weierstrass. Plus g\'en\'eralement on peut consid\'erer la champ $\overline{\cal H}^j_g$ form\'e des courbes hyperelliptique stables avec $j \, (1\leq j\leq 2g+2)$ points de Weierstrass marqu\'es. Il est ais\'e de justifier son existence en adaptant les arguments de la section 6.2. Noter que si $j\leq 2g$ il est possible que sur une courbe hyperelliptique stable appartenant \`a $\overline{\cal H}^j_g$ deux points de Weierstrass libres s'effondrent. Par oubli du marquage par le point d'indice $j$ on obtient un morphisme  $\pi_j: \overline{\cal H}^j_g\to \overline{\cal H}^{j-1}_g$ de degr\'e $2g+3-j$.   On notera aussi   $\overline H_g \subset \overline{\Cal M}_g$   l'image stable non marqu\'ee de $\overline{\Cal H}_g$, i.e le lieu hyperelliptique. On a donc la factorisation    de $\imath: \overline{\Cal H}_g \rightarrow \overline{\Cal M}_g$, oubli du marquage par les points de Weirstrass,  en
   $$\overline{\Cal H}_g\buildrel\pi\over \rightarrow \overline{{\Cal H}^1}_g \rightarrow \overline{H}_g \subset \overline{\Cal M}_g $$
  
Le morphisme $\pi$ est l'oubli du marquage par les points de num\'eros $2g+2,\cdots,2$ dans l'ordre indiqu\'e, c'est  \`a dire la composition $\pi_2\cdots\pi_{2g+2}$, son degr\'e est $(2g+1)!$. Pour ramener  par le principe de correspondance le calcul \`a un calcul dans $A^\bullet(\overline{\Cal M}_{0,2g+2})$, on note que 
$$\int_{\overline{{\Cal H}^1}_g}\!\! \xi = {1\over (2g+1)!} \int_{\overline{\Cal H}_g}\!\! \pi^*(\xi)\tag (10.52)$$ 
Il faut   expliciter l'image r\'eciproque par $\pi$ des classes $\mu_1$ et $\kappa_a$.  Par utilisation it\'er\'ee de la relation g\'en\'erale $\pi_{n+1}^*(\psi_1) = \psi_1 - [D_{1,n+1}]$ (10.25), on trouve que, avec nos notations 
$$\pi_j^*(\mu_1) = \mu_1 - [D_{1,j}] \cap [\overline{{\Cal H}^j}_g]$$
 (l'indice $j$ signifie que $j$ points de Weierstrass sont marqu\'es). Comme $D_{1,j}$ est le lieu des courbes obtenues en collant une ''bulle'' ($\cong \Bbb P^1$)  contenant les points de Weierstrass de num\'eros $1$ et $j$  \`a une courbe de genre $g$ marqu\'ee par $j-1$ points,  l'intersection avec le lieu hyperelliptique est vide. On obtient donc que $\pi^*(\mu_1) = \mu_1$, ce qui nous autorise \`a noter cette classe $\mu_1$ sans pr\'ecision suppl\'ementaire. Le m\^eme argument montre que pour la classe $\kappa_a$ on a 
$$\pi^*(\kappa_a) = \kappa_a - \sum_{\alpha = 2}^{2g+2} {\mu_\alpha^{a}} \tag (10.53)$$

Traitons  en premier  la seconde int\'egrale, qui est imm\'ediate.   Elle se ram\`ene  au calcul de 
 $\int_{\overline{\Cal H}_g}\!\! \mu_1^{2g-1}$. 
  Vu que $2\mu_1 = \delta^*(\psi_1)$,  elle   s'obtient imm\'ediatement partant  du r\'esultat connu 
$$\int_{\overline{\Cal M}_{0,n}} \psi_1^{\alpha_1}\cdots \psi_n^{\alpha_n} = {(n-3)!\over \alpha_1!\cdots\alpha_n!}  \quad(\alpha_1+\cdots+\alpha_n = n-3)\tag (10.54)$$
on trouve en tenant compte du fait que le degr\'e de $\delta$ est ${1\over 2}$
$$\int_{\overline{{\Cal H}^1}_g} \!\! \mu_1^{2g-1}  = {1\over 2^{2g}} \int_{\overline{\Cal M}_{0,2g+2}}\!\! \psi_1^{2g-1} = {1\over 2^{2g}(2g+1)!} $$
 
Pour la premi\`ere int\'egrale,   on montre  qu'on a plus g\'en\'eralement:
\proclaim Th\'eor\`eme  10.22. Avec les notations pr\'ec\'edentes, si $0\leq a \leq 2g-1$, on a  
$$\int_{\overline{{\Cal H}^1}_g} \!\! \kappa_a\mu_1^{2g-1-a}  \, = \, {1\over 2^{2g-1-a} (2g+1)!} \left({2g\choose a+1} - {2g+1\over 2^{a+1}} {2g-1\choose a}\right) \tag (10.55)$$

\dem En utilisant la relation (10.52) on voit que l'int\'egrale (10.55) se ram\`ene \`a  une int\'egrale sur le champ $\overline {\Cal H}_g$ classifiant les courbes hyperelliptiques avec points de Weierstrass marqu\'es
$$(2g+1)! \int_{\overline{{\Cal H}^1}_g} \!\! \kappa_a\mu_1^{2g-1-a}  = \int_{\overline{\Cal H}_g}\!\! (\kappa_a - \sum_2^{2g+2} \mu_\alpha^{a})\mu_1^{2g-1-a} = \int_{\overline{\Cal H}_g}\!\! \kappa_a\mu_1^{2g-1-a} - \sum_2^{2g+2} \int_{\overline{\Cal H}_g}\!\!   
\mu_\alpha^{a}\mu_1^{2g-1-a}$$
En tenant compte de la relation de Riemann-Hurwitz \`a l'ordre $a$, qui dans le cas hyperelliptique s'\'ecrit sous la forme simple $\kappa_l = 2\delta^*(\kappa'_l)$, et de la formule de projection, on a  
$$\int_{\overline{\Cal H}_g}\!\!  \kappa_a{\mu_1}^{2g-1-a} = {1\over 2^{2g-2-a}} \int_{\overline{\Cal H}_g} \!\! \delta^*(\kappa'_a {\psi_1}^{2g-1-a}) = {1\over 2^{2g-2-a}} \int_{ \overline{\Cal M}_{0,2g+2}}\!\! \kappa'_a{\psi_1}^{2g-1-a} \delta_* [1]$$
Mais $\delta_*[1] =   \delta_*([\overline{\Cal H}]_Q)  = {1\over 2}\delta_*([\overline{\Cal H}_g]) = {1\over 2}[\overline{M}_{0,2g+2}]$,  de sorte que   l'int\'egrale se r\'eduit \`a 
$$  {1\over 2^{2g-1-a}}\int_{\overline{\Cal M}_{0,2g+2}}\!\! \kappa'_a{\psi_1}^{2g-1-a}$$
De la  m\^eme mani\`ere  on obtient 
$$\int_{\overline{\Cal H}_g}\!\!    \mu_\alpha^{a}\mu_1^{2g-1-a} =  {1\over 2^{2g-1}}\int_{\overline{\Cal M}_{0,2g+2}}\!\! \psi_\alpha^{a}{\psi_1}^{2g-1-a} = {1\over 2^{2g}} {2g- 1 \choose a}$$
  Pour conclure, notons    le r\'esultat \'el\'ementaire suivant, cons\'equence facile de l'\'equation des cordes:
\proclaim Lemme 10.23. Soit pour $4\leq a+3\leq n,$ $\tau_{a,n} = \int_{\overline{\Cal M}_{0,n}} \!\! \kappa_a \psi_1^{n-3-a}$.  On a $ \tau_{a,n} =     {n-2\choose  a+1}$.
 
\dem On proc\`ede  par r\'ecurrence. Soit le morphisme   $\pi_{n+1}: \overline {M}_{0,n+1} \rightarrow {\overline M}_{0,n}$ oubli du point $x_{n+1}$. On a $\pi_{n+1}^*(\kappa_a) = \kappa_a - \psi_{n+1}^{a}$, donc 
$$\pi_{n+1}^*(\kappa_a)\psi_1^{n-a-2} = \kappa_a\psi_1^{n-a-2} - \psi_1^{n-a-2}\psi_{n+1}^{a}$$
En int\'egrant, et par la formule de projection,  on obtient
$$  \int_{\overline{\Cal M}_{0,n+1}}\!\!  \pi_{n+1}^*(\kappa_a)\psi_1^{n-a-2} = \tau_{a,n+1} - \int_{\overline{\Cal M}_{0,n+1}}\!\! \psi_1^{n-a-2}\psi_{n+1}^{a} =
 \int_{\overline{\Cal M}_{0,n}}\!\!  \kappa_a \pi_{n+1,*}(\psi_1^{n-a-2})$$
Mais par utilisation de l'\'equation des cordes \cite{4} 
$${\pi_n}_*(\psi_1^{a_1}\cdots \psi_{n-1}^{a_{n-1}})  = \sum_{j, a_j>0} \, \psi_1^{a_1}\cdots \psi_j^{a_j-1}\cdots \psi_{n-1}^{a_{n-1}} \tag (10.56)$$
on a en particulier  $\pi_{n+1,*}(\psi_1^{n-a-2}) = \psi_1^{n-a-3}$, d'o\`u  finalement la relation pour $a+3\leq n,$
$\tau_{a,n+1} = \tau_{a,n} + {n-2\choose a},  $
donc finalement $\tau_{a,n} = {n-3\choose a}+ {n-4\choose a} +  + \cdots + {a\choose a} = {n-2\choose a+1}$.\hfill\break\qed

Pour conclure la preuve de (10.55), les r\'esultats qui pr\'ec\`edent   donnent apr\`es division par le facteur $(2g+1)!$
 $$\int_{\overline{{\Cal H}^1}_g} \!\! \kappa_a\mu_1^{2g-1-a}  ={1\over (2g+1)!}\left( {1\over 2^{2g-1-a} } {2g\choose a+1} - {2g+1\over 2^{2g}} {2g-1\choose a}\right)$$ qui est exactement l'expression (10.55).  \hfill\break\qed  
 
 Ce calcul pour $a = 0$, a pour cons\'equence (Faber-Pandharipande \cite{30}, prop 4) la relation
 $$1+\sum_{g\geq 0} t^{2g} \int_{\overline{\Cal M}_{g,1}} \,   \left(\sum_{i=1}^g (-1)^{i+1}(2^{g+1-i} - 1)\psi_1^{g-i}\lambda_{i-1}\right) \psi_1^{2g-1} = {\sin (t/2)\over t/2} \tag (10.57)$$
 
 Pour $a=1$, A.J. Benne \cite {10} obtient  de mani\`ere analogue, via la relation de Mumford qui d\'ecrit la classe fondamentale $[{\cal H}_g]$, une relation (plus compliqu\'ee) reliant certaines  int\'egrales de Hodge
 $$  \int_{\overline{\Cal M}_{g,1}} \,   \left(\sum_{i=1}^g (-1)^{i+1}(2^{g+1-i} - 1)\psi_1^{g-i}\lambda_{i-1}\right) \kappa_1\psi_1^{2g-2} ={14g^2 - 11g+3\over 3.2^{2g}(2g+1)!} \tag (10.58)$$
On peut s'interroger   sur la validit\'e d'une relation analogue pour $a\geq 2$.  Si $a = 2g-1$, un calcul similaire donne l'int\'egrale de $\kappa_{2g-1}$  \'evalu\'ee sur le lieu hyperelliptique  $\overline H_g$, on trouve
$$\int_{\overline H_g} \! \kappa_{2g-1} = {1\over (2g+2)!} - {1\over 2^{2g} (2g+1)!} \tag (10.59)$$
Si  $a=2$, c'est le calcul de Mumford $\int_{\overline M_2} \! \kappa_3 = {1\over 1152}$.

Soit maintenant   $\overline {\Cal H}_{g,p,\xi}$ le champ  des rev\^etements de degr\'e $p$, et de genre $g$,  de $\Bbb P^1$. On suppose que la donn\'ee de ramification, versus le diviseur de branchement, est $B = \sum_{i=1}^b \nu_i Q_i \,\, (1\leq \nu_i<p)$.   On a  donc $b = {2g-2+2p\over p-1}$. On va observer que la relation (10.46) permet d'\'evaluer,   au prix de relations compliq\'ees,  certaines int\'egrales de Hodge sur $\overline {\Cal H}_{g,p,\xi}$, celles contenant $\lambda$. Le r\'esultat sera en g\'en\'eral sensible  au choix de la composante $\overline{\Cal H}_{g,p,\xi}$, c'est \`a dire de $\xi$, \`a la diff\'erence du calcul  par  Bryan-Graber-Pandharipande\footnote{The orbifold quantum cohomology of $\Bbb C^2/\Bbb Z_3$ and Hurwitz-Hodge integrals arXiv:math.AG/0510335.}  de $\int_{\overline{\Cal H}_{g,3,\xi}} \lambda_{g-1}$  le (cas $p=3$). Posons
$$B_{g,\xi} = \int_{\overline {\Cal H}_{g,p,\xi}} \!\! \lambda^{b-3}$$
Le cas $p=2$ (hyperelliptique) est un calcul de Faber. Pour initialiser cette suite, il faut calculer $B_{g,\xi}$ lorsque  $b = 3$, c'est \`a dire $g = {p-1\over 2}$. Dans ce cas le champ est ponctuel.  Soit   $e = \vert \Aut (C_{a,b,c}\to \Bbb P^1)\vert $, o\`u  $C_{a,b,c} $ est la courbe $y^p = x^{a}(1-x)^b$, avec $1\leq a,b,c< p, \, a+b+c\equiv 0 \pmod p$, et $C_{a,b,c} \to \Bbb P^1$, le rev\^etement  $(x,y) \mapsto x$. On trouve imm\'ediatement   
$$ e = \cases{p & si $a,b,c$ distincts\cr 2p & si deux entre les trois sont \'egaux\cr 18 & si $a=b=c, \, p=3$\cr}\tag (10.60)$$

On a a $[\overline {\Cal H}_{g,p,\xi}]_Q = {1\over e}[\overline {\Cal H}_{g,p,\xi}]$, et  $B_{g,\xi} = {1\over e}$.  
La relation (10.46) jointe \`a (10.49) montre que  qu'il faut au prŽalable \'evaluer  les int\'egrales du bord $\int_{\Delta_\pi} \! \lambda^{b-4}$.
Il y a deux cas \`a consid\'erer, selon que l'unique orbite de point double est  du type R, ou NS. Le premier cas correspond \`a une partition $\pi$ de $B$ en $B = B_1 + B_2, \,\, B_i = \sum_{i\in I_i} \nu_iQ_i$, et $\sum_{i\in I_i} \equiv 0 \pmod p$. Avec des notations \'evidentes, on a $g = g_1+g_2+p-1$, et la proposition 7.13 dit que la composante $\Delta_\pi$  est image du morphisme $\imath = \xi_\Gamma$,  d\'ecrit par le diagramme
$$\overline {\Cal H}_{g_1,p}\times \overline{\Cal H}_{g_2,p}\leftarrow \overline {\Cal C}_{g_1,p}\times \overline{\Cal C}_{g_2,p}\buildrel \imath\over \rightarrow  \overline{\Cal H}_{g,p} 
\tag (10.61)$$
dans lequel $\overline {\Cal C}_{g_1,p}\to \overline {\Cal H}_{g_1,p}$ d\'esigne le rev\^etement universel  (\S \ 6.4.2).  Le morphisme $\imath$ identifie deux par deux les points des orbites des points s\'electionn\'es \`a isotropie triviale. Dans le cas NS, on a $g = g_1+g_2+1$, et (10.60) se r\'eduit \`a 
$$\overline {\Cal H}_{g_1,p}\times \overline{\Cal H}_{g_2,p}\buildrel\imath\over  \to \overline{\Cal H}_{g,p} \tag (10.62)$$
Si $\pi: B = B_1+B_2$ est  une partition de type NS, donc $\sum_{j\in I_1} \nu_j \equiv \eta \pmod p$, on notera $\xi_1$ la donn\'ee de Hurwitz d\'efinie par $B_1+\eta P_1$, et   $\xi_2$ celle d\'efinie par $B_2 + (p-\eta)P_2$ (voir \S \ 7.4). 

\proclaim Lemme 10.24.   On a $\int_{\Delta_\pi}\! \lambda^{b-4} = \cases {0 & si $\pi $ est R\cr 
{b-4\choose b_1-2} B_{g_1,\xi_1}B_{g_2,\xi_2} & si $\pi$ est NS\cr}$. 

\dem Dans le cas NS, c'est imm\'ediat du fait que $\imath$ \'etant le ''clutching'' morphisme, on sait  \cite {48},  qu'avec des notation \'evidentes
$$ \imath^*(\lambda) = \lambda_1+\lambda_2$$
donc  en notant que le $b$ de $\overline{\Cal H}_{g_i,p}$ est $b_i+1$ 
$$\int_{\Delta_\pi} \! \imath^*(\lambda^{b-4}) = \int_{\Delta_\pi}\! (\lambda_1+\lambda_2)^{b-4} = { b-4 \choose  b_1-2} \int_{\overline {\Cal H}_{g_1,p}\times \overline{\Cal H}_{g_2,p}}\lambda_1^{b_1-2}\lambda_2^{b_2-2} =  { b-4 \choose b_1-2} B_{g_1,\xi_1}B_{g_2,\xi_2}$$
Dans le cas R le morphisme $\imath: \overline {\Cal C}_{g_1,p}\times \overline{\Cal C}_{g_2,p} \rightarrow  \overline{\Cal H}_{g,p} $ identifie les paires de points $(\sigma^j(x_1),\sigma^j(x_2))$ si $x_1$ (resp. $x_2$) est la section de $p_1: \overline{\Cal C}_{g_1,p}\to \overline{\Cal H}_{g_1,p}$ (resp. 
$p_2: \overline{\Cal C}_{g_2,p}\to \overline{\Cal H}_{g_2,p}$).  Dans ce cas   le fibr\'e vectoriel $\imath^*(\Bbb E)$ est d\'ecrit par une extension
$$0 \to p_1^*(\Bbb E_1) \bigoplus p_2^*(\Bbb E_2) \to \imath^*(\Bbb E) \buildrel {\rm res}\over \to {\Cal O}^{p-1} \to 0 \tag (10.63)$$
l'application res \'etant induite par le r\'esidu.  En effet, soit un rev\^etement $\pi: C \to D$  au-dessus de $k$ dans l'image de $\imath$, avec $C = C_1\cup C_2$, et une orbite de $p$ points doubles, l'orbite du point $x = \imath (x_1,x_2)$. Dans la d\'ecomposition en sous-espaces propres (proposition 10.3),   ${\rm H}^0(C,\omega_C) = \bigoplus_{j=1}^{p-1} {\rm H}^0(C,\omega_C) ^j$, il est imm\'ediat de voir qu'on a 
$$\dim {\rm H}^0(C,\omega_C) ^j = \dim {\rm H}^0(C_1,\omega_{C_1}) ^j + \dim {\rm H}^0(C_2,\omega_{C_2}) ^j + 1$$
D'autre part le r\'esidu en $x$ donne pour tout $1\leq j\leq p-1$ une application  
${\rm res}_x: {\rm H}^0(C,\omega_C) ^j \to k$, qui  conduit \`a  la suite exacte
$$0\to  {\rm H}^0(C_1,\omega_{C_1}) ^j \oplus  {\rm H}^0(C_2,\omega_{C_2}) ^j \to  {\rm H}^0(C,\omega_{C}) ^j \buildrel {{\rm res}_x}\over \to k \to 0 \tag (10.64)$$
la suite exacte  (10.62) en d\'ecoule. On a donc
$$\int_{\overline {\Cal C}_{g_1,p}\times \overline{\Cal C}_{g_2,p}} \imath^*(\lambda^{b-4}) = \int_{\overline {\Cal C}_{g_1,p}\times \overline{\Cal C}_{g_2,p}} (\lambda_1+\lambda_2)^{b-4} = {b-4\choose b_1-2} \int_{\overline {\Cal C}_{g_1,p}\times \overline{\Cal C}_{g_2,p}} \lambda_1^{b_1-2}\lambda_2^{b_2-2}$$
Mais  dans cette configuration les classes $\lambda_i$  sur $\overline{\Cal C}_{g_i,p}$ sont les images inverses des  classes de m\^eme nom de $\overline{\Cal H}_{g_i,p}$, donc pour une raison de dimension, l'int\'egrale est nulle. \hfill\break \qed

On peut ainsi n\'egliger les composantes de type R.  Les relations (10.46) et (10.49)  en tenant compte du lemme 10.24,  conduisent facilement \`a la relation  de r\'ecursion
$${24(b-1)\over p^2-1} B_{g,\xi} = \sum_{\pi = (I_1,I_2) } (b_1b_2 - 2(b-1)){b-4\choose b_1-2} B_{g_1,\xi_1}B_{g_2,\xi_2} \tag (10.65)$$
somme portant sur les partitions de type NS.

\vskip.6cm
 
\centerline{\bf R\'EFERENCES}
\vskip.3cm


\ref{1}   D. Abramovich,  A. Vistoli,  {\it  Complete  moduli 
 for  families  over  semistable  curves, }  \hfill\break
[arXiv: math.AG/9811059]

 \ref{2}   D. Abramovich, A. Corti, A. Vistoli,  {\it  Twisted
bundles and Admissible covers,}  Commun. Algebra  {\bf 31}, No.8
 (2003)   3547--3618

\ref{3} D. Abramovich, T. Jarvis,  {\it  Moduli of twisted spin curves,}
   Proc. AMS {\bf 131}, No.3 \  (2003) \ 685--699

\ref{4} A. Arsie, A.Vistoli, { \it    Stacks of cyclic covers
of  projective spaces,}      
[arXiv: math.AG/0301008]

\ref{5} M. Asada, M Matsumoto, T. Oda, {\it Local monodromy
on the fundamental groups
of algebraic curves along a denerate stable curve},   J. Pure and
Applied Algebra,  {\bf 103} (1995)   235--283

\ref{6}  H. Bass, {\paper Covering theory for graphs of groups},  J. Pure and Applied Algebra, {\bf 89}  (1993)   3--47 

\ref{7} A. Beauville,  {\paper    Prym varieties and the Schottky
problem},   Invent.Math.  (1977)  149--196  

\ref{8} S. Beckmann,  {\paper Ramified primes in the field of
moduli of branched
coverings of curves},   J. of Algebra  {\bf 125}   (1989)  236--255 

\ref{9} A.A. Beilinson, Y.I. Manin,  {\paper The Mumford Form
and the Polyakov Measure in String theory,}   Commun. Math. Phys.
 {\bf 107}   (1986)  359--376  
 
\ref{10} A.J. Bene, { \paper  Combinatorial classes, Hyperelliptic loci, and Hodge integrals,}    \hfill\break
[arXiv math.GT/0610603]  

\ref{11}  J. Bertin, {\paper Compactification des sch\'emas de Hurwitz,}
  C. R. Acad. Sci., Paris, S\'er. I  {\bf 322}  No.11  (1996)
  1063--1066  

\ref{12}   J. Bertin, A. M\'ezard,  { \paper  D\'eformations
 formelles des  rev\^etements  sauvagement  ramifi\'es  de
courbes alg\'ebriques,}  
Invent. Math.  {\bf141} (2000)  195--238  

\ref{13}  M. Boggi, M. Pikaart, { \paper Galois covers of moduli
of curves,}   Comp. Math.  {\bf 120}, No.2  (2000)  171--191

\ref{14}   J.L. Brylinski,  { \paper  Propri\'et\'es de ramification
\`a l'infini du groupe de Teichm\"{u}ller,}  Ann. Sci. Ec. Norm. Sup.
{\bf 12}  (1979)  295--333  

\ref{15}  K. Combes, D. Harbater,   {\paper Hurwitz families and arithmetic
Galois groups,}  Duke Math. J.  {\bf 52}  (1985)   821--839 

\ref{16}  P. D\`ebes,   {\paper Arithm\'etique et espaces de modules de
rev\^etements,}  Number theory in progress. Proc. Int. Conf. in honor
of Andrzej Schinzel, Zakopane, Poland, June 30-July 9 1997
(Gy\"ory, K\'alm\'an, et al., eds.),  de Gruyter (1999)   75--102

\ref{17}   A.J. De Jong, M. Pikaart, { \paper Moduli of curves
with non abelian level structures,}   The moduli space of curves
(R. Dijkgraaf, C. Faber, and G. Van der Geer eds.),  Prog. in Math.
{\bf 129}   Birkh\"{a}user  (1995)  483--510

\ref{18} P. Deligne, { \paper Le d\'eterminant de la
cohomologie,}    Contemp. Math. {\bf 67}  (1987)  
93--178  

\ref{19} P. Deligne, { \paper Le lemme de Gabber,}  
in  Sem. sur les pinceaux arithm\'etiques: la conjecture de Mordell,
Ast\'erisque  {\bf 127}  (1985)   131--150  

\ref{20}    P. Deligne, D. Mumford, { \paper The irreducibility of
the space of curves of a given genus,}  Publ. Math. IHES.
{\bf 36}  (1969)   75--100

\ref{21}   M. Demazure, P.Gabriel, { \paper Groupe Alg\'ebriques,}
  North-Holland  (Amsterdam) (1970)  

\ref{22}  S. Diaz, D. Edidin,   {\paper Towards the homology of
Hurwitz spaces,} J. Diff. Geom.  {\bf 43}  (1996)   66--97 

\ref{23}   R. Dijkgraaf,   {\paper Mirror symmetry and elliptic
curves,}  The Moduli space of Curves (R.Dijkgraaf, C.Faber, and G.van der Geer, eds.),
Prog. in Math.  {\bf 129}   Birkh\"{a}user  (1995)   149--163  

\ref{24}  C.J. Earle,   {\paper On the moduli of closed Riemann
surfaces with symmetries,}  Math. Studies  {\bf 66 }  (1971)   119--130 

\ref{25}  A. Edmonds,  { \paper Surface symmetry,}  
Mich. Math. J. {\bf 29}    (1982)   171--183  

\ref{26}   T. Ekedahl,  {\paper Boundary behaviour of Hurwitz schemes,}
  The Moduli space of Curves (R. Dijkgraaf, C. Faber, and G. van der Geer,
eds.), Prog. in Math. {\bf 129}   Birkh\"{a}user  (1995)  173-198

\ref{27}   T. Ekedahl, S. Lando, M. Shapiro, A. Vainshtein,
{\paper Hurwitz numbers and intersections on moduli spaces of curves,}
  Invent. Math.  {\bf 146}, No.2  (2001)  297--327  

\ref{28}  G. Ellingsrud, K. Lonsted,   {\paper An equivariant
Lefschetz formula for finite reductive groups,}   Math. Ann.
{\bf 251}  (1980)    253--261  

\ref{29}   M. Emsalem,   {\paper Familles de rev\^etements de la
droite projective,}  Bull. Soc. Math. Fr. {\bf 123}  (1995)
 47--85 

\ref{30}   C. Faber, R. Pandharipande,  { \paper  Hodge integrals and Gromov-Witten  theory}
  Invent. Math.  {\bf 139}  (2000)  109--129  

\ref{31}   G. Faltings, { \paper Moduli stacks for bundles on
semistable curves,}   Math. Ann.  {\bf 304} (1996)  
489--515

\ref{32}  M. Fried,   {\paper Fields of definition of function
fields and Hurwitz families, groups as Galois groups,}   Comm. in Alg.,
 {\bf 5}  (1977)  17--81  

\ref{33}   M. Fried, { \paper Introduction to modular towers,}  
Comtemp. Math  {\bf 186}  (1995)   111--171  

\ref{34}   M. Fried, H. V\"{o}lklein,   { \paper  The inverse Galois
problem and rational points on moduli spaces,}   Math. Annn.
 {\bf 290}  (1991)   771--800

\ref{35}    W. Fulton,   {\paper Hurwitz schemes and irreducibility
of moduli of algebraic curves,}  Ann. of Math.  {\bf 90}  (1969)
 542--575  

\ref{36}   B. van Geemen, F. Oort,   { \paper A compactification of
a fine moduli space of curves,}   Moduli of Abelian varieties, 
Proc. Conf. Texel Island, Netherlands, April 1999 (C. Faber, G. van der Geer,
F. Oort, eds.), Prog. in Math. {\bf 195} Birkh\"{a}user (2001) 

\ref{37}  D. Gieseker,  { \paper Lectures on Moduli of curves,}  
Tata Institute lectures notes  {1982}  

\ref{38}  T. Graber, R. Vakil,   {\paper Hodge integrals and Hurwitz
numbers via virtual localization,}   Compos. Math. {\bf 135}, No.1
 {2003}   25--36  

\ref{39}  A. Grothendieck, et al. {\paper rev\^etements \'etales et groupe
fondamental,}  SGA1, Lectures Notes in Math.  {\bf 224} (1971)

\ref{40}   J. Harris, I. Morrison,  { \paper Moduli of Curves,}  
Graduate Texts in Math.  {\bf 187}   Springer  (1998)  

\ref{41}  J. Harris, T. Graber, J. Starr,  {\paper A note on
Hurwitz schemes of covers of a positive genus curve, }  
[arXiv:math.AG/0205056]    

\ref{42}  J. Harris, D. Mumford, { \paper On the Kodaira dimension
of the moduli space of
curves,}  Invent.Math., {\bf 67}  (1982)   23--86  

\ref{43}   T. J. Jarvis,  {\paper Geometry of the moduli of higher
spin curves,}   Int. J. Math.  {\bf 11}  No.5  (2000)   637--663 

\ref{44}  T. J. Jarvis, {\paper Torsion-free sheaves and moduli of
generalized spin curves,}   Comp. Math. {\bf 110}   (1998)
  291--333

\ref{45}   T. J. Jarvis, T. Kimura, { \paper Orbifold quantum
cohomology of the classifying space of a finite group,}  
Orbifolds in mathematics and physics. Proc. Conf. on mathematical
aspects of orbifold string theory, Madison, WI, USA, May 4-8, 2001.
(A. Adem et al. eds.), Contemp. Math. {\bf 310}  (2002)  
123--134  

\ref{46}   T. J. Jarvis, R. Kaufmann, T. Kimura, {\paper Pointed
admissible $G$-covers and $G$-equivariant cohomological field theories,}
 [arXiv:math.AG/0302316] 
 
\ref{47}   G. Katz,  { \paper How tangents solve algebraic equations,
or a remarkable geometry of discriminant varieties,}
 [arXiv:math.AG/021128]

\ref{48}   F. F. Knudsen,  { \paper The projectivity of the moduli
space of stable curves (II),}   Math. Scand.  {\bf 52}  (1983)
 161--189  

\ref{49}   M. Kontsevich,  { \paper Enumeration of rational curves
via torus actions,}  The Moduli space of Curves (R. Dijkgraaf,
C. Faber, and G. van der Geer, eds), Prog. in Math.  {\bf 129}  
Birkh\"{a}user  (1995)   335--368  

\ref{50}   S. K. Lando,  { \paper Ramified coverings of the
two-dimensional sphere and intersection theory in spaces of meromorphic
functions on algebraic curves,}   Russ. Math. Surv. {\bf 57}   No.3
 (2002)   463--533  

\ref{51}  G. Laumon, L. Moret-Bailly, { \paper Champs Alg\'ebriques
Ergebnisse der Mathematik und ihrer Grenzgebiete, }  3. Folge. {\bf 39}
 Springer-Verlag  (2000)

\ref{52}   E. Looijenga, { \paper Smooth Deligne-Mumford
compactifications by means of Prym level structures}
  J. Alg. Geom., {\bf 3},  No.2  (1994)   283--293  

\ref{53}   Y. Manin, { \paper Frobenius manifolds, quantum
cohomology, and moduli spaces,}  l Colloquium Publications AMS
 (1999)  

\ref{54}   L. Moret-Bailly,  {\paper Construction de
rev\^etements de courbes point\'ees,}  
J. Alg. {\bf 240}, No.2  (2001)   505--534  

\ref{55}   S. Mochizuki,  { \paper The geometry of the
compactification of the Hurwitz scheme,}  Publ. RIMS.
Kyoto.Univ.  {\bf 31}, No.3  (1995)  355--441 

\ref{56}  I. Morisson, H. Pinkham,   {\paper Galois Weirstra\ss\
points and Hurwitz characters,}  Ann. of Math.   {\bf 124}  
(1986)  591--625  

\ref{57}  D. Mumford,  { \paper Geometric invariant theory, second ed.
Springer-Verlag}  (1981)  

\ref{58}  D. Mumford, {\paper Towards an enumerative geometry of the
moduli space of curves,  Arithmetic and geometry, Pap. dedic. I. R.
Shafarevich,} Vol. II, Prog. Math. {\bf 36}  Birka\"{u}ser
 (1983)   271--328  

\ref{59}  A. Okounkov, R. Pandharipande, { \paper Gromov-Witten
theory, Hurwitz numbers, and matrix models  I,}  
[arXiv:math.AG/0101147]

\ref{60} M. Raynaud,  { \paper $p$-groupes et r\'eduction
semi-stables des courbes,}  The Grothendieck Festschrift, Vol III,
Prog. Math.  {\bf 88}  Birkh\"{a}user  (1990) 179--197 

\ref{61}   M. Romagny,  { \paper Group actions on stacks and applications}
  Pr\'epublication Max Planck Institut f\"ur Mathematik
No.118, Bonn (2003)  

\ref{62}  M. Romagny, { \paper The stack of Potts curves and its fibre
at a prime of wild ramification,}  J. Alg.  {\bf 274}, No.2  (2004)
 772--803 

\ref{63}    M. Sa\"\i di, { \paper Rev\^etements mod\'er\'es et
groupe fondamental de graphe de groupes,}  Comp. Math.
 {\bf 107}, No.3  (1997)   319--338  

\ref{64}  J.-P. Serre, { \paper Repr\'esentations lin\'eaires des
groupes finis,} 5\`eme \'ed.   Hermann, Paris  (1988)  

\ref{65}  J. -P.  Serre, {\paper Topics in Galois theory,}  
Research Notes in Math.  {\bf 1}  Bartlett and Jones  (1992)  

\ref{66}   A. Vistoli,  { \paper Intersection theory on algebraic stacks and on their moduli spaces,}   Inventiones math.  {\bf 97} (1989) 613--670 

\ref{67}  B. Wajnryb,  { \paper Orbits of Hurwitz action for
coverings of a sphere with two special fibers,}   Indag. Math., New Ser.
 {\bf 7}, No.4  (1996)  549--558  

\ref{68}   S. Wewers,  { \paper Deformation of tame admissible covers
of curves,}   Aspects of Galois theory, London Math. Soc.
Lecture Note Ser. {\bf 256}   Cambridge Univ. Press  (1999)  

\ref{69}  S. Wewers, { \paper Construction of Hurwitz spaces,}
  Thesis, Preprint No. 21 of the IEM, Essen  (1998)

\end

\bye